\documentclass{article}
\usepackage{latexsym,amsmath,amssymb,amscd,bm,fancyheadings, pifont,}
\usepackage{mathrsfs}
\usepackage[dvips]{color}
\usepackage[dvips]{graphicx}
\begin{document}

\title{Linearized Analysis of Adiabatic Oscillations of Rotating Gaseous Stars  }
\author{Tetu Makino \footnote{Professor Emeritus at Yamaguchi University, Japan;  E-mail: makino@yamaguchi-u.ac.jp}}
\date{\today}
\maketitle

\newtheorem{Lemma}{Lemma}
\newtheorem{Proposition}{Proposition}
\newtheorem{Theorem}{Theorem}
\newtheorem{Definition}{Definition}
\newtheorem{Remark}{Remark}
\newtheorem{Corollary}{Corollary}
\newtheorem{Notation}{Notation}
\newtheorem{Assumption}{Assumption}
\newtheorem{Question}{Question}
\newtheorem{Situation}{Situation}

\numberwithin{equation}{section}

\begin{abstract}
We study adiabatic oscillations of rotating self-gravitating gaseous stars in mathematically rigorous manner. The internal motion of the star is supposed to be governed by the Euler-Poisson equations with rotation of constant angular velocity under the equation of state of the ideal gas. The motion is supposed to be adiabatic, but not to be barotropic in general. This causes a free boundary problem to gas-vacuum interface. Existence of solutions to the linearized equation in the Lagrangian coordinates of the perturbations around a fixed stationary solution, the eigenvalue problem with concept of quadratic pencil of operators, and the stability problem with a new concept of stability introduced in this article are discussed. 



\end{abstract}

\section{Introduction}

We consider the motion of a rotating gaseous star governed by the Euler-Poisson equation:

\begin{subequations}
\begin{align}
&\frac{D\rho}{Dt}+\rho\mathrm{div}_{\bm{x}}\bm{v}=0, \label{EPa}\\
&\rho\Big[\frac{D\bm{v}}{Dt}+2\bm{\Omega}\times\bm{v}
+\bm{\Omega}\times(\bm{\Omega}\times\bm{x})\Big]+{\nabla}_{\bm{x}}P+\rho{\nabla}_{\bm{x}}
\Phi 
=0, \label{EPb}\\
&\rho\frac{DS}{Dt}=0, \label{EPc} \\
&\triangle \Phi=4\pi\mathsf{G}\rho \label{1c}
\end{align}
\end{subequations}
on $(t,\bm{x})=(t,x^1,x^2,x^3) \in [0,+\infty[\times \mathbb{R}^3$.
We are denoting
\begin{equation}
 \nabla_{\bm{x}}=
\begin{bmatrix}
\frac{\partial}{\partial x^1} \\
\\
\frac{\partial}{\partial x^2} \\
\\
\frac{\partial}{\partial x^3}
\end{bmatrix}, \quad
\frac{D}{Dt}=\frac{\partial}{\partial t}+\sum_{k=1}^3 v^k\frac{\partial}{\partial x^k}, \quad \triangle=\sum_{k=1}^{3}\Big(\frac{\partial}{\partial x^k}\Big)^2.
\end{equation}
The unknown variables $\rho, P, S, \bm{v}=
(v^1,v^2,v^3)^{\top}$ are density, the pressure, the specific entropy, and 
 the velocity field, while
 \begin{equation}
\bm{\Omega}=\Omega\frac{\partial}{\partial x^3},
\end{equation}
$\Omega$ being a constant, the angular velocity of the rotation around the $x^3$-axis. 
Hereafer we denote
\begin{equation}
B=2\Omega
\begin{bmatrix}
0 & -1 & 0 \\
1 & 0 & 0 \\
0 & 0 & 0
\end{bmatrix}
\end{equation}
so that
\begin{equation}
B\bm{v}=2\bm{\Omega} \times \bm{v} =2\Omega
\begin{bmatrix}
-v^2 \\
v^1 \\
0
\end{bmatrix}.
\end{equation}

On the other hand $\Phi$ is the 
gravitational potential and $\mathsf{G}$ is a positive constant, while we suppose that the support of $\rho(t,\cdot)$ is compact for $\forall t $ and replace the Poisson equation \eqref{1c} by the Newton potential

\begin{equation}
\Phi(t,\mbox{\boldmath$x$})=-4\pi\mathsf{G}\mathcal{K}[\rho(t,\cdot)](\bm{x}), \label{NPot}
\end{equation}
where
\begin{equation}
\mathcal{K}[f](\bm{x}):=\frac{1}{4\pi}\int\frac{f(\mbox{\boldmath$x$}')}{
\|\mbox{\boldmath$x$}-\mbox{\boldmath$x$}'\|}d\mbox{\boldmath$x$}'. \label{defK}
\end{equation}

The pressure $P$ is supposed to be a prescribed function $\mathscr{P}(\rho, S)$ of $\rho, S$, where $\mathscr{P} \in C^1([0,+\infty[\times \mathbb{R})$ and
$\displaystyle \mathscr{P} >0,  \frac{\partial\mathscr{P}}{\partial\rho} >0$ for $\rho >0$ and $\mathscr{P}(0, S)=0$. So we read
\begin{equation}
P=\mathscr{P}(\rho, S) \label{1.8}
\end{equation}
in \eqref{EPb}.

We consider the initial data $\rho^0, S^0, \bm{v}^0$ for which the motion enjoys
\begin{equation}
\rho(0,\bm{x})=\rho^0(\bm{x}), \quad S(0,\bm{x})=S^0(\bm{x}), \quad \bm{v}(0,\bm{x})=\bm{v}^0(\bm{x}).
\end{equation}
We denote
\begin{equation}
P^0(\bm{x})=\mathscr{P}(\rho^0(\bm{x}), S^0(\bm{x})), \quad
\Phi^0(\bm{x})=-4\pi\mathsf{G}\mathcal{K}[\rho^0](\bm{x}).
\end{equation}\\

We suppose that there is fixed an axially symmetric stationary solution
 $(\rho, S, \bm{v})=(\rho_b(\bm{x}),S_b(\bm{x}), \bm{v}_b(\bm{x})) $, which satisfy \eqref{EPa}, \eqref{EPb}, \eqref{EPc}, \eqref{NPot},
such that $\mathfrak{R}_b=\{ \bm{x} \in \mathbb{R}^3 | \rho_b(\bm{x}) >0\}$ is a bounded domain.
The existence of stationary solutions is discussed in Section 2.\\

We consider small perturbations  at this fixed stationary solution.
The linearized equation which governs the perturbations, which will be called ``
{\bf ELASO }( the equation of linear adiabatic stellar oscillations )``, 
 turns out to be

\begin{equation}
\frac{\partial^2\bm{u}}{\partial t^2}+
B\frac{\partial\bm{u}}{\partial t}+\mathcal{L}\bm{u}=0
\quad\mbox{on}\quad [0,+\infty[\times \mathfrak{R}_b,
\label{Eqxi}
\end{equation}
where
\begin{equation}
\mathcal{L}\bm{u}
=\frac{1}{\rho_b}{\nabla}\mathfrak{d} P
-\frac{\mathfrak{d}\rho}{\rho_b^2}{\nabla}P_b+
{\nabla} \mathfrak{d}\Phi, \label{L}
\end{equation}
with
\begin{subequations}
\begin{align}
&\mathfrak{d}\rho=-\mathrm{div}(\rho_b\bm{u}),  \label{LASOa}\\
&\mathfrak{d} S=-(\bm{u}|{\nabla} S_b), \label{LASOb}\\
&\mathfrak{d} P=
\Big(\frac{\partial P}{\partial \rho}\Big)_b\mathfrak{d}\rho+
\Big(\frac{\partial P}{\partial S}\Big)_b\mathfrak{d} S, \label{LASOc}\\
&\mathfrak{d}\Phi=-4\pi\mathsf{G}\mathcal{K}[\mathfrak{d}\rho]. \label{LASOd}
\end{align}
\end{subequations}
Here we are denoting
$$
\Big(\frac{\partial P}{\partial \rho}\Big)_b=\partial_{\rho}\mathscr{P}(\rho_b, S_b), \quad
\Big(\frac{\partial P}{\partial S}\Big)_b=\partial_S\mathscr{P}(\rho_b, S_b). $$

The unknown variable $\bm{u}$ means
\begin{equation}
\bm{u}(t,\bm{x})=\bm{\varphi}(t,\bm{x})-\bm{\varphi}_b(t,\bm{x}) + \bm{u}^0 \label{1.16}
\end{equation}
where $\bm{u}^0$ is a given vector field on $\mathfrak{R}_b$ and $\bm{\varphi}(t,\bm{x})$ is the solution of 

\begin{equation}
\frac{\partial}{\partial t}\bm{\varphi}(t,{\bm{x}})=\bm{v}(t,\bm{\varphi}(t,\bm{x})), \quad\bm{\varphi}(0,\bm{x})=\bm{x},
\end{equation}
while $\bm{\varphi}_b(t,\bm{x})$ is the solution of
\begin{equation}
\frac{\partial}{\partial t}\bm{\varphi}_b(t,\bm{x})=\bm{v}_b(\bm{\varphi}_b(t,\bm{x})), \quad
\bm{\varphi}_b(0,\bm{x})=\bm{x}. \label{PHb}
\end{equation}

Here we suppose

\begin{Assumption} \label{Ass.LASO}
It holds that
\begin{equation}
\{ \bm{x} \  |\  \rho^0(\bm{x}) > 0 \}=\mathfrak{R}_b
\end{equation}
and $\bm{u}^0$ enjoys
\begin{subequations}
\begin{align}
&\rho^0(\bm{x})-\rho_b(\bm{x})=-\mathrm{div}\Big(\rho_b(\bm{x})\bm{u}^0(\bm{x})\Big), \\
&S^0(\bm{x})-S_b(\bm{x})=-\Big(\bm{u}^0(\bm{x})\Big|{\nabla} S_b(\bm{x})\Big).
\label{1.18b}
\end{align}
\end{subequations}
\end{Assumption}

The equation \eqref{Eqxi} is considered with the initial condition
\begin{equation}
\bm{u}|_{t=0}=\bm{u}^0,\quad \frac{\partial \bm{u}}{\partial t}\Big|_{t=0}=\bm{v}^0
\quad\mbox{on}\quad \mathfrak{R}_b. \label{IVxi}
\end{equation}\\

In Section 3  the derivation of this {\bf ELASO} \eqref{Eqxi} is discribed
under the situation which enjoys the assumption

\begin{Assumption} \label{Ass.ELASO}
$\rho-\rho_b, S-S_b, \bm{v}-\bm{v}_b$ are of $\mathscr{O}(\varepsilon)$ and moreover 
$\bm{v}_b$ itself is of $\mathscr{O}(\varepsilon)$.
\end{Assumption}

Here we use

\begin{Notation}
Considering small $\varepsilon$, a quantity $Q$ will denoted by $\mathscr{O}(\varepsilon)$ if $Q$ and its derivatives are of order $O(\varepsilon)$ uniformly on $t \in [0, T]$ for $\forall$ fixed $T \in [0, +\infty[$. 
\end{Notation}

In Section 4 we give the basic existence result to the initial value problem
\eqref{Eqxi}, \eqref{IVxi}, by realizing the integrodifferential operator $\mathcal{L}$ as a
self-adjoint operator $\bm{L}$ in the Hilbert space $\mathfrak{H}=L^2(\mathfrak{R}_b,\rho_bd\bm{x} ; \mathbb{C}^3)$. This is done by considring the first order system, which will be called `{\bf ELASO \ddag}', 
\begin{equation}
\frac{dU}{dt}+\bm{A}U=\bm{0},
\end{equation}
where
\begin{equation}
U=\begin{bmatrix}
\bm{u} \\
\\
\frac{d\bm{u}}{dt}
\end{bmatrix},
\quad
\bm{A}=
\begin{bmatrix}
O & -I \\
\\
\bm{L} & \bm{B}
\end{bmatrix},
\end{equation}
which is equivalent to {\bf ELASO}, and
 applying the Hille-Yosida theory. The applied theorem is formulated and proved in {\bf Appendix A} for the sake of selfcontainedness.\\

In Section 5 we discuss about the eigenvalue problem to the equation
 {\bf ELASO} \eqref{Eqxi}.  The issue is addressed to the problem of possibility of eigenmode expansions associated with {\bf ELASO}. The `eigenvalue' of {\bf ELASO} \eqref{Eqxi}  is 
$\pm \sqrt{-\lambda}$ for the 
 usual eigenvalue $\lambda$ of $\bm{L}$ only when $\Omega=0, B=O$, But, when $\Omega\not=0$, we need the concept `quadratic pencil' (after Bognar \cite{Bognar} ) and its `spectrum'. The result of the structure of the `spectrum' of the quadratic pencil by J. Dyson and B. F. Schultz \cite{DysonS} is justified. \\

In Section 6 we discuss about the stability of solutions of {\bf ELASO}\eqref{Eqxi}\eqref{IVxi}. We propose a new concept of stability based on seminorms on the space of the values $\bm{u}(t,\cdot)$ of the solution at instant $t$, taking into account that the magnitude of $\bm{u}(t,\cdot)$ itself is not essential but the magnitude of the perturbation $\mathfrak{d} \rho = -\mathrm{div}(\rho_b\bm{u}), \mathfrak{d}\bm{v} =\partial \bm{u}/\partial t -\bm{v}_b^L$ is essential in the discussion of stability. Examples of the seminorms are presented
with an open question. \\

Traditionally astrophysicists used to
discuss a linearized equation of adiabatic stellar oscillations which is slightly different from {\bf ELASO}\eqref{Eqxi}. In Section 7 we introduce the equation, which will be called `{\bf ELASO($\blacklozenge$)}', of astrophysists:

\begin{equation}
\frac{\partial^2\bm{\xi}}{\partial t^2}
+\bm{B}_{\blacklozenge}\frac{\partial\bm{\xi}}{\partial t}+\bm{L}_{\blacklozenge}\bm{\xi}=\bm{0}.
\end{equation}
Here $\bm{B}_{\blacklozenge}$ and $\bm{L}_{\blacklozenge}$ are slightly different from $\bm{B}$ and $\bm{L}$,
but  $\bm{B}_{\blacklozenge}=\bm{B}$ and $\bm{L}_{\blacklozenge}=\bm{L}$ when $\bm{v}_b=\bm{0}$. The
unknown variable $\bm{\xi}$ means
\begin{equation}
\bm{\xi}(t.\bm{x})=\bm{\varphi}(t\Big|0, \bm{\varphi}(t | 0 \| \bm{v}_b)\Big\| \bm{v} )-\bm{x}.
\end{equation}
Here $\bm{\varphi}(t|s,\bm{y})\| \bm{V})$ stands for the solution of
\begin{equation}
\frac{\partial}{\partial t}\bm{\varphi}(t|s,\bm{y}\|\bm{V})=\bm{V}(t,\bm{\varphi}(t|s,y\|\bm{V})),\quad
\bm{\varphi}(s|s,\bm{y}\|\bm{V}))=\bm{y},
\end{equation}
while $\bm{V}=\bm{v}$ or $\bm{v}_b$. The derivation of {\bf ELASO($\blacklozenge$)}
is performed inder the situation which enjoys the assumption

\begin{Assumption}\label{Ass.0ELASObl}
The initial configurations of $\rho, S$ are the same to those of the background, that is,
\begin{subequations}
\begin{align}
\rho^0&=\rho(0,\bm{x})=\rho_b(\bm{x}), \\
S^0&=S(0,\bm{x})=S_b(\bm{x}).
\end{align}
\end{subequations}
\end{Assumption}
This is stronger than Assumption \ref{Ass.LASO}, which reduces this when $\bm{u}^0=\bm{0}$.
But we suppose, instead of Assumption \ref{Ass.ELASO}, the weaker
\begin{Assumption}\label{Ass.ELASObl}
$\rho-\rho_b, S-S_b, \bm{v}-\bm{v}_b$ are of $\mathscr{O}(\varepsilon)$,
\end{Assumption}
without supposing that $\bm{v}_b= \mathscr{O}(\varepsilon)$. \\

Let us note the following relation between {\bf ELASO} and
{\bf ELASO($\blacklozenge$) }:
The solution $\bm{\xi}$ of {\bf ELASO($\blacklozenge$)} should satisfy
$\bm{\xi}(0, \cdot)=\bm{0}$; Consider the solution $\bm{u}$ of {\bf ELASO }
such that $\bm{u}^0=\bm{u}(0, \cdot)=\bm{0}$;  Then we have the relation
\begin{equation}
 \bm{\xi}(t,\bm{x})=\bm{u}(t, \bm{\varphi}(0|t,\bm{x}||\bm{v}_b)),
 \end{equation}
since
the Lagrange displacement
$\bm{u}$ in {\bf ELASO} \eqref{Eqxi} was defined by
 \eqref{1.16}, which reads
$$ \bm{u}(t, {\bm{x}})=\bm{\varphi}(t|0,{\bm{x}}|| \bm{v})
-\bm{\varphi}(t|0,{\bm{x}} ||\bm{v}_b)+
\bm{u}^0({\bm{x}}); $$
Therefore $ \bm{\xi}=\bm{u}$ when $\bm{v}_b=\bm{0}$ so that 
$\bm{\varphi}(t|s, \bm{x}||\bm{v}_b)=\bm{x} $, but $ \bm{\xi}\not=\bm{u}$ in general.\\

In Section 8 some properties of {\bf ELASO($\blacklozenge$) } are discussed in view of functional analysis, keeping in mind the analysis of {\bf ELASO} in Section 4.

\section{Existence of stationary solutions --Distorted Lane-Emden functions}
 
In this section we establish the existence of stationary solutions which enjoy good properties used in the following consideration, supposing

\begin{Assumption} \label{Ass.0}
There are a function $\mathscr{F} =\mathscr{F}(\upsilon, S) \in C^{\infty}([0,+\infty[\times\mathbb{R})$ and a number $\gamma$ such that $1<\gamma <2$, 
$\displaystyle \mathscr{F}(\upsilon, S) >0$ for $\upsilon \geq 0, |S|<\infty$, $\mathscr{F}+\frac{\gamma-1}{\gamma}\upsilon\frac{\partial\mathscr{F}}{\partial\upsilon} >0 $ for $\upsilon >0, |S| <\infty$, $\displaystyle \Big|\frac{\partial}{\partial\upsilon}\log\mathscr{F}\Big| \leq C$, for which
\begin{equation}
\mathscr{P}(\rho, S)=\rho^{\gamma}\mathscr{F}(\rho^{\gamma-1}, S) \quad \mbox{for}\quad \rho \geq 0, |S|<\infty.
\end{equation}
\end{Assumption}

We keep in mind the case of $\mathscr{F}(\upsilon, S)=\exp(S/\mathsf{C}_V)$, $\mathsf{C}_V$ being a positive constant, namely, the case when the equation of state \eqref{1.8} turns out to be
\begin{equation}
P=\rho^{\gamma}\exp\Big(\frac{S}{\mathsf{C}_V}\Big). \label{EOS.g}
\end{equation}
\\

Let us put the following 

\begin{Definition}\label{Def.2}
A triple of $t$-independent axially symmetric functions $(\rho_b, S_b, \bm{v}_b)
\in C_0^1(\mathbb{R}^3; [0,+\infty[)\times C^1(\mathbb{R}^3; \mathbb{R})
\times C^1(\mathbb{R}^3; \mathbb{R}^3)$
which satisfies \eqref{EPa}\eqref{EPb}\eqref{EPc} with $\Phi=\Phi_b$, $P=P_b$ determined by \eqref{NPot} \eqref{1.8} with $\rho=\rho_b, S=S_b$
 is called an {\bf admissible stationary solution}, if 
 
 0) there is a function $\omega_b(\varpi, z)$ of $(\varpi, z) \in [0,+\infty[\times \mathbb{R}$  such that
 \begin{equation}
 \bm{v}_b(\bm{x})=\omega_b(\varpi, z)\frac{\partial}{\partial x^3}\times \bm{x}
 =\omega_b(\varpi, z)
 \begin{bmatrix}
 -x^2 \\
 x^1 \\
 0
 \end{bmatrix},
 \end{equation}
 where 
 \begin{equation}
 \varpi=\sqrt{ (x^1)^2+(x^2)^2 }, \quad z=x^3 \qquad \mbox{for}\quad\bm{x}=(x^1,x^2,x^3),
 \end{equation}
 and the function $\bm{x} \mapsto \omega_b(\varpi, z)$ belongs to
 $C^1(\mathbb{R}^3)$;

 1) $\rho_b(\bm{x}), S_b(\bm{x})$ are functions of $(r,\zeta)$, and 
 $\mathfrak{R}_b=\{ \bm{x} | \rho_b(\bm{x})>0\}$ is of the form $\{ \bm{x} | r <R_b(\zeta) \}$, where
\begin{equation}
r=\|\bm{x}\|=\sqrt{(x^1)^2+(x^2)^2+(x^3)^2}, \quad \zeta=\frac{x^3}{r}
\qquad \mbox{for}\quad\bm{x}=(x^1,x^2,x^3),
\end{equation}
and $R_b(\zeta)$ is a continuous function of $\zeta \in [-1,1]$ such that
$R_b(\zeta) >0 \  \forall \zeta \in [-1,1]$;

2) The boundary $\partial\mathfrak{R}_b$ is of $C^{3,\alpha}$-class and the functions
$\rho_b^{\gamma-1}, S_b, \Phi_b=-4\pi\mathsf{G}\mathcal{K}[\rho_b] \in C^{\infty}(\mathfrak{R}_b) \cap
C^{3,\alpha}(\mathfrak{R}_b\cup\partial\mathfrak{R}_b)$, $\alpha$ being a positive number such that
$\displaystyle 0<\alpha <\Big(\frac{1}{\gamma-1}-1 \Big)\wedge 1 \Big(:=\min
\Big\{\frac{1}{\gamma-1}-1, 1\Big\}\Big)$ ;

3) $\displaystyle \frac{\partial\rho_b}{\partial r}, \frac{\partial P_b}{\partial r} <0$ in  $\mathfrak{R}_b$ and
\begin{equation}
\frac{\partial\rho_b}{\partial r}\leq -\frac{r}{C}, \quad
\frac{\partial P}{\partial r} \leq  -\frac{r}{C}
\quad\mbox{for}\quad 0<r \ll1; \label{PON}
\end{equation}

4) The boundary $\partial \mathfrak{R}_b$, on which $\rho_b=0$, is a physical vacuum boundary, that is,
\begin{equation} 
0 < \frac{\partial}{\partial \bm{n}}\sqrt{\Big(\frac{\partial P}{\partial \rho}\Big)_b } <+\infty\quad\mbox{on}
\quad \partial\mathfrak{R}_b,
\end{equation}
where $\bm{n}$ stands for the unit normal vector at the boundary point directed inward to $\mathfrak{R}_b$.
Here we note that $\partial P/\partial\rho$ is the square of the sound speed.

\end{Definition}

Note that the condition 2) 
requires $\|\nabla \rho_b\|, \|\nabla P_b\|\not=0$ for $r\not=0$.

\begin{Definition}

The admissible stationary solution $(\rho_b, S_b, \bm{v}_b)$ is said to be` {\bf barotropic}', if 
$\omega_b$ is independent of $z$, that is, 
$$\bm{v}_b(\bm{x})=\omega_b(\varpi)
\begin{bmatrix}
-x^2 \\
x^1 \\
0
\end{bmatrix},
$$ and
there is a function $\Pi \in C^{1,\alpha}(\mathbb{R})$ such that
$d\Pi(\upsilon)/\upsilon >0$ and it holds
\begin{equation}
P_b(\bm{x})=\Pi(\rho_b(\bm{x})^{\gamma-1})\quad\mbox{for}\quad \forall \bm{x} \in \mathfrak{R}_b.
\end{equation}

\end{Definition}

Note that, then, 
\begin{equation}
(\nabla P_b|\nabla\rho_b)=(\gamma-1)\upsilon\frac{d\Pi}{d\upsilon}\Big|_{\upsilon=\rho_b^{\gamma-1}}\|\nabla\rho_b\|^2 >0
\quad\mbox{for}\quad r\not=0.
\end{equation}\\

\begin{Remark}
It is well-known as the Poincar\'{e}-Wavre theorem (\cite [CHAPITRE II, 20, p.33]{Wavre} ) that 
$\displaystyle \frac{\partial \Omega_b}{\partial z}=0$ if and only if
$P_b=f(\rho_b)$, where $\Omega_b(\varpi, z):=\Omega+\omega_b(\varpi, z)$.
\end{Remark}

For admissible stationary solutions, let us put

\begin{Definition} \label{Def.G}
We denote
\begin{subequations}
\begin{align}
\Gamma &:=\frac{\rho}{P}\frac{\partial P}{\partial\rho}={\rho}\frac{\partial_{\rho}\mathscr{P}(\rho, S)}{\mathscr{P}(\rho, S)}, \\
\mathfrak{a}&:=
-\Big(\rho\frac{\partial P}{\partial \rho}\Big)^{-1}\Big(\frac{\partial P}{\partial S}\Big)\nabla S=
-\Big(\frac{1}{\Gamma P}\frac{\partial P}{\partial S}\Big)\nabla S \nonumber \\
&=-\frac{1}{\Gamma P}\nabla P +\frac{1}{\rho}\nabla \rho, \label{Def.aa}\\
\mathscr{A}&:= 
\frac{1}{\Gamma P}\frac{(\nabla P|\nabla\rho)}{\|\nabla\rho\|}
-\frac{\|\nabla\rho\|}{\rho} \nonumber \\
&=
-(\mathfrak{a}|\bm{n}), \quad\mbox{with}\quad 
\bm{n}:=\frac{\nabla\rho}{\|\nabla \rho\|} \quad
\mbox{wherever}\quad \|\nabla\rho\|\not=0.
\end{align}
\end{subequations}
and
\begin{subequations}
\begin{align}
\Gamma_b &:=\Big(\frac{\rho}{P}\frac{\partial P}{\partial\rho}\Big)_b={\rho_b}\frac{\partial_{\rho}\mathscr{P}(\rho_b, S_b)}{\mathscr{P}(\rho_b, S_b)}, \\
\mathfrak{a}_b&:=
-\Big(\rho\frac{\partial P}{\partial \rho}\Big)_b^{-1}\Big(\frac{\partial P}{\partial S}\Big)_b\nabla S_b=
-\Big(\frac{1}{\Gamma P}\frac{\partial P}{\partial S}\Big)_b\nabla S_b \nonumber \\
&=-\frac{1}{\Gamma_b P_b}\nabla P_b +\frac{1}{\rho_b}\nabla \rho_b, \\
\mathscr{A}_b&:=
\frac{1}{\Gamma_b P_b}\frac{(\nabla P_b|\nabla\rho_b)}{\|\nabla \rho_b\|}
-\frac{\|\nabla\rho_b\|}{\rho_b} \nonumber \\
 &=-(\mathfrak{a}_b|\bm{n}), \quad\mbox{with}\quad 
\bm{n}:=\frac{\nabla\rho_b}{\|\nabla \rho_b\|} \quad
\mbox{on}\quad   \mathfrak{R}_b\setminus\{O\}.
\end{align}
\end{subequations}

\end{Definition}

Since
\begin{equation}
\Gamma =\gamma +(\gamma-1)\upsilon\frac{\partial}{\partial\upsilon} \log \mathscr{F},
\end{equation}
we have $\Gamma _b \in C^{\infty}(\mathfrak{R}_b) \cap
C^{3,\alpha}(\mathfrak{R}_b\cup\partial\mathfrak{R}_b)$. Since
\begin{equation}
\frac{1}{\Gamma P}\frac{\partial P}{\partial S}=
\frac{\partial_S\log\mathscr{F} }{\gamma+(\gamma-1)\upsilon\frac{\partial}{\partial\upsilon}\log\mathscr{F}}, \label{*2.12}
\end{equation}
we see
$\mathfrak{a}_b \in  C^{\infty}(\mathfrak{R}_b) \cap
C^{2,\alpha}(\mathfrak{R}_b\cup\partial\mathfrak{R}_b)$. And
we see, at least, $$\mathscr{A}_b \in \bigcap_{0<\delta\ll 1}C^{1,\alpha}(\mathfrak{R}_b\cup \partial\mathfrak{R}_b \setminus \{ r < \delta\}) \cap L^{\infty}(\mathfrak{R}_b).$$ When $(\rho_b, S_b)$ is equatorially symmetric moreover, we can claim $\mathscr{A}_b \in C(\mathfrak{R}\cup \partial\mathfrak{R}_b)$
by evaluating $\mathscr{A}_b(O)=0$.\\

If the admissible stationary soltion is barotropic, then 
we have

\begin{equation}
\mathfrak{a}_b=-\mathscr{A}_b\bm{n} \quad\mbox{in}
\quad \mathfrak{R}_b, \label{2.5}
\end{equation}
and $ \mathscr{A}_b \in C^{\alpha}(\mathfrak{R}\cup \partial\mathfrak{R})$, since
$$
\mathscr{A}_b=
\frac{\gamma-1}{\Gamma P}\upsilon\frac{d\Pi}{d\upsilon}\Big|_{\upsilon=\rho_b^{\gamma-1}}\|\nabla\rho_b\|-\frac{\|\nabla\rho_b\|}{\rho_b}. $$

\begin{Remark}
When the stationary solution is barotropic, this does not mean that the perturbed motion to be considered should be barotropic. The perturbed state $\rho=\rho_b+\mathfrak{d}\rho, P=P_b+\mathfrak{d} P$ is approximately
$$\rho=\rho_b-\mathrm{div}(\rho_b\bm{u}), \quad P=P_b-
\Big(\frac{\Gamma P}{\rho}\Big)_b\mathrm{div}(\rho_b\bm{u})
+\Gamma_b P_b(\bm{u}|\mathfrak{a}_b),$$
 and ${\nabla}\rho,
{\nabla}P $ 
are not necessarily parallel.
\end{Remark}

\begin{Remark}
After \cite{DysonS} we can call $\mathfrak{a}_b$ the ``{\bf vector Schwarzschild discriminens}', and $\mathscr{A}_b$ the ``{\bf generalized Schwarzschild discriminant}''. But, in spite of the definition \cite[(3.7)]{DysonS}, we would like to define
the square of the ``{\bf generalized Brunt-V\"{a}is\"{a}l\"{a}  frequency (local buoyancy frequency)}'' 
$\mathscr{N}^2$ by
\begin{equation}
\mathscr{N}^2:=-\mathscr{A}_b\Big(\frac{1}{\rho_b}{\nabla} P_b\Big|\bm{n}\Big)
=-\mathscr{A}_b\frac{({\nabla} P_b|{\nabla} \rho_b)}{\rho_b\|{\nabla}\rho_b\|}. \label{DefN}
\end{equation}
Then, at $\|\bm{x}\|\not=0$,
$\mathscr{N}_b(\bm{x}) \in \mathbb{R}$ if and only if
$\mathscr{A}_b(\bm{x}) \leq 0$
, since $(\nabla P_b|\nabla\rho_b) >0$. It seems curious that the Lagrangian displacemnt of perturbation $\bm{\xi}$ is involved in the definition \cite[(3.7)]{DysonS}. 
Anyway, if $\Omega=\omega_b=0$ and the background is spherically symmetric, then it turns out to be $\displaystyle \bm{n}=-\bm{e}_r=-\frac{\bm{x}}{\|\bm{x}\|}$ and the definitions \eqref{2.5}, \eqref{DefN} of $\mathscr{A}$, $\mathscr{N}^2$ coincide with those \cite[(1.11), (1.12)]{TM2023}. In {\bf Appendix B} we derive the definition \eqref{DefN} of the  Brunt-V\"{a}is\"{a}l\"{a}  frequency $\mathscr{N}$ here by generalizing that found in meteorological texts, e.g. \cite[Chapter II, Section 21]{Brunt}, \cite[Chapter 2, Section 2.7.3]{Holton}.
\end{Remark} 

We are going to construct a barotropic admissible stationary solution.

Suppose 
\begin{equation}
\bm{v}_b(\bm{x})=\omega_b(\varpi)\frac{\partial}{\partial x^3}\times \bm{x}
=\omega_b(\varpi)
\begin{bmatrix}
-x^2 \\
x^1 \\
0
\end{bmatrix}, \label{vb}
\end{equation}
where $\omega_b \in C^1([0,+\infty[)$.

The equation \eqref{EPa} is reduced to
\begin{equation}
(\bm{v}_b|{\nabla} \rho_b)+ \rho_b\mathrm{div}\bm{v}_b=0. \label{**a}
\end{equation}
But 
$$(\bm{v}_b|{\nabla} \rho_b)=\omega(\varpi)\Big(-x^2\frac{\partial  \rho_b}{\partial  x^1}
+x^1\frac{\partial  \rho_b}{\partial  x^2}\Big)=
\omega(\varpi)\frac{\partial  \rho_b}{\partial  \phi}=0,
$$
if $\rho_b$ is axially symmetric. Moreover we have $\mathrm{div}\bm{v}_b=0$
for \eqref{vb}. Thus \eqref{**a} holds for axially symmetric $\rho_b$.

The equation \eqref{EPc} is reduced to
\begin{equation}
(\bm{v}_b|{\nabla} S_b)=0. \label{**2}
\end{equation}
This holds provided \eqref{vb} and axial symmetricity of $\rho_b$.

Thus the equations are reduced to
\begin{equation}
(\bm{v}_b|{\nabla})\bm{v}_b+2\bm{\Omega}\times \bm{v}_b +
\bm{\Omega}\times(\bm{\Omega}\times \bm{x})+
\frac{1}{\rho_b}{\nabla} P_b +{\nabla} \Phi_b=0 \label{1.1bb}
\end{equation}
with
$$\Phi_b=-4\pi\mathsf{G}\mathcal{K}[\rho_b],\quad
P_b=\mathscr{P}(\rho_b, S_b).
$$\\

We claim

\begin{Theorem}\label{Th.1}
Let a smooth function $\omega_b$ on $[0,+\infty[$, 
  a smooth function $\Sigma$ on $\mathbb{R}$ and a positive number $\rho_O$ be given.
  Supoose that
 \begin{align}
& \frac{d}{d\rho}\mathscr{P}(\rho, \Sigma(\rho^{\gamma-1}))= \nonumber \\
&=\upsilon\Big[ 
 \gamma +(\gamma-1)\Big(
\upsilon\frac{\partial}{\partial\upsilon}+
\upsilon\frac{d\Sigma}{d\upsilon}\frac{\partial}{\partial S}\Big) \Big]\mathscr{F}
 >0
\quad\mbox{for}\quad \upsilon =\rho^{\gamma-1}>0. \label{pDP}
\end{align}
If $\displaystyle \frac{6}{5}<\gamma <2$, and 
$\frac{1}{\rho_O}
\|\Omega +\omega_b\|_{L^{\infty}}^2$
is sufficiently small, 
 then there uniquely exists a
 barotropic admissible  stationary solution $(\rho_b, S_b, \bm{v}_b))$ such that $S_b=\Sigma(\rho_b^{\gamma-1})$  
$\rho_b(O)=\rho_{O}$, and $\bm{v}_b=\omega_b(\varpi)\frac{\partial}{\partial x^3}\times\bm{x}$.
\end{Theorem}

Proof . Consider the functions $f^P, f^{\varUpsilon}$ defined by 
\begin{align}
f^P(\rho)=f^P(\rho; \Sigma)&:=\mathscr{P}(\rho, \Sigma(\rho^{\gamma-1})) \nonumber \\
&=\rho^{\gamma}
\mathscr{F}(\rho^{\gamma-1}, \Sigma(\rho^{\gamma-1})), \\
f^{\varUpsilon}(\rho)=f^{\varUpsilon}(\rho; \Sigma)
&:=\int_0^{\rho}
\frac{Df^P(\rho')}{\rho'}d\rho'
\end{align}
for $\rho >0$. Thanks to the assumption \eqref{pDP} we have
$Df^P(\rho)
 >0$
for $\rho >0$, and there exists a smooth function $\Lambda$ on $\mathbb{R}$ such that
$\Lambda(0)=0$ and
\begin{equation}
f^P(\rho)=\mathsf{A}\rho^{\gamma}(1+\Lambda(\rho^{\gamma-1}))
\end{equation}
for $\rho >0$. Here 
\begin{equation}
\mathsf{A}=\mathsf{A}(\Sigma):=\mathscr{F}(0, \Sigma(0)) 
\end{equation}
 is a positive constant. 
Then we have
\begin{equation}
\varUpsilon=f^{\varUpsilon}(\rho)=\frac{\gamma\mathsf{A}}{\gamma-1}\rho^{\gamma-1}(1+\Lambda_1(\rho^{\gamma-1}))
\end{equation}
for $\rho >0$, where $\Lambda_1$ is a smooth function on $\mathbb{R}$ such that 
$\Lambda_1(0)=0$, and  the inverse function $f^{\rho}$ of $f^{\varUpsilon}$ 
\begin{equation}
f^{\rho}(\varUpsilon)=\Big(\frac{\gamma-1}{\gamma\mathsf{A}}\Big)^{\frac{1}{\gamma-1}}
(\varUpsilon \vee 0)^{\frac{1}{\gamma-1}}(1+\Lambda_2(\varUpsilon))
\end{equation}
is given so that
$ \rho=f^{\rho}(\varUpsilon) \Leftrightarrow \varUpsilon=f^{\varUpsilon}(\rho)$ for $\varUpsilon>0 (\rho >0)$. Here $\varUpsilon\vee 0$ stands for $\max( \varUpsilon, 0)$ and  $ \Lambda_2$ is a smooth functions on
$\mathbb{R}$ such that $\Lambda_2(0)=0$.

As shown in \cite{JJTM2019} with the inertial coordinate system, 
the problem is reduced to :
\begin{equation}
\varUpsilon_b(\bm{x})+\Phi_b(\bm{x})
-\mathfrak{B}(\varpi)
=
\varUpsilon_O+\Phi_b(O) \quad\mbox{on}\quad \mathfrak{D}_0 \label{2.27}
\end{equation}
with
\begin{align}
\Phi_b(\bm{x})&=-4\pi\mathsf{G} \mathcal{K}[\rho_b], \nonumber \\
 \rho_b(\bm{x})&=f^{\rho}(\varUpsilon_b(\bm{x}))=\Big(\frac{\gamma-1}{\gamma\mathsf{A}}\Big)^{\frac{1}{\gamma-1}}
(\varUpsilon_b(\bm{x}) \vee 0)^{\frac{1}{\gamma-1}}(1+\Lambda_2(\varUpsilon_b(\bm{x})).
\end{align}
and
\begin{equation}
\mathfrak{B}(\varpi)=\mathfrak{B}(\varpi; \Omega+\omega_b):=\int_0^{\varpi}(\Omega+\omega_b(\varpi'))^2\varpi'd\varpi'.
\end{equation}
Here $\mathfrak{D}_0$ is a domain such
that
$\varUpsilon_b \in C^{2,\alpha}(\mathfrak{D}_0 \cup \partial \mathfrak{D}_0)$ and
$\{ \varUpsilon_b \geq 0 \} \subset \mathfrak{D}_0$.

Put
\begin{equation}
\varUpsilon_O=\varUpsilon_O(\rho_O, \Sigma):=f^{\varUpsilon}(\rho_O)=
\frac{\gamma\mathsf{A}}{\gamma-1}\rho_O^{\gamma-1}
(1+\Lambda_1(\rho_O^{\gamma-1}))
\end{equation}
by the given $\rho_O$, and put 
\begin{align}
&\mathsf{a}=\mathsf{a}(\rho_O, \Sigma):=
\frac{1}{\sqrt{4\pi\mathsf{G}}}
\Big(\frac{\mathsf{A}\gamma}{\gamma-1}\Big)^{ \frac{1}{2(\gamma-1)} }
\varUpsilon_O^{ -\frac{2-\gamma}{2(\gamma-1)} },\\
&\mathfrak{b}(\eta)=\mathfrak{b}(\eta; \rho_O, \Sigma, \Omega+\omega_b):=
\varUpsilon_O^{-1}\mathfrak{B}(\mathsf{a}\eta).
\end{align}\\

We put 
$$ \|\mathfrak{b}\|_1=\sup_{\eta \leq r_{\infty}/\mathsf{a}} |\mathfrak{b}| +\sup_{\eta \leq r_{\infty}/\mathsf{a}} \Big|\frac{d\mathfrak{b}}{d\eta}\Big|,$$
where $\mathfrak{D}_0 \subset \{ \|\bm{x}\| \leq r_{\infty} \}$,  
and
note
$$
\|\mathfrak{b}\|_1 
\leq \frac{1}{4\pi\mathsf{G}}\Big(\frac{\mathsf{A}\gamma }{\gamma-1}\Big)^{\frac{1}{\gamma-1}}
\varUpsilon_O^{-\frac{1}{\gamma-1}}
\|\Omega+\omega_b\|_{L^{\infty}(0,r_{\infty})}
\Big(\frac{r_{\infty}^2}{2}+r_{\infty}\Big).
$$

We suppose\\

{\bf Assumption}: {\it
 $\frac{6}{5} < \gamma <2$, and
$\mathfrak{b}$ is sufficiently small, say,  $\|\mathfrak{b}\|_1 \leq \beta^0$, $\beta^0$ being a positive number depending on $\gamma,  \Lambda$.} \\

Then we put
\begin{align}
&\rho_b(\bm{x})=f^{\rho}(\varUpsilon_O\Theta)=\Big(\frac{\gamma-1}{\mathsf{A}\gamma}\Big)^{\frac{1}{\gamma-1}}(\varUpsilon_O\Theta \vee 0)^{\frac{1}{\gamma-1}}
(1+\Lambda_2(\varUpsilon_O \Theta)) \nonumber \\
\mbox{with} & \nonumber \\
& \Theta=\Theta\Big(\frac{r}{\mathsf{a}},\zeta \Big),\nonumber \\
\mbox{and}& \nonumber \\
&S_b(\bm{x})=\Sigma(\rho_b(\bm{x})^{\gamma-1}).
\end{align}
Here
 $\Theta(\xi, \zeta)$
is the "{\bf distorted Lane-Emden function}",
which is the solution of the integral equation
$$
\Theta(\xi,\zeta)=1+\mathfrak{b}(\xi \sqrt{1-\zeta^2})+
\mathscr{K}g(\xi,\zeta)-\mathscr{K}g(0,0)
$$
with
\begin{align*}
\mathscr{K}g(\xi,\zeta)&=\frac{1}{4\pi}\int_{-1}^1\int_0^{+\infty}
K(\xi,\zeta, \xi', \zeta')
g(\xi', \zeta')\xi'^2d\xi'd\zeta', \\
K(\xi, \zeta, \xi', \zeta')&=\int_0^{2\pi}
\frac{d\phi}{ \sqrt{\xi^2+\xi'^2-2\xi\xi'(\sqrt{1-\zeta^2}\sqrt{1-\zeta'^2}\cos\phi +\zeta\zeta')} }
, \\
g(\xi,\zeta)&=\Big(\frac{\gamma \mathsf{A}}{\gamma-1}\Big)^{\frac{1}{\gamma-1}}f^{\rho}\Big(\Upsilon_O\Theta(\xi, \zeta) \Big)\\
&=
(\Theta(\xi,\zeta)\vee 0)^{\frac{1}{\gamma-1}}
\Big(1+\Lambda_2(\varUpsilon_O\Theta(\xi,\zeta)) \Big),
\end{align*}
and enjoys the following properties:

1)\  The function $\bm{\xi}\mapsto \Theta(\|\bm{\xi}\|, \frac{\xi^3}{\|\bm{\xi}\|}
)$ belongs to $C^{2,\alpha}(\{ \|\bm{\xi}\| \leq {\Xi}_0\})$
with $0 <\alpha <(\frac{1}{\gamma-1}-1)\wedge 1$,
${\Xi}_0=4\xi_1(\frac{1}{\gamma-1})$, $\xi_1(\frac{1}{\gamma-1})$ being the zero of the Lane-Emden function
$\theta(\xi; \frac{1}{\gamma-1})$ of the index $\frac{1}{\gamma-1}$, that is, the solution
of
$$
-\frac{1}{\xi^2}\frac{d}{d\xi}\xi^2\frac{d\theta}{d\xi}=(\theta\vee 0)^{\frac{1}{\gamma-1}},\quad
\theta|_{\xi=0}=1. $$

2)\   $\Theta(0,\zeta)=1$ and there is a curve
$\zeta \in [-1,1]\mapsto {\Xi}_1(\zeta)$ such that
$$\xi_1(\frac{1}{\gamma-1})\leq {\Xi}_1(\zeta) <\frac{1}{2}{\Xi}_0=2
\xi_1(\frac{1}{\gamma-1})$$
and, for $0 \leq \xi <{\Xi}_0$,
\begin{equation}
 0<\Theta(\xi, \zeta)  \Leftrightarrow
0\leq \xi <{\Xi}_1(\zeta). \label{*}
\end{equation}

3)\  $\displaystyle \frac{\partial\Theta}{\partial \xi} \leq -\frac{\xi}{C} $ with a positive number $C$.\\

For proof of the existence of the distorded Lane-Emden function $\Theta$, see \cite{JJTM2019}.\\

For the sake of simplicity we consider $\Theta(\xi, \zeta) <0 $ for $\xi \geq {\Xi}_0$
by modyfing the values of $\Theta$ on $\xi > {\Xi}_0$,
 so that \eqref{*} is valid for $0 \leq \xi <\infty, |\zeta|\leq 1$.
 Let us denote
 \begin{equation}
 \mathfrak{D}_0=\{ \bm{x} \in \mathbb{R}^3 \Big| \|\bm{x}\| < R_{b0} \}, \quad R_{b0}=\mathsf{a}{\Xi}_0=4\mathsf{a}\xi_1(\frac{1}{\gamma-1}).
 \end{equation}
 Then 
 \begin{equation}
\varUpsilon_b(\bm{x})=\varUpsilon_O\Theta\Big(\frac{\|\bm{x}\|}{\mathsf{a}}, \frac{x^3}{\|\bm{x}\|}\Big),
\end{equation}
which belongs to $C^{2,\alpha}(\mathfrak{D}_0\cup\partial \mathfrak{D}_0)$, satisfies
\eqref{2.27}. \\

We  denote 
\begin{subequations}
\begin{align}
\mathfrak{R}_b&=\{ \rho_b >0\}
=\{ \Theta >0\} =\{ r < R_{1b}(\zeta) \}, \\
\partial\mathfrak{R}_b&=\{ \Theta=0\}=\{ r=R_{1b}(\zeta) \},
\end{align}
\end{subequations}
where we put
$R_{1b}(\zeta)
=\mathsf{a}{\Xi}_1(\zeta) 
$.

Since 
\begin{equation}
P_b=\Pi(\rho_b^{\gamma-1}),
\end{equation}
where
\begin{equation}
\Pi(\upsilon)=(\upsilon \vee 0)^{\frac{\gamma}{\gamma-1}}
\mathscr{F}(\upsilon, \Sigma(\upsilon)),
\end{equation}
  the stationary solution is barotropic. Here $\Pi \in C^{1,\alpha}$ for 
$\displaystyle \frac{\gamma}{\gamma-1} >1+\alpha$.

As for 
$\displaystyle \mathfrak{a}_b=-\Big(\frac{1}{\Gamma P}\frac{\partial P}{\partial S}\Big)_b\nabla S_b $,
 since $S_b=\Sigma(\rho_b^{\gamma-1})$ with $\Sigma \in C^{\infty}$, we see
\begin{equation}
\mathfrak{a}_b=-
\Big(\frac{1}{\Gamma P}\frac{\partial P}{\partial S}\Big)_b
{\nabla}S_b
=-
\Big(\frac{1}{\Gamma P}\frac{\partial P}{\partial S}\Big)_b
\frac{d\Sigma}{d\upsilon}{\nabla}\upsilon\Big|_{\upsilon=\rho_b^{\gamma-1}},
\end{equation}
and
\begin{equation}
\mathscr{A}_b=\Big(\frac{1}{\Gamma P}\frac{\partial P}{\partial S}\Big)_b
\upsilon\frac{d\Sigma}{d\upsilon}
\|{\nabla}\upsilon\|\Big|_{\upsilon=\rho_b^{\gamma-1}}. \label{*2.41}
\end{equation}
Since $\rho_b^{\gamma-1} \in C^{2,\alpha}(\mathfrak{R}_b \cup \partial\mathfrak{R}_b)$, we have
$\mathscr{A}_b \in C^{1,\alpha}(\mathfrak{R}_b \cup \partial \mathfrak{R}_b)$.
Recall \eqref{*2.12}.
\\

Summing up, we get an admissible barotropic stationary solution $(\rho_b, S_b, \bm{v}_b)$. $\square$\\

\begin{Remark} \label{Rem.2}
Let us note that, by \eqref{*2.41}, we see
$$
\Big(\frac{d\varUpsilon}{d\rho}\Big)_b=Df^{\varUpsilon}(\rho_b)=
=\Big(\frac{\Gamma P}{\rho^2}\Big)_b
\Big(1+\mathscr{A}_b
\frac{\rho_b}{\|\rho_b\|}\Big).
$$
Therefore, if $\displaystyle \frac{d\Sigma}{d\upsilon}=0\quad \forall \upsilon$, that is, $\Sigma$ is constant,
then $\mathscr{A}_b=0$ and it holds
\begin{equation}
\Big(\frac{d\varUpsilon}{d\rho}\Big)_b=\Big(\frac{\Gamma P}{\rho^2}\Big)_b.
\end{equation}
We shall say that the background is `{\bf isentropic}' when $\mathfrak{a}_b=0$ as in this case. Then, if
$\mathscr{F}(\upsilon, S)=\exp(S/\mathsf{C}_V)$, we have $\Gamma=\gamma$ and
$$ P=\mathsf{A}\rho^{\gamma},\quad \varUpsilon=\frac{\gamma \mathsf{A}}{\gamma-1}\rho^{\gamma-1}; \quad S=\mathsf{C}_V\log \mathsf{A}.
$$
\end{Remark}

\section{Derivatin of {\bf ELASO}  }

Let us fix a compactly supported axially symmetric stationary solution $(\rho_b, S_b, \bm{v}_b)$ with
$\displaystyle \bm{v}_b(\bm{x})=\omega_b(\varpi,z)\frac{\partial}{\partial x^3}\times \bm{x}$.
In order to investigate solutions $(\rho, S, \bm{v})$ near the fixed stationary solution, for which the domain
$\mathfrak{R}(t)=\{ \bm{x} \Big| \rho(t,\bm{x}) >0 \}$ should move with the free boundary
$\Sigma(t)=\partial \mathfrak{R}(t)$, we introduce the Lagrangian co-ordinate to describe the equations.

Let us look at a solution $(\rho, S, \bm{v})$ which belongs to
$C^1([0, T[\times \mathbb{R}^3)$ with
$\mathfrak{R}(t)=\{ \bm{x} \Big| \rho(t,\bm{x}) >0 \}$,
provided that it exists. Suppose that $(\mathfrak{R}(t), \partial \mathfrak{R}(t))$ is differmorphic onto
$(\mathfrak{R}_b, \partial \mathfrak{R}_b)$. We are supposing that
\begin{equation}
\mathfrak{R}(0)=\{ \bm{x} | \rho^0(\bm{x}) >0 \}=\mathfrak{R}_b. \label{3.1}
\end{equation}
However we do not suppose that $\rho^0=\rho_b$.

We can consider the flow mapping
$\bm{\varphi} \in C^1([0,T[\times \mathbb{R}^3)$ defined by
\begin{equation}
\frac{\partial}{\partial t}\bm{\varphi}(t,\bar{\bm{x}})=\bm{v}(t,\bm{\varphi}(t,\bar{\bm{x}})), \quad\bm{\varphi}(0,\bar{\bm{x}})=\bar{\bm{x}},
\end{equation}
and $\bm{\varphi}_b(t, \bar{\bm{x}}) \in C^1([0,+\infty[; \mathbb{R}^3)$ defined by 
\begin{equation}
\frac{\partial}{\partial t}\bm{\varphi}_b(t,\bar{\bm{x}})=\bm{v}_b(\bm{\varphi}_b(t,\bar{\bm{x}})), \quad\bm{\varphi}_b(0,\bar{\bm{x}})=\bar{\bm{x}},
\end{equation}
which can be written explicitely as
\begin{equation}
\bm{\varphi}_b(t,\bar{\bm{x}})=
\begin{bmatrix}
\cos(\omega_bt) & -\sin(\omega_bt) & 0 \\
\sin(\omega_bt) & \cos(\omega_bt) & 0 \\
0 & 0 & 1
\end{bmatrix}
\bar{\bm{x}},
\end{equation}
where $\omega_b=\omega_b(\bar{\varpi}, \bar{z}),
\bar{\varpi}=\sqrt{ (\bar{x}^1)^2+(\bar{x}^2)^2 }, \bar{z}=\bar{x}^3$.

Through the mapping $\bm{x}=\bm{\varphi}(t, \bar{\bm{x}})$ the equations on $(t,\bm{x})$ can be written as  equations on $(t,\bar{\bm{x}})$. This is the Lagrangian description. \\

Hereafter we denote
\begin{equation}
(D_{\bar{\bm{x}}}{\bm{x}})(t, \bar{\bm{x}})=D\bm{\varphi} (t,\bar{\bm{x}})=\Big(\frac{\partial}{\partial \bar{\bm{x}}^{\alpha}}\Big(\bm{\varphi}(t,\bar{\bm{x}})\Big)^j\Big)_{j,\alpha}.
\end{equation}
Then it hold that $D_{\bar{\bm{x}}}\bm{x}(0,\bar{\bm{x}})=I$ and
\begin{equation}
D_{\bar{\bm{x}}}\bm{x}(t,\bar{\bm{x}})=D\bm{\varphi}(t,\bar{\bm{x}})=
\exp\Big[
\int_0^t(D_{\bm{x}}\bm{v})(t',\bm{\varphi}(t',\bar{\bm{x}}))dt'\Big] \label{1.29}
\end{equation}
 Here
$D_{\bm{x}}\bm{v}$ stands for the matrix field
$\displaystyle \Big(\frac{\partial v^k}{\partial x^j}\Big)_{k,j}$. 

Let us denote
\begin{align}
&\rho^L(t,\bar{\bm{x}})=\rho(t,\bm{\varphi}(t,\bar{\bm{x}})), \quad
 P^L(t,\bar{\bm{x}})=P(t,\bm{\varphi}(t,\bar{\bm{x}})), \nonumber \\
& S^L(t,\bar{\bm{x}})=S(t,\bm{\varphi}(t,\bar{\bm{x}})), \quad
\bm{v}^L(t,\bar{\bm{x}})=\bm{v}(t,\bm{\varphi}(t,\bar{\bm{x}}))
\end{align}
for $(t,\bar{\bm{x}})\in [0, +\infty[ \times \mathfrak{R}_b$.

The equations \eqref{EPa}, \eqref{EPb}, \eqref{EPc}
read
\begin{subequations}
\begin{align}
&\frac{\partial\rho^L}{\partial t}+(\mathrm{div}_{\bm{x}}\bm{v})(t,\bm{\varphi}(t,\bar{\bm{x}}))=0, \label{LEa} \\
&\frac{\partial }{\partial t}
(\bm{v}^L-\bm{v}_b^L)+
B(\bm{v}^L-\bm{v}_b^L)+
\Big(
\frac{1}{\rho}({\nabla}_{\bm{x}}P)-
\frac{1}{\rho_b}({\nabla}_{\bm{x}}P_b) + \nonumber \\
&+({\nabla}_x(\Phi-\Phi_b))
\Big)^L=0, \label{LEb}\\
&\frac{\partial S^L}{\partial t}=0. \label{LEc}
\end{align}
\end{subequations}
Here, for function $Q(t,\bm{x})$, we denote 
\begin{equation}
Q^L(t,\bar{\bm{x}})=Q(t, \bm{\varphi}(t, \bar{\bm{x}}))
\end{equation}
 in general.\\

Integrating \eqref{LEc}, we have
\begin{equation}
S^L(t,\bar{\bm{x}})=S^0(\bar{\bm{x}}).
\end{equation}

Integrating \eqref{LEa}, we have
\begin{equation}
\rho^L(t,\bar{\bm{x}})=\rho^0(\bar{\bm{x}})
\exp\Big[-\int_0^t
(\mathrm{div}_{\bm{x}}\bm{v})(s,\bm{\varphi}(s,\bar{\bm{x}}))ds \Big] \label{Eqrho}
\end{equation}
for $(t,\bar{\bm{x}}) \in [0,+\infty[\times \mathfrak{R}_b$.\\

Let us denote
\begin{subequations}
\begin{align}
\mathfrak{d}\rho(t,\bar{\bm{x}})&=\rho^L(t,\bar{\bm{x}})-\rho_b^L(t,\bar{\bm{x}})=
\rho(t,\bm{\varphi}(t,\bar{\bm{x}}))-\rho_b(\bm{\varphi}(t,\bar{\bm{x}})), \\
\mathfrak{d} P(t,\bar{\bm{x}})&=P^L(t,\bar{\bm{x}})-P_b^L(t,\bar{\bm{x}})=
P(t,\bm{\varphi}(t,\bar{\bm{x}}))-P_b(\bm{\varphi}(t,\bar{\bm{x}})), \\
\mathfrak{d} S(t,\bar{\bm{x}})&=S^L(t,\bar{\bm{x}})-S_b^L(t,\bar{\bm{x}})=
S(t,\bm{\varphi}(t,\bar{\bm{x}}))-S_b(\bm{\varphi}(t,\bar{\bm{x}}))= \nonumber \\
&=S^0(\bar{\bm{x}})-S_b(\bm{\varphi}(t, \bar{\bm{x}})), \\
\mathfrak{d}\Phi(t,\bar{\bm{x}})&=\Phi^L(t,\bar{\bm{x}})-\Phi_b^L(t,\bar{\bm{x}})=
\Phi(t,\bm{\varphi}(t,\bar{\bm{x}}))-\Phi_b(\bm{\varphi}(t,\bar{\bm{x}})), \\
\mathfrak{d}\bm{v}(t,\bar{\bm{x}})&= \bm{v}^L(t,\bar{\bm{x}})-\bm{v}_b^L(t,\bar{\bm{x}}),
\end{align}
\end{subequations}
and
\begin{subequations}
\begin{align}
{\delta}\rho(t,\bm{x})&=\rho(t,\bm{x})-\rho_b(t,\bm{x}), \\
{\delta} P(t,{\bm{x}})&=P(t,{\bm{x}})-P_b(t,{\bm{x}}), \\
{\delta} S(t,{\bm{x}})&=S(t,{\bm{x}})-S_b(t,{\bm{x}}), \\
{\delta}\Phi(t,{\bm{x}})&=\Phi(t,{\bm{x}})-\Phi_b(t,{\bm{x}}), \\
{\delta}\bm{v}(t,{\bm{x}})&= \bm{v}(t,{\bm{x}})-\bm{v}_b(t,\bm{x}).
\end{align}
\end{subequations}

Then the equations 
\eqref{LEb} and \eqref{1.1bb} read
\begin{align}
&\frac{\partial}{\partial t}\mathfrak{d}\bm{v}+2\bm{\Omega}\times\mathfrak{d}\bm{v} + \nonumber \\
&+\Big(
\frac{1}{\rho_b}{\nabla}_{\bm{x}}P -\frac{1}{\rho_b}{\nabla}_{\bm{x}}P_b +
{\nabla}_{\bm{x}}(\Phi-\Phi_b)
\Big)^L
=
-\Big((D_{\bm{x}}\bm{v}_b)\delta \bm{v}
\Big)^L. \label{*1.1}
\end{align} \\

\textbullet \  
Now let us derive the linearized approximation of \eqref{Eqrho},
\eqref{*1.1},
 supposing Assumption \ref{Ass.ELASO}: {\it 
 $\rho-\rho_b, S-S_b, \bm{v}_b, \bm{v}-\bm{v}_b$ are of $\mathscr{O}(\varepsilon)$.} \\

First of all we note that \eqref{1.29} shows

\begin{align}
D_{\bar{\bm{x}}}\bm{x}(t,\bar{\bm{x}})&=I+\int_0^t(D_{\bm{x}}\bm{v})(s,\bm{\varphi}(s,\bar{\bm{x}}))ds +
\mathscr{O}(\varepsilon^2) \\
&=I+\mathscr{O}(\varepsilon), \label{3.16}
\end{align}
since we are assuming $\bm{v}=\bm{v}_b+(\bm{v} -\bm{v}_b) =\mathscr{O}(\varepsilon)$.
Of course
\begin{align}
\bm{\varphi}(t,\bar{\bm{x}})&=\bar{\bm{x}}+\int_0^t\bm{v}(s, \bm{\varphi}(s, \bar{\bm{x}}))ds \nonumber \\
&=\bar{\bm{x}} + \mathscr{O}(\varepsilon).
\end{align}

Therefore, if $\bm{Q}=\bm{Q}(t,\bm{x}) =\mathscr{O}(\varepsilon)$, then
\begin{equation}
\mathrm{div}_{\bm{x}}\bm{Q}(t,\bm{x})=\mathrm{div}_{\bar{\bm{x}}}\bm{Q}^L(t,\bar{\bm{x}})+\mathscr{O}(\varepsilon^2),
\end{equation}
and, if $Q=Q(t,\bm{x})=\mathscr{O}(\varepsilon)$, then
\begin{equation}
{\nabla}_{\bm{x}}Q(t,\bm{x})={\nabla}_{\bar{\bm{x}}} Q^L(t,\bar{\bm{x}})+\mathscr{O}(\varepsilon^2),
\end{equation}
 while $\bm{x}=\bm{\varphi}(t,\bar{\bm{x}})$. In other words
\begin{align}
&\Big(\mathrm{div}_{\bm{x}}\bm{Q}\Big)^L=\mathrm{div}_{\bar{\bm{x}}}\bm{Q}^L+\mathscr{O}(\varepsilon^2), \label{3.20}\\
&\Big({\nabla}_{\bm{x}}Q\Big)^L={\nabla}_{\bar{\bm{x}}}Q^L +\mathscr{O}(\varepsilon^2). \label{3.21}
\end{align}\\

\begin{Remark}
If we assume only $\bm{v}-\bm{v}_b =\mathscr{O}(\varepsilon)$, but do not assume  $\bm{v}_b=\mathscr{O}(\varepsilon)$, then \eqref{3.16} can be not the case, but we have
\begin{align*}
D_{\bar{\bm{x}}}\bm{\varphi}(t,\bar{\bm{x}})&=D_{\bar{\bm{x}}}\bm{\varphi}_b(t, \bar{\bm{x}})+\mathscr{O}(\varepsilon), \\
D_{\bar{\bm{x}}}\bm{\varphi}_b(t,\bar{\bm{x}})&=
\begin{bmatrix}
\cos(\omega_bt) & -\sin(\omega_bt) & 0 \\
\sin(\omega_bt) & \cos(\omega_bt) & 0 \\
0 & 0 & 1
\end{bmatrix}+
\frac{1}{\varpi}\frac{\partial\omega_b}{\partial\varpi}t
\begin{bmatrix}
-x_b^2\bar{x}^1 & -x_b^1\bar{x}^2 & 0 \\
x_b^1\bar{x}^1 & x_b^1\bar{x}^2 & 0 \\
0 & 0 &  0
\end{bmatrix}
\\
&+ \frac{\partial\omega_b}{\partial z}t
\begin{bmatrix}
0 & 0 & -x_b^2 \\
0 & 0 & x_b^1 \\
0 & 0 & 0
\end{bmatrix}
,
\end{align*}
where $\bm{x}_b=\bm{\varphi}_b(t,\bar{\bm{x}})$, 
$\omega_b=\omega_b(\bar{\varpi},\bar{z})=\omega_b(\varpi, z)$,
$\varpi=\sqrt{(x_b^1)^2+(x_b^2)^2) }, \bar{\varpi}=\sqrt{(\bar{x}^1)^2+(\bar{x}^2)^2) }$,
$z=x_b^3, \bar{z}=\bar{x}^3$.
Thus, when $\omega_b \not= 0$, we have
$D_{\bar{\bm{x}}}\bm{x} \not= I +\mathscr{O}(\varepsilon)$, and the formulae \eqref{3.20},\eqref{3.21}  do  not work.
\end{Remark}

Let us look at \eqref{Eqrho}, which can be written as 
\begin{equation}
\rho(t,\bm{\varphi}(t,\bar{\bm{x}}))=
\rho^0(\bar{\bm{x}})
\exp\Big[-\int_0^t(\mathrm{div}_{\bm{x}}{\delta}\bm{v})(s,\bm{\varphi}(s,\bar{\bm{x}}))ds\Big]. \label{*3.1}
\end{equation}
 Let us verify it. Look at
 \begin{align*}
 &\exp\Big[-\int_0^t(\mathrm{div}_{\bm{x}}\bm{v})(s,\bm{\varphi}(s,\bar{\bm{x}}))ds\Big]= \\
 &=
 \exp\Big[-\int_0^t(\mathrm{div}_{\bm{x}}\bm{v}_b)(s,\bm{\varphi}(s,\bar{\bm{x}}))ds\Big]\cdot
  \exp\Big[-\int_0^t(\mathrm{div}_{\bm{x}}{\delta}\bm{v})(s,\bm{\varphi}(s,\bar{\bm{x}}))ds\Big].
  \end{align*}
 Since
 $$\bm{v}_b(\bm{x})=\omega_b(\varpi, z)
 \begin{bmatrix}
 -x^2 \\
 x^1 \\
 0
 \end{bmatrix},
 $$
 we see $\mathrm{div}_{\bm{x}}\bm{v}_b=0$. Therefore
 $$
 \exp\Big[-\int_0^t(\mathrm{div}_{\bm{x}}\bm{v}_b)(s,\bm{\varphi}(s,\bar{\bm{x}}))ds\Big]= I
 $$
 and \eqref{*3.1} follows.
 
 Recalling 
 $$
 \frac{\partial}{\partial t}\bm{\varphi}_b(t,\bar{\bm{x}})=\bm{v}_b(t,\bm{\varphi}_b(t,\bar{\bm{x}})),\quad
 \bm{\varphi}_b(0,\bar{\bm{x}})=\bar{\bm{x}},
 $$
 we have
 \begin{equation}
 \bm{\varphi}(t,\bar{\bm{x}})-\bm{\varphi}_b(t,\bar{\bm{x}})=
 \int_0^t({\delta}\bm{v})(s,\bm{\varphi}(s,\bar{\bm{x}}))ds +\mathscr{O}(\varepsilon^2), \label{*3.3}
 \end{equation}
 since
 \begin{align*}
 &\bm{v}(t,\bm{\varphi}(t,\bar{\bm{x}}))-\bm{v}_b(\bm{\varphi}_b(t,\bar{\bm{x}})) = \\
&= ({\delta}\bm{v})(t,\bm{\varphi}(t,\bar{\bm{x}}))+
 \Big(\bm{v}_b(\bm{\varphi}(t,\bar{\bm{x}}))-\bm{v}_b(\bm{\varphi}_b(t,\bar{\bm{x}})) \Big)  \\
 &=({\delta}\bm{v})(t,\bm{\varphi}(t,\bar{\bm{x}}))+\mathscr{O}(\varepsilon^2).
 \end{align*}
 Then we have
 \begin{align*}
 \mathrm{div}_{\bar{\bm{x}}}(\bm{\varphi}(t,\bar{\bm{x}})-\bm{\varphi}_b(t,\bar{\bm{x}}))&=
 \int_0^t\mathrm{div}_{\bar{\bm{x}}}\Big({\delta}\bm{v}(s,\bm{\varphi}(s,\bar{\bm{x}}))\Big)ds +\mathscr{O}(\varepsilon^2) \\
 &=
 \int_0^t\Big(\mathrm{div}_{\bm{x}}({\delta}\bm{v})\Big) (s,\bm{\varphi}(s,\bar{\bm{x}}))ds +\mathscr{O}(\varepsilon^2).
 \end{align*}
Thus \eqref{*3.1}
reads
$$
\rho(t,\bm{\varphi}(t,\bar{\bm{x}}))=\rho^0(\bar{\bm{x}})
\exp\Big[-\mathrm{div}_{\bar{\bm{x}}}(\bm{\varphi}(t,\bar{\bm{x}})-\bm{\varphi}_b(t,\bar{\bm{x}}))\Big]
+\mathscr{O}(\varepsilon^2),$$
or
we can claim
\begin{equation}
\rho(t,\bm{\varphi}(t,\bar{\bm{x}}))=\rho^0(\bar{\bm{x}})
-\rho^0(\bar{\bm{x}})
\mathrm{div}_{\bar{\bm{x}}}(\bm{\varphi}(t,\bar{\bm{x}})-\bm{\varphi}_b(t,\bar{\bm{x}}))+
\mathscr{O}(\varepsilon^2). \label{*3.4}
\end{equation} 

Further we have, in modulo $\mathscr{O}(\varepsilon^2)$ ,
\begin{align*}
\mathfrak{d} \rho(t,\bar{\bm{x}})&=\rho(t,\bm{\varphi}(t,\bar{\bm{x}}))-\rho_b(\bm{\varphi}(t,\bar{\bm{x}})) \\
&\equiv \rho^0(\bar{\bm{x}})-\rho_b(\bm{\varphi}(t,\bar{\bm{x}}))
-\rho^0(\bar{\bm{x}})\mathrm{div}_{\bar{\bm{x}}}(\bm{\varphi}(t,\bar{\bm{x}})-\bm{\varphi}_b(t,\bar{\bm{x}})) \\
&\equiv
\rho^0(\bar{\bm{x}})-\rho_b(\bm{\varphi}_b(t,\bar{\bm{x}}))
-\Big({\nabla} \rho_b(\bm{\varphi}_b(t,\bar{\bm{x}})-\bm{\varphi}_b(t,\bar{\bm{x}}))\Big|\bm{\varphi}(t, \bar{\bm{x}})\Big) \\
&
-\mathrm{div}_{\bar{\bm{x}}}\rho^0(\bar{\bm{x}})(\bm{\varphi}(t,\bar{\bm{x}})-\bm{\varphi}_b(t,\bar{\bm{x}})) 
+\Big({\nabla} \rho^0(\bar{\bm{x}})\Big|\bm{\varphi}(t,\bar{\bm{x}})-\bm{\varphi}_b(t,\bar{\bm{x}})\Big) \\
&= \rho^0(\bar{\bm{x}})-\rho_b(\bm{\varphi}_b(t,\bar{\bm{x}}))
-\mathrm{div}_{\bar{\bm{x}}}\Big(\rho^0(\bar{\bm{x}})(\bm{\varphi}(t,\bar{\bm{x}})-\bm{\varphi}_b(t,\bar{\bm{x}}))\Big).
\end{align*}

Since it holds that
\begin{equation}
\rho_b(\bm{\varphi}_b(t,\bar{\bm{x}}))=\rho_b(\bar{\bm{x}}), \label{*3.5}
\end{equation}
we can claim
\begin{equation}
\mathfrak{d} \rho(t,\bar{\bm{x}})=
\rho^0(\bar{\bm{x}})-\rho_b(\bar{\bm{x}})
-\mathrm{div}_{\bar{\bm{x}}}\Big(\rho_b(\bar{\bm{x}})(\bm{\varphi}(t,\bar{\bm{x}})-\bm{\varphi}_b(t,\bar{\bm{x}})\Big)
+\mathscr{O}(\varepsilon^2), \label{*3.6}
\end{equation}
since we are assuming that
$\rho_b-\rho^0=\mathscr{O}(\varepsilon)$.

Now let us introduce $\bm{u}$ by
\begin{equation}
\bm{u}(t,\bar{\bm{x}})=\bm{\varphi}(t,\bar{\bm{x}})-\bm{\varphi}_b(t,\bar{\bm{x}})+\bm{u}^0(\bar{\bm{x}}), \label{3.8}
\end{equation}
where $\bm{u}^0=\mathscr{O}(\varepsilon)$. Then we can claim
\begin{equation}
\mathfrak{d} \rho(t,\bar{\bm{x}})=-\mathrm{div}_{\bar{\bm{x}}} (\rho_b(\bar{\bm{x}})\bm{u}(t,\bar{\bm{x}}))
+\mathscr{O}(\varepsilon^2), \label{*3.9}
\end{equation}
provided that
\begin{equation}
\rho^0(\bar{\bm{x}})-\rho_b(\bar{\bm{x}})=-\mathrm{div}_{\bar{\bm{x}}}(\rho_b(\bar{\bm{x}})\bm{u}^0(\bar{\bm{x}})). \label{*3.10}
\end{equation}
 We note that
 \begin{equation}
 \frac{\partial \bm{u}}{\partial t}(t,\bar{\bm{x}})=\mathfrak{d} \bm{v}(t,\bar{\bm{x}})+\mathscr{O}(\varepsilon^2), \label{3.11}
 \end{equation}
 where we recall
 $\mathfrak{d}\bm{v}(t,\bar{\bm{x}})=({\delta}\bm{v})(t,\bm{\varphi}(t,\bar{\bm{x}}))$
 and $\bm{v}(t,\bm{\varphi}(t,\bar{\bm{x}}))-
 \bm{v}_b(\bm{\varphi}_b(t,\bar{\bm{x}}))=\mathfrak{d} \bm{v}(t,\bar{\bm{x}})+\mathscr{O}(\varepsilon^2)$. \\
 
 Let us look at $\mathfrak{d} S$. We can claim
 
 \begin{equation}
 \mathfrak{d} S(t,\bar{\bm{x}})=-\Big(\bm{u}(t,\bar{\bm{x}})\Big|{\nabla}_{\bar{\bm{x}}}S_b(\bar{\bm{x}})\Big)+\mathscr{O}(
 \varepsilon^2), \label{*3.34}
 \end{equation}
 provided that
 \begin{equation}
S^0(\bar{\bm{x}})-S_b(\bar{\bm{x}})=-(\bm{u}^0(\bar{\bm{x}})\Big|{\nabla}_{\bar{\bm{x}}}S_b(\bar{\bm{x}})). \label{S0}
\end{equation}

In fact, we see
\begin{align*}
\mathfrak{d} S(t,\bar{\bm{x}})&=S^0(\bar{\bm{x}})-S_b(\bm{\varphi}(t,\bar{\bm{x}})) \\
&=S^0(\bar{\bm{x}})-\Big[S_b(\bm{\varphi}_b(t,\bar{\bm{x}}))+\Big(\bm{\varphi}(t,\bar{\bm{x}})-
\bm{\varphi}_b(t,\bar{\bm{x}}))\Big|{\nabla}_{\bar{\bm{x}}}S_b(\bm{\varphi}_b(t,\bar{\bm{x}}))\Big)
+\mathscr{O}(\varepsilon^2)\Big] \\
&=S^0(\bar{\bm{x}})-S_b(\bar{\bm{x}})-
\Big(\bm{u}-\bm{u}^0\Big|{\nabla}_{\bar{\bm{x}}}S_b(\bar{\bm{x}})\Big)+\mathscr{O}(\varepsilon^2)
\\
&=-(\bm{u}\Big|{\nabla}_{\bar{\bm{x}}}S_b(\bar{\bm{x}})) +\mathscr{O}(\varepsilon^2) ,
\end{align*}
provided \eqref{S0}. Here we have used the indentity
\begin{equation}
S_b(\bm{\varphi}_b(t,\bar{\bm{x}}))=S_b(\bar{\bm{x}}), 
\end{equation}
which holds thanks to the axial symmetricity of $S_b$.

Moreover we can replace \eqref{S0} by
\begin{equation}
\rho^0(\bar{\bm{x}})S^0(\bar{\bm{x}})-\rho_b(\bar{\bm{x}})S_b(\bar{\bm{x}})=-\mathrm{div}_{\bar{\bm{x}}}\Big(
\rho_b(\bar{\bm{x}})S_b(\bar{\bm{x}})\bm{u}^0(\bar{\bm{x}})\Big) \label{*3.36}
\end{equation}
in modulo $\mathscr{O}(\varepsilon^2)$, since

\begin{align*}
\rho^0S^0 & -\rho_bS_b+\mathrm{div}(\rho_bS_b\bm{u}^0)=
\rho^0\Big(S^0-S_b+(\bm{u}^0|{\nabla}S_b)\Big) + \\
&+S_b\mathrm{div}\Big((\rho^0-\rho_b)\bm{u}^0\Big)+
(\rho^0-\rho_b)(S^0-S_b),
\end{align*}
provided \eqref{*3.10}.
Note that \eqref{*3.10} and \eqref{*3.36} imply
\begin{align}
&\int_{\mathfrak{R}_b}\rho^0(\bm{x})d\bm{x}=\int_{\mathfrak{R}_b}\rho_b(\bm{x})d\bm{x}, \\
&\int_{\mathfrak{R}_b}\rho^0(\bm{x})S^0(\bm{x})d\bm{x}=\int_{\mathfrak{R}_b}\rho_b(\bm{x})S_b(\bm{x})d\bm{x},
\end{align}
which mean that the total mass and the total entropy are the same for the perturbed configulation and the unperturbed background. \\

Let us look at $\mathfrak{d} P$.
We have
\begin{align*}
\mathfrak{d} P(t,\bar{\bm{x}})&=\Big(\frac{\partial  P}{\partial \rho}\Big)\Big|_{(\rho, S)=(\rho_b^L, S_b^L)(t,\bar{\bm{x}})}
\mathfrak{d} \rho(t,\bar{\bm{x}}) + \\
&+\Big(\frac{\partial P}{\partial S}\Big)\Big|_{(\rho, S)=(\rho_b^L,
S_b^L(t,\bar{\bm{x}}) }\mathfrak{d} S(t,\bar{\bm{x}})+
\mathscr{O}(\varepsilon^2).
\end{align*}

Since
\begin{align*}
&\Big(\frac{\partial  P}{\partial \rho}\Big)\Big|_{(\rho, S)=(\rho_b^L, S_b^L)(t,\bar{\bm{x}})}
=\Big(\frac{\partial  P}{\partial \rho}\Big)\Big|_{(\rho, S)=(\rho_b, S_b)(\bar{\bm{x}})} +\mathscr{O}(\varepsilon), \\
&\Big(\frac{\partial  P}{\partial S}\Big)\Big|_{(\rho, S)=(\rho_b, S_b)(\bar{\bm{x}})}
=\Big(\frac{\partial  P}{\partial S}\Big)\Big|_{(\rho, S)=(\rho_b, S_b)(\bar{\bm{x}})}+\mathscr{O}(\varepsilon),
\end{align*}
we see
\begin{align*}
\mathfrak{d} P(t,\bar{\bm{x}})&=\Big(\frac{\partial  P}{\partial \rho}\Big)\Big|_{(\rho, S)=(\rho_b, S_b)(\bar{\bm{x}})}
\mathfrak{d} \rho(t,\bar{\bm{x}}) + \\
&+\Big(\frac{\partial P}{\partial S}\Big)\Big|_{(\rho,S)=(\rho_b, S_b)(\bar{\bm{x}})} \mathfrak{d} S(t,\bar{\bm{x}})+
\mathscr{O}(\varepsilon^2).
\end{align*}

Let us look at
\begin{align}
\Phi(t,\bm{x})-\Phi_b(\bm{x})&=
-\mathsf{G}\int \frac{\rho(t,\bm{x}')-\rho_b(\bm{x}')}{\|\bm{x}-\bm{x}'\|}d\bm{x}' \nonumber \\
&=-\mathsf{G}\int_{\mathfrak{R}_b}
\frac{\mathfrak{d}\rho(t,\bar{\bm{x}}')}{\|\bar{\bm{x}}-\bar{\bm{x}}'\|}d\bar{\bm{x}}' +\mathscr{O}(\varepsilon^2) \label{1.51}
\end{align}

In fact, applying the Gronwall argument to the identity
$$
\bm{\varphi}(t,\bar{\bm{x}})-\bm{\varphi}(t,\bar{\bm{x}}')=
\bar{\bm{x}}-\bar{\bm{x}}'+
\int_0^t\Big[
\bm{v}(s,\bm{\varphi}(s,\bar{\bm{x}}))-\bm{v}(s,\bm{\varphi}(s,\bar{\bm{x}}'))\Big]ds,
$$
we can derive the estimate
$$
\frac{1}{1+C_1|\varepsilon| T}
\|\bm{\varphi}(t,\bar{\bm{x}})-\bm{\varphi}(t,\bar{\bm{x}}')\|
\leq \|\bar{\bm{x}}-\bar{\bm{x}}'\| \leq
\frac{1}{1-C_1|\varepsilon| T}
\|\bm{\varphi}(t,\bar{\bm{x}})-\bm{\varphi}(t,\bar{\bm{x}}')\|
$$
for $t \in [0,T], \bar{\bm{x}}, \bar{\bm{x}}' \in \mathfrak{R}_b$, provided that $C_1|\varepsilon| T <1$,
where
$$
C_1=\frac{1}{|\varepsilon|}\sup
\Big\{\|D_{\bm{x}}\bm{v}(t,\bm{x})\|\Big| 0\leq t\leq T, \bm{x} \in \mathfrak{D}_0\Big\}
$$
so that
$$
\Big|
\int_0^t\Big[
\bm{v}(s,\bm{\varphi}(s,\bar{\bm{x}}))-\bm{v}(s,\bm{\varphi}(s,\bar{\bm{x}}'))\Big]ds
\Big|\leq
C_1|\varepsilon|
\int_0^t\|\bm{\varphi}(s, \bar{\bm{x}})-\bm{\varphi}(s,\bar{\bm{x}}')\|ds,
$$
and we have
$$
d\bm{x}'=\Big|\mathrm{det}D_{\bar{\bm{x}}}\bm{x}(t,\bar{\bm{x}}')\Big|d\bar{\bm{x}}' =(1+\mathscr{O}(\varepsilon))d\bar{\bm{x}}'.
$$
Hence
$$
\frac{\rho(t,\bm{x}')-\rho_b(\bm{x}')}{\|\bm{x}-\bm{x}'\|}d\bm{x}'=
\frac{\mathfrak{d}\rho(t,\bar{\bm{x}}')}{\|\bar{\bm{x}}-\bar{\bm{x}}'\|}d\bar{\bm{x}}' +\mathscr{O}(\varepsilon^2)
$$
For $\bm{x}=\bm{\varphi}(1,\bar{\bm{x}}), \bm{x}'=\bm{\varphi}
(t,\bar{\bm{x}}')$, therefore 
\eqref{1.51} holds.

Combining \eqref{*3.9} and \eqref{1.51}, we see
\begin{equation}
\Phi(t,\bm{x})-\Phi_b(\bm{x})=
\mathsf{G}\int_{\mathfrak{R}_b}
\frac{\mathrm{div}_{\bar{\bm{x}}}(\rho_b\bm{u}(t,\bar{\bm{x}}'))}{\|\bar{\bm{x}}-\bar{\bm{x}}'\|}d\bar{\bm{x}}' +\mathscr{O}(\varepsilon^2).
\end{equation}\\

We note
\begin{align}
\frac{1}{\rho}{\nabla}_{\bm{x}}P-\frac{1}{\rho_b}{\nabla}_{\bm{x}}P_b 
&=\frac{1}{\rho_b}{\nabla}_{\bm{x}}(P-P_b)
-\frac{\rho-\rho_b}{\rho_b^2}{\nabla}_{\bm{x}}P_b
\nonumber \\
&=\frac{1}{\rho_b^L}{\nabla}_{\bar{\bm{x}}}(\mathfrak{d} P)
-\frac{\mathfrak{d} \rho}{(\rho_b^L)^2}{\nabla}_{\bar{\bm{x}}}P_b^L 
+\mathscr{O}(\varepsilon^2)
\end{align}
with $\bm{x}=\bm{\varphi}(t,\bar{\bm{x}})$. \\

Summing up, we get the linearized approximation of \eqref{*1.1}:
\begin{equation}
\frac{\partial\mathfrak{d}\bm{v}}{\partial t}+B{\mathfrak{d}\bm{v}}+
\mathcal{L}\bm{u}=0, \label{LE.1}
\end{equation}
with
\begin{align}
&\mathcal{L}\bm{u}=\frac{1}{\rho_b}{\nabla}_{\bar{\bm{x}}} (\mathfrak{d} P)
-\frac{{\nabla}_{\bar{\bm{x}}} P_b}{\rho_b^2}\mathfrak{d} \rho +
{\nabla}_{\bar{\bm{x}}}\Big(4\pi\mathsf{G}\mathcal{K}[-\mathfrak{d} \rho]\Big), \label{LE.1L}\\
&\mathfrak{d}\rho =-\mathrm{div}_{\bar{\bm{x}}}(\rho_b\bm{u}), \\
&\mathfrak{d} S=-(\bm{u}|\nabla_{\bar{\bm{x}}} S_b(\bar{\bm{x}})), \\
&\mathfrak{d} P(t,\bar{\bm{x}})=\Big(\frac{\partial  P}{\partial \rho}\Big)\Big|_{(\rho, S)=(\rho_b, S_b)(\bar{\bm{x}})}
\mathfrak{d} \rho + \\
&+\Big(\frac{\partial P}{\partial S}\Big)\Big|_{(\rho,S)=(\rho_b, S_b)(\bar{\bm{x}})} \mathfrak{d} S,
\end{align}

where we read $\rho_b=\rho_b(\bar{\bm{x}}), P_b=P_b(\bar{\bm{x}}),
S_b=S_b(\bar{\bm{x}})$.\\

Since $\mathfrak{d}\bm{v}=\partial\bm{u}/\partial t +\mathscr{O}(\varepsilon^2)$ ( \eqref{3.11} ), \eqref{LE.1}
reads
\begin{equation}
\frac{\partial^2\bm{u}}{\partial t^2}+B\frac{\partial \bm{u}}{\partial t}+
\mathcal{L}\bm{u}=0.
 \label{LLEq}
\end{equation}

This equation \eqref{LLEq} is the linearized equation for perturbations described by the Lagrangian coordinate to be analyzed.\\

Let us recall 
that we are assuming \eqref{3.1}, \eqref{*3.10} and \eqref{S0}.\\

\begin{Remark} \label{Rem.X}
Strictly speaking, we should introduce $\lceil \mathfrak{d}\rho \rfloor, \lceil \mathfrak{d} S \rfloor,
\lceil \mathfrak{d} P \rfloor, \lceil \mathfrak{d} \Phi \rfloor $ by defining
\begin{subequations}
\begin{align}
 \lceil \mathfrak{d}\rho \rfloor(t,\bar{\bm{x}})&=-\mathrm{div}_{\bar{\bm{x}}}(\rho_b(\bar{\bm{x}})\bm{u}(t,\bar{\bm{x}})), \\
\lceil \mathfrak{d} S \rfloor(t,\bar{\bm{x}})&=-(\bm{u} (t,\bar{\bm{x}})|{\nabla}_{\bar{\bm{x}}}S_b(\bar{\bm{x}})), \\
\lceil \mathfrak{d} P \rfloor(t,\bar{\bm{x}})&=
\Big(\frac{\partial  P}{\partial \rho}\Big)_b(\bar{\bm{x}})
\lceil \mathfrak{d} \rho \rfloor (t,\bar{\bm{x}})+\Big(\frac{\partial P}{\partial S}\Big)_b(\bar{\bm{x}}) \lceil \mathfrak{d} S \rfloor (t,\bar{\bm{x}}), \\
\lceil \mathfrak{d} \Phi \rfloor(t,\bar{\bm{x}})&=4\pi\mathsf{G}\mathcal{K}[-
\lceil \mathfrak{d} \rho \rfloor](t,\bar{\bm{x}}).
\end{align}
\end{subequations}

Then \eqref{LE.1L} should read
\begin{equation}
\bm{L}\bm{u}= \frac{1}{\rho_b}{\nabla}_{\bar{\bm{x}}}\lceil \mathfrak{d} P \rfloor
-\frac{{\nabla}_{\bar{\bm{x}}}P_b}{\rho_b^2}\lceil \mathfrak{d}\rho \rfloor 
+{\nabla}_{\bar{\bm{x}}}\lceil \mathfrak{d} \Phi \rfloor,
\end{equation}
while
\begin{align*}
\mathfrak{d} \rho &= \lceil \mathfrak{d}\rho \rfloor +\mathscr{O}(\varepsilon^2) ,\\
\mathfrak{d} S &=\lceil \mathfrak{d} S \rfloor +\mathscr{O}(\varepsilon^2) ,\\
\mathfrak{d} P &=\lceil \mathfrak{d} P \rfloor+\mathscr{O}(\varepsilon^2) ,\\
\mathfrak{d} \Phi &=
\lceil \mathfrak{d} \Phi \rfloor+\mathscr{O}(\varepsilon^2) .
\end{align*}
\end{Remark}

\begin{Remark}  

 The derivation of the linearized approximation of the equations in Lagrangian co-ordinate system can be found
 \cite[Sect. 56]{LedouxW}, \cite[pp. 139-140.]{Batchelor}, \cite[p.11, (A)]{Bjerknes}, \cite{LyndenBO}, \cite[p.500, (1)]{Lebovitz} and so on. But there was considered only the case of $\rho^0=\rho_b, S^0=S_b$. Moreover, the definition of Lagrange displacement and ELASO, the linearized equation of adiabatic perturbations, derived by them and traditionally used by astrophysicists are slightly different from those formulated here, when $\bm{v}_b=\mathscr{O}(\varepsilon)$ is not supposed. For the details see Section 7.
 
 \end{Remark}

\section{Basic existence theorem for {\bf ELASO} --Application of the Hille-Yosida theory}

We are going to discuss on the existence of solutions to the 
equation of {\bf ELASO}  \eqref{Eqxi}.
We formulate the initial value problem to be considered as:\\

\begin{subequations}
\begin{align}
&\frac{\partial ^2\bm{u}}{\partial  t^2} +B\frac{\partial  \bm{u}}{\partial  t}+ \mathcal{L}\bm{u} = \bm{f}(t,\bm{x}) \quad\mbox{on}\quad [0,+\infty[\times \mathfrak{R}_b, \label{IVPa}\\
&\bm{u}|_{t=0}=\bm{u}^0, \quad \frac{\partial \bm{u}}{\partial  t}\Big|_{t=0}=\bm{v}^0 \quad\mbox{on}\quad 
\mathfrak{R}_b, \label{IVPb}
\end{align}
\end{subequations}
where $\bm{u}^0, \bm{v}^0, \bm{f}$ are given data.

Let us rewrite
$$\mathcal{L}_0\bm{u}:=\frac{1}{\rho_b}\nabla \mathfrak{d} P-\frac{\nabla P_b}{\rho_b^2}\mathfrak{d}\rho,
$$
using $\Gamma_b, \mathfrak{a}_b$ defined in Definition \ref{Def.G}. Since
$$
\mathfrak{d} P=\Big(\frac{\Gamma P}{\rho}\Big)_b\mathfrak{d}\rho+
(\Gamma P)_b(\bm{u}|\mathfrak{a}_b),
$$
we have
\begin{align*}
\mathcal{L}_0\bm{u}&=\nabla\Big[-\Big(\frac{\Gamma P}{\rho^2}\Big)_b\mathrm{div}(\rho_b\bm{u})+
\Big(\frac{\Gamma P}{\rho}\Big)_b(\bm{u}|\mathfrak{a}_b)\Big] + \\
&+\Big(\frac{\Gamma P}{\rho^2}\Big)_b
\Big[-\mathfrak{a}_b\mathrm{div}(\rho_b\bm{u})
+(\bm{u}|\mathfrak{a}_b)\nabla\rho_b\Big].
\end{align*}
Consequently we can write
\begin{equation}
\mathcal{L}\bm{u}=\mathcal{L}_0\bm{u}+4\pi\mathsf{G}\mathcal{L}_1
\end{equation}
with
\begin{align}
\mathcal{L}_0\bm{u}&=
\nabla\Big[-\Big(\frac{\Gamma P}{\rho^2}\Big)_b\mathrm{div}(\rho_b\bm{u})+
\Big(\frac{\Gamma P}{\rho}\Big)_b(\bm{u}|\mathfrak{a}_b)\Big] +
\nonumber \\
&+\Big(\frac{\Gamma P}{\rho^2}\Big)_b
\Big[-\mathfrak{a}_b\mathrm{div}(\rho_b\bm{u})
+(\bm{u}|\mathfrak{a}_b)\nabla\rho_b\Big] \\
\mathcal{L}_1\bm{u}&={\nabla}\mathcal{K}(\mathrm{div}(\rho_b \bm{u})).
\end{align}\\

\subsection{Case of barotropic backgrounds}
Now we suppose \\

{\bf  Assumption} : {\it The equation of state enjoys {\bf Assumption} \ref{Ass.0} and the background $(\rho_b, S_b, \bm{v}_b)$ is a barotropic admissible stationary solution.}\\

Hereafter we  denote 
$
\mathfrak{R}_b
$, 
$ \rho_b(\bm{x}),
P_b(\bm{x}), S_b(\bm{x})
$, $\mathfrak{a}_b, \mathscr{A}_b$
by 
$
\mathfrak{R}
$,  
$\rho(\bm{x}), P(\bm{x}), S(\bm{x})
$, $\mathfrak{a}, \mathscr{A}$
for the sake of simplicity of symbols.\\

Let us consider the integro-differential operator $\mathcal{L}$ in the Hilbert space
$\mathfrak{H}$:
\begin{Definition}
We denote
$\mathfrak{H}=L^2(\mathfrak{R},\rho d\bm{x}; \mathbb{C}^3)$ endowed with the inner product
\begin{equation}
(\bm{u}_1|\bm{u}_2)_{\mathfrak{H}}=\int_{\mathfrak{R}}(\bm{u}_1(\bm{x})|\bm{u}_2(\bm{x}))\rho(\bm{x})d\bm{x}.
\end{equation}
\end{Definition}

Of course
$$(\bm{u}_1(\bm{x})|\bm{u}_2(\bm{x})):=\sum_k u_1^k(\bm{x})u_2^k(\bm{x})^*
\quad\mbox{for}\quad 
\bm{u}_{\mu}(\bm{x})=
\begin{bmatrix}
u_{\mu}^1(\bm{x}) \\
\\
u_{\mu}^2(\bm{x})\\
\\
 u_{\mu}^ 3(\bm{x})
\end{bmatrix}, \mu=1,2.$$\\
 Here and hereafter $Z^*$ denotes the complex conjugate $X-\mathrm{i}Y$of $Z=X+\mathrm{i}Y$, while $\mathrm{i}$ stands for the imaginary unit, $\sqrt{-1}$.\\
 
 First we observe $\mathcal{L}$ restricted on $C_0^{\infty}(\mathfrak{R};  \mathbb{C}^3)$. 
 
 We have the following formula by integration by parts:
\begin{align*}
(\mathcal{L}_0\bm{u}_1|\bm{u}_2)&=\int \frac{\Gamma P}{\rho^2}\mathrm{div}(\rho \bm{u}_1)\mathrm{div}(\rho \bm{u}_2)^* + \\
&+\int \frac{\Gamma P}{\rho}\Big(-(\bm{u}_1|\mathfrak{a})\mathrm{div}(\rho\bm{u}_2)^* 
-\mathrm{div}(\rho \bm{u}_1)(\bm{u}_2|\mathfrak{a})^*\Big) + \\
&+\int \frac{\Gamma P}{\rho}(\bm{u}_1|\mathfrak{a})({\nabla}\rho|\bm{u}_2) \\
&=\int \frac{\Gamma P}{\rho^2}\mathrm{div}(\rho \bm{u}_1)\mathrm{div}(\rho \bm{u}_2)^* + \\
&+\int \frac{\Gamma P \mathscr{A}}{\rho}
\Big((\bm{u}_1|\bm{n})\mathrm{div}(\rho \bm{u}_2)^*
+(\mathrm{div}(\rho \bm{u}_1)(\bm{u}_2|\bm{n})^*\Big) \\
&-\int \frac{\Gamma P\mathscr{A}}{\rho}\|{\nabla}\rho\|^2(\bm{u}_1|\bm{n})(\bm{u}_2|\bm{n})^*,
\end{align*}
where
we recall
\begin{equation}
\mathfrak{a}=-\mathscr{A}\bm{n}, \quad \bm{n}=\frac{{\nabla}\rho}{\|{\nabla}\rho\|},
\end{equation}
and
\begin{equation}
(\mathcal{L}_1\bm{u}_1|\bm{u}_2)=-
\int \mathcal{K}[\mathrm{div}(\rho \bm{u}_1)]\mathrm{div}(\rho \bm{u}_2)^*.
\end{equation}

Hence
\begin{align}
(\mathcal{L}\bm{u}_1|\bm{u}_2)&=
\int \frac{\Gamma P}{\rho^2}\mathrm{div}(\rho \bm{u}_1)\mathrm{div}(\rho \bm{u}_2)^*  +2\mathfrak{Re}\Big[\int \frac{\Gamma P\mathscr{A}}{\rho}(\bm{u}_1|\bm{n})\mathrm{div}(\rho \bm{u}_2)^*\Big] \nonumber \\
&-\int \frac{\Gamma P\mathscr{A}}{\rho}\|{\nabla}\rho\|^2(\bm{u}_1|\bm{n})(\bm{u}_2|\bm{n})^* + \nonumber \\
&-4\pi\mathsf{G}
\int \mathcal{K}[\mathrm{div}(\rho \bm{u}_1)]\mathrm{div}(\rho \bm{u}_2)^*. \label{*4.15}
\end{align}

Then it is clear that

\begin{equation}
(\mathcal{L}\bm{u}_1|\bm{u}_2)=(\bm{u}_1|\mathcal{L}\bm{u}_2)
\quad\mbox{for}\quad \forall \bm{u}_1,\bm{u}_2 \in C_0^{\infty}(\mathfrak{R}),
\end{equation} 
that is, $\mathcal{L}$ restricted on $C_0^{\infty}(\mathfrak{R})$ is a symmetric operator.\\

Moreover we have
\begin{align}
(\mathcal{L}_0\bm{u}|\bm{u})&=\int \frac{\Gamma P}{\rho^2}|\mathrm{div}(\rho \bm{u})|^2d\bm{x} + 2\mathfrak{Re}\Big[
\int \frac{\Gamma P\mathscr{A}}{\rho}(u|\bm{n})\mathrm{div}(\rho \bm{u})^*\Big] + \nonumber \\
&-\int \frac{\Gamma P \mathscr{A} }{\rho}\|{\nabla}\rho\|^2
|(\bm{u}|\bm{n})|^2.
\end{align}

Since $\mathscr{A} \in C^{\alpha}(\mathfrak{R} \cup \partial\mathfrak{R}))$, we have
\begin{equation}
\kappa_1:=\sup_{\mathfrak{R}}\Big\{\   |\mathscr{A}|^2\frac{\Gamma P}{\rho}
\  \Big\} <\infty.
\end{equation}
Therefore
$$
\Big|2\int \frac{\Gamma P \mathscr{A}}{\rho}(\bm{u}|\bm{n})\mathrm{div}(\rho \bm{u})^* \Big|\leq
\frac{1}{2}\int\frac{\Gamma P}{\rho^2}|\mathrm{div}(\rho \bm{u})|^2 d\bm{x}
+2\kappa_1\|\bm{u}\|_{\mathfrak{H}}^2.
$$

On the other hand, we have

\begin{equation}
\kappa_2:=\sup_{\mathfrak{R}}\Big\{\frac{\Gamma P |\mathscr{A}|}{\rho}\|{\nabla}\rho\|^2 \Big\}
<\infty.
\end{equation}
Then
$$
\Big| \int \frac{\Gamma P \mathscr{A} }{\rho}\|{\nabla}\rho\|^2
 |(\bm{u}|\bm{n})|^2 \Big|
\leq \kappa_2 \|\bm{u}\|_{\mathfrak{H}}^2.
$$

Thus we have

\begin{Proposition}\label{Prop.8v1}
It holds that
\begin{align} 
&\frac{1}{2}\int\frac{\Gamma P}{\rho^2}|\mathrm{div}(\rho \bm{u})|^2 d\bm{x}
-(2\kappa_1+\kappa_2)\|\bm{u}\|_{\mathfrak{H}}^2 \leq
(\mathcal{L}_0\bm{u}|\bm{u}) \leq \nonumber \\
&\leq \frac{3}{2}\int\frac{\Gamma P}{\rho^2}|\mathrm{div}(\rho \bm{u})|^2 d\bm{x}
+(2\kappa_1+\kappa_2)\|\bm{u}\|_{\mathfrak{H}}^2.
\end{align}
\end{Proposition}

As for $\mathcal{L}_1$, on the other hand, we have
\begin{Proposition}\label{Prop.8v2}
It holds that
 \begin{equation}
 -\rho_O \|\bm{u}\|_{\mathfrak{H}}^2\leq
 (\mathcal{L}_1\bm{u}|\bm{u})_{\mathfrak{H}} \leq 0.
 \end{equation}
 \end{Proposition}

Proof.   Look at
 \begin{equation}
  \Psi:=-\mathcal{K}[g], \quad \bm{C}:=\rho \bm{u} -{\nabla}\Psi.
 \end{equation}
 Since $\triangle \Psi=g$, we have
 $$
 \mathrm{div}\bm{C}=0.
 $$
 Keeping in mind that $\Psi=O(\frac{1}{r}), {\nabla}\Psi, \bm{C}=O(\frac{1}{r^2})$
 as $ r\rightarrow +\infty$, we derive from this that
 \begin{equation}
 \int_{\mathbb{R}^3}({\nabla}\Psi|\bm{C})d\bm{x}=0  \label{4.12}
 \end{equation} 
Then we see
\begin{align*}
\int_{\mathfrak{R}}\mathcal{K}[g]g^* d\bm{x} &=\int_{\mathfrak{R}}({\nabla}\Psi|\rho\bm{u}) \\
&=\int_{\mathbb{R}^3} ({\nabla}\Psi|\bm{C}+{\nabla}\Psi)d\bm{x}  \\
&=\int_{\mathbb{R}^3} \|{\nabla}\Psi\|^2d\bm{x} \quad \mbox{(by \eqref{4.12})}\\
&=\int_{\mathfrak{R}} \|\rho\bm{u}\|^2d\bm{x} -\int_{\mathbb{R}^3} \|\bm{C}\|^2d\bm{x} \quad \mbox{(again by \eqref{4.12})}\\
&\leq \int_{\mathfrak{R}} \|\rho\bm{u}\|^2d\bm{x} \\
&\leq \rho_O\int_{\mathfrak{R}} \|u\|^2\rho d\bm{x} =\rho_O\|\bm{u}\|_{\mathfrak{H}}^2.
\end{align*}
and
$$\int_{\mathfrak{R}}\mathcal{K}[g]g^* d\bm{x} =\int_{\mathbb{R}^3} \|{\nabla}\Psi\|^2d\bm{x} \geq 0.
$$ 
$\square$

This proof is due to Juhi Jang. See \cite[Proposition 2]{JJTM2020}.\\

Summing up, we have
\begin{align} 
&\frac{1}{2}\int\frac{\Gamma P}{\rho^2}|\mathrm{div}(\rho \bm{u})|^2 d\bm{x}
-(2\kappa_1+\kappa_2+4\pi\mathsf{G}\rho_O)\|\bm{u}\|_{\mathfrak{H}}^2 \leq
(\mathcal{L}\bm{u}|\bm{u}) \leq \nonumber \\
&\leq \frac{3}{2}\int\frac{\Gamma P}{\rho^2}|\mathrm{div}(\rho \bm{u})|^2 d\bm{x}
+(2\kappa_1+\kappa_2)\|\bm{u}\|_{\mathfrak{H}}^2.
\end{align}\\

Therefore the operator $\mathcal{L}$ restricted on $C_0^{\infty}(\mathfrak{R})$ is 
symmetric, and is bounded from below as an operator in $\mathfrak{H}$. Applying the Stone-Friedrichs theorem, (e.g., \cite[Chapter VI, Section 2.3]{Kato}),
we can claim:

\begin{Theorem} \label{Th.*2}
 $\mathcal{L}\restriction C_0^{\infty}(\mathfrak{R})$ admits the Friedrichs extension
$\bm{L}$, which is a self-adjoint operator
in $\mathfrak{H}$, whose domain is
\begin{equation}
\mathsf{D}(\bm{L})=\{ \bm{u} \in \mathfrak{G}_0\  \Big|\  \mathcal{L}\bm{u} \in \mathfrak{H} \}.
\end{equation}
\end{Theorem}

Here we use
\begin{Definition}
$\mathfrak{G}$ stands for the Hilbert space 
of all $\bm{u} \in \mathfrak{H}$ such that
$\displaystyle \int_{\mathfrak{R}} \frac{\Gamma P}{\rho^2}|\mathrm{div}(\rho \bm{u})|^2
d\bm{x} <\infty $ endowed with the inner product
\begin{equation}
(\bm{u}_1|\bm{u}_2)_{\mathfrak{G}}=
(\bm{u}_1|\bm{u}_2)_{\mathfrak{H}}+
\int_{\mathfrak{R}} \frac{\Gamma P}{\rho^2}(\mathrm{div}(\rho \bm{u}_1)(\mathrm{div}(\rho\bm{u}_2)^*
d\bm{x}, 
\end{equation}
and $\mathfrak{G}_0$ is the closure of $C_0^{\infty}(\mathfrak{R})$ in $\mathfrak{G}$. 
\end{Definition}

In fact, we define the quadratic form $Q_0$ by
\begin{align}
Q_0[\bm{u}]=Q_0(\bm{u},\bm{u})&= \int \frac{\Gamma P}{\rho^2}|\mathrm{div}(\rho \bm{u})|^2 d\bm{x} +2\mathfrak{Re}\Big[
\int \frac{\Gamma P \mathscr{A} }{\rho}(\bm{u}|\bm{n})\mathrm{div}(\rho \bm{u})^*d\bm{x} \Big] + \nonumber \\
&-\int \frac{\Gamma P \mathscr{A} }{\rho}|(\bm{u}|\bm{n})|^2\|{\nabla}\rho\|d\bm{x}
-4\pi\mathsf{G}\int \mathcal{K}[\mathrm{div}(\rho\bm{u})]\mathrm{div}(\rho \bm{u})^*d\bm{x}.
\end{align}
Then
it holds that
\begin{equation}
(\mathcal{L}\bm{u}_1|\bm{u}_2)_{\mathfrak{H}}=Q_0(\bm{u}_1,\bm{u}_2)
\end{equation}
for $\forall \bm{u}_1 , \bm{u}_2 \in C_0^{\infty}(\mathfrak{R})$. 
We define
\begin{equation}
Q_a(\bm{u}_1, \bm{u}_2)=Q_0(\bm{u}_1, \bm{u}_2)+a(\bm{u}_1 | \bm{u}_2)_{\mathfrak{H}}
\end{equation}
with
\begin{equation}
a := 2\kappa_1+\kappa_2+4\pi\mathsf{G}\rho_O. \label{7.4.26}
\end{equation}

Since
\begin{equation} 
 \frac{1}{2}\int\frac{\Gamma P}{\rho^2}|\mathrm{div}(\rho \bm{u})|^2 d\bm{x}
\leq
Q_a[\bm{u}]
\leq
\frac{3}{2}\int\frac{\Gamma P}{\rho^2}|\mathrm{div}(\rho \bm{u})|^2 d\bm{x}
+(2\kappa_1+\kappa_2+a)\|\bm{u}\|_{\mathfrak{H}}^2,
\end{equation}
we have
\begin{equation}
0 \leq Q_a[\bm{u}] \quad\mbox{for}\quad \bm{u} \in \mathfrak{G}. \label{4.26}
\end{equation}

Summing up, we define the alternative inner product $(\cdot|\cdot)_{\breve{\mathfrak{G}}}$ of $\mathfrak{G}$ by
\begin{align}
(\bm{u}_1|\bm{u}_2)_{\breve{\mathfrak{G}}}& =Q_a(\bm{u}_1|\bm{u}_2)+(\bm{u}_1|\bm{u}_2)_{\mathfrak{H}} \nonumber \\
&=Q_0(\bm{u}_1|\bm{u}_2)+(a+1)(\bm{u}_1|\bm{u}_2)_{\mathfrak{H}}
\end{align}
with the norm
\begin{equation}
\|\bm{u}\|_{\breve{\mathfrak{G}}}=\sqrt{ (\bm{u}|\bm{u})_{\breve{\mathfrak{G}}} }.
\end{equation}
Then $\|\cdot\|_{\breve{\mathfrak{G}}}$ is equivalent to $\|\cdot \|_{\mathfrak{G}}$, that is, there is a constant $C$
such that
\begin{equation}
\frac{1}{C}\|\bm{u}\|_{\mathfrak{G}} \leq \|\bm{u}\|_{\breve{\mathfrak{G}}} \leq C \|\bm{u}\|_{\mathfrak{G}} \quad \mbox{for}\quad \forall 
\bm{u} \in \mathfrak{G}.
\end{equation}
It hold that
\begin{equation}
\|\bm{u}\|_{\breve{\mathfrak{G}}} \geq \|\bm{u}\|_{\mathfrak{H}} \quad\mbox{for}\quad \forall \bm{u} \in \mathfrak{G},
\end{equation}
and
\begin{equation}
(\bm{u}|\bm{v})_{\breve{\mathfrak{G}}}=
((\bm{L}+a+1)\bm{u}|\bm{v})_{\mathfrak{H}} \quad\mbox{for}\quad \forall \bm{u} \in \mathsf{D}(\bm{L}), \bm{v} \in \mathfrak{G}_0.
\end{equation}\\

In order to avoid confusions, hereafter we use 

\begin{Notation}
We denote the Hilbert space of all funstions in $\mathfrak{G}$ $ [(\mathfrak{G}_0 )]$ endowed with the inner product $(\cdot |\cdot)_{\breve{\mathfrak{G}}}$ and norm $\|\cdot\|_{\breve{\mathfrak{G}}}$ by $\breve{\mathfrak{G}}$ $[(\breve{\mathfrak{G}_0})]$ respectively.
\end{Notation}

In this situation we can apply the Hille-Yosida theory to the first order evoltion equation {\bf ELASO\ddag} in $[0,+\infty[\times \mathfrak{E}$:
\begin{align}
& \frac{dU}{dt}+\bm{A}U=\bm{0}, \\
& \bm{A}=
\begin{bmatrix}
O & -I \\
\\
\bm{L} & \bm{B}
\end{bmatrix},
\quad
\mathsf{D}(\bm{A})=\mathsf{D}(\bm{L})\times \breve{\mathfrak{G}_0}, \\
& \mathfrak{E}=\breve{\mathfrak{G}_0} \times \mathfrak{H}
\end{align}
which is equivalent to {\bf ELASO} through the interpretation
\begin{equation}
U=\begin{bmatrix}
\bm{u} \\
\\
\frac{d\bm{u}}{dt}
\end{bmatrix}.
\end{equation}
This application iss described in {\bf Appendix A}. The conclusion is:\\

\begin{Theorem} \label{Th.*3}

Suppose $\bm{u}^0 \in \mathsf{D}(\bm{L})$, 
$\bm{v}^0 \in \mathfrak{G}_0$
and $ \bm{f} \in C([0,+\infty[; \mathfrak{H})$. Then the initial value problem 
\eqref{IVPa} \eqref{IVPb} admits a unique solution
$\bm{u}=\bm{u}(t,\bm{x})$ in
$$ C^2([0,+\infty[, \mathfrak{H})\cap
C^1([0,+\infty[,\mathfrak{G}_0)
\cap C([0,+\infty[, \mathsf{D}(\bm{L}))$$
and the energy
\begin{align}
E(t, \bm{u})&:=\||\bm{u}\|_{\breve{\mathfrak{G}}}^2
+\Big\|\frac{\partial\bm{u}}{\partial  t}\Big\|_{\mathfrak{H}}^2 \nonumber \\
&=
\|\bm{u}||_{\mathfrak{H}}^2+Q_a[\bm{u}]+\Big\|\frac{\partial\bm{u}}{\partial  t}\Big\|_{\mathfrak{H}}^2
=(1+a)\|\bm{u}\|_{\mathfrak{H}}^2+
Q_0[\bm{u}]+\Big\|\frac{\partial\bm{u}}{\partial  t}\Big\|_{\mathfrak{H}}^2
\end{align}
enjoys the estimate
\begin{equation}
\sqrt{E(t,\bm{u})}\leq e^{\kappa t}\sqrt{E(0, \bm{u})}+\int_0^te^{\kappa(t-s)}
\|\bm{f}(s)\|_{\mathfrak{H}}ds
\end{equation}
for 
$\kappa=1+a $. 
\end{Theorem}

Let us claim the following

\begin{Theorem} \label{Th.*4}
$0$ is an eigenvalue of $\bm{L}$ and $\mathrm{dim}\mathsf{N}(\bm{L}) = \infty$.
\end{Theorem}

Proof. 1) Suppose $\mathscr{A}\not= 0$, that is, ${\nabla}S_b(\bm{x})\not= 0$ for $\exists \bm{x} \in \mathfrak{R}$. Then
\begin{equation}
\bm{u}^0(\bm{x})=\frac{1}{\rho_b(\bm{x})}{\nabla}S_b(\bm{x})\times {\nabla}f(\bm{x}),
\end{equation}
where $f \in C_0^{\infty}(\mathfrak{R}; \mathbb{R})$ is arbitrary, enjoys
\begin{align*}
&\lceil \mathfrak{d} \rho; \bm{u}^0 \rfloor=-\mathrm{div}(\rho_b\bm{u}^0)=0, \\
&\lceil \mathfrak{d} S; \bm{u}^0 \rfloor=-(\bm{u}^0|{\nabla}S_b)=0, \\
& \lceil \mathfrak{d} P; \bm{u}^0 \rfloor=-\Big(\frac{\Gamma P}{\rho}\Big)_b\lceil \mathfrak{d} \rho; \bm{u}^0 \rfloor+
\Big(\frac{\partial P}{\partial S}\Big)_b \lceil \mathfrak{d} S; \bm{u}^0) \rfloor=0, \\
&\lceil \mathfrak{d}\Phi; \bm{u}^0 \rfloor=-4\pi\mathsf{G}\mathcal{K}[\lceil \mathfrak{d}\rho; \bm{u}^0\rfloor ]=0,  
\end{align*}
therefore $\bm{L}\bm{u}^0=\bf{0}$. (Recall the notations of Remark \ref{Rem.X}.) If ${\nabla}f (\bm{x}_1)\not=0$
and ${\nabla}S_b(\bm{x}_1)\not= \bf{0}$ at $\exists \bm{x}_1 \in \mathfrak{R}$,  then $\bm{u}^0
\not=0$ and $\bm{u}^0$ is an eigenvector of the eigenvalue $0$.

2) Suppose $\mathscr{A}=0$, that is, ${\nabla}S_b(\bm{x})=0$ for $\forall \bm{x} \in \mathfrak{R}$,
or, the background is isentropic. Then, 
for any vector field $\bm{f} \in C_0^{\infty}(\mathfrak{R}; \mathbb{R}^3)$, the vector field 
$\bm{u}^0$ defined by
\begin{equation}
\bm{u}^0(\bm{x})=\frac{1}{\rho_b(\bm{x})}\mathrm{rot}\bm{f}(\bm{x})
\end{equation}
enjoys
$ \mathfrak{d} \rho=0, \mathfrak{d} S=0, \mathfrak{d} P=0, \mathfrak{d} \Phi=0$, so that $\bm{L}\bm{u}^0=\bf{0}$.
If $\mathrm{rot}\bm{f}\not=\bf{0}$ somewhere, then $\bm{u}^0 \not=\bf{0}$ and it is an eigenvector. $\square$\\

Roughly or formally speaking, the boundary condition for $\bm{u} \in \mathsf{D}(\bm{L})$ is:

$$
(\rho_b\bm{u}|\bm{N})=0 \quad\mbox{on}\quad \partial\mathfrak{R}_b,
$$
where $\bm{N}$ is the normal vector at the boundary point, namely,
$$ \Big(\rho_b(\bm{x})\bm{u}(\bm{x})\Big| \bm{n}(\bm{x})\Big) \rightarrow 0\quad\mbox{as}\quad
\bm{x} \in \mathfrak{R}_b \rightarrow \bm{z} \in \partial\mathfrak{R}_b,$$
where $\bm{n}(\bm{x})=\nabla\rho_b(\bm{x})/\|\nabla\rho_b(\bm{x})\|$
which is in $C^{1,\alpha}(\mathfrak{R}_b\cup\partial\mathfrak{R}_b)$ and enjoys $$\bm{n}(\bm{z})=
\lim_{\bm{x}\in \mathfrak{R}_b\rightarrow\bm{z}}\bm{n}(\bm{x})=\bm{N}(\bm{z})$$
for $\bm{z}\in \partial\mathfrak{R}_b$. Note that, if $\bm{u}$ is bounded near the boundary point, this condition is satisfied, since $\rho_b(\bm{z})=0$,
and that the boundary condition does not require $(\bm{N}|\bm{u})=0$
on the boundary.

More precisely speaking, the boundary condition means as follows. Let $\bm{u} \in \mathfrak{G}_0$. Put
$\hat{\bm{u}}:=\rho_b\bm{u}$. Then $\hat{\bm{u}} \in H_0(\mathrm{div})$, where
\begin{align*}
&H(\mathrm{div}):=\{ \bm{w} \in L^2(\mathfrak{R}_b; \mathbb{C}^3) \  |\   \mathrm{div}\bm{w} \in L^2(\mathfrak{R}_b; \mathbb{C} \}, \\
&H_0(\mathrm{div}):=\mbox{the closure of}\quad C_0^{\infty}(\mathfrak{R}_b; \mathbb{C}^3)\quad
\mbox{in}\quad H(\mathrm{div}),
\end{align*}
since
$$\inf_{\mathfrak{R}_b}\Big(\frac{\Gamma P}{\rho^2}\Big)_b >0, \quad
\inf_{\mathfrak{R}_b}\frac{1}{\rho_b} =\frac{1}{\rho_O} >0. $$
It is known that the boundary value operator
$\gamma_{\bm{N}}: \mathcal{D}(\mathfrak{R}_b) \rightarrow C(\partial\mathfrak{R}_b; \mathbb{C}):
\bm{w} \mapsto (\bm{w}|\bm{N}) $, where
$$\mathcal{D}(\mathfrak{R}_b):=C^{\infty}(\mathbb{R}^3; \mathbb{C}^3)\restriction R_b,$$
can be extended to a continous linear mapping
$\tilde{\gamma}_{\bm{N}}$ from $H(\mathrm{div})$ into
$H^{-\frac{1}{2}}(\partial \mathfrak{R}_b; \mathbb{C})$ (\cite[Theorem I.2.5]{GiraultR}), 
and $H_0(\mathrm{div})=\mathsf{N}(\tilde{\gamma}_{\bm{N}})$ (\cite[Theorem I.2.6]{GiraultR} ).
The boundary condition means 
$$
(\hat{\bm{u}}|\bm{N})=(\rho_b\bm{u}|\bm{N})=0 \quad\mbox{on}\quad \partial\mathfrak{R}_b$$ 
in this distribution sense. \\

\subsection{Case of baloclinic backgrounds}  

Let us consider the case when the background is not barotropic, but `{\bf baloclinic}', that is, $\nabla P$ and $\nabla \rho$ are not parallel somewhere in $\mathfrak{R}$.\\

 In this case $\mathcal{L}\restriction C_0^{\infty}(\mathfrak{R}; \mathbb{C}^3)$ is not symmetric.
In fact, we consider
\begin{align*}
(\mathcal{L}_0\bm{u}_1|\bm{u}_2)_{\mathfrak{H}}-(\bm{u}_1|\mathcal{L}_0\bm{u}_2)_{\mathfrak{H}}&=
\int\frac{\Gamma P}{\rho}\Big(
(\bm{u}_1|\mathfrak{a})(\nabla \rho|\bm{u}_2)
-(\bm{u}_1|\nabla \rho )(\mathfrak{a}|\bm{u}_2)\Big)d\bm{x} \\
&=\int \frac{1}{\rho}\Big(-(\bm{u}_1|\nabla P)(\nabla \rho|\bm{u}_2)+(\bm{u}_1|\nabla\rho)(\nabla P|\bm{u}_2)\Big) d\bm{x}\\
&=: \clubsuit
\end{align*}
for $\bm{u}_1,\bm{u}_2 \in C_0^{\infty}(\mathfrak{R})$. Put
$$
c(\bm{x}):=\frac{|(\nabla \rho (\bm{x})|\nabla P(\bm{x}))|}{\|\nabla \rho(\bm{x})\|\|\nabla P(\bm{x})\|}.
$$
If $\nabla P$ is not parallel to $\nabla\rho$ at a point $\bm{x}_1 \in \mathfrak{R}$, then there is a neighborhood $\mathscr{V}$ of $\bm{x}_1$, say a small ball with center $\bm{x}_1$ and radius $\delta$, such that $$\mathscr{V} \Subset \mathfrak{R}, \qquad c(\bm{x})^2 \leq 1-\kappa \quad \forall \bm{x} \in \mathscr{V}$$
with $0<\kappa <1$. Consider $\bm{u}_1,\bm{u}_2 \in C_0^{\infty}(\mathfrak{R})$ such that
\begin{equation}
\bm{u}_1 =\nabla \rho, \qquad \bm{u}_2=\nabla P\quad \mbox{on}\quad \mathscr{V}. \label{*V}
\end{equation}
Then
\begin{align*}
\clubsuit&=\int_{\mathscr{V}}\frac{1}{\rho}\|\nabla\rho\|^2\|\nabla P\|^2(1-c(\bm{x})^2)d\bm{x}
+I \\
&\geq \varepsilon + I,
\end{align*}
where
\begin{align*}
&\varepsilon=\kappa \inf_{\mathscr{V}}\frac{1}{\rho}\|\nabla \rho\|^2\|\nabla P\|^2\cdot \frac{4\pi}{3}\delta^3 > 0, \\
&I=\int_{\mathfrak{R}\setminus \mathscr{V}}
 \frac{1}{\rho}\Big(-(\bm{u}_1|\nabla P)(\nabla \rho|\bm{u}_2)+(\bm{u}_1|\nabla\rho)(\nabla P|\bm{u}_2)\Big)d\bm{x}.
 \end{align*}
 Clearly we can find $\bm{u}_1, \bm{u}_2 \in C_0^{\infty}(\mathfrak{R})$ which enjoy \eqref{*V} and
 $|I|\leq \frac{\varepsilon}{2}$.  Then $\clubsuit \geq \frac{\varepsilon}{2} >0$, that is,
 $(\mathcal{L}_0\bm{u}_1|\bm{u}_2)_{\mathfrak{H}}\not=(\bm{u}_1|\mathcal{L}_0\bm{u}_2)_{\mathfrak{H}}$.
 
Therefor we cannot suppose that $(\mathcal{L}\bm{u}|\bm{u})_{\mathfrak{H}}$ is real for any $\bm{u} \in C_0^{\infty}(\mathfrak{R}; \mathbb{C}^3)$. 
In fact we have a following counterexample. Look at
\begin{align*}
\mathfrak{Im}\Big[(\mathcal{L}\bm{u}|\bm{u})_{\mathfrak{H}}\Big]&=\int_{\mathfrak{R}}\frac{\Gamma P}{\rho}
\mathfrak{Im}[(\bm{u}|\mathfrak{a})(\nabla\rho|\bm{u})]d\bm{x} \\
&=-\int_{\mathfrak{R}}\frac{1}{\rho}\mathfrak{Im}[(\bm{u}|\nabla P)(\nabla\rho|\bm{u})]d\bm{x}.
\end{align*}
Noting that
$$\nabla P \parallel \nabla\rho \quad \Leftrightarrow\quad |(\nabla P|\nabla\rho)|=\|\nabla P\|\|\nabla \rho\|,
$$
we suppose that there is a small subdomain $\mathcal{V}$ of $\mathfrak{R}$ and a small positive number $\delta$ such that
$$|(\nabla P|\nabla\rho)|\leq (1-\delta)\|\nabla P\|\|\nabla \rho\|
\quad\mbox{on}\quad \mathcal{V}.
$$
Let us consider $\bm{u} \in C_0^{\infty}(\mathfrak{R}; \mathbb{C}^3)$ such that
$$\bm{u}=\nabla P+\mathrm{i}\nabla \rho \quad
\mbox{on}\quad \mathcal{V}.
$$
Look at
$$\int_{\mathfrak{R}}\frac{1}{\rho}\mathfrak{Im}[(\bm{u}|\nabla P)(\nabla\rho|\bm{u})]d\bm{x}
=J_1+J_2,
$$
where
\begin{align*}
J_1&=\int_{\mathcal{V}}\frac{1}{\rho}\mathfrak{Im}[(\bm{u}|\nabla P)(\nabla\rho|\bm{u})]d\bm{x} \\
J_2&=\int_{\mathfrak{R}\setminus \mathcal{V}}\frac{1}{\rho}\mathfrak{Im}
[(\bm{u}|\nabla P)(\nabla\rho|\bm{u})]d\bm{x}.
\end{align*}
Then we have 
\begin{align*}
J_1&=\int_{\mathcal{V}}\Big(\|\nabla P\|^2\|\nabla\rho\|^2-
|(\nabla P|\nabla\rho)|^2\Big)d\bm{x} \\
&\geq \delta(2-\delta)\int_{\mathcal{V}}\frac{1}{\rho}\|\nabla P\|^2\|\nabla\rho\|^2d\bm{x} \\
&=:\delta' >0.
\end{align*}
We can suppose $\displaystyle |J_2|\leq \frac{\delta'}{2}$ by evaluating $\bm{u}$ outside $\mathcal{V}$ appropriately, we have
$$
\int_{\mathfrak{R}}\frac{1}{\rho}\mathfrak{Im}[(\bm{u}|\nabla P)(\nabla\rho|\bm{u})]d\bm{x}
\geq \frac{\delta'}{2}.
$$
Then $\displaystyle \mathfrak{Im}[(\mathcal{L}\bm{u}|\bm{u})_{\mathfrak{H}} ] \leq -\frac{\delta'}{2}<0$ and
$(\mathcal{L}\bm{u}|\bm{u})_{\mathfrak{H}}  \not\in \mathbb{R}$.
This is a counterexample against $(\mathcal{L}\bm{u}|\bm{u})_{\mathfrak{H}}$ to be real.
Due to this trouble we discuss the problem as follows.\\

Now, supposing only \\

{\bf  Assumption} : {\it The equation of state enjoys {\bf Assumption} \ref{Ass.0} },\\

\noindent we consider the
operator
\begin{align}
\mathcal{L}\bm{u}&=
\nabla\Big[-\frac{\Gamma P}{\rho^2}\mathrm{div}(\rho\bm{u})+\frac{\Gamma P}{\rho}(\bm{u}|\mathfrak{a})\Big]
+\frac{\Gamma P}{\rho^2}\Big[
-\mathrm{div}(\rho\bm{u})\mathfrak{a}+(\bm{u}|\mathfrak{a})\nabla\rho\Big] + \nonumber \\
&+4\pi\mathsf{G}\nabla\mathcal{K}[\mathrm{div}(\rho\bm{u})].
\end{align}

Let us consider the formal adjoint operator $\mathcal{L}^{\dagger}$ of $\mathcal{L}$ defined by
\begin{align}
\mathcal{L}^{\dagger}\bm{u}&:=
\nabla\Big[-\frac{\Gamma P}{\rho^2}\mathrm{div}(\rho\bm{u})+\frac{\Gamma P}{\rho}(\bm{u}|\mathfrak{a})\Big]
+\frac{\Gamma P}{\rho^2}\Big[
-\mathrm{div}(\rho\bm{u})\mathfrak{a}+(\bm{u}|\nabla\rho)\mathfrak{a}\Big] + \nonumber \\
&+4\pi\mathsf{G}\nabla\mathcal{K}[\mathrm{div}(\rho\bm{u})].
\end{align}
Then, for $\bm{u}_1,\bm{u}_2 \in C_0^{\infty}$, we have
\begin{equation}
(\mathcal{L}\bm{u}_1|\bm{u}_2)_{\mathfrak{H}}=(\bm{u}_1|\mathcal{L}^{\dagger}\bm{u}_2)_{\mathfrak{H}}.
\end{equation}

Thus `` $\mathcal{L}\bm{u}=\bm{f}$ in the distribution sense'' means
$$ (\bm{u}|\mathcal{L}^{\dagger}\bm{\psi})_{\mathfrak{H}}=
(\bm{f}|\bm{\psi})_{\mathfrak{H}} \quad\mbox{for}\quad \forall
\bm{\psi} \in C_0^{\infty}(\mathfrak{R}).
$$

Put
\begin{align}
\mathcal{L}^S\bm{u}&:=\frac{1}{2}(\mathcal{L}\bm{u}+\mathcal{L}^{\dagger}\bm{u}) = \nonumber \\
&=
\nabla\Big[-\frac{\Gamma P}{\rho^2}\mathrm{div}(\rho\bm{u})+\frac{\Gamma P}{\rho}(\bm{u}|\mathfrak{a})\Big] + \nonumber \\
&+\frac{\Gamma P}{\rho^2}\Big[
-\mathrm{div}(\rho\bm{u})\mathfrak{a}+
\frac{1}{2}\Big((\bm{u}|\mathfrak{a})\nabla\rho
+(\bm{u}|\nabla\rho)\mathfrak{a}\Big)\Big] + \nonumber \\
&+4\pi\mathsf{G}\nabla\mathcal{K}[\mathrm{div}(\rho\bm{u})],
\end{align}
and
\begin{align}
\mathcal{N}\bm{u}&:=\frac{1}{2}(\mathcal{L}\bm{u}-\mathcal{L}^{\dagger}\bm{u})= \nonumber \\
&=\frac{1}{2}\frac{\Gamma P}{\rho^2}\Big[
(\bm{u}|\mathfrak{a})\nabla\rho
-(\bm{u}|\nabla\rho)\mathfrak{a}\Big].
\end{align}
Then
$\mathcal{L}^S$ is symmetric, $\mathcal{N}$ is skew-symmetric, 
that is, $(\bm{u}_1|\mathcal{N}\bm{u}_2)_{\mathfrak{H}}=-(\bm{u}_1|\mathcal{N}\bm{u}_2)_{\mathfrak{H}}$
for $\bm{u}_1,\bm{u}_2 \in C_0^{\infty}$, and
\begin{equation}
\mathcal{L}=\mathcal{L}^S+\mathcal{N},\qquad \mathcal{L}^{\dagger}=\mathcal{L}^S-\mathcal{N}.
\end{equation}

We consider the sesquilinear form $Q_0^S$ defined by
\begin{align}
Q_0^S(\bm{u}_1, \bm{u}_2)&:=
\int_{\mathfrak{R}}\frac{\Gamma P}{\rho^2}(\mathrm{div}(\rho \bm{u}_1)(\mathrm{div}(\rho \bm{u}_2)^* + \nonumber \\
&-\int_{\mathfrak{R}}\frac{\Gamma P}{\rho^2}
\Big( (\mathfrak{a}|\bm{u}_1)\mathrm{div}(\rho \bm{u}_2)
+(\bm{u}_2|\mathfrak{a})^*\mathrm{div}(\rho \bm{u}_1)\Big)d\bm{x} + \nonumber \\
&+\int_{\mathfrak{R}}\frac{\Gamma P}{\rho}(\bm{u}_1|\mathfrak{a})(\nabla\rho|\bm{u}_2)d\bm{x}+ \nonumber \\
&-4\pi\mathsf{G}\int_{\mathfrak{R}}\mathcal{K}[\mathrm{div}(\rho \bm{u}_1)]\mathrm{div}(\rho \bm{u}_2)^*
d\bm{x}
\end{align}
 on $\mathfrak{G}$. 
 Then we have
 \begin{equation}
 Q_0^S(\bm{u}_1,\bm{u}_2)=(\mathcal{L}^S\bm{u}_1|\bm{u}_2)_{\mathfrak{H}}
 \quad\mbox{for}\quad \bm{u}_1, \bm{u}_2 \in C_0^{\infty}.
 \end{equation}
 Moreover we have $Q_0^S[\bm{u}]=Q_0^S(\bm{u},\bm{u}) \in \mathbb{R}$ for $\forall \bm{u} \in \mathfrak{G}$. 
 
 Now we see that $\mathcal{L}^S\restriction C_0^{\infty}$ is bounded from below.
 Actually we look at
 \begin{align}
Q_0^S[\bm{u}]&:=
\int_{\mathfrak{R}}\frac{\Gamma P}{\rho^2}|(\mathrm{div}(\rho \bm{u})|^2d\bm{x} + \nonumber \\
&-2\int_{\mathfrak{R}}\frac{\Gamma P}{\rho^2}
\mathfrak{Re}[
 (\mathfrak{a}|\bm{u})\mathrm{div}(\rho \bm{u})^*]
d\bm{x} + \nonumber \\
&+\int_{\mathfrak{R}}\frac{\Gamma P}{\rho}\mathfrak{Re}[(\bm{u}|\mathfrak{a})(\nabla\rho|\bm{u})]d\bm{x}+ \nonumber \\
&-4\pi\mathsf{G}\int_{\mathfrak{R}}\mathcal{K}[\mathrm{div}(\rho \bm{u})]\mathrm{div}(\rho \bm{u})^*
d\bm{x}.
\end{align}

Since
$$ \frac{\Gamma P}{\rho} \in \rho^{\gamma-1}C^{3,\alpha}, \mathfrak{a}\in C^{3,\alpha}, \frac{\nabla \rho}{\rho}\in \rho^{\gamma-1}
C^{3,\alpha}(\mathfrak{R}\cup \partial\mathfrak{R}),
$$
the numbers 
\begin{align}
&\kappa_1:=\sup_{\mathfrak{R}}\Big\{\frac{\Gamma P}{\rho}\|\mathfrak{a}\|^2\Big\}, \\
&\kappa_2:=\sup_{\mathfrak{R}}\Big\{\frac{\Gamma P}{\rho^2}\|\mathfrak{a}\|\|\nabla\rho\|\Big\},
\end{align}
are finite, and, putting
\begin{align}
& a:=2\kappa_1+\kappa_2+4\pi\mathsf{G}\rho_O, \\
&Q_{a}^S(\bm{u}_1,\bm{u}_2):=Q_0^S(\bm{u}_1,\bm{u}_2)+a(\bm{u}_1,\bm{u}_2)_{\mathfrak{H}},
\end{align}
we see
\begin{align}
&\frac{1}{2}\int_{\mathfrak{R}}\frac{\Gamma P}{\rho^2}|\mathrm{div}(\rho\bm{u})|^2d\bm{x} \leq Q_a^S[\bm{u}] \leq
\nonumber  \\
&\leq 
\frac{3}{2}\int_{\mathfrak{R}}\frac{\Gamma P}{\rho^2}|\mathrm{div}(\rho\bm{u})|^2d\bm{x}
+2a\|\bm{u}\|_{\mathfrak{H}}^2.
\end{align}

Therefore $\mathcal{L}^S$ is symmetric and bounded from below, namely,
for $\bm{u} \in C_0^{\infty}$,
$$
(\mathcal{L}^S\bm{u}|\bm{u})_{\mathfrak{H}}
=Q_a^S[\bm{u}]-a\|\bm{u}\|_{\mathfrak{H}}^2 \geq -a \|\bm{u}\|_{\mathfrak{H}}^2.
$$
The Stone-Friedrichs therem says that $\mathcal{L}^S\restriction C_0^{\infty}$ admits the self-adjoint extension $\bm{L}^S$ defined on
$\mathsf{D}(\bm{L}^S)=\{ \bm{u} \in \mathfrak{G}_0 | \mathcal{L}^S \bm{u} \in \mathfrak{H} \quad
\mbox{in distribution sense} \}$. Let us put
\begin{equation}
\bm{L}:=\bm{L}^S+ \bm{N},
\end{equation}
where $\bm{N}$ is the bounded operator on $\mathfrak{H}$, $\bm{N}\bm{u}=\mathcal{N}\bm{u}$, namely,
\begin{equation}
\|\bm{N}\bm{u}\|_{\mathfrak{H}} \leq \kappa_2\|\bm{u}\|_{\mathfrak{H}},
\end{equation}
so that
$ \mathcal{L}\bm{u} \in \mathfrak{H} $ $\Leftrightarrow $ $\mathcal{L}^S\bm{u} \in \mathfrak{H}$. Hence $\bm{L}$ is a closed extension  of $\mathcal{L}\restriction C_0^{\infty}$
with domain
$$\mathsf{D}(\bm{L})=\mathsf{D}(\bm{L}^S)=\{ \bm{u} \in \mathfrak{G}_0|\mathcal{L}\bm{u} \in \mathfrak{H} \},
$$
and its adjoint operator $\bm{L}^*$ turns out to be
$$
\bm{L}^*=(\bm{L}^S)^*+ \bm{N}^*=\bm{L}^S -\bm{N}.
$$
In fact we can claim that $(\bm{T}+\bm{A})^*=\bm{T}^*+\bm{A}^*$ when $\bm{A}$
is a bounded linear operator defined on the whole space.

Note that 
\begin{equation}
(\bm{L}\bm{u}_1|\bm{u}_2)_{\mathfrak{H}}-(\bm{u}_1|\bm{L}\bm{u}_2)_{\mathfrak{H}}=-2(\bm{u}_1|\bm{N}\bm{u}_2)_{\mathfrak{H}},
\end{equation}
therefore
\begin{equation}
|(\bm{L}\bm{u}_1|\bm{u}_2)_{\mathfrak{H}}-(\bm{u}_1|\bm{L}\bm{u}_2)_{\mathfrak{H}}|\leq 2\kappa_2\|\bm{u}_1\|_{\mathfrak{H}}\|\bm{u}_2\|_{\mathfrak{H}}.
\end{equation}\\

Then we define
\begin{equation}
(\bm{u}_1|\bm{u}_2)_{\breve{\mathfrak{G}}}=Q_a^S(\bm{u}_1,\bm{u}_2)+(\bm{u}_1|\bm{u}_2)_{\mathfrak{H}}
\end{equation}
and
\begin{equation}
\|\bm{u}\|_{\breve{\mathfrak{G}}}=\sqrt{ (\bm{u}|\bm{u})_{\breve{\mathfrak{G}}} }
\end{equation}
 Then $\|\cdot\|_{\breve{\mathfrak{G}}}$ is a equivalent norm to $\|\cdot\|_{\mathfrak{G}}$ and
it hold that
\begin{equation}
\|\bm{u}\|_{\breve{\mathfrak{G}}} \geq \|\bm{u}\|_{\mathfrak{H}} \quad\mbox{for}\quad \forall \bm{u} \in \mathfrak{G},
\end{equation}
and
\begin{equation}
\mathfrak{Re}\Big[(\bm{u}|\bm{v})_{\breve{\mathfrak{G}}}\Big]=
\mathfrak{Re}\Big[((\bm{L}+a+1)\bm{u}|\bm{v})_{\mathfrak{H}}\Big] \quad\mbox{for}\quad \forall \bm{u} \in \mathsf{D}(\bm{L}), \bm{v} \in \mathfrak{G}_0.
\end{equation}\\

Anyway we can claim

\begin{Proposition}
Let $\lambda >a$. Then the equation

\begin{equation}
(\bm{L}+\lambda)\bm{u}=\bm{f}, \qquad \bm{f} \in \mathfrak{H} \label{+4.1}
\end{equation}
admits a unique solution $\bm{u}$ in $\mathsf{D}(\bm{L})$ and the solution enjoys
\begin{equation}
\|\bm{u}\|_{\mathfrak{H}}\leq \frac{1}{\lambda-a}\|\bm{f}\|_{\mathfrak{H}}.   \label{+4.2}
\end{equation}
\end{Proposition}

Proof. We consider the operator $\bm{L}+\lambda$ defined on $\mathsf{D}(\bm{L})$. If $\bm{u} \in \mathsf{D}(\bm{L})$ satisfies \eqref{+4.1}, then
$$
\mathfrak{Re}[((\bm{L}+\lambda)\bm{u}|\bm{u})_{\mathfrak{H}}] =Q_a^S[\bm{u}]+(\lambda-a)\|\bm{u}\|_{\mathfrak{H}}^2 \geq
(\lambda-a)\|\bm{u}\|_{\mathfrak{H}}^2, 
$$
since
$\mathfrak{Re}[(\bm{N}\bm{u}|\bm{u})_{\mathfrak{H}}]=0$, and
$$
\mathfrak{Re}[((\bm{L}+\lambda)\bm{u}|\bm{u})_{\mathfrak{H}} ]\leq \|\bm{f}\|_{\mathfrak{H}}\|\bm{u}\|_{\mathfrak{H}},
$$
therefore
$$
\|\bm{u}\|_{\mathfrak{H}} \leq \frac{1}{\lambda-a}\|\bm{f}\|_{\mathfrak{H}}.
$$
Hence $\bm{L}+\lambda$ is one-to-one. The range
$\mathsf{R}(\bm{L}+\lambda)=\mathsf{D}((\bm{L}+\lambda)^{-1})$ is 
dense in $\mathfrak{H}$. In fact, if
$\bm{h} \in \mathfrak{H}$ satisfies
$$((\bm{L}+\lambda)\bm{u}|\bm{h})_{\mathfrak{H}}=0 \quad \forall \bm{u} \in \mathsf{D}(\bm{L}),
$$
then $\bm{h} \in \mathsf{D}(\bm{L}^*)$ and
$$(\bm{u}|(\bm{L}^*+\lambda)\bm{h})_{\mathfrak{H}}=0 \quad \forall \bm{u} \in \mathsf{D}(\bm{L}).$$
Since $\mathsf{D}(\bm{L}) \supset C_0^{\infty}$ is dense in $\mathfrak{H}$,
we have $(\bm{L}^*+\lambda)\bm{h}=0$. Then
$$0=\mathfrak{Re}[(\bm{L}^*+\lambda)\bm{h}|\bm{h})_{\mathfrak{H}}]
\geq (\lambda-a)\|\bm{h}|_{\mathfrak{H}}^2,
$$
therefore $\bm{h}=\bm{0}$.
Since $\bm{L}$ is closed, we can claim that $(\bm{L}+\lambda)^{-1} \in \mathcal{B}(\mathfrak{H})$ and the estimate
 \eqref{+4.2} holds. $\square$ \\

Then we can apply the Hille-Yosida theory in the version described in {\bf Appendix C} .
The conclusion is

\begin{Theorem}
Suppose $\bm{u}^0 \in \mathsf{D}(\bm{L}), \bm{v}^0 \in \mathfrak{G}_0$ and
$\bm{f} \in C([0,+\infty[; \mathfrak{H})$. Then the initial value problem
\eqref{IVPa}\eqref{IVPb} admits a unique solution $\bm{u}=\bm{u}(t,\bm{x})$
in 
$$C^2([0,+\infty[; \mathfrak{H})\cap
C^1([0,+\infty[; \mathfrak{G}_0 )
\cap
C([0,+\infty[; \mathsf{D}(\bm{L})),$$
and the energy
\begin{align}
E(t, \bm{u})&:=\|\bm{u}\|_{\breve{\mathfrak{G}}}^2
+\Big\|\frac{\partial\bm{u}}{\partial  t}\Big\|_{\mathfrak{H}}^2 \nonumber \\
&=
\|\bm{u}||_{\mathfrak{H}}^2+Q_a^S[\bm{u}]+\Big\|\frac{\partial\bm{u}}{\partial  t}\Big\|_{\mathfrak{H}}^2
=(1+a)\|\bm{u}\|_{\mathfrak{H}}^2+
Q_0^S[\bm{u}]+\Big\|\frac{\partial\bm{u}}{\partial  t}\Big\|_{\mathfrak{H}}^2
\end{align}
enjoys the estimate
\begin{equation}
\sqrt{E(t,\bm{u})}\leq e^{\kappa t}\sqrt{E(0, \bm{u})}+\int_0^te^{\kappa(t-s)}
\|\bm{f}(s)\|_{\mathfrak{H}}ds
\end{equation}
for 
$$\kappa=\frac{5}{4}+\frac{a+\kappa_2}{2}=
\frac{5}{4}+\kappa_1+\kappa_2+2\pi\mathsf{G}\rho_{O}.
$$
\end{Theorem}

\section{Eigenvalue problem --Possibility of eigenmode expansions}

We are considering the equation {\bf ELASO}:
\begin{equation}
\frac{\partial^2\bm{u}}{\partial t^2}+\bm{B}\frac{\partial\bm{u}}{\partial t}+\bm{L}\bm{u} =0. \label{5.1}
\end{equation}

Since 
$\displaystyle \bm{B}\bm{v}=2\Omega
\begin{bmatrix}
0 & -1 & 0 \\
1 & 0 & 0 \\
0 & 0 & 0
\end{bmatrix}
\bm{v}$,
 $\bm{B}$ is a skew-symmetric bounded linear operator in $\mathcal{B}(\mathfrak{H})$,
 $\bm{B}^*=-\bm{B}$, and
  $\mathfrak{Re}[(\bm{B}\bm{v}|\bm{v})]=0 \quad \forall \bm{v}$. Put
\begin{equation}
 \beta^{\star} := |\|\bm{B}\||_{\mathcal{B}(\mathfrak{H})} = \sup \frac{\|\bm{B}\bm{v}\|_{\mathfrak{H}}}{\|\bm{v}\|_{\mathfrak{H}}} =2|\Omega|.
 \end{equation}
 $\bm{L}$ is a self-adjoint operator in $\mathfrak{H}$
 and
 $$((\bm{L}+a)\bm{u}|\bm{u})_{\mathfrak{H}} \geq 0 \quad\mbox{for}\quad \bm{u} \in \mathsf{D}(\bm{L}).
 $$
Put
 \begin{equation}
 \mu_*:=\Big(-\inf_{\bm{u} \in \mathsf{D}(\bm{L})} \frac{(\bm{L}\bm{u}|\bm{u})_{\mathfrak{H}}}{\|\bm{u}\|_{\mathfrak{H}}^2}\Big)\vee 0.
 \end{equation}
 Then  
  \begin{equation}
 (\bm{L}\bm{u}|\bm{u})_{\mathfrak{H}} \geq -\mu_*\|\bm{u}\|_{\mathfrak{H}}^2
 \quad \forall \bm{u} \in \mathsf{D}(\bm{L}),
 \end{equation}
 and
 $ \mu_*\leq a$.\\

 Now we look for nontrivial solution of \eqref{5.1} of the form
 $ \bm{u}(t,\bm{x})=e^{\lambda t}\bm{u}^0(\bm{x})$, where $\lambda \in \mathbb{C}, \bm{u}^0 \in \mathsf{D}(\bm{L})$. Namely
 
 \begin{Definition}
 If there is a $\bm{u}^0 \in \mathsf{D}(\bm{L}) $ such that $\bm{u}^0 \not= \bm{0}$, and
 \begin{equation}
 \Big[\lambda^2+\lambda \bm{B} +\bm{L}\Big]\bm{u}^0=\bm{0}
 \end{equation}
 holds, then $\lambda$ is called an eigenvalue of {\bf ELASO}\eqref{5.1} and $\bm{u}^0$ is called an eigenvector associated with the eigenvalue $\lambda$.
 \end{Definition}
 
 If $\lambda$ is an eigenvalue of {\bf ELASO}\eqref{5.1} and $\bm{u}^0$ is an associated eigenvector such that $\|\bm{u}^0\|_{\mathfrak{H}}=1$, then \\
 
 {\it 
\noindent  1)
 $\bm{u}: (t,\bm{x}) \mapsto e^{\lambda t}\bm{u}^0(\bm{x})$ is a solution of \eqref{5.1}
 with the initial condition
 \begin{equation}
 \bm{u}|_{t=0}=\bm{u}^0(\bm{x}),\quad
 \frac{\partial\bm{u}}{\partial t}\Big|_{t=0}=\lambda\bm{u}^0(\bm{x}),
 \end{equation} \\
 2) $\bm{u}^{\mathfrak{Re}}: (t,\bm{x}) \mapsto \mathfrak{Re}[e^{\lambda t}\bm{u}^0(\bm{x})]$
 is a solution of \eqref{5.1} with the initial conditon
 \begin{equation}
 \bm{u}^{\mathfrak{Re}}|_{t=0}=\mathfrak{Re}[\bm{u}^0(\bm{x})],
 \quad
 \frac{\partial \bm{u}^{\mathfrak{Re}}}{\partial t}\Big|_{t=0}=
 \mathfrak{Re}[\lambda]\mathfrak{Re}[\bm{u}^0(\bm{x})]-
 \mathfrak{Im}[\lambda]\mathfrak{Im}[\bm{u}^0(\bm{x})],
 \end{equation} \\
 3) $\lambda$ satisfies the quadratic equation
 \begin{equation}
  \lambda^2+\mathrm{i} b \lambda + c =0
 \end{equation}
 where the coefficents $b,c$ are given by
 \begin{equation}
  b=-\mathrm{i} (B\bm{u}^0|\bm{u}^0)_{\mathfrak{H}},
 \quad c=(\bm{L}\bm{u}^0|\bm{u}^0)_{\mathfrak{H}},
 \end{equation}
 and $b,c$ are real numbers.
  }\\
 
By Theorem \ref{Th.*4} we see that
 $0$ is an eigenvalue of {\bf ELASO}\eqref{5.1}.
The solution $\bm{u}(t,\bm{x})=\bm{u}^0(\bm{x})$ with $\bm{u}^0 \in \mathsf{N}(\bm{L})$ may be called a `trivial perturbation' after \cite{FriedmanS}.\\

Recall that by the operator $\bm{A}$ in $\mathfrak{E}=\mathsf{G}_0 \times \mathfrak{H}$ defined by
\begin{equation}
\bm{A} U=
\begin{bmatrix}
O & -I \\
\\
\bm{L} & \bm{B}
\end{bmatrix}
U,
\quad \mathsf{D}(\bm{A})=\mathsf{D}(\bm{L}) \times \mathsf{G}_0,
\end{equation}
the equation \eqref{5.1} reads {\bf ELASO\ddag}:
\begin{equation}
\frac{d U}{dt} +\bm{A}U=\bm{0} \label{E5.11}
\end{equation}
with $$U=
\begin{bmatrix}
\bm{u} \\
\\
\displaystyle  \frac{d \bm{u}}{d t}
 \end{bmatrix}.
 $$

\begin{Proposition}
$\lambda \in \mathbb{C}$ is an eigenvalue of {\bf ELASO}\eqref{5.1} if and only if
$\lambda$ is an eigenvalue of the operator $-\bm{A}$.
\end{Proposition}

In fact, if $\bm{u}^0$ is an eigenvector of {\bf ELASO}\eqref{5.1} associated with eigenvector $\lambda$, then 
$U^0=\begin{bmatrix}
\bm{u}^0 \\
\\
\lambda \bm{u}^0
\end{bmatrix}
$
is an eigenvector of $\bm{A}$ with eigenvalue $-\lambda$; if
$U^0=\begin{bmatrix}
\bm{u}^0 \\
\\
\bm{v}^0
\end{bmatrix}
(\not= 
\begin{bmatrix}
\bm{0} \\
\\
\bm{0}
\end{bmatrix})$
is an eigenvector of $\bm{A}$ with eigenvalue $-\lambda$, then
$(\lambda^2+\lambda\bm{B}+\bm{L})\bm{u}^0=\bm{0}, \bm{u}^0\not=\bm{0}$,
since $\bm{v}^0=\lambda \bm{u}^0$, that is, $\lambda$ is an eigenvalue of {\bf ELASO}\eqref{5.1}.\\

Following \cite{DysonS}, we use the terminology defined by
\begin{Definition}
The family of operators $\mathfrak{L}=(\mathfrak{L}(\lambda))_{\lambda \in \mathbb{C}}$,
where
\begin{equation}
\mathfrak{L}(\lambda)=\lambda^2+\lambda \bm{B}+\bm{L},
\end{equation}
is called the `quadratic pencil' of the equation \eqref{5.1}.

The resolvent set $\mbox{\Pisymbol{psy}{82}}(\mathfrak{L})$ is the set of
$\lambda \in \mathbb{C}$ such that
$\mathfrak{L}(\lambda)$ has the bounded incverse in $\mathcal{B}(\mathfrak{H})$.
$\mathbb{C} \setminus \mbox{\Pisymbol{psy}{82}}(\mathfrak{L})$ is the spectrum of $\mathfrak{L}$, and is denoted by $\mbox{\Pisymbol{psy}{83}}(\mathfrak{L})$.

The set of all eigenvalues of {\bf ELASO}\eqref{5.1} is denoted by $\mbox{\Pisymbol{psy}{83}}_{\mathrm{p}}(\mathfrak{L})$.
\end{Definition}

If $\lambda $ is an eigenvalue of \eqref{5.1}, then it belongs to the spectrum of $\mathfrak{L}$, that is, 
$\mbox{\Pisymbol{psy}{83}}_{\mathrm{p}}(\mathfrak{L})
\subset \mbox{\Pisymbol{psy}{83}}(\mathfrak{L})$. However we cannot say that any $\lambda \in \mbox{\Pisymbol{psy}{83}}(\mathfrak{L})$ is an eigenvalue of \eqref{5.1},
that is,  $\mbox{\Pisymbol{psy}{83}}_{\mathrm{p}}(\mathfrak{L})
= \mbox{\Pisymbol{psy}{83}}(\mathfrak{L})$, a priori.
We are going to study the structure of $\mbox{\Pisymbol{psy}{83}}(\mathfrak{L})$
as much as possible at the moment.\\

%%%%%%%%%%%%%%%%
%%%%%%%%%%%%%%%%%%%

\begin{Proposition}
$\mbox{\Pisymbol{psy}{82}}(\mathfrak{L})$ is an open subset of $\mathbb{C}$ and $\mbox{\Pisymbol{psy}{83}}(\mathfrak{L})$ is  closed .
\end{Proposition}

Proof. Let us consider $\lambda \in \mbox{\Pisymbol{psy}{82}}(\mathfrak{L})$. Then 
$$\mathfrak{L}(\lambda+\Delta\lambda)=
\mathfrak{L}(\lambda)\Big[
I+ \mathfrak{L}(\lambda)^{-1}(\Delta\lambda(2\lambda+\bm{B}+\Delta\lambda)\Big]
$$
admits the bounded inverse and $\lambda+\Delta\lambda \in \mbox{\Pisymbol{psy}{82}}(\mathfrak{L})$, if
$$|\|\mathfrak{L}(\lambda)^{-1}\Delta\lambda(2\lambda+\bm{B}+\Delta\lambda)\||_{\mathcal{B}(\mathfrak{H})} <1.$$
For this inequality, it is sufficient that
$$|\|\mathfrak{L}(\lambda)^{-1}\||_{\mathcal{B}(\mathfrak{H})}\cdot|\Delta\lambda|\cdot|(2|\lambda|
+\beta^{\star}+|\Delta\lambda| ) < 1,$$
or
$$|\Delta\lambda|<-\Big(|\lambda|+\frac{\beta^{\star}}{2}\Big)+\sqrt{\Big(|\lambda|+\frac{\beta^{\star}}{2}\Big)^2+
\frac{1}{|\|\mathfrak{L}(\lambda)^{-1}\||_{\mathcal{B}(\mathfrak{H})}}}.$$
This means $\mbox{\Pisymbol{psy}{82}}(\mathfrak{L})$ is open.
$\square$\\

\begin{Proposition}
$\mbox{\Pisymbol{psy}{82}}(\mathfrak{L})$ and $\mbox{\Pisymbol{psy}{83}}(\mathfrak{L})$ are symmetric about the imaginary axis $\mathrm{i} \mathbb{R}$. 
\end{Proposition}

In fact we see $\mathfrak{L}(\lambda)^*=\mathfrak{L}(-\lambda^*)$.

\begin{Proposition}
It hold
$$ \mbox{\Pisymbol{psy}{82}}(\mathfrak{L})=\mbox{\Pisymbol{psy}{82}}(-\bm{A}),\quad \mbox{\Pisymbol{psy}{83}}(\mathfrak{L})=\mbox{\Pisymbol{psy}{83}}(-\bm{A}),
\quad
\mbox{\Pisymbol{psy}{83}}_{\mathrm{p}}(\mathfrak{L})
\subset \mbox{\Pisymbol{psy}{83}}_{\mathrm{p}}(-\bm{A}).$$
Here $\mbox{\Pisymbol{psy}{82}}(-\bm{A}), \mbox{\Pisymbol{psy}{83}}(-\bm{A}),
\mbox{\Pisymbol{psy}{83}}_{\mathrm{p}}(-\bm{A})$
stand for the usual resolvent set, spectrum, point spectrum (the set of all eigenvalues) of the operator
$-\bm{A}$ in $ \breve{\mathfrak{G}_0} \times \mathfrak{H}$.
\end{Proposition}

Proof. Let $\lambda \in \mbox{\Pisymbol{psy}{82}}(\mathfrak{L})$. Consider the equation
$$\mathbf{A}U+\lambda U=F=
\begin{bmatrix}
\bm{g} \\
\\
\bm{f}
\end{bmatrix}
\in \breve{\mathfrak{G}_0}\times \mathfrak{H}, $$
or,
$$
\begin{cases}
&-\bm{v}+\lambda \bm{u}=\bm{g} \\
& \\
& \bm{L}\bm{u} +\bm{B}\bm{u}+\lambda \bm{v}=\bm{f}
\end{cases}
$$
This system of equations can be solved
as
$$
\begin{cases}
\bm{u}&=
\mathfrak{L}(\lambda)^{-1}(\bm{f} +\bm{B}\bm{g} -\lambda\bm{g})
 \\
& \\
\bm{v}&=\lambda \mathfrak{L}(\lambda)^{-1}(\bm{f}+\bm{B}\bm{g}-\lambda\bm{g})-\bm{f}
\end{cases}
$$
since $\lambda \in \mbox{\Pisymbol{psy}{82}}(\mathfrak{L})$
, while
$F \mapsto U$ is continuous. Here we note that $\mathfrak{L}(\lambda)^{-1} \in \mathcal{B}(\mathfrak{H}, \mathfrak{H})$ belongs to $ \mathcal{B}(\mathfrak{H}, \breve{\mathfrak{G}_0})$ thanks to the following Proposition \ref{PropBHBG}.
 Therefore $\lambda \in 
\mbox{\Pisymbol{psy}{82}}(-\mathbf{A})$.

Inversely let $\lambda \in \mbox{\Pisymbol{psy}{82}}(-\mathbf{A})$. Consider the equation
$$(\lambda^2+\lambda\bm{B}+\bm{L})\bm{u}=\bm{f}
\in \mathfrak{H},$$
which is equivalent to the system of equations
$$
\begin{cases}
&  \lambda \bm{v} +\bm{B}\bm{v} +\bm{L}\bm{u}=\bm{f},\\
& \\
&\bm{v}=\lambda\bm{u}
\end{cases}
$$
But this is nothing but
$$\mathbf{A}U+\Lambda U=
\begin{bmatrix}
\mathbf{0} \\
\\
\bm{f}
\end{bmatrix}
$$
for $U=(\bm{u}, \bm{v})^{\top}$.
Since $\lambda \in \mbox{\Pisymbol{psy}{82}}(-\mathbf{A})$ is supposed, this admits the solution
$$U=
\begin{bmatrix}
\bm{u} \\
\\
\bm{v}
\end{bmatrix}=
(\mathbf{A}+\lambda)^{-1}
\begin{bmatrix}
\mathbf{0} \\
\\
\bm{f}
\end{bmatrix},
$$
and $\bm{f} \mapsto \bm{u}$ is continuous, that is,
$\lambda  \in \mbox{\Pisymbol{psy}{82}}
(\mathfrak{L})$. $\square$\\

\begin{Proposition}\label{PropBHBG}
$\lambda \in \mbox{\Pisymbol{psy}{82}}(\mathfrak{L})$, that is, $\mathfrak{L}(\lambda)^{-1} \in \mathcal{B}(\mathfrak{H})$ if and only if
$\mathfrak{L}(\lambda)^{-1} \in \mathcal{B}(\mathfrak{H}, \breve{\mathfrak{G}_0})$.
\end{Proposition}

Proof. Suppose $\mathfrak{L}(\lambda)^{-1} \in \mathcal{B}(\mathfrak{H})$. We consider
$$
\mathfrak{L}(\lambda)\bm{u}=(\lambda^2+\lambda \bm{B}+\bm{L})\bm{u}=\bm{f} \in\mathfrak{H},\quad \bm{u} \in \mathsf{D}(\bm{L}).
$$
Then 
$$ \lambda^2\|\bm{u}\|_{\mathfrak{H}}^2+\lambda (\bm{B}\bm{u}|\bm{u})_{\mathfrak{H}})_{\mathfrak{H}}+(\bm{L}\bm{u}|\bm{u})_{\mathfrak{H}}=(\bm{f}|\bm{u})_{\mathfrak{H}},
$$
while
$$ \|\bm{u}\|_{\breve{\mathfrak{G}}}^2=(\bm{L}\bm{u}|\bm{u})_{\mathfrak{H}}+(1+a)\|\bm{u}\|_{\mathfrak{H}}^2.
$$
Therefore
we have 
\begin{align*}
\|\bm{u}\|_{\breve{\mathfrak{G}}}^2 & \leq \|\bm{f}\|_{\mathfrak{H}}\|\bm{u}\|_{\mathfrak{H}}+
(|1+a-\lambda^2|+|\lambda|\beta^{\star})\|\bm{u}\|_{\mathfrak{H}}^2 \\&\leq \Big( |\|\mathfrak{L}(\lambda)^{-1}\||_{\mathcal{B}(\mathfrak{H})}
+|1+a-\lambda^2|+|\lambda|\beta^{\star} \Big)\|\bm{u}\|_{\mathfrak{H}}^2,
\end{align*}
so that
\begin{align*}
\|\bm{u}\|_{\breve{\mathfrak{G}}} & \leq 
\Big[ |\|\mathfrak{L}(\lambda)^{-1}\||_{\mathcal{B}(\mathfrak{H})}
+|1+a-\lambda^2|+|\lambda|\beta^{\star} \Big]^{\frac{1}{2}}\|\bm{u}\|_{\mathfrak{H}} \\
& \leq
|\|\mathfrak{L}(\lambda)^{-1}\||_{\mathcal{B}(\mathfrak{H})}
\Big[ |\|\mathfrak{L}(\lambda)^{-1}\||_{\mathcal{B}(\mathfrak{H})}
+|1+a-\lambda^2|+|\lambda|\beta^{\star} \Big]^{\frac{1}{2}}\|\bm{f}\|_{\mathfrak{H}},
\end{align*}
Consequently we see $\mathfrak{L}(\lambda)^{-1} \in \mathcal{B}(\mathfrak{H}, \breve{\mathfrak{G}_0})$ and
$$
|\| \mathfrak{L}(\lambda)^{-1}|\|_{\mathcal{B}(\mathfrak{H},\breve{\mathfrak{G}_0})}
\leq 
|\|\mathfrak{L}(\lambda)^{-1}\||_{\mathcal{B}(\mathfrak{H})}
\Big[ |\|\mathfrak{L}(\lambda)^{-1}\||_{\mathcal{B}(\mathfrak{H})}
+|1+a-\lambda^2|+|\lambda|\beta^{\star} \Big]^{\frac{1}{2}}.
$$

Suppose $\mathfrak{L}(\lambda)^{-1} \in \mathcal{B}(\mathfrak{H},\breve{\mathfrak{G}_0})$.
We consider
$$
\mathfrak{L}(\lambda)\bm{u}=(\lambda^2+\lambda \bm{B}+\bm{L})\bm{u}=\bm{f} \in\mathfrak{H},\quad \bm{u} \in \mathsf{D}(\bm{L}).
$$
We have
\begin{align*}
\| \bm{u}\|_{\mathfrak{H}}^2&=\| \bm{u} |_{\breve{\mathfrak{G}}}^2-(\bm{L}+a)\bm{u}|\bm{u})_{\mathfrak{H}} \\
&\leq  \|\bm{u}\|_{\breve{\mathfrak{G}}}^2 \\
&\leq \Big[ {C}|\|\mathfrak{L}(\lambda)^{-1}\||_{\mathcal{B}(\mathfrak{H},\breve{\mathfrak{G}_0})} \|\bm{f}\|_{\mathfrak{H}}^2
\Big]^2.
\end{align*}

Consequently $\mathfrak{L}(\lambda)^{-1} \in \mathcal{B}(\mathfrak{H})$ and
$$
|\|\mathfrak{L}(\lambda)^{-1}\||_{\mathcal{B}(\mathfrak{H})}
\leq  |\|\mathfrak{L}(\lambda)^{-1}\||_{\mathcal{B}(\mathfrak{H},\breve{\mathfrak{G}})}. $$
$\square$.\\

We consider the question:

 \begin{Question}\label{Q.1}
 Is $\mbox{\Pisymbol{psy}{83}}_{\mathrm{p}}(\mathfrak{L})=\mbox{\Pisymbol{psy}{83}}(\mathfrak{L})$?
 \end{Question}
 
\noindent  in the following view point.\\

 Let us consider the case of $\bm{B}=O$, or $\Omega=0$. Then, we see that, if
$\mbox{\Pisymbol{psy}{83}}_{\mathrm{p}}(\bm{L})=\mbox{\Pisymbol{psy}{83}}(\bm{L})$,  
there is an orthonormal system $(\bm{\phi}_n)_n$ which consists of eigenvectors of 
$\bm{L}$, say, $\bm{L}\bm{\phi}_n=\lambda_n\bm{\phi}_n, \lambda_n \in \mathbb{R}$, and is complete in $\mathfrak{H}$, according to
\begin{quote}

{\bf Theorem}(\cite[Theorem X.3.4]{DunfordS}): {\it If the spectrum of the normal operator $T$ in the Hilbert space $\mathsf{X}$ is countable, then there is an orthonormal basis $\mathscr{B}$ for $\mathsf{X}$ consisting of eigenvectors of $T$;
$$x=\sum_{y \in \mathscr{B}}(x|y)_{\mathsf{X}}y \quad \forall x \in \mathsf{X},
$$
and, for each $x$, all but a countable number of the coefficients $(x|y)_{\mathsf{X}}$ are zero. }
\end{quote}

Namely, if $\Omega=0$ and if
 $\mbox{\Pisymbol{psy}{83}}_{\mathrm{p}}(\bm{L})=\mbox{\Pisymbol{psy}{83}}(\bm{L})$ is the case, then 
 all solutions $\bm{u}(t, \bm{x})$ of \eqref{5.1}
 in $
C([0,+\infty[; \mathfrak{H}) \cap
C^1([0,+\infty[; \mathfrak{G}_0) \cap
C^2([0,+\infty[; \mathsf{D}(\bm{L}))$
can be expanded by 
 the eigenmodes 
$$
\bm{u}_n^{\pm}(t,\bm{x})=
\begin{cases}
e^{\pm\sqrt{\lambda_n}\mathrm{i}t}\bm{\phi}_n(\bm{x}) \quad (\lambda_n \geq 0) \\
e^{\pm\sqrt{|\lambda_n|}t}\bm{\phi}_n(\bm{x}) \quad (\lambda_n <0)
\end{cases}
$$
as
$$
\bm{u}(t, \bm{x})=\sum_n(c_n^+\bm{u}_n^+(t,\bm{x})+c_n^-\bm{u}_n^-(t,\bm{x})),
$$
where
$$
c_n^{\pm}=
\begin{cases}
\frac{1}{2}
\Big((\bm{u}^0|\bm{\phi}_n)_{\mathfrak{H}}\pm\frac{1}{\sqrt{\lambda_n}\mathrm{i}}(\bm{v}^0|\bm{\phi}_n)_{\mathfrak{H}}\Big) \quad (\lambda_n >0) \\
\frac{1}{2}\Big((\bm{u}^0|\bm{\phi}_n)_{\mathfrak{H}}\pm\frac{1}{\sqrt{|\lambda_n|}}(\bm{v}^0|\bm{\phi}_n)_{\mathfrak{H}}\Big) \quad (\lambda_n <0) \\
\frac{1}{2}(\bm{u}^0|\bm{\phi}_n)_{\mathfrak{H}} \quad (\lambda_n=0),
\end{cases}
$$
with  $\bm{u}^0(\bm{x})=\bm{u}(0,\bm{x}), \bm{v}^0(\bm{x})=\frac{\partial \bm{u}}{\partial t}(0,\bm{x})$.

But, if $\mbox{\Pisymbol{psy}{83}}_{\mathrm{p}}(\bm{L})=\mbox{\Pisymbol{psy}{83}}(\bm{L})$ is not the case, we cannot claim the completeness of eigenmodes, the possibility of
the eigenmode expansion. In this sense we we are concerned with whether  $\mbox{\Pisymbol{psy}{83}}_{\mathrm{p}}(\mathfrak{L})=\mbox{\Pisymbol{psy}{83}}(\mathfrak{L})$ or not. \\

Ler us note the following observations:\\

1)

Suppose that $\Omega\not=0, \bm{B}\not=O$. Then, even if $\mbox{\Pisymbol{psy}{83}}(\mathfrak{L})=
\mbox{\Pisymbol{psy}{83}}_p(\mathfrak{L})$, we know neither whether $\mbox{\Pisymbol{psy}{83}}(\bm{L})=\mbox{\Pisymbol{psy}{83}}_p(\bm{L})$
nor whether eigenmode expansion is available or not. Actually, for an eigenvalue $\mu\not=0$ of $\bm{L}$ with eigenvector $\bm{\phi}$ such that 
$\overset{=}{\bm{\phi}}=
\begin{bmatrix}
\phi^1 \\
\phi^2 \\
0
\end{bmatrix}
\not=\bm{0}$, we have $\mathfrak{L}(\pm\sqrt{-\mu})
\bm{\phi}=\pm\sqrt{-\mu}\bm{B}\bm{\phi}\not=\bm{0}$. The relation between 
$\mbox{\Pisymbol{psy}{83}}(\mathfrak{L})$ and $\mbox{\Pisymbol{psy}{83}}(\bm{L})$, and between 
$\mbox{\Pisymbol{psy}{83}}_p(\mathfrak{L})$ and $\mbox{\Pisymbol{psy}{83}}_p(\bm{L})$ are not clear.\\

2)

 Supoose
 $\mbox{\Pisymbol{psy}{83}}(\mathfrak{L})=\mbox{\Pisymbol{psy}{83}}_p(\mathfrak{L})$. Then,  although $\mbox{\Pisymbol{psy}{83}}(\bm{A})=\mbox{\Pisymbol{psy}{83}}_p(\bm{A})$, we cannot claim that there is a complete orthonormal basis in $\mathfrak{E}=\breve{\mathfrak{G}_0} \times \mathfrak{H}$ which consistes of eigenvectors of $\bm{A}$, since $\bm{A}$ is  not normal. In fact we have
 \begin{Proposition}
 The adjoint of $\bm{A}$ turns out to be
 \begin{align}
 &\mathsf{D}(\bm{A}^*)=\mathsf{D}(\bm{L})\times \breve{\mathfrak{G}_0}, \\
 &\bm{A}^*=
 \begin{bmatrix}
 O & I-(a+1)\bm{L}_a^{-1} \\
 \\
 -\bm{L}_a & -B
 \end{bmatrix},
 \end{align}
 where $\bm{L}_a=\bm{L}+a+1$.
 \end{Proposition}
 
 Proof. Recall that
 $ \|\bm{u}\|_{\breve{\mathfrak{G}}}^2=
 ((\bm{L}+a+1)\bm{u}|\bm{u})_{\mathfrak{H}}$ for $\forall  \bm{u} \in \mathsf{D}(\bm{L})$, $a \geq 0$
 and that
 $\bm{L}_a^{-1}$ belongs to $\mathcal{B}(\mathfrak{H})$ and is self-adjoint.
 
 For $\bm{F}=\begin{bmatrix}
 \bm{g} \\
\bm{f}
 \end{bmatrix}, 
 \bm{F}^{\dag}=\begin{bmatrix}
 \bm{g}^{\dag} \\
 \bm{f}^{\dag}
 \end{bmatrix}
 $,
 `` $\bm{F} \in \mathsf{D}(\bm{A}^*), \bm{A}^*\bm{F}=\bm{F}^{\dag}$
 ``
 means
 \begin{equation}
 (\bm{A}U|\bm{F})_{\mathfrak{E}}=(U|\bm{F}^{\dag})_{\mathfrak{E}} \quad\mbox{for}\quad \forall U \in \mathsf{D}(\bm{A})
 \end{equation}
 or
 \begin{align}
 -(\bm{v}|\bm{g})_{\breve{\mathfrak{G}}}+(\bm{L}\bm{u}+\bm{B}\bm{v}|\bm{f})_{\mathfrak{H}}&=(\bm{u}|\bm{g}^{\dag})_{\breve{\mathfrak{G}}}+(\bm{v}|\bm{f}^{\dag})_{\mathfrak{H}} \nonumber \\
 &\mbox{for}\quad \forall \bm{u} \in \mathsf{D}(\bm{A}) \forall \bm{v} \in \breve{\mathfrak{G}_0}. \label{**A.2}
 \end{align}
 We want to show that \eqref{**A.2} holds if and only if
 $\bm{g} \in \mathsf{D}(\bm{L}), \bm{f} \in \breve{\mathfrak{G}_0}$ and
 \begin{equation}
 \bm{g}^{\dag}=\bm{f}-(a+1)\bm{L}_a^{-1}\bm{f},\quad
 \bm{f}^{\dag}=-\bm{L}_a\bm{g} -\bm{B}\bm{f}. \label{**A.3}
 \end{equation}
 
 Suppose \eqref{**A.2}. Taking $\bm{u}=\bm{0}$, we see
 $$ -(\bm{L}_a\bm{v}|\bm{g})_{\mathfrak{H}}+(\bm{B}\bm{v}|\bm{f})_{\mathfrak{H}}=(\bm{v}|\bm{f}^{\dag})_{\mathfrak{H}}\quad \forall \bm{v} \in \mathsf{D}(\bm{L}).
 $$
 Consider $\bm{h}=\bm{L}_a\bm{v}$, while $\forall \bm{v} \in \mathsf{D}(\bm{L}) \Leftrightarrow
 \forall \bm{h} \in \mathfrak{H}$, we have
 $$ -(\bm{h}|\bm{g})_{\mathfrak{H}}+(\bm{B}\bm{L}_a^{-1}\bm{h}|\bm{f})_{\mathfrak{H}}=
 (\bm{L}_a^{-1}\bm{h}|\bm{f}^{\dag})_{\mathfrak{H}}\quad \forall \bm{h} \in \mathfrak{H}.$$
 Therefore
 $$-\bm{g}-\bm{L}_a^{-1}\bm{B}f=\bm{L}_a^{-1}\bm{f}^{\dag}.
 $$
 Here we know that $\bm{L}_a^{-1} \in \mathcal{B}(\mathfrak{H})$ is self-adjoint. Hence
 $$\bm{f}^{\dag}=-\bm{L}_a \bm{g}-\bm{B}\bm{f}.$$
 So \eqref{**A.2} reduces
 to
 $$ (\bm{L}\bm{u}|\bm{f})_{\mathfrak{H}}=(\bm{u}|\bm{g}^{\dag})_{\breve{\mathfrak{G}}} \quad \forall \bm{u} \in \mathsf{D}(\bm{L}).
 $$
 Recall $\bm{B}^*=-\bm{B}$. Consider $\bm{k}=\bm{L}_a\bm{u}$, while $\forall \bm{u} \in \mathsf{D}(\bm{L})
 \Leftrightarrow \bm{k} \in \mathfrak{H}$. Then 
 $$
 (\bm{L}\bm{L}_a^{-1}\bm{k}|\bm{f})_{\mathfrak{H}}=(\bm{k}|\bm{g}^{\dag})_{\mathfrak{H}} \quad \forall \bm{k} \in \mathfrak{H}.$$
 Since $\bm{L}\bm{L}_a^{-1}=I- (a+1)\bm{L}_a^{-1}$, this reads
  $$
 (\bm{k}|\bm{f})_{\mathfrak{H}}-(a+1)(\bm{k}|\bm{L}_a^{-1}\bm{f})_{\mathfrak{H}}
 =(\bm{k}|\bm{g}^{\dag})_{\mathfrak{H}} \quad \forall \bm{k} \in \mathfrak{H},$$
 so that
 $$\bm{f} -(a+1)\bm{L}_a^{-1}\bm{f}=\bm{g}^{\dag}.
 $$
 Since $\bm{g}^{\dag} \in \breve{\mathfrak{G}_0}$, we see $\bm{f} = \bm{g}^{\dag}+(a+1)\bm{L}_a^{-1}\bm{f} $  should belong to $\breve{\mathfrak{G}_0}$. Hence $\bm{g} \in \mathsf{D}(\bm{L}), \bm{f} \in \breve{\mathfrak{G}_0}$ and \eqref{**A.3} should holds.
  
  Inversely suppose that $\bm{g} \in \mathsf{D}(\bm{L}), \bm{f} \in \breve{\mathfrak{G}_0}$ and \eqref{**A.3} holds. Then
  we see, for $\bm{u} \in \mathsf{D}(\bm{L}), \bm{v} \in \breve{\mathfrak{G}_0}$, 
  \begin{align*}
  -(\bm{v}|\bm{g})_{\breve{\mathfrak{G}}}+(\bm{L}\bm{u}+\bm{B}\bm{v}|\bm{f})_{\mathfrak{H}}&=
  -(\bm{v}|\bm{L}_a\bm{g})_{\mathfrak{H}}+
  (\bm{L}\bm{u}+\bm{B}\bm{v}|\bm{f})_{\mathfrak{H}} \\
  &=(\bm{v}|-\bm{L}_a\bm{f}-\bm{B}\bm{f})_{\mathfrak{H}}+(\bm{L}\bm{u}|\bm{f})_{\mathfrak{H}} \\
  &=(\bm{v}|\bm{f}^{\dag})_{\mathfrak{H}}+(\bm{L}_a\bm{u}|\bm{g}^{\dag})_{\mathfrak{H}} \\
  &=(\bm{v}|\bm{f}^{\dag})_{\mathfrak{H}}+(\bm{u}|\bm{g}^{\dag})_{\mathfrak{H}}.
  \end{align*}
  That is, \eqref{**A.2} holds. $\square$\\
  
  3)
  
  Suppose $\Omega=0, \bm{B}=O$.   
  If $\phi$ is an eigenvector of $\bm{L}$ associated with a positive eigenvalue $\mu$, then $\displaystyle \Phi_{\pm}=
\begin{bmatrix}
\phi \\
\\
\pm\sqrt{\mu}\mathrm{i}\phi
\end{bmatrix}
$ are two different eigenvectors of $\bm{A}$ associated with the different eigenvalues $\pm\sqrt{\mu}\mathrm{i}$. But $$(\Phi_+|\Phi_-)_{\mathfrak{E}}=\|\phi\|_{\breve{\mathfrak{G}}}^2
-\mu \|\phi\|_{\mathfrak{H}}^2=
(a+1)\|\phi\|_{\mathfrak{H}}^2 \not=0, $$
that is, they are not orthogonal.\\

%%%%%%%%%%%%%%%%%%%%%%%%%%%%%%%%%%%%%%

\textbullet\   Here we note the following fact:

\begin{Theorem} \label{Th.X}
Suppose that the background is isentropic, that is, $\mathfrak{a}_b=0$ everywhere.
Let $\bm{L}^G$ be te Friedrichs extension of the operator $\mathcal{L}\restriction C_0^{\infty}(\mathfrak{R}_b,; \mathbb{C}^3)$ in the functional space
$\mathfrak{G}=\{ \bm{u} \in \mathfrak{H} | \mathrm{div}(\rho_b\bm{u}) \in L^2(\mathfrak{R}_b, \Big(\frac{\Gamma P}{\rho^2}\Big)_bd\bm{x}; \mathbb{C})\}$. Then
 $ \mbox{\Pisymbol{psy}{83}}(\bm{L}^G)=
 \mbox{\Pisymbol{psy}{83}}_{\mathrm{p}}(\bm{L}^G)$, and
$ \mbox{\Pisymbol{psy}{83}}(\bm{L}^G)$ consists of $\{0\}$ and a sequence of eigenvalues of finite multiplicities which accumates $+\infty$.
\end{Theorem}

The proof is the same as that of \cite[Sections III, IV ]{JJTM2020} for the case of spherically symmetric background for $\Omega=\omega_b=0$, but we describe a sketch of proof for the sake of self-containedness. First of all we note that, since $\mathfrak{a}_b=0$, it is the case
$$\mathcal{L}\bm{u}=
\mathrm{grad}\Big(
-\frac{\Gamma P}{\rho^2}\mathrm{div}(\rho\bm{u})+
4\pi\mathsf{G}\mathcal{K}[\mathrm{div}(\rho\bm{u})]\Big).
$$

Let us set the Hilbert spaces $ \hat{\mathfrak{H}}, \hat{\mathfrak{G}}$, 
$\mathfrak{g}, \mathfrak{g}^1,  \mathfrak{g}_0^1$ defined by

\begin{align}
\hat{\mathfrak{H}}&=\rho \mathfrak{H}=L^2\Big(\mathfrak{R}, \frac{1}{\rho}d\bm{x}; \mathbb{C}^3\Big)  \nonumber \\
&\|\hat{\bm{u}}\|_{\hat{\mathfrak{H}}}^2=\int_{\mathfrak{R}}\|\hat{\bm{u}}\|^2\frac{1}{\rho}d\bm{x}, \\
\hat{\mathfrak{G}}&=\Big\{ \hat{\bm{u}} \in \hat{\mathfrak{H}} \quad \Big| \mathrm{div}\hat{\bm{u}} \in
L^2\Big(\mathfrak{R}, \frac{\Gamma P}{\rho^2}d\bm{x}; \mathbb{C}\Big) \Big\}, \nonumber \\
&\|\hat{\bm{u}}\|_{\hat{\mathfrak{G}}}^2=\|\hat{\bm{u}}\|_{L^2(\frac{1}{\rho})}^2+
\|\mathrm{div}\hat{\bm{u}}\|_{L^2(\frac{\Gamma P}{\rho^2})}^2, \\
\mathfrak{g}&=L^2\Big(\mathfrak{R}, \frac{\Gamma P}{\rho^2}d\bm{x}; \mathbb{C} \Big), \nonumber \\
&\|g\|_{\mathfrak{g}}^2= \int_{\mathfrak{R}}|g|^2\frac{\Gamma P}{\rho^2}d\bm{x}, \\
\mathfrak{g}^1&=\Big\{ g \in L^2\Big(\mathfrak{R}, \frac{\Gamma P}{\rho^2}d\bm{x}; \mathbb{C} \Big) \quad \Big| \nonumber \\
&\mathrm{grad}\Big(\frac{\Gamma P}{\rho^2}g \Big) \in L^2(\mathfrak{R}, \rho d\bm{x}; \mathbb{C}^3),  \int_{\mathfrak{R}}gd \bm{x} =0 \Big\}, \nonumber \\
&\|g\|_{\mathfrak{g}^1}^2=\|g\|_{L^2(\frac{\Gamma P}{\rho^2})}^2 +
\Big\|\mathrm{grad}\Big(\frac{\Gamma P}{\rho^2}g\Big)\Big\|_{L^2(\rho)}^2, \\
\mathfrak{g}_0^1&=\mbox{the closure of}\quad C_0^{\infty}(\mathfrak{R}; \mathbb{C})\quad\mbox{in}\quad \mathfrak{g}^1.
\end{align}
The operators $\hat{\bm{L}}$ in $\hat{\mathfrak{G}}$, $\hat{\bm{M}}$ from $\mathsf{D}(\hat{\bm{M}} (\subset \mathfrak{g} ) $ into $\hat{\mathfrak{G}}$, and $\bm{N}$ in $\mathfrak{g}$ are defined by
\begin{align}
\mathsf{D}(\hat{\bm{L}}&=\{ \hat{\bm{u}} \in \hat{\mathfrak{G}}_0 \quad\Big| \hat{\mathcal{L}}\hat{\bm{u}} \in \hat{\mathfrak{G}} \} \nonumber \\
\hat{\bm{L}}&=\hat{\mathcal{L}}\hat{\bm{u}}=\rho \mathcal{L}\Big(\frac{\hat{\bm{u}}}{\rho}\Big) 
, \\
\mathsf{D}(\hat{\bm{M}})&=\{ g \in \mathfrak{g}_0^1  \quad\Big| \hat{\mathcal{M}}g \in \hat{\mathfrak{G}} \} \nonumber \\
\hat{\bm{M}}g&=\hat{\mathcal{M}}g
=\rho\mathrm{grad}\Big(
-\frac{\Gamma P}{\rho} g +4\pi\mathsf{G}[g] \Big), \\
\mathsf{D}(\bm{N})&=\{ g \in \mathfrak{g}_0^1 \quad\Big| \mathcal{N}g \in \mathfrak{g} \} \nonumber \\
\bm{N}g&=\mathcal{N}g=\mathrm{div}\hat{\bm{M}}g \nonumber \\
&=\mathrm{div}\Big(
\rho\mathrm{grad}\Big(
-\frac{\Gamma P}{\rho} g +4\pi\mathsf{G}[g] \Big),.
\end{align}

Then
$$\bm{N}g=\mathrm{div}(\hat{\bm{M}} g)\quad\mbox{for}\quad g \in \mathsf{D}(\bm{N})
$$ and
$$\bm{L}\bm{u}=\frac{1}{\rho}\hat{\bm{M}}\mathrm{div}\hat{\bm{u}}=
\frac{1}{\rho}\hat{\bm{L}}(\rho \bm{u}) \quad\mbox{for}\quad \bm{u} \in 
C_0^{\infty}(\mathfrak{R}; \mathbb{R}^3).
$$

It follows that $\bm{N}$ is a self-adjoint operator bounded from below in $\mathfrak{g}$. In fact, $\mathcal{N}$ is symmetric and we have
$$(\mathcal{N}g|g)_{\mathfrak{g}}\geq 
(1-\tau)\|g\|_{\mathfrak{g}^1}^2
-C_{\tau}\|g\|_{\mathfrak{g}}^2, \quad \forall g \in C_0^{\infty}(\mathfrak{R}; \mathbb{C}),
$$
where
$0<\tau <1$,
$$C_{\tau}=\frac{1}{4\tau}(4\pi\mathsf{G}\kappa)^2\frac{\rho_O}{\nu}+\tau
$$
with
\begin{align*}
\kappa&=\sup\Big\{\|\mathcal{K}[g]\|_{L^{\infty}(\mathfrak{R})}\quad \Big|\quad \|g\|_{L^2(\mathfrak{R})}=1 \Big\}
<+\infty, \\
\nu&=\inf_{\mathfrak{R}}\frac{\Gamma P}{\rho^2} >0.
\end{align*}
On the other hand, the imbedding of $\mathfrak{g}_0^1$ into $\mathfrak{g}$ is compact. (See \cite[Proposition 4]{JJTM2020}.) Consequently the resovent $(\lambda -\bm{N})^{-1}, \lambda > C_{\tau},$ is a compact operator in $\mathfrak{g}$. Hence

{\it The spectrum $ \mbox{\Pisymbol{psy}{83}}(\bm{N})$ consists of a sequence of eigenvalues $\lambda_n$, $ -C_{\tau} <\lambda_1 < \cdots < \lambda_n<\lambda_{n+1}< \cdots \rightarrow +\infty$,
of finite multiplicities. }

We can claim

{\it If $\lambda \in  \mbox{\Pisymbol{psy}{82}}(\bm{N}), \lambda\not=0$, then
$\lambda \in  \mbox{\Pisymbol{psy}{82}}(\hat{\bm{L}})$. }

Proof. In fact the resolvent $(\lambda-\hat{\bm{L}})^{-1}$ is given by
$$
(\lambda-\hat{\bm{L}})^{-1}\hat{\bm{f}}=\frac{1}{\lambda}
\Big(\hat{\bm{f}}+
\hat{\bm{M}}(\lambda-\bm{N})^{-1}\mathrm{div}\hat{\bm{f}} \Big)
$$
for $\hat{\bm{f}} \in \hat{\mathfrak{G}}$. Here $\hat{\bm{M}}(\lambda-\bm{N})^{-1}$ turns out to be a bounded operator from $\mathfrak{g}$ into $\hat{\mathfrak{G}}$. $\square$

{\it $\hat{\bm{L}}$ is a self-adjoint operator in $\hat{\mathfrak{G}}$. }

Proof. Since it can be shown that $\hat{\bm{L}}$ is symmetric and closable, \cite[Th. 3.16]{Kato} can be applied. $\square$.

{\it If $\lambda \in  \mbox{\Pisymbol{psy}{83}}(\bm{N}), \lambda\not=0$, then $\lambda$ is an eigenvalue of
$\hat{\bm{L}}$ of finite multiplicity. }

Proof. Let $\bm{N}g=\lambda g, \|g\|_{\mathfrak{g}}=1$. Then
$\hat{\bm{u}}:=\frac{1}{\lambda}\hat{\bm{M}}g \in \hat{\mathfrak{G}}$, therefore
$\lambda\mathrm{div}\hat{\bm{u}}=\bm{N}g=\lambda g$. Since $\lambda\not=0$, we have $\mathrm{div}\hat{\bm{u}}=g$, therefore, $
\lambda\hat{\bm{u}}=\hat{\bm{M}}g=\hat{\bm{M}}\mathrm{div}\hat{\bm{u}}=\hat{\bm{L}}\hat{\bm{u}}$. Since
$$\hat{\bm{u}}\in \mathsf{N}(\lambda-\hat{\bm{L}})
\quad \Leftrightarrow \quad
\hat{\bm{u}}=\frac{1}{\lambda}\hat{\bm{M}}g\quad\mbox{for}\quad
\exists g \in \mathsf{N}(\lambda-\bm{N}),
$$
we have
$$\mathrm{dim}.\mathsf{N}(\lambda-\hat{\bm{L}})\leq
\mathrm{dim}. \mathsf{N}(\lambda-\bm{N}).
$$ $\square$.\\

Summing up, we can claim
{\it The spectrum of the self-adjoint operator $\hat{\bm{L}}$ coincides with
$\{0\}\cup  \mbox{\Pisymbol{psy}{83}}(\bm{N})$ and
$\mathrm{dim}\mathsf{N}(\lambda-\hat{\bm{L}}) <\infty$ for $\lambda \in  \mbox{\Pisymbol{psy}{83}}(\hat{\bm{L}})\setminus \{0\}$. }

This completes the proof of Theorem \ref{Th.X}. \\

It is obvious that $\mbox{\Pisymbol{psy}{83}}_p(\bm{L})=\mbox{\Pisymbol{psy}{83}}_p(\bm{L}^G)$, that is, $\lambda \in \mathbb{C}$ is an eigenvalue of $\bm{L}$ if and only if 
$\lambda $ is an eigenvalue of $\bm{L}^G$. However the relation between 
$\mbox{\Pisymbol{psy}{83}}(\bm{L})$ and $\mbox{\Pisymbol{psy}{83}}(\bm{L}^G)$ is not clear. Namely is 
$\mbox{\Pisymbol{psy}{83}}(\bm{L}) \subset \mbox{\Pisymbol{psy}{83}}(\bm{L}^G)$ the case? In other words, if $\lambda \in 
\mbox{\Pisymbol{psy}{82}}(\bm{L}^G)$, then $\lambda \not\in \mbox{\Pisymbol{psy}{83}}_p(\bm{L}) $
$\Big(=\mbox{\Pisymbol{psy}{83}}_p(\bm{L}^G) \Big)$ and so $(\lambda-\bm{L})^{-1}$ exists as an operator in $\mathfrak{H}$; Although $(\lambda-\bm{L}^G)^{-1} \in \mathcal{B}(\mathfrak{G})$, we do not know whether this $(\lambda-\bm{L})^{-1} \Big(\supset (\lambda-\bm{L}^G)^{-1} \Big)$ is bounded in $\mathfrak{H}$ or not. This is the {\bf Question} \ref{Q.1}.\\

\textbullet\  When $\Omega=0,\bm{B}=O$, 
$$
\mbox{\Pisymbol{psy}{83}}(\mathfrak{L})=\Big\{ \lambda \in \mathbb{C}\ \Big|\  -\lambda^2 \in \mbox{\Pisymbol{psy}{83}}(\bm{L}) \Big\},
$$
where $\mbox{\Pisymbol{psy}{83}}(\bm{L})$ is the spectrum of the operator $\bm{L}$ in $\mathfrak{H}$.
We know
$\mbox{\Pisymbol{psy}{83}}(\bm{L}) \subset [-\mu_*, +\infty[$. 
 Therefore
 we see
 $$\mbox{\Pisymbol{psy}{83}}(\mathfrak{L}) \subset
 \mathrm{i} \mathbb{R} \cup
 [-\sqrt{\mu_*}, \sqrt{\mu_*}].
 $$

 Generaly speaking, at least
 we can claim
 
 \begin{Proposition}\label{roughR}
 It holds
 \begin{equation}
 ]-\infty, -\sqrt{\mu_*}[\quad\cup\quad
 ]\sqrt{\mu_*}, +\infty[ \quad \subset \quad \mbox{\Pisymbol{psy}{82}}(\mathfrak{L})=\mbox{\Pisymbol{psy}{82}}(-\bm{A}).
 \end{equation}
  \end{Proposition}
 
 For proof see the proof of Propositions {\bf A1},  {\bf A2} of {\bf Appendix A}.\\
 
 \begin{Proposition}\label{Rbelow}
 For $\lambda \in \mbox{\Pisymbol{psy}{82}}(\mathfrak{L})$, we have
 \begin{equation}
 |\|\mathfrak{L}(\lambda)^{-1}\||_{\mathcal{B}(\mathfrak{H})} \geq
 \frac{1}{d(2|\lambda|+\beta^{\star} +d)}
 \end{equation}
 where
 $d:=\mathrm{dist.}(\lambda, \mbox{\Pisymbol{psy}{83}}(\mathfrak{L}))$.
 \end{Proposition}
 
 Proof. Let $\lambda \in \mbox{\Pisymbol{psy}{82}}(\mathfrak{L})$. Then, for $\Delta\lambda \in \mathbb{C}$, the operator
$$\mathfrak{L}(\lambda+\Delta\lambda)=
\mathfrak{L}(\lambda)\Big[I+
\mathfrak{L}(\lambda)^{-1}\Delta\lambda(2\lambda+\bm{B}+\Delta\lambda)\Big]$$
admits the bounded inverse in $\mathcal{B}(\mathfrak{H})$ and $\lambda+\Delta\lambda \in
\mbox{\Pisymbol{psy}{82}}(\mathfrak{L})$, if
$$|\|\mathfrak{L}(\lambda)^{-1}\Delta\lambda(2\lambda+\bm{B}+\Delta\lambda)\||_{\mathcal{B}(\mathfrak{H})} <1.$$
For this inequality, it is sufficient that
$$|\|\mathfrak{L}(\lambda)^{-1}\||_{\mathcal{B}(\mathfrak{H})}\cdot|\Delta\lambda|\cdot(2|\lambda|
+\beta^{\star}+|\Delta\lambda| ) <1.$$
In other words, if  $\lambda+\Delta\lambda \in \mbox{\Pisymbol{psy}{83}}(\mathfrak{L})$, then it should hold
$$|\|\mathfrak{L}(\lambda)^{-1}\||_{\mathcal{B}(\mathfrak{H})}\cdot|\Delta\lambda|\cdot(2|\lambda|
+\beta^{\star}+|\Delta\lambda| ) \geq 1.
$$

If $d <+\infty$, then there is a sequence $\lambda+(\Delta\lambda)_n
\in \mbox{\Pisymbol{psy}{83}}(\mathfrak{L}) $
 such that $|(\Delta\lambda)_n| \rightarrow d$, and the assertion follows. $\square$\\
 
 \begin{Theorem}(\cite[Theorem 1]{DysonS})\label{Th.4}
 $\mbox{\Pisymbol{psy}{83}}(\mathfrak{L}) $ is a subset of 
 \begin{equation}
 S:=\mathrm{i}\mathbb{R}\cup \Big\{ \lambda \in \mathbb{C}\  \Big|\  |\lambda|\leq \sqrt{\mu_*}, \quad |\mathfrak{Im}[\lambda]|\leq \frac{\beta^{\star}}{2} \Big\}.
 \end{equation}
 \end{Theorem}
 
 Proof. First we consider $\lambda_{\infty} \in \partial \mbox{\Pisymbol{psy}{83}}(\mathfrak{L}) $. 
 
  Let us take a sequance $(\lambda_n)_n$ such that $\lambda_n \in \mbox{\Pisymbol{psy}{82}}(\mathfrak{L})$ and
$\lambda_n \rightarrow \lambda_{\infty}$ as $ n \rightarrow \infty$. By Proposition \ref{Rbelow} 
we have $|\|(\mathfrak{L}(\lambda_n)^{-1}\||_{\mathcal{B}(\mathfrak{H})} \rightarrow +\infty$, therefore
there are $\bm{f}_n \in \mathfrak{H}$ such that $\|\bm{f}_n\|_{\mathfrak{H}}=1$ and $\|\mathfrak{L}(\lambda_n)^{-1}\bm{f}_n\|_{\mathfrak{H}}
 \rightarrow +\infty$ as $n \rightarrow \infty$. Put $\bm{u}_n=\mathfrak{L}(\lambda_n)^{-1}\bm{f}_n (\in \mathsf{D}(\bm{L}))$ and $\bm{u}_n=\bm{u}_n/\|\bm{u}_n\|_{\mathfrak{H}}$. Then $\|\bm{u}_n\|_{\mathfrak{H}}=1$ and $$
 \Big|(\mathfrak{L}(\lambda_n)\bm{u}_n|\bm{u}_n)_{\mathfrak{H}}\Big|
 =\Big|\frac{1}{\|\bm{u}_n\|_{\mathfrak{H}}^2}(\bm{f}_n|\bm{u}_n)_{\mathfrak{H}}\Big|
 \leq \frac{1}{\|\bm{u}_n\|_{\mathfrak{H}}} \rightarrow 0.$$
 But we see
 $$
 (\mathfrak{L}(\lambda_n)\bm{u}_n|\bm{u}_n)_{\mathfrak{H}}=\lambda_n^2+\mathrm{i}\lambda_nb_n+c_n,$$
 where
 $$b_n:=-\mathrm{i}(\bm{B}\bm{u}_n|\bm{u}_n)_{\mathfrak{H}},
 \quad c_n=(\bm{L}\bm{u}_n|\bm{u}_n)_{\mathfrak{H}}.$$
 Here $b_n, c_n \in \mathbb{R}$  and 
 $|b_n|\leq \beta^{\star}, c_n \geq -\mu_*$.
 Hence, by taking a subsequence if necessary, we can suppose that $b_n$ tends to a limit $b_{\infty}$ such that $b_{\infty} \in \mathbb{R}, |b_{\infty}|\leq \beta^{\star}$. Put
 $c_{\infty}:=-\lambda_{\infty}^2-\mathrm{i} \lambda_{\infty}b_{\infty}$. Then we see 
 $c_n \rightarrow c_{\infty}$. Hence $c_{\infty}\in \mathbb{R}$  and $ c_{\infty} \geq -\mu_*$,
 and $\lambda_{\infty}$ turns out to enjoy the quadratic equation
 $$\lambda_{\infty}^2+\mathrm{i} b_{\infty}\lambda_{\infty}+c_{\infty}=0.$$
 Consequently, 
 $$\lambda_{\infty}=\mathrm{i}\Big[
 -\frac{b_{\infty}}{2}\pm\sqrt{\frac{b_{\infty}^2}{4}+c_\infty}\Big].
 $$
 If $\displaystyle\frac{b_{\infty}^2}{4}+c_{\infty} \geq 0$, then $\lambda_{\infty} \in \mathrm{i}\mathbb{R}$. If $\displaystyle \frac{b_{\infty}^2}{4}+c_{\infty} \leq 0$, then 
 $|\lambda_{\infty}|^2 =-c_{\infty} \leq \mu_*$ and
 $\displaystyle \Big|\mathfrak{Im}[\lambda_{\infty}]\Big|= \Big|-\frac{b_{\infty}}{2}\Big|\leq \frac{\beta^{\star}}{2}$.  Hence $\lambda_{\infty} \in S$.

 Let us consider $\lambda_0 \in \mbox{\Pisymbol{psy}{83}}(\mathfrak{L})$. We claim $\lambda_0 \in S$. Suppose
 $\lambda_0 \not\in S$. By the symmetricity, we suppose $\mathfrak{Re}[\lambda_0] >0$.
 Then there would exist a curve $\Gamma: t \in [0,1] \mapsto \lambda(t) \in\mathbb{C} \setminus S$ such that $\lambda(0)=\lambda_0$ and
 $\lambda(1) \in ] \sqrt{\mu}_*,
  +\infty[$. Recall $\lambda(1) \in \mbox{\Pisymbol{psy}{82}}(\mathfrak{L})$ by Proposition
\ref{roughR}. The time
 $$\bar{t}=\sup\Big\{ t \in [0,1]\quad\Big|\quad \lambda(t) \in \mbox{\Pisymbol{psy}{83}}(\mathfrak{L}) \Big\}
 $$
 would enjoy
 $\bar{t}\in [0,1[, \lambda(\bar{t}) \in \partial\mbox{\Pisymbol{psy}{83}}(\mathfrak{L}), \lambda(\bar{t})\not\in S$,
 a contradiction to $\partial\mbox{\Pisymbol{psy}{83}}(\mathfrak{L}) \subset S$. Therefore we can claim $\lambda_0 \in S$.
 $\square$\\

We want to clarify the structure of $\mbox{\Pisymbol{psy}{83}}(\mathfrak{L})$ more concretely. In \cite[p.405]{DysonS} J. Dyson and B. F. Schutz say that the following assumption about $\mbox{\Pisymbol{psy}{83}}(\mathfrak{L})$ is very reasonable and provable:

\begin{quote}
{\bf Dyson-Schutz' Assumption}: {\it
Let $S=S_{(i)}\cup S_{(ii)}\cup S_{(iii)} $, where
\begin{align*}
&S_{(i)}=\{ \lambda \in \mathrm{i}\mathbb{R} \  |\  |\mathfrak{Im}[\lambda]| >\frac{\beta^{\star}}{2}\} \\
&S_{(ii)}=\{ \lambda \in S \  |\  \lambda \not\in \mathrm{i}\mathbb{R} \}, \\
&S_{(iii)}=\{ \lambda \in \mathrm{i}\mathbb{R} \  |\  |\mathfrak{Im}[\lambda]| \leq\frac{\beta^{\star}}{2}\} .
\end{align*}
Then 
\begin{align*}
\mbox{(i)}\quad & \mbox{\Pisymbol{psy}{83}}(\mathfrak{L}) \cap S_{(i)} \subset \mbox{\Pisymbol{psy}{83}}_p(\mathfrak{L}), \\
\mbox{(ii)}\quad & \mbox{\Pisymbol{psy}{83}}(\mathfrak{L}) \cap S_{(ii)}\subset \mbox{\Pisymbol{psy}{83}}_p(\mathfrak{L}),\quad \mbox{and}\quad
\partial(\mbox{\Pisymbol{psy}{83}}_p(\mathfrak{L})\cap S_{(ii)} ) \cap (\mathrm{i}\mathbb{R})=\emptyset, \\
\mbox{(iii)}\quad & \Big(\mbox{\Pisymbol{psy}{83}}(\mathfrak{L}) \cap S_{(iii)}\Big) \setminus \mbox{\Pisymbol{psy}{83}}_p(\mathfrak{L}) \not=\emptyset.
\end{align*}
}
\end{quote}

However it seems that this assumption has not yet been justified mathematically.
Note that this assumption imlies that, when $\Omega=0$, it holds that
$\mbox{\Pisymbol{psy}{83}}(\mathfrak{L})=\mbox{\Pisymbol{psy}{83}}_p(\mathfrak{L})$ and that, when $\Omega\not=0$, it holds that $\mbox{\Pisymbol{psy}{83}}(\mathfrak{L})\setminus \mbox{\Pisymbol{psy}{83}}_p(\mathfrak{L}) \not=\emptyset$.\\

 \section{Stability problem --New concepts of stability}
 
 Let us discuss on stability of solutions of {\bf ELASO}, supposing \\
 
 {\bf  Assumption} : {\it The equation of state enjoys {\bf Assumption} \ref{Ass.0} and the background $(\rho_b, S_b, \bm{v}_b)$ is a barotropic admissible stationary solution.}\\
  
  D. Lynden-Bell and J. P. Ostriker, \cite[p.301, line 18]{LyndenBO}, say:
\begin{quote}
Equation (36)  shows that the system is stable if $c$ is positive for each eigen $\bm{\xi}$ [ read $\bm{u}^0$ ] . This assured if $\mathbf{C}$ [ read $\mbox{\boldmath$L$}$] is positive definite. Thus:

A sufficient condition for stability is that $\mathbf{C}$ [read $\mbox{\boldmath$L$}$ ] is positive definite. This is {\it the } condition for secular stability.
\end{quote}

We should be careful to understand the meaning of the saying `the system is stable' and `$\bm{L}$ is positve definite'. As for `stability' 
 C. Hunter \cite{Hunter} says:

\begin{quote}
A general system is said to be ordinarily or dynamically unstable if the amplitude of some mode grows exponentially in time, but ordinarily stable if every mode is oscillatory in time. An ordinarily stable system can be said to be secularly unstable if small additional dissipative forces can cause some perturbation to grow. Otherwise, the system is said to be secularly stable.
\end{quote}

So  `(ordinary) stability' means the condition:\\

{\bf (ST.1)}\  {\it For any eigenvalue $\lambda \in \mbox{\Pisymbol{psy}{83}}_{\mathrm{p}}(\mathfrak{L})$ of {\bf ELASO} \eqref{5.1} and its associated solution $\bm{u}(t,\bm{x})=e^{\lambda t}\bm{u}^0(\bm{x}), \bm{u}^0 \in \mathsf{D}(\bm{L}), \bm{u}^0 \not= \bm{0}, \mathfrak{L}(\lambda)\bm{u}^0=\bm{0}$, it holds}
\begin{equation}
\|{U}(t,\cdot)\|_{\mathfrak{G}\times\mathfrak{H}} \leq C\quad \mbox{for}\quad \forall t \in [0,+\infty[,
\end{equation}\\

Here and hereafter we use the notation
\begin{equation}
{U}=
\begin{bmatrix}
\bm{u} \\
\\
\dot{\bm{u}}
\end{bmatrix},
\quad \dot{\bm{u}}=\frac{\partial \bm{u}}{\partial t}.
\end{equation}

Since
$$ \|\bm{u}(t,\cdot)\|_{\mathfrak{H}}=e^{\mathfrak{Re}[\lambda]t}\|\bm{u}^0\|_{\mathfrak{H}},
\quad \|\dot{\bm{u}}(t,\cdot)\|_{\mathfrak{H}}=|\lambda|e^{\mathfrak{Re}[\lambda]t}\|\bm{u}^0\|_{\mathfrak{H}}, $$
the following condition is necessary and sufficient for the stability in the sense of {\bf (ST.1)}: \\

{\bf (*) }\  {\it For any eigenvalue $\lambda \in \mbox{\Pisymbol{psy}{83}}_{\mathrm{p}}(\mathfrak{L})$ it hlds that $\mathfrak{Re}[\lambda] =0$, that is,
$\mbox{\Pisymbol{psy}{83}}_{\mathrm{p}}(\mathfrak{L}) \subset \mathrm{i}\mathbb{R}$. }\\

On the other hand, any $\lambda \in \mbox{\Pisymbol{psy}{83}}_{\mathrm{p}}(\mathfrak{L})$ satisfies the quadratic equation
 \begin{equation}
  \lambda^2+\mathrm{i} b \lambda + c =0
 \end{equation}
 where the coefficents $b,c$ are given by
 \begin{equation}
  b=-\mathrm{i} (\bm{B}\bm{u}^0|\bm{u}^0)_{\mathfrak{H}},
 \quad c=(\bm{L}\bm{u}^0|\bm{u}^0)_{\mathfrak{H}},
 \end{equation}
 $\bm{u}^0$ being an associated eigenvector such that $\|\bm{u}^0\|_{\mathfrak{H}}=1$,
  and $b,c$ being real numbers,
  therefore
  \begin{equation}
  \lambda=\mathrm{i}\Big[\frac{b}{2}\pm \sqrt{ \frac{b^2}{4}+c}\Big].
  \end{equation}
  Hence
  \begin{equation}
  \mathfrak{Re}[\lambda]=0 \quad\Leftrightarrow\quad
  \frac{b^2}{4}+c \geq 0
  \quad \Leftarrow c \geq 0.
  \end{equation}
  
  Therefore, considering the condition\\
  
  {\bf (PD.1)}: \  {\it It holds}
  \begin{equation}
  (\bm{L}\bm{u}|\bm{u})_{\mathfrak{H}} \geq 0\quad \forall \bm{u} \in \mathsf{D}(\bm{L}),
  \end{equation}
 \\
 
\noindent  we can claim
 
 \begin{Proposition}
 If $\bm{L}$ is `positive definite' in the sense of {\bf (PD.1)}, then the background $(\rho, S, \bm{v})=(\rho_b,S_b, \bm{v}_b)$ is `stable' in the sense of {\bf (ST.1)}.
 \end{Proposition}
 
 Thus the saying of \cite{LyndenBO} can be interpreted as this Proposition.\\
 
 Note that {\bf (PD.1)} implies
 $$ \inf_{\bm{u} \in \mathsf{D}(\bm{L})}\frac{(\bm{L}\bm{u}|\bm{u})_{\mathfrak{H}}}{\|\bm{u}\|_{\mathfrak{H}}^2} \geq 0,
 $$
 so that $$\mu_*=\Big(\inf_{\bm{u} \in \mathsf{D}(\bm{L})}\frac{(\bm{L}\bm{u}|\bm{u})_{\mathfrak{H}}}{\|\bm{u}\|_{\mathfrak{H}}^2}           \Big)\vee 0=0,$$ and 
 $\mbox{\Pisymbol{psy}{83}}(\mathfrak{L}) \subset \mathrm{i}\mathbb{R}$ by Theorem \ref{Th.4}.\\
 
 The concept of `stability' in the sense of {\bf (ST.1)} is concerned with only `every modes', say, waves associated with eigenvalues. We may consider all solutions 
 $\bm{u}$ of \eqref{5.1} and the wider stability concept:\\
 
 {\bf (ST.2) }: \  {\it For all solution $\bm{u} \in C^2([0,+\infty[, \mathfrak{H})\cap
 C^1([0,+\infty[; \mathfrak{G}_0) \cap C([0,+\infty[; \mathsf{D}(\bm{L}))$ of
{\bf ELASO}  \eqref{5.1} it holds}
 \begin{equation}
 \|{U}(t,\cdot)\|_{\mathfrak{G}\times \mathfrak{H}}\leq C\quad \mbox{for}\quad \forall t \in [0,+\infty[.
 \end{equation}
 \\
 
 There is a gap between {\bf (ST.1)} and {\bf (ST.2)}. Of course {\bf (ST.2)} $\Rightarrow$
 {\bf (ST.1)}, but the inverse is not obvious, since we do not have an answer to the 
 {\bf Question 1}

 Let us note that, if $\Omega=\omega_b=0$ and the background is spherically symmetric and isentropic, then it holds that $\mbox{\Pisymbol{psy}{83}}(\bm{L}^{G})=\mbox{\Pisymbol{psy}{83}}_{\mathrm{p}}(\bm{L}^{G})$, where $\bm{L}^{G}$ is the Friedrichs extention of
 $\mathcal{L} \restriction C_0^{\infty} $ in the Hilbert space
 $\displaystyle \mathfrak{G}=\Big\{ \bm{u} \in \mathfrak{H}\ \Big|\ \mathrm{div}(\rho\bm{u}) \in L^2\Big(\frac{\Gamma P}{\rho^2}d\bm{x}\Big)\Big\}$. See 
 Theorem \ref{Th.X}. But, even in this non-rotating case, we do not know whether $\mbox{\Pisymbol{psy}{83}}(\bm{L})=\mbox{\Pisymbol{psy}{83}}_{\mathrm{p}}(\bm{L})$ or not, $\bm{L}$ being the Friedrics extension in $\mathfrak{H}$. Moreover, when $\Omega\not=0$, 
 $\bm{B}^G=\bm{B}\restriction \mathfrak{G}$
   is not a bounded linear operator on $\mathfrak{G}$ so that the application of Hille-Yosida theory for the basic existence proof in Section 4 does not work if we take $\mathfrak{G}$ as the basic space instead of $\mathfrak{H}$. \\
 
 Also let us note that, if $\Omega=0, \bm{B}=\bf{0}$, when {\bf ELASO} \eqref{5.1} reduces
 \begin{equation}
 \frac{\partial^2\bm{u}}{\partial t^2}+\bm{L}\bm{u} =0, \label{*6.9}
 \end{equation}
 { \bf (PD.1)} implies that
 \begin{equation}
 \bm{u}=\int_0^{+\infty}\cos(\sqrt{\sigma}t)dE(\sigma)\bm{a} +
 \int_0^{+\infty}\sin(\sqrt{\sigma}t)dE(\sigma)\bm{b} 
 \end{equation}
 solves \eqref{*6.9} with $\bm{u}^0=\bm{a}, \bm{v}^0=\sqrt{\bm{L}}\bm{b}$,
 provided that $\bm{a}, \bm{b} \in \mathsf{D}(\bm{L})$. Here
  $(E(\sigma))_{\sigma \in\mathbb{R}}$ is the spectral decomposition of the self-adjoint operator $\bm{L}$,
 which enjoys $E(-0)=O$ thanks to {\bf (PD.1)}. In this case we have
 $$ \|{U}\|_{\mathfrak{G}\times \mathfrak{H}}
 \leq C\Big(\|\bm{a}\|_{\mathsf{D}(\sqrt{\bm{L}})}
 +\|\bm{b}\|_{\mathsf{D}(\sqrt{\bm{L}})}\Big).
 $$  
 Here
 $$ \sqrt{\bm{L}}=\int_0^{+\infty}\sqrt{\sigma}dE(\sigma) $$
 so that $(\sqrt{\bm{L}})^2=\bm{L}$ and
 $$ \|\bm{u}\|_{\mathsf{D}(\sqrt{\bm{L}})}
 =\Big[ \|\bm{u}\|_{\mathfrak{H}}^2+ \|\sqrt{\bm{L}}\bm{u}\|_{\mathfrak{H}}^2 \Big]^{\frac{1}{2}}. $$\\

 {\bf \textbullet\  Proposal of new concepts of stability}\\

Now let us introduce more general concepts of stability. Namely, let
 $\mathfrak{N}$ be a semi-norm on $\mathsf{D}(\bm{L}) \times \mathfrak{G}_0$, and consider the condition:\\
 
  {\bf (ST.3)}:\  {\it For any solution $\bm{u}$ of {\bf ELASO} \eqref{5.1} it holds

 \begin{equation}
 \mathfrak{N}({U}(t,\cdot)) \leq C \quad\mbox{for}\quad \forall t\in [0,+\infty[.
 \end{equation}
 Here $C$ is a bounce deending on ${U}^0={U}(0,\cdot)=(\bm{u}^0,\bm{v}^0)^{\top}$.
 }\\
 
 Recall
 $$
 \mathfrak{N}_0({U}(t,\cdot)) \leq e^{\Lambda t}\mathfrak{N}_0({U}^0)
 $$
 for
 $$\mathfrak{N}_0(U):=\|U\|_{\mathfrak{E}}=\Big[\|\bm{u}\|_{\breve{\mathfrak{G}}}^2+\|\dot{\bm{u}}\|_{\mathfrak{H}}^2\Big]^{\frac{1}{2}}
 \quad\mbox{for}\quad
 U=\begin{bmatrix}
 \bm{u} \\
 \\
 \dot{\bm{u}}
 \end{bmatrix}
 $$
 Here $\Lambda =1+\mu_* \geq 1$ so this does not give the stability {\bf (ST.3)}. We are keeping in mind the following situation. 
 
 Suppose a seminorm $\mathfrak{n}$ on $\mathsf{D}(\bm{L})$ satisfy\\
 
 {\bf (PD.2)}:\  {\it There is a positive number ${\delta}$ such that}
 \begin{equation}
 (\bm{L}\bm{u}|\bm{u}) \geq {\delta} \mathfrak{n}(\bm{u})^2
 \quad \forall \bm{u} \in \mathsf{D}(\bm{L}).
 \end{equation}\\
 
 Put
 \begin{equation}
 \mathfrak{N}(U):=\Big[ {\delta} \mathfrak{n}(\bm{u})^2+\|\dot{\bm{u}}\|_{\mathfrak{H}}^2 \Big]^{\frac{1}{2}}.
 \end{equation}
 
 Then {\bf (ST.3)} holds with
 $$ C=\Big[\|\bm{v}^0\|_{\mathfrak{H}}^2+(\bm{L}\bm{u}^0|\bm{u}^0)_{\mathfrak{H}}\Big]^{\frac{1}{2}}.
 $$
 
 In fact, multipling \eqref{5.1} by $\dot{\bm{u}}$ by the $\mathfrak{H}$-inner product, taking the real part, using $\mathfrak{Re}[(B\dot{\bm{u}}|\dot{\bm{u}})_{\mathfrak{H}}]=0$, we have
 $$\frac{d}{dt}\Big[ \|\dot{\bm{u}}\|_{\mathfrak{H}}^2+(\bm{L}\bm{u}|\bm{u})_{\mathfrak{H}}\Big]=0.
 $$
 
 Of course {\bf (PD.2)} implies {\bf (PD.1)} but is much stronger concept of `positive definiteness'.\\
 
 Now we are going to find a seminorm $\mathfrak{n}$ which enjoys {\bf (PD.2)}, under the following situation:
 
 \begin{Situation} \label{As.3}
1)  The equation of state is \eqref{EOS.g}: $P=\rho^{\gamma}\exp(S/\mathsf{C}_V)$;

2) $\Omega =\omega_b=0$ and the background $(\rho_b, S_b, \bm{0})$ is spherically symmetric, and is isentropic, that is, $S_b=S_0$, $S_0$ being a constant; 

3) $\displaystyle \frac{4}{3}<\gamma <2$.
 \end{Situation}
 
 Then
 $\displaystyle \varUpsilon=\varUpsilon_O\theta\big(\frac{r}{\mathsf{a}}; \frac{1}{\gamma-1}\Big) $  and  
  $\mathfrak{R}=\{ \bm{x} | r=\|\bm{x}\|<R\}$ with $R:=\mathsf{a}\xi_1\Big(\frac{1}{\gamma-1}\Big)$, where $\theta(\cdot, \frac{1}{\gamma-1})$ is the Lane-Emden function and $\xi_1(\frac{1}{\gamma-1})$ is its zero.
    It holds
  $$
  P=\mathsf{A}\rho^{\gamma},\quad
  \varUpsilon=\frac{\mathsf{A}\gamma}{\gamma-1}\rho^{\gamma-1},
  \quad
  \frac{\gamma P}{\rho^2}=\frac{d\varUpsilon}{d\rho}=\mathsf{A}\rho^{-(2-\gamma)},
  $$ and
 
 \begin{equation}
 (\bm{L}\bm{u}|\bm{u})_{\mathfrak{H}} =\int_{\mathfrak{R}}\Big(
 \frac{d\Upsilon}{d\rho}|g|^2-4\pi\mathsf{G}\mathcal{K}[g]g^*\Big)d\bm{x},
 \end{equation}
 with $g=\mathrm{div}(\rho \bm{u})$. ( Recall Remark \ref{Rem.2}. )
 
 Denote
 \begin{equation}
 \mathfrak{g}_{L}:= \Big\{ g \in L^2\Big(\frac{d\Upsilon}{d\rho}d\bm{x}\Big) \Big| \exists \bm{u} \in \mathsf{D}(\bm{L}) : g=\mathrm{div}(\rho \bm{u}) \Big\},
 \end{equation}
 and put
 \begin{equation}
 \Psi=-\mathcal{K}[g]
 \end{equation}
 for $g \in \mathfrak{g}_L$. Then, for $g \in \mathfrak{g}_L$, $g \in C(\mathfrak{R}), \Psi \in C(\mathbb{R}^3)$, and the expansion with respect to
 spherical harmonics
 \begin{subequations}
 \begin{align}
 &g(\bm{x})=\sum _{0\leq \ell, |m|\leq \ell}
 g_{\ell m}(r)Y_{\ell  m}(\vartheta, \phi), \\
& \Psi(\bm{x})=\sum _{0\leq \ell, |m|\leq \ell}
 \Psi_{\ell m}(r)Y_{\ell  m}(\vartheta, \phi)
 \end{align}
 \end{subequations}
 can be used.
 Here
 $$
 \bm{x}=(x^1,x^2, x^3)=(r\sin\vartheta\cos\phi, r\sin\vartheta\sin\phi,r \cos\vartheta)
 $$
 and
 \begin{align}
 &Y_{\ell  m}(\vartheta, \phi)=\sqrt{ \frac{2\ell  +1}{4\pi}\frac{(\ell  -m)!}{(\ell  +m)!} }
 P_{\ell  }^m(\cos\vartheta )e^{\mathrm{i}m\phi}, \nonumber \\
 &Y_{\ell, -m}=(-1)^mY_{\ell m}^* \qquad(0\leq m\leq \ell),
 \end{align}
 \begin{align}
 &g_{\ell m}(r)=\int_0^{2\pi}\int_0^{\pi}g(\bm{x})Y_{\ell m}(\vartheta,\phi)^*\sin\vartheta d\vartheta d\phi, \\
 &\Psi_{\ell m}(r)=\int_0^{2\pi}\int_0^{\pi}\Psi(\bm{x})Y_{\ell m}(\vartheta,\phi)^*\sin\vartheta d\vartheta d\phi.
 \end{align} 
 See \cite[Lemma 5]{JJTM2020}. Note that
 \begin{equation}
 \triangle^{(\ell)}\Psi_{\ell m}=g_{\ell m},  
 \quad \Psi_{\ell m}=-\mathcal{H}_{\ell}[g_{\ell m}],
\end{equation}
where
\begin{align}
&\triangle^{(\ell)}w=\frac{1}{r^2}\frac{d}{dr}r^2\frac{dw}{dr}-\frac{\ell(\ell+1)}{r^2}w, \\
&\mathcal{H}_{\ell}[\check{g}](r)=
\frac{1}{2\ell +1}\Big[\int_0^r\check{g}(s)\Big(\frac{r}{s}\Big)^{-\ell -1}sds
+\int_r^{+\infty}\check{g}(s)\Big(\frac{r}{s}\Big)^{\ell}sds
\Big],
\end{align}
and
\begin{equation}
\|{\nabla}\Psi\|^2=\sum_{\ell m}\|\nabla^{(\ell)}\Psi_{\ell m}\|^2,
\end{equation}
where
\begin{equation}
\|\nabla^{(\ell)}w\|^2=\Big|\frac{dw}{dr}\Big|^2+\frac{\ell (\ell +1)}{r^2}|w|^2.
 \end{equation}\\
 
 We are observing
 \begin{equation} 
 (\bm{L}\bm{u}|\bm{u})_{\mathfrak{H}}=4\pi\sum_{\ell, m}Q_{\ell m},
 \end{equation}
 where
 \begin{equation}
Q_{\ell m}=\int_0^{+\infty}
 \Big(\frac{d\varUpsilon}{d\rho}|\triangle^{(\ell)}\Psi_{\ell m}|^2
 -4\pi\mathsf{G}\|\nabla^{(\ell )}\Psi_{\ell m}\|^2\Big)r^2dr.
 \end{equation}\\

 First we consider $\ell=0$. Then we have
 \begin{equation}
 Q_{00}=\int_0^R(\mathcal{L}^{\mathsf{ss}}y)y^*\rho r^4dr,
 \end{equation}
 where
 
 \begin{align}
 & y=-\frac{1}{r\rho}\frac{d}{dr}\Psi_{00},
 \quad
 g_{00}=\frac{1}{r^2}\frac{d}{dr}(r^3\rho y), \\
 &\mathcal{L}^{\mathsf{ss}}y=-\frac{1}{\rho r^4}\frac{d}{dr}\Big(\gamma r^4\rho\frac{dy}{dr}\Big)
 -(3\gamma -4)\frac{1}{r}\frac{d\varUpsilon}{dr}y, \\
 &\int_0^R(\mathcal{L}^{\mathsf{ss}}y)y^* \rho r^4dr
 =
 \int_0^R
 \Big(\gamma \Big|\frac{dy}{dr}\Big|^2
 -(3\gamma -4)\frac{1}{r}\frac{d\varUpsilon}{dr}|y|^2\Big)\rho r^4dr.
 \end{align}
 
 Since we are supposing $\frac{4}{3}<\gamma$, we have
 \begin{equation}
 {\delta}_*:=\inf_{0<r<R} -(3\gamma-4)\frac{1}{r}\frac{d\varUpsilon}{dr} >0,
 \end{equation}
 and
 \begin{equation}
 Q_{00} \geq {\delta}_*\int_0^R|y|^2\rho r^4dr ={\delta}_*\int_0^R\frac{1}{\rho}\Big|\frac{d\Psi_{00}}{dr}\Big|^2 r^2dr. 
 \end{equation}
 
 But 
 $$ \frac{d\Psi_{00}}{dr}=\frac{1}{r^2}\int_0^r g_{00}(s)s^2ds =0 \quad\mbox{for}\quad r >R,
 $$
 since
 $g_{00}=0 $ on $]R,+\infty[$ and
 
 $$\int_0^Rg_{00}(r)r^2dr=\frac{1}{\sqrt{4\pi}}\int_{\mathfrak{R}}g(\bm{x})d\bm{x}=0$$
 for $g \in \mathfrak{g}_L$. 
 On the other hand, $$
 \frac{1}{\rho}\Big|\frac{d\Psi_{00}}{dr}\Big|^2 \in L^{\infty}(0,R),$$
  since 
 $$\frac{d\Psi_{00}}{dr}=-\frac{1}{r^2}\int_r^R g_{00}(s)s^2ds \quad \mbox{for}\quad 0<r <R
 $$
 enjoys
 \begin{align*}
 \Big| \frac{d\Psi_{00}}{dr} \Big|
 &\leq
  \frac{1}{r^2}\sqrt{ \int_r^R\frac{d\rho}{d\varUpsilon}s^2ds }
 \sqrt{ \int_0^R\frac{d\varUpsilon}{d\rho}|g_{00}(s)|^2s^2ds } \\
 &
 \leq C (R-r)^{\frac{1}{2(\gamma-1)}} \quad\mbox{for}\quad \frac{R}{2}\leq r <R,
 \end{align*}
 for $\displaystyle \frac{d\rho}{d\varUpsilon} =O((R-r)^{\frac{2-\gamma}{\gamma-1}})$.

 Therefore, under the convention that 
 $\displaystyle \frac{1}{\rho}\Big|\frac{\Psi_{00}}{dr}\Big|^2 $ means $0$ for $r >R$,  we can write
 \begin{equation}
 Q_{00}\geq {\delta}_*\int_0^{+\infty}\frac{1}{\rho}\Big|\frac{d\Psi_{00}}{dr}\Big|^2r^2dr.
 \end{equation}
 
 In other words, if we adopt the decomposition
 \begin{align}
 &u^r(\bm{x})=\sum_{\ell m}u^r_{\ell m}(r)Y_{\ell m}(\vartheta,\phi), \\
 &u^r_{\ell m}(r)=\int_0^{2\pi}\int_0^{\pi}u^r(\bm{x})Y_{\ell m}(\vartheta, \phi)^*\sin\vartheta d\vartheta d\phi,
 \end{align}
 while
 \begin{align}
& \bm{u} (\bm{x})= u^r (\bm{x})\bm{e}_r +u^{\vartheta}(\bm{x})\bm{e}_{\vartheta}+ u^{\phi}(\bm{x})\bm{e}_{\phi}, \\
\mbox{where}& \nonumber \\
& \bm{e}_r=\frac{\partial}{\partial r},\quad
\bm{e}_{\vartheta}=\frac{1}{r}\frac{\partial }{\partial \vartheta},\quad \bm{e}_{\phi}=\frac{1}{r\sin\vartheta}\frac{\partial}{\partial\phi}.
 \end{align}
 We see
 \begin{equation}
 \frac{1}{r^2}\frac{d}{dr}(r^2\rho u^r_{00})=\frac{1}{r^2}\frac{d}{dr}(r^2g_{00}),
 \end{equation}
 for
 $$ g=\mathrm{div}(\rho \bm{u})=
 \frac{1}{r^2}\frac{\partial}{\partial r}(r^2\rho u^r)+
 \frac{1}{r\sin\vartheta}\frac{\partial}{\partial \vartheta}(\sin\vartheta u^{\vartheta})+
 \frac{1}{r\sin\vartheta}\frac{\partial u^{\phi}}{\partial \phi}.
 $$
 Hence we have
 $$\rho(r) u^r_{00}(r)=\int_0^rg_{00}(s)s^2ds=
 -\frac{d}{dr}\Psi_{00}(r)
 \quad\mbox{for}\quad 0 \leq r <R,
 $$
 therefore we can write
 $$ \int_0^{+\infty}\frac{1}{\rho}\Big|\frac{d\Psi_{00}}{dr}\Big|^2r^2dr =
 \int_0^{+\infty}|u^r_{00}|^2 \rho(r)r^2dr.
 $$
  Thus
  \begin{equation}
  Q_{00}\geq {\delta}_*\int_0^{+\infty}|u^r_{00}|^2 \rho(r)r^2dr.
  \end{equation}\\

Moreover there is a positive number $\mu_0$ such that
\begin{equation}
\int_0^R\frac{1}{\rho}\Big|\frac{d\Psi_{00}}{dr}\Big|^2 r^2dr \geq
\mu_0\int_0^R|\Psi_{00}(r)|^2r^2dr.
\end{equation}\\

Proof. Looking at
\begin{equation}
\mu:=\inf_{\psi \in \mathcal{F}}
\frac{ \int_0^R\frac{1}{\rho}\Big|\frac{d\Psi_{00}}{dr}\Big|^2 r^2dr }{
\int_0^R|\Psi_{00}(r)|^2r^2dr},
\end{equation}
where
$$\mathcal{F}=\Big\{ \psi \in C^1[0,R]\  \Big|\  
\psi(R)=0, |\psi(r)| \leq \|\Psi_{00}\|_{L^{\infty}} \Big\},
$$
we prove $\mu >0$.

Suppose $\mu=0$ for reductio ad absurdum. Then there is a sequence $(\psi_n)_{n=1,2,\cdots,}$ in $\mathcal{F}$
such that
$\displaystyle \int_0^R|\psi_n|^2r^2dr =1$ and
$$M_n:=\int_0^R\frac{1}{\rho}\Big|\frac{d\psi_n}{dr}\Big|^2r^2dr \searrow 0
\quad\mbox{as}\quad n \to \infty.
$$ 
For $0<r \leq r' \leq R$ it hold
\begin{align*}
|\psi_n(r)|&=\Big|-\int_r^R D\psi_n(s)ds\Big| \\
&\leq
\sqrt{ \int_r^R\frac{\rho(s)}{s^2}ds }
\sqrt{ \int_r^R\frac{1}{\rho(s)}|D\psi_n(s)|^2s^2ds } \\
&\leq \sqrt{\rho_O}\sqrt{ \frac{R-r}{Rr} } M_n \leq \frac{\sqrt{\rho_O}}{r}M_1, \\
\mbox{and} & \\
|\psi_n(r')-\psi_n(r)|&
\leq \sqrt{\rho_O}\sqrt{ \frac{r'-r}{r'r} }M_n
\leq \sqrt{\rho_O}\sqrt{ \frac{r'-r}{r'r} }M_1.
\end{align*}
Therefore by the Ascoli-Arzel\`{a} theorem we can suppose $\psi_n$ tends to a limit 
$\psi_{\infty}$ uniformly on any compact subinterval of $]0, R]$, by taking a subsequence. Since $M_n \to 0$, we have $\psi_{\infty}=0$.
For any $0<r_* \ll 1$ we have
$$
\int_0^{r_*}|\psi_n(r)|^2r^2dr \leq \|\Psi_{00}\|_{L^{\infty}}^2\frac{r_*^3}{3}.
$$
Therefore
\begin{align*}
1=&\int_0^R|\psi_n(r)^2r^2dr \\
  &\leq  \|\Psi_{00}\|_{L^{\infty}}^2\frac{r_*^3}{3}+
  \int_{r_*}^R|\psi_n(r)|^2r^2dr \\
  &\to  \|\Psi_{00}\|_{L^{\infty}}^2\frac{r_*^3}{3}+
  \int_{r_*}^R|\psi_{\infty}(r)|^2r^2dr \\
  &=\|\Psi_{00}\|_{L^{\infty}}^2\frac{r_*^3}{3}.
  \end{align*}
  Taking $r_*$ so small that $\displaystyle \|\Psi_{00}\|_{L^{\infty}}^2\frac{r_*^3}{3} <1$, we see a contradiction.
  $\square$ \\

 %%%%%%%%%%%%%%%%%%%%%%%%
 
 Next we consider $\ell \geq 1$. We claim
 \begin{equation}
 Q_{\ell m} \geq
 4\pi\mathsf{G}\int_0^{+\infty}
 \Big(\|\nabla^{(\ell)}\Psi_{\ell m}\|^2-
 -4\pi\mathsf{G}\frac{d\rho}{d\varUpsilon}|\Psi_{\ell m}|^2\Big)
 r^2dr. \label{6.27}
 \end{equation}
 
 For $0<\epsilon \ll 1$, we put
 $\displaystyle \Big(\frac{d\varUpsilon}{d\rho}\Big)_{\epsilon}=\gamma \mathsf{A}(\rho+\epsilon)^{\gamma-2}$. Then, as $\epsilon \searrow +0$, 
 $\displaystyle \Big(\frac{d\varUpsilon}{d\rho}\Big)_{\epsilon} \nearrow 
 \frac{d\varUpsilon}{d\rho}$ on $0<r<R$, and
 
 \begin{align*}
 & \int_0^{+\infty}\Big(
 \Big(\frac{d\varUpsilon}{d\rho}\Big)_{\epsilon}|g_{\ell m}|^2-4\pi\mathsf{G}
 \|\nabla^{(\ell )}\Psi_{\ell m}\|^2\Big)r^2dr \\
 \nearrow 
 &\quad Q_{\ell m}=
 \int_0^{+\infty}\Big(
 \Big(\frac{d\varUpsilon}{d\rho}\Big)|g_{\ell m}|^2-4\pi\mathsf{G}
 \|\nabla^{(\ell )}\Psi_{\ell m}\|^2\Big)r^2dr.
 \end{align*}
 But we see
 \begin{align*}
 & \int_0^{+\infty}\Big(
 \Big(\frac{d\varUpsilon}{d\rho}\Big)_{\epsilon}|g_{\ell m}|^2-4\pi\mathsf{G}
 \|\nabla^{(\ell )}\Psi_{\ell m}\|^2\Big)r^2dr \\
 &=
 \int_0^{+\infty}\Big|\sqrt{\Big(\frac{d\varUpsilon}{d\rho}\Big)_{\epsilon} }\triangle^{(\ell)}\Psi_{\ell m}+
 4\pi\mathsf{G}\sqrt{ \Big( \Big(\frac{d\varUpsilon}{d\rho}\Big)_{\epsilon} \Big)^{-1} }
 \Psi_{\ell m}\Big|^2r^2dr \\
 &-4\pi\mathsf{G}
 \int_0^{+\infty}
 2\mathfrak{Re}[(\triangle^{(\ell)}\Psi_{\ell m})\Psi_{\ell m}^*]
 +\|\nabla^{(\ell)}\Psi_{\ell m}\|^2
 +4\pi\mathsf{G}\Big( \Big(\frac{d\varUpsilon}{d\rho}\Big)_{\epsilon} \Big)^{-1}
 |\Psi_{\ell m}|^2 \Big)r^2dr \\
 &\geq 4\pi\mathsf{G}
 \int_0^{+\infty}\Big(\|\nabla^{(\ell)}\Psi_{\ell m}\|^2
 -4\pi\mathsf{G}\Big( \Big(\frac{d\varUpsilon}{d\rho}\Big)_{\epsilon} \Big)^{-1}
 |\Psi_{\ell m}|^2
 \Big)r^2dr,
 \end{align*}
 for
 $$
 \int_0^{+\infty}(\triangle^{(\ell)}\Psi_{\ell m})\Psi_{\ell m}^* r^2dr =
 -\int_0^{+\infty}\|\nabla^{(\ell)}\Psi_{\ell m}\|^2r^2dr,
 $$
 since
 $\displaystyle r^2\frac{d\Psi_{\ell m}}{dr}\Psi_{\ell m}^* =O\Big(\frac{1}{r}\Big)$ as $r \rightarrow +\infty$.
  Since $\displaystyle \Big( \Big(\frac{d\varUpsilon}{d\rho}\Big)_{\epsilon} \Big)^{-1}=\frac{(\rho+\epsilon)^{2-\gamma}}{\gamma \mathsf{A}} \searrow 
  \frac{d\rho}{d\varUpsilon}=\frac{\rho^{2-\gamma}}{\gamma \mathsf{A}}$ as $\epsilon \searrow 0$, we have \eqref{6.27}.\\
  
  Let us show that
  \begin{equation}
  \mathcal{Q}_1[w]:=\int_0^{+\infty}\Big(\|\nabla^{(1)}w\|^2-
  4\pi\mathsf{G}\frac{d\rho}{d\varUpsilon}|w|^2 \Big)r^2dr \geq 0 \label{6.28}
  \end{equation}
  for $w=\Psi_{\ell m}$.
  
  Fixing ${N} \geq R$, put
  \begin{equation}
  \mathcal{Q}_1^{N}[w]:= \int_0^{N}
  \Big(\|\nabla^{(1)}w\|^2-
  4\pi\mathsf{G}\frac{d\rho}{d\varUpsilon}|w|^2 \Big)r^2dr.
  \end{equation}
  This is the quadratic form associated with the operator
  \begin{equation}
  \mathcal{P} w=\Big[-\triangle^{(1)}-4\pi\mathsf{G}\frac{d\rho}{d\varUpsilon}\Big] w.
  \end{equation}
  We consider $\mathcal{P}$ in $\mathfrak{X}_0^{N}=L^2([0,{N}]; r^2dr)$. We have the Friedrichs extension $\bm{P}^{N}$ of $\mathcal{P}\restriction C_0^{\infty}(]0,{N}[)$,
  keeping in mind that
  $$ 4\pi\mathsf{G}\frac{d\rho}{d\varUpsilon} \leq \frac{4\pi\mathsf{G}}{\gamma \mathsf{A}}\rho_O^{2-\gamma}.
  $$
  We see that $\bm{P}^{N}$ is of the Strum-Liuville type. the boundary condition at $r={N}$ is the Dirichlet boundary condition $w|_{r=S}=0$.
  
  Let $\mu_1$ be the least eigenvalue of $\bm{P}^{N}$ . We are going to show $\mu_1 \geq 0$. Let $\phi_1$ be an eigenfunction. We can suppose $\phi_1(r)>0$ for $0<r<{N}$.
  Consider $\displaystyle \varUpsilon':= \frac{d\varUpsilon}{dr}$, which satisfies
  $\mathcal{P} \varUpsilon'=0$ on $[0, +\infty[$. Here we consider
  $$
  \varUpsilon(r)=-K\Big(\frac{1}{R}-\frac{1}{r}\Big),
  \quad \varUpsilon'(r)=-\frac{K}{r^2} \quad\mbox{for}\quad R<r
  $$
with $\displaystyle K=-\Big(r^2\frac{d\varUpsilon}{dr}\Big)_{r=R-0} (>0)$. 
we look at
\begin{align*}
0=(\mathcal{P} \varUpsilon'|\phi_1)_{\mathfrak{X}_0^{N}}&=(\varUpsilon'|\mathcal{P}\phi_1)_{\mathfrak{X}_0^{N}}+
r^2\varUpsilon'\frac{d\phi_1}{dr}\Big|_{r={N}-0} \\
&=\mu_1(\varUpsilon'|\phi_1)+r^2\varUpsilon'\frac{d\phi_1}{dr}\Big|_{r={N}-0}
\end{align*}
Since $\varUpsilon' <0$ on $]0,{N}]$, we have $(\varUpsilon'|\phi_1)_{\mathfrak{X}_0^{N}} <0$. But $\displaystyle r^2\varUpsilon'\frac{d\phi_1}{dr} \Big|_{r={N}-0} \geq 0$, since
$\displaystyle \frac{d\phi_1}{dr}\Big|_{r={N}-0}\leq 0$. Consequently $\mu_1 \geq 0$ and
$$\mathcal{Q}_1^{N}[w] \geq 0 \quad \mbox{for}\quad \forall w \in \mathsf{D}(\bm{P}^{N}). $$

 Since we cannot say $w=\Psi_{\ell m} \in \mathsf{D}(\bm{P}^{N})$, for maybe $\Psi_{\ell m}({N})\not=0$,
 we decompose it as
 \begin{equation}
 w=\Psi_{\ell m}=W^{N}+\Pi^{N},
 \end{equation}
 where
 \begin{equation}
 \Pi^{N}=-C_{\ell m}{N}^{-\ell -1}\Big(\frac{r}{{N}}\Big)^{\ell},
 \quad
 C_{\ell m}=\frac{1}{2\ell +1}\int_0^Rg_{\ell m}(s)s^{\ell +2}ds.
 \end{equation}
 Then
 \begin{equation}
 \triangle^{(\ell)}W^{N}=\triangle^{(\ell)}w,\quad W^{N}({N})=0.
 \end{equation}
 Note
 $$ \mathcal{P} W^{N}=-g_{\ell m}-\Big(\frac{\ell(\ell+1)-2}{r^2}+4\pi\mathsf{G}\frac{d\rho}{d\varUpsilon}\Big)(w-\Pi^{N}) \in \mathfrak{X}_0^{N}, $$
 since $g_{\ell m}\in C([0, R[, w=\Psi_{\ell m} =O(r^2)$ for $\ell \geq 2$, $\Pi^{N}=O(r^{\ell})$. Therefore $W^{N} \in \mathsf{D}(\bm{P}^{N})$ and
 $$ \mathcal{Q}_1^{N}[W^{N}] \geq 0.$$
 Since
 \begin{align*}
 &\|\nabla^{(1)}\Pi^{N}\|^2 \leq 3C_{\ell m}{N}^{-2\ell -4}, \\
 &4\pi\mathsf{G}\frac{d\rho}{d\varUpsilon}|\Pi^{N}|^2 \leq 4\pi\mathsf{G}\frac{\rho_O^{2-\gamma}}{\gamma \mathsf{A}}C_{\ell m}^2{N}^{-2\ell -2}
 \end{align*}
 on $[0,{N}]$, we see
 \begin{align*}
 \mathcal{Q}_1^{N}[w] &\geq \mathcal{Q}_1^S[W^{N}] +O({N}^{-2\ell +1}) \\
 &\geq O({N}^{-2\ell +1}),
 \end{align*}
 Taking the limit as ${N} \rightarrow +\infty$, we get \eqref{6.28} for $w=\Psi_{\ell m}$.\\

Now 
\begin{equation}
Q_{\ell m}\geq 4\pi\mathsf{G}\mathcal{Q}_1[\Psi_{\ell m}] \geq 0,
\end{equation}
and, for $\ell \geq 2$, we have
\begin{align}
Q_{\ell m}&\geq 4\pi\mathsf{G}\int_0^{+\infty}
\Big(\|\nabla^{(\ell)}\Psi_{\ell m}\|^2
-4\pi\mathsf{G}\frac{d\rho}{d\varUpsilon}|\Psi_{\ell m}|^2
\Big)r^2dr \nonumber \\
&=4\pi\mathsf{G}
\Big( \mathcal{Q}_1[\Psi_{\ell m}]+
\int_0^{+\infty}(\ell(\ell +1)-2)|\Psi_{\ell m}|^2dr \Big) \nonumber \\
&\geq 4\pi\mathsf{G}
 \int_0^{+\infty}(\ell(\ell +1)-2)|\Psi_{\ell m}|^2dr .
 \end{align}
 
 Note that
\begin{equation}
(\ell(\ell+1)-2)\int_0^{+\infty}|\Psi_{\ell m}(r)|^2dr
\geq
\frac{4}{R^2}\int_0^R|\Psi_{\ell m}(r)|^2r^2dr
\end{equation}
for $\ell \geq 2$.\\

Consequently
\begin{align}
(\bm{L}\bm{u}|\bm{u})_{\mathfrak{H}}&=4\pi\sum_{0\leq \ell, |m|\leq \ell }Q_{\ell m} \nonumber \\
&\geq \mathfrak{d}\Big(
\int_0^{+\infty}\frac{1}{\rho}\Big|\frac{d\Psi_{00}}{dr}\Big|^2 r^2dr +
\sum_{2\leq \ell, |m|\leq \ell} (\ell(\ell+1)-2)
\int_0^{+\infty}|\Psi_{\ell m}|^2dr
\Big),
\end{align}
with
$\displaystyle {\delta} := 4\pi ({\delta}_* \wedge \mathsf{G})$. \\

Therefore

\begin{Theorem} Under  {\bf Situation \ref{As.3} },  {\bf (PD.2)} holds
for $\mathfrak{n}=\mathfrak{n}_1, \mathfrak{n}_{10}, \mathfrak{n}_1 \geq \mathfrak{n}_{10},$ defined by
\begin{subequations}
\begin{align}
\mathfrak{n}_1(\bm{u})&=\Big[
\int_0^{+\infty}\frac{1}{\rho}\Big|\frac{d\Psi_{00}}{dr}\Big|^2 r^2dr +
\sum_{2\leq \ell, |m|\leq \ell}
(\ell(\ell+1)-2)\int_0^{+\infty}|\Psi_{\ell m}|^2dr
\Big]^{\frac{1}{2}} \nonumber \\
&=\Big[
\int_0^{+\infty}|u^r_{00}|^2 \rho(r)r^2dr
+
\sum_{2\leq \ell, |m|\leq \ell}
(\ell(\ell+1)-2)\int_0^{+\infty}|\Psi_{\ell m}|^2dr
\Big]^{\frac{1}{2}} , \\
\mathfrak{n}_{10}(\bm{u})&=\sqrt{\mu_1}\Big[
\int_0^{R}|\Psi_{00}(r)|^2r^2dr
+\sum_{2 \leq \ell, |m|\leq\ell}
\int_0^{R}|\Psi_{\ell m}(r)|^2 r^2dr
\Big]^{\frac{1}{2}}  \\
&\mbox{with} \quad
\mu_1=\mu_0 \wedge \frac{4}{R^2}. \nonumber
\end{align}
\end{subequations}
Here
\begin{align}
\Psi_{\ell m}(r)&=-
\int_0^{2\pi}\int_0^{\pi}
\mathcal{K}[g](\bm{x})Y_{\ell m}(\vartheta, \phi)^*
\sin\vartheta d\vartheta d\phi, \\
u^r_{00}(r)&=\sqrt{4\pi}\int_0^{2\pi}\int_0^{\pi}\Big(\bm{u}(\bm{x})\Big|\frac{\bm{x}}{\|\bm{x}\|}\Big)\sin\vartheta d\vartheta d\phi,
\end{align}
in which
$$
 \bm{x}=(x^1,x^2, x^3)=(r\sin\vartheta\cos\phi, r\sin\vartheta\sin\phi,r \cos\vartheta)
 $$
and $g=\mathrm{div}(\rho \bm{u})$.
\end{Theorem}

Note that the seminorm $\mathfrak{n}=\mathfrak{n}_1, \mathfrak{n}_{10}$ and the associated $\mathfrak{N}$ are not norm, that is, $\mathfrak{N}(U)=0$ does not imply $U=\bm{0}$, say, does not imply $\bm{u}=0$.
But $\mathfrak{N}({U})$ can control the amplitude of 
$\Psi=\mathcal{K}[\mathrm{div}(\rho \bm{u})]$, except for the $\ell=1$ components, and $\dot{\bm{u}}$. Here $
-4\pi\mathsf{G}\Psi=\mathfrak{d} \Phi =\Phi_{\rho_b+\mathfrak{d}\rho}-\Phi_b,
\mathfrak{d}\rho=-\mathrm{div}(\rho \bm{u})$ and $\bm{v}=\dot{\bm{u}}$ are the  quantities essential for the perturbation, and
the control of the magnitude of $\bm{u}=\bm{\varphi}(t,\bm{x})-\bm{x}+\bm{u}^0$ itelf is not essential in the discussion of the stability of the background. 
( Recall $\bm{\varphi}_b(t,\bar{\bm{x}})=\bar{\bm{x}}$ for $\omega_b=0$.)
In fact
$\mathfrak{n}_1(\bm{u})=0$ implies that $\Psi_{\ell m}=0, g_{\ell m}=0$ on $[0,+\infty[$ for $\ell \not=1$, and
$\mathfrak{n}_{10}(\bm{u})=0$ implies that $\Psi_{\ell m}=0, g_{\ell m}=0$ on $[0,+R]$ for $\ell \not=1$,
and
$$g(\bm{x})=
\sqrt{\frac{3}{4\pi}}g_{1,0}(r)\cos\vartheta 
-\sqrt{\frac{3}{8\pi}}g_{1,1}(r)\sin\vartheta e^{\mathrm{i}\phi}
+\sqrt{\frac{3}{8\pi}}g_{1,-1}(r)\sin\vartheta e^{-\mathrm{i}\phi}.
$$

\begin{Remark}
It seems impossible to control the magnitude of $\Psi_{1
 m}, m=0,\pm1$ by the following reason. The differential operator $\mathcal{P}\restriction C_0^{\infty}(]0,+\infty[)$ admits the Friedrichs extention $\bm{P}$ in $\mathfrak{X}=L^2([0,+\infty[, r^2dr)$, and $\mathcal{Q}_1$ is the quadratic form associated with $\bm{P}$. But $\bm{P}$ is not of the Sturm-Liouville type. In fact, let $\phi \in \mathsf{D}(\bm{P})$ and $\mathcal{P}\phi=\lambda \phi$ with $\lambda >0$. Since $-\triangle^{(1)}\phi=0$ on $]R,+\infty[$,
there are constants $C_{\pm}$ such that
$$
\phi=C_+\phi_+ + C_-\phi_- \quad\mbox{on}\quad ]R,+\infty[,$$
where
$$\phi_{\pm}(r)=\sqrt{r}J_{\pm\frac{3}{2}}(\sqrt{\lambda}r),\quad J_{\pm\frac{3}{2}}\quad\mbox{being the Bessel function}.
$$
Since
$$J_{\pm\frac{3}{2}}(r) \sim \sqrt{ \frac{2}{\pi r} }\cos
\Big(r\pm\frac{3\pi}{4} -\frac{\pi}{4}\Big)\quad\mbox{as}\quad r \rightarrow +\infty,
$$
we see $\phi_{\pm} \not\in L^2([R,+\infty[, r^2dr)$, and
$\phi \in L^2([R, +\infty[, r^2dr)$ requires $C_{\pm}=0$, $\phi=0$ on $]R, +\infty[$.
By the uniqueness of solutions of ODE, it follows that $\phi=0$ on $[0, +\infty[$.
In other words, any positive real number cannot be an eigenvalue of $\bm{P}$, so $\bm{P}$ is not of the Sturm-Liouville type.

However, if we restrict ourselves to axially and equatorially symmetric perturbations, functions of $(r,|\zeta|)$, we have $\Psi_{\ell m}=0$ for odd $\ell$ a priori, since
$\int_0^{2\pi}e^{\mathrm{i}m\phi}d\phi=0$ for $m\not=0$ and $P_{\ell}^0(-\zeta)=-P_{\ell}^0(\zeta)$ for odd $\ell$, we need not control the magnitudes of $\Psi_{1 m}, m=0,\pm1$. In this situation $\mathfrak{n}(\bm{u})=0$ implies $\Psi=0, g=0$, or,
$\mathfrak{d}\Phi=0, \mathfrak{d}\rho=0$, or, more precisely, we have
\begin{equation}
\mathfrak{n}_{10}(\bm{u}) \geq
\sqrt{\frac{\mu_1}{4\pi}}\|\Psi\|_{L^2(\mathfrak{R})}
=\sqrt{\frac{\mu_1}{4\pi}}\|\mathcal{K}[\mathfrak{d} \rho]\|_{L^2(\mathfrak{R})},
\end{equation}
since
$$
\|\Psi\|_{L^2(\mathfrak{R})}^2=
4\pi
\sum_{\substack{\ell : \mbox{even} \\
|m| \leq \ell }}
\|\Psi_{\ell m}\|_{L^2([0, R],r^2dr)}^2
$$
in this case.
\end{Remark}

\textbullet\  Now we are going to try to derive {\bf (PD.2)} with another seminorm
 \begin{equation}
 \mathfrak{n}(\bm{u}):=
 \|\mathrm{div}(\rho \bm{u})\|_{L^2\Big(\frac{\Gamma P}{\rho^2}d\bm{x}\Big)}
 =
 \Big[
 \int_{\mathfrak{R}}\frac{\Gamma P}{\rho^2}|\mathrm{div}(\rho \bm{u})|^2d\bm{x} \Big]^{\frac{1}{2}}
 \end{equation}
 for 
 \begin{align}
(\mathcal{L}\bm{u}|\bm{u})_{\mathfrak{H}}&=
\int_{\mathfrak{R}} \frac{\Gamma P}{\rho^2}|\mathrm{div}(\rho \bm{u})|^2 d\bm{x} +2\mathfrak{Re}
\Big[\int_{\mathfrak{R}} \frac{\Gamma P\mathscr{A}}{\rho}(\bm{u}|\bm{n})\mathrm{div}(\rho \bm{u})^*\Big]d\bm{x} \nonumber \\
&-\int_{\mathfrak{R}} \frac{\Gamma P\mathscr{A}}{\rho}\|{\nabla}\rho\|^2|(\bm{u}|\bm{n})|^2 d\bm{x}+ \nonumber \\
&-4\pi\mathsf{G}
\int_{\mathfrak{R}} \mathcal{K}[\mathrm{div}(\rho \bm{u})]\mathrm{div}(\rho \bm{u})^*d\bm{x}. 
\end{align}\\

First we claim:\\

\begin{Proposition}\label{Prop.8v11}
Let $\varepsilon$ be a positive number. If
\begin{equation}
-\varepsilon \frac{\|{\nabla}\rho\|^2}{\rho}
\leq \mathscr{A} \leq 0 \quad\mbox{on}\quad \mathfrak{R}, \label{6.14}
\end{equation}
then
\begin{equation}
(\bm{L}_0\bm{u}|\bm{u})\geq 
(1-\varepsilon)
\int_{\mathfrak{R}}\frac{\Gamma P}{\rho^2}|g|^2d\bm{x}.
 \label{6.15}
\end{equation}
Here $g=\mathrm{div}(\rho\bm{u})$.
\end{Proposition}

Proof. We see
\begin{align*}
&\Big|2\mathfrak{Re}\Big[\int_{\mathfrak{R}}\frac{\Gamma P\mathscr{A}}{\rho}(\bm{u}|\bm{n})g^*d\bm{x}\Big]\Big|
\leq \varepsilon\int_{\mathfrak{R}}\int_{\mathfrak{R}}\frac{\Gamma P}{\rho^2}|g|^2d\bm{x}+
\frac{1}{\varepsilon}\int_{\mathfrak{R}}\Gamma P \mathscr{A}^2|(\bm{u}|\bm{n}))|^2d\bm{x} \\
&\leq
\varepsilon\int_{\mathfrak{R}}\frac{\Gamma P}{\rho^2}
|g|^2d\bm{x}+
\int_{\mathfrak{R}}\frac{\Gamma P |\mathscr{A}|}{\rho}\|{\nabla} \rho\|^2|(\bm{u}|\bm{n}))|^2d\bm{x} 
\end{align*}
by \eqref{6.14}, where  $g=\mathrm{div}(\rho\bm{u})$, since
$-|\mathscr{A}|-\mathscr{A} =0$ for $\mathscr{A} \leq 0$, we have \eqref{6.15}.
$\square$

Therefore, if we adopt the Cowling approximation, which neglect the perturbed self-gravitation term
$-4\pi\mathsf{G}\int_{\mathfrak{R}}\mathcal{K}[g]g^*d\bm{x}$,
we have {\bf (PD.2)} with $0<\varepsilon <1$. 
Namely we consider the following 

\begin{Situation} \label{S.2}
The equation of state is \eqref{EOS.g}: $P=\rho^{\gamma}\exp(S/\mathsf{C}_V)$.
Instead of $\bm{L}$, we consider $\bm{L}_0$:

\begin{align}
\bm{L}_0\bm{u}&={\nabla}\Big(
-\frac{\gamma P}{\rho^2}\mathrm{div}(\rho\bm{u})+\frac{\gamma P}{\rho}(\bm{u}|\mathfrak{a})\Big)+ \nonumber\\
&+\frac{\gamma P}{\rho}
\Big(-\mathfrak{a}\mathrm{div}(\rho \bm{u})+(\bm{u}|\mathfrak{a}){\nabla}\rho \Big).
\end{align}

\end{Situation}

We can claim

\begin{Theorem}
Under the {\bf Situation \ref{S.2} }suppose \eqref{6.14} with $0<\varepsilon <1$. Put $\delta:=1-\varepsilon$. Then {\bf (PD.2)} holds for $\bm{L}$ replaced by $\bm{L}_0$ with 
$\delta$ and $\mathfrak{n}=\mathfrak{n}_2$ defined by
\begin{equation}
 \mathfrak{n}_2(\bm{u})=
 \Big[
 \int_{\mathfrak{R}}\frac{\gamma P}{\rho^2}|\mathrm{div}(\rho \bm{u})|^2d\bm{x} \Big]^{\frac{1}{2}}.
 \end{equation}
 \end{Theorem}

So let us look at the perturbed self-gravitation term in order to justify the Cowling approximation in this context.\\

Put
\begin{align}
&K^{\spadesuit}:=
\sup_{g \in \mathfrak{g}_L} \frac{4\pi\mathsf{G}\int_{\mathfrak{R}}\mathcal{K}[g]g^*d\bm{x}}{\int_{\mathfrak{R}}\frac{\gamma P}{\rho^2}|g|^2d\bm{x}},  \\
\mbox{where}& \nonumber \\
&\mathfrak{g}_L:=\Big\{ g \in L^2\Big(\frac{\gamma P}{\rho^2}d\bm{x}, \mathfrak{R}\Big) 
  \Big|\  g=\mathrm{div}(\rho \bm{u}) \quad \mbox{for}\quad \exists\bm{u} \in \mathsf{D}(\bm{L})
  \Big\}. 
\end{align}
We can claim $K^{\spadesuit} <\infty$ by the following\\

\begin{Proposition}
There exists a constant $C$ such that
\begin{equation}
0\leq 
4\pi\mathsf{G}\int_{\mathfrak{R}}\mathcal{K}[g]g^*d\bm{x}
\leq C
\int_{\mathfrak{R}}\frac{\gamma P}{\rho^2}|g|^2d\bm{x}. \label{6.17}
\end{equation}
\end{Proposition}

Proof. We have
\begin{align*}
|\mathcal{K}[g](\bm{x})|&\leq
\frac{1}{4\pi}\sqrt{ \int_{\mathfrak{R}}\frac{d\bm{x}'}{\|\bm{x}-\bm{x}'\|} }
\sqrt{ \int_{\mathfrak{R}}|g|^2d\bm{x} } \\
&\leq \sqrt{\frac{R}{2\pi}}\|g\|_{L^2(\mathfrak{R})},
\end{align*}
where
$R:=\sup \{ \|\bm{x}\|\  |\  \bm{x} \in \mathfrak{R} \}$. Therefore
\begin{align*}
\int_{\mathfrak{R}} \mathcal{K}[g]g^*d\bm{x} &\leq
\sqrt{\frac{R}{2\pi}}\|g\|_{L^2}\int_{\mathfrak{R}}|g|d\bm{x} \\
&\leq \sqrt{\frac{2}{3}}R^2\|g\|_{L^2}^2
\end{align*}
On the other hand, 
\begin{equation}
\nu_*:=\inf_{\mathfrak{R}}\frac{\gamma P}{\rho^2} >0,
\end{equation}
since
$$
\frac{1}{C} \mathsf{d}^{-\frac{2-\gamma}{\gamma-1}}
\leq \frac{\gamma P}{\rho^2} \leq 
C \mathsf{d}^{-\frac{2-\gamma}{\gamma-1}} \quad\mbox{on}\quad \mathfrak{R},
$$
where $\mathsf{d}:=\mathrm{dist} (\bm{x}, \partial \mathfrak{R})$ and $1<\gamma <2$.
Then 
$$\int_{\mathfrak{R}} |g|^2d\bm{x} \leq \frac{1}{\nu_*}\int_{\mathfrak{R}}\frac{\gamma P}{\rho^2}|g|^2d\bm{x},
$$
therefore \eqref{6.17} holds with
$\displaystyle C=\sqrt{\frac{2}{3}}\frac{R^2}{\nu_*}$. $\square$\\

Consequently we can claim \\

{\it 
Suppose $K^{\spadesuit} <1$ and let $0<\varepsilon < 1-K^{\spadesuit}$ so that
$\delta:=1-\varepsilon-K^{\spadesuit} >0$. If \eqref{6.14} holds with $\varepsilon$,
then {\bf (PD.2)} holds with ${\delta}$ and}
$$
 \mathfrak{n}(\bm{u}):=
 \|\mathrm{div}(\rho \bm{u})\|_{L^2\Big(\frac{\gamma P}{\rho^2}d\bm{x}\Big)}
 =
 \Big[
 \int_{\mathfrak{R}}\frac{\gamma P}{\rho^2}|\mathrm{div}(\rho \bm{u})|^2d\bm{x} \Big]^{\frac{1}{2}}.
$$\\

Thus we have the following question, which is still open:\\

\begin{Question}
When the background  realizes the condition $K^{\spadesuit} <1$ ?
\end{Question}

Let us note the nondimensionalization of the estimate of $K^{\spadesuit}$. We have
$$ \frac{\gamma P}{\rho^2}=4\pi\mathsf{G} \mathsf{a}^2 (\gamma-1)\Theta^{-\frac{2-\gamma}{\gamma-1}}(1+\omega_1)(1+\omega_2),
$$
where
\begin{align*}
&\Theta=\Theta\Big(\frac{\|\bm{x}\|}{\mathsf{a}},\frac{x^3}{\|\bm{x}\|}; \frac{1}{\gamma-1},\mathfrak{b}\Big), \\
& \omega_1=O(\varUpsilon_O),\quad \omega_2=O\Big(\frac{\gamma-1}{\gamma \mathsf{A}}
\varUpsilon_O\Big),
\end{align*}
and
$$\int_{\mathfrak{R}}\frac{\gamma P}{\rho^2}|g|^2d\bm{x} =
4\pi\mathsf{G} \mathsf{a}^5\int_{\mathfrak{R}_1}
(\gamma-1)\Theta^{-\frac{2-\gamma}{\gamma-1}}(1+\omega_1)(1+\omega_2)|\hat{g}(\underline{\bm{x}})|^2d\underline{\bm{x}},
$$
where
\begin{align*}
&\mathfrak{R}_1=\Big\{ \underline{\bm{x}} \quad \Big| \quad \|\underline{\bm{x}}\| < {U}_1\Big(\frac{\underline{x}^3}{\|\underline{\bm{x}}\|}\Big) \quad \Big\}
, \\
&\Theta=\Theta\Big(\|\underline{\bm{x}}\|,\frac{\underline{x}^3}{\|\underline{\bm{x}}\|}; \frac{1}{\gamma-1},\mathfrak{b}\Big),
\quad \hat{g}(\underline{\bm{x}})=g(\mathsf{a}\underline{\bm{x}}).
\end{align*}
On the other hand, we see
$$4\pi\mathsf{G}\int_{\mathfrak{R}}\mathcal{K}[g]g^*d\bm{x} =4\pi\mathsf{G} \mathsf{a}^5
\int_{\mathfrak{R}_1}\mathcal{K}[\hat{g}](\underline{\bm{x}})\hat{g}^*(\underline{\bm{x}})d\underline{\bm{x}}.
$$
Therefore we have the reduction
\begin{equation}
K^{\spadesuit}=\sup_{\hat{g}}
\frac{ \int_{\mathfrak{R}_1} \mathcal{K}[\hat{g}]\hat{g}^*d\underline{\bm{x}} }{ \int_{\mathfrak{R}_1}(\gamma-1)\Theta^{-\frac{2-\gamma}{\gamma-1}}(1+\omega_1)(1+\omega_2)|\hat{g}|^2d\underline{\bm{x}} }.
\end{equation}
If the background is isentropic, we have $\omega_1=\omega_2=0$.
If $\Omega=0$ and the background is isentropic and spherically symmetric, then
$$
K^{\spadesuit}=\sup_{\hat{g}}
\frac{ \int_{\mathfrak{R}_1} \mathcal{K}[\hat{g}]\hat{g}^*d\underline{\bm{x}} }{ \int_{\mathfrak{R}_1}(\gamma-1)\theta^{-\frac{2-\gamma}{\gamma-1}}|\hat{g}|^2d\underline{\bm{x}} }.
$$
where $\theta=\theta(\|\underline{\bm{x}}\|,; \frac{1}{\gamma-1})$ is the Lane-Emden function.

\section{{\bf  ELASO($\blacklozenge$)} of astrophysicists}

%\subsection{}

The equation of linearized adiabatic stellar oscillations traditionally used by astrophysicists is little bit different from {\bf ELASO}
\eqref{Eqxi} with \eqref{L}, \eqref{LASOa}, \eqref{LASOb}, \eqref{LASOc}, \eqref{LASOd}
 derived in Section 3. In this section we derive the equation {\bf ELASO($\blacklozenge$) } used 
 by astrophysicists after V. R. Lebovitz 1961 \cite{Lebovitz}, D. Lynden-Bell and J. P. Ostriker 1967 \cite{LyndenBO} and J. L. Friedman and B. F. Schutz 1978 \cite{FriedmanS}.\\

We consider a pair of solutions $(\rho, S, \bm{v})$ and $(\rho_b, S_b, \bm{v}_b)$ of the Euler-Poisson equations described in a rotating frame \eqref{EPa},\eqref{EPb},\eqref{EPc},\eqref{NPot} with the equation of state \eqref{1.8}.
Here $(\rho_b, S_b, \bm{v}_b)$ is a compactly supported axially symmetric stationary solution with
\begin{equation}
\bm{v}_b(t,\bm{x})=\omega_b(\varpi, z)
\begin{bmatrix}
0 \\
0 \\
1
\end{bmatrix}
\times \bm{x} =
\omega_b(\varpi, z)
\begin{bmatrix}
-x^2 \\
x^1 \\
0
\end{bmatrix}.
\end{equation}

We define the Lagrange displacement $ \bm{\xi}$ in the Lebovitz sense (\cite{Lebovitz}) by
\begin{equation}
 \bm{\xi}(t,\bm{x})=\bm{\varphi}(t|0,\bm{\varphi}(0|t,\bm{x}||\bm{v}_b)||\bm{v})-\bm{x}. \label{z7}
\end{equation}
Here, for $\bm{V}=\bm{v},\bm{v}_b$, $\bm{\varphi}(t|s,\bm{y}||\bm{V})$ stands for the solution $\bm{x}=\bm{\varphi}(t|s,\bm{y}||\bm{V})$ of the problem
\begin{equation}
\frac{d\bm{x}}{d t}=\bm{V}(t,\bm{x}), \quad \bm{x}|_{t=s}=\bm{y}.
\end{equation}\\

Let us suppose Assumption \ref{Ass.0ELASObl}: 

{\it The initial configurations of $\rho, S$ are the same to those of the background, that is,
\begin{subequations}
\begin{align}
&\rho^0=\rho(0,\cdot)=\rho_b, \\
&S^0=S(0,\cdot)=S_b.
\end{align}
\end{subequations}
} \\

Then the linearized equation for $\bm{\xi}$ turns out to be

\begin{equation}
\frac{\partial^2 \bm{\xi}}{\partial t^2}+\mathcal{B}_{\blacklozenge}
\frac{\partial \bm{\xi}}{\partial t} + \mathcal{L}_{\blacklozenge} \bm{\xi}=0
\quad\mbox{in}\quad [0,+\infty[\times \mathfrak{R}_b, \label{ZLASO}
\end{equation}
where
\begin{equation}
\mathcal{B}_{\blacklozenge}\dot{ \bm{\xi}}=2\bm{\Omega} \times \dot{ \bm{\xi}}+
 \mathcal{B}_{\blacklozenge 1}\dot{ \bm{\xi}} \label{7.73}
 \end{equation}
 with
 \begin{equation}
   \mathcal{B}_{\blacklozenge 1}\dot{ \bm{\xi}}=2(D\dot{ \bm{\xi}})\bm{v}_b
   =2\omega_b\frac{\partial \dot{\bm{\xi}}}{\partial \phi},\label{7.74}
  \end{equation}
and
\begin{equation}
\mathcal{L}_{\blacklozenge}=\mathcal{L}_{\blacklozenge 2}+\mathcal{L}_{\blacklozenge1}+\mathcal{L}_{\blacklozenge 0}+
\mathcal{L}, \label{7.75}
\end{equation}

with
 \begin{subequations}
 \begin{align} 
 \mathcal{L}_{\blacklozenge 2} \bm{\xi}&=D((D \bm{\xi})\bm{v}_b)\bm{v}_b
 =\omega_b^2\frac{\partial^2\bm{\xi}}{\partial \phi^2},   \label{7.76a}\\
 \mathcal{L}_{\blacklozenge 1} \bm{\xi}&=2\bm{\Omega} \times ((D \bm{\xi})\bm{v}_b) =2\Omega \omega_b\frac{\partial}{\partial \phi}
\begin{bmatrix}
- \xi^2 \\
\\
 \xi^1 \\
\\
0
\end{bmatrix}, \\
 \mathcal{L}_{\blacklozenge 0} \bm{\xi}&=
-\Big( D((D\bm{v}_b)\bm{v}_b)
+D(2\bm{\Omega} \times \bm{v}_b) \Big)
 \bm{\xi} =W_b\bm{\xi}, \\
 \mbox{with}& \nonumber \\
 W_b&:=
-D((D\bm{v}_b)\bm{v}_b)-D(2\bm{\Omega}\times \bm{v}_b) \nonumber \\
&=(2\Omega+\omega_b)\omega_b
\begin{bmatrix}
1 & 0 & 0 \\
0 & 1 & 0 \\
0 & 0 & 0
\end{bmatrix}
+2(\Omega+\omega_b)
\begin{bmatrix}
\frac{1}{\varpi}\frac{\partial\omega_b}{\partial\varpi}(x^1)^2 & \frac{1}{\varpi}\frac{\partial\omega_b}{\partial\varpi} x^1x^2 & \frac{\partial \omega_b}{\partial z}x^1 \\
\\
\frac{1}{\varpi}\frac{\partial\omega_b}{\partial\varpi} x^1x^2 & \frac{1}{\varpi}\frac{\partial\omega_b}{\partial\varpi} (x^2)^2 & \frac{\partial \omega_b}{\partial z}x^2 \\
\\
0 & 0 & 0
\end{bmatrix},  \\
\mathcal{L} \bm{\xi}&=
\frac{1}{\rho_b}\nabla \lceil \delta P \rfloor
-\frac{\nabla P_b}{\rho_b^2}\lceil \delta\rho \rfloor
-4\pi\mathsf{G}\nabla\mathcal{K}[ \lceil \delta\rho \rfloor   ], \\
\mbox{with}& \nonumber \\
 \lceil \delta\rho \rfloor&=-\mathrm{div}(\rho_b \bm{\xi}), \quad \lceil \delta S \rfloor =-(\bm{\xi}|\nabla S_b), \nonumber \\
 \lceil \delta P \rfloor &=\Big(\frac{\Gamma P}{\rho}\Big)_b
\lceil \delta\rho \rfloor
+\Big(\frac{\partial P}{\partial S}\Big)_b\lceil \delta S \rfloor
=\Big(\frac{\Gamma P}{\rho}\Big)_b \lceil \delta\rho \rfloor+ 
( \Gamma P)_b( \bm{\xi} |\mathfrak{a}_b). \label{ZLASO7.9f}
\end{align}
\end{subequations}\\

Here and hereafter $\nabla, D$ stand for $\nabla_{\bm{x}}, D_{\bm{x}}$,
namely,
\begin{align*}
\nabla f &=
\begin{bmatrix}
\frac{\partial f}{\partial x^1} \\
\\
\frac{\partial f}{\partial x^2} \\
\\
\frac{\partial f}{\partial x^3}
\end{bmatrix}, \\
D \bm{Q} &=
\begin{bmatrix}
\frac{\partial Q^1}{\partial x^1} & \frac{\partial Q^1}{\partial x^2} & \frac{\partial Q^1}{\partial x^3} \\
\\
\frac{\partial Q^2}{\partial x^1} & \frac{\partial Q^2}{\partial x^2} & \frac{\partial Q^2}{\partial x^3} \\
\\
\frac{\partial Q^3}{\partial x^1} & \frac{\partial Q^3}{\partial x^2} & \frac{\partial Q^3}{\partial x^3} 
\end{bmatrix}, 
\end{align*}
and
$\displaystyle \frac{\partial \bm{u}}{\partial\phi}$ means
$\displaystyle -x^2\frac{\partial\bm{u}}{\partial  x^1}+x^1\frac{\partial \bm{u}}{\partial x^2}$.\\

The equation \eqref{ZLASO} with  \eqref{7.73} \eqref{7.74} \eqref{7.75} with \eqref{7.76a} - \eqref{ZLASO7.9f} is the equation used by astrophysicists after \cite{Lebovitz} \cite{FriedmanS} \cite{LyndenBO}.
 We call this equation the `{\bf ELASO($\blacklozenge$)} `.\\

Let us explain how to derive the equation {\bf ELASO($\blacklozenge$) }\eqref{ZLASO} following 
\cite{Lebovitz}, \cite{FriedmanS} and \cite{LyndenBO}.\\

We shall use the identities
\begin{align}
&D_{\bm{y}}\bm{\varphi}(t|s,\bm{y}||\bm{V})=\Big(D_{\bm{y}}\bm{\varphi}(s|t, \bm{\varphi}(t|s,\bm{y}||\bm{V})||\bm{V}) \Big)^{-1}, \label{Z4.1} \\
&\partial_{s}\bm{\varphi}(t|s,\bm{y}||\bm{V})=-D_{\bm{y}}\bm{\varphi}(t|s,\bm{y}||\bm{V})\bm{V}(s,\bm{y}). \label{Z4.2}
\end{align}\\

In fact, \eqref{Z4.1} can be deduced by differentiating the identity
\begin{equation}
\bm{y}=\bm{\varphi}(s|t,\bm{\varphi}(t|s,y||\bm{V})||\bm{V} )\label{Z4.3}
\end{equation}
by $\bm{y}$,
and \eqref{Z4.2} can be deduced by differentiating \eqref{Z4.3} by $s$
and using \eqref{Z4.1}.\\

Note that, for $\bm{x}=\bm{x}(t,\bar{\bm{x}})=\bm{\varphi}(t|0,\bar{\bm{x}}||\bm{V})$, it holds that
\begin{align}
&D_{\bar{\bm{x}}}\bm{x}(t,\bar{\bm{x}})=
\exp \Big[
\int_0^t D_{\bm{x}}\bm{V}(t', \bm{\varphi}(t'|0,\bar{\bm{x}}||\bm{V}))dt' \Big], \label{Z4.4} \\
& \mathrm{det}
D_{\bar{\bm{x}}}\bm{x}(t,\bar{\bm{x}})=
\exp \Big[
\int_0^t \mathrm{div}_{\bm{x}}\bm{V}(t', \bm{\varphi}(t'|0,\bar{\bm{x}}||\bm{V}))dt' \Big], \label{Z4.5} 
\end{align}\\

\begin{Definition}
For a pair $(Q, Q_b)$ of scalar fields, we define the scaler field $\Delta Q$, the Lagrangian perturbation, and $\delta Q$, the Eulerian perturbation, by
\begin{align}
& \Delta Q(t,\bm{x})=Q(t,\bm{x}+ \bm{\xi}(t,\bm{x}))-Q_b(t,\bm{x}), \\
& \delta Q(t,\bm{x})=Q(t,\bm{x})-Q_b(t,\bm{x}). 
\end{align}
For a pair $(\bm{Q}, \bm{Q}_b)$ of vector fields, we define
\begin{align}
& \Delta \bm{Q}(t,\bm{x})=\bm{Q}(t,\bm{x}+ \bm{\xi}(t,\bm{x}))-\bm{Q}_b(t,\bm{x}), \\
&\delta \bm{Q}(t,\bm{x})=\bm{Q}(t,\bm{x})-\bm{Q}_b(t,\bm{x}).
\end{align}
\end{Definition}

We suppose that $\delta \rho, \delta S, \delta\bm{v} = \mathscr{O}(\varepsilon)$.\\

Then we have $ \bm{\xi}=\mathscr{O}(\varepsilon)$.\\

In fact, writing
$$ \bm{\xi}=\underline{\bm{x}} -\bm{x}, \quad \underline{\bm{x}}=\bm{\varphi}(t|0,\bar{\bm{x}}|| \bm{v}), \quad \bm{x}=\bm{\varphi} (t|0,\bar{\bm{x}}||\bm{v}_b), $$
we see that 
$\bm{\psi}(t)= \bm{\xi}(t,\bm{x})$ satisfies
$$\frac{d\bm{\psi}}{dt}=\bm{v}(t, \bm{x}+\bm{\psi})-\bm{v}(t,\bm{x})+\bm{v}(t,\bm{x})-\bm{v}_b(t,\bm{x}).$$
Hence
$$\Big\|\frac{d\bm{\psi}}{dt}\Big\| \leq C_0+C_1\|\bm{\psi}\|, \quad \bm{\psi}(0)=\bm{0},$$
where 
$$C_0=\|\delta \bm{v}\|_{L^{\infty}},\quad C_1=\| D_{\bm{x}}\bm{v}\|_{L^{\infty}}.
$$
The Gronwall argument imlies 
$$ \|\bm{\psi} (t)\| \leq C_0(e^{C_1t}-1). $$
Thus $ \bm{\xi}={O}(\varepsilon)$.
As for the derivatives $\displaystyle \psi_j(t)=\frac{\partial}{\partial x^j} \bm{\xi}(t,\bm{x})$, we see 
$\psi_j=O(\varepsilon)$ by the same argument,  so we have 
$\displaystyle \frac{\partial}{\partial x^j} \bm{\xi}=O(\varepsilon)$.\\

Consequently, when $\delta Q=\mathscr{O}(\varepsilon)$ for scalar fileds paired, we can claim
\begin{equation}
\Delta Q=\delta Q+(\bm{\xi}|\nabla_{\bm{x}}Q_b)+\mathscr{O}(\varepsilon^2). \label{*7.20}
\end{equation}
and when $\delta\bm{Q}=\mathscr{O}(\varepsilon)$ for vector fields, we have
\begin{equation}
\Delta\bm{Q}=\delta\bm{Q}+(D\bm{Q}_b) \bm{\xi}+\mathscr{O}(\varepsilon^2). \label{z19}
\end{equation}\\

\textbullet\  Note that it hold that
\begin{equation}
\mathrm{div}\bm{v}_b=0,
\quad
(\bm{v}_b|\nabla \rho_b)=0,
\quad
(\bm{v}_b|\nabla S_b)=0.
\end{equation}\\

\textbullet\  Integration of the equation \eqref{EPc} gives
\begin{equation}
S(t.\bm{\varphi}(t|0,\bar{\bm{x}}||\bm{v}))=S^0(\bar{\bm{x}}):=S(0,\bar{\bm{x}}).
\end{equation}
Then 
$$\Delta S(t,\bm{x})=S^0(\bar{\bm{x}})-S_b(\bm{x}). $$
Here $\bar{\bm{x}}=\bm{\varphi}(0|t,\bm{x}||\bm{v}_b)$. 

We suppose
\begin{equation}
S^0=S_b. \label{Z14.1}
\end{equation}
Since $S_b$ is supposed to be axially symmetric, and 
the values of coordinates $\varpi, z$ of $\bm{\varphi}(t|s,\bm{y}||\bm{v}_b)$ are invariant in $t$, we see $S_b(\bar{\bm{x}})=S_b(\bm{x})$. Therefore
\begin{equation}
\Delta S=0,
\end{equation}
provided \eqref{Z14.1}. 
Thanks to \eqref{*7.20}, this implies
\begin{equation}
\delta S=-(\bm{\xi}|\nabla S_b)+ \mathscr{O}(\varepsilon^2). \label{7*dS}
\end{equation}\\

\textbullet\  Integration of the equation of continuity \eqref{EPa} gives
\begin{equation}
\rho(t,\bm{\varphi}(t|0, \bar{\bm{x}}||\bm{v}))=
\rho^0(\bar{\bm{x}})\exp
\Big[-\int_0^t(\mathrm{div}_{\bm{x}}\bm{v})(t', \bm{\varphi}(t'|0,\bar{\bm{x}}||\bm{v}) dt' \Big], \label{Z8.0}
\end{equation}
where $\rho^0=\rho(0,\cdot)$. 

We denote 
\begin{equation}
\mathfrak{R}(t)=\{ \bm{x} \quad|\quad  \rho(t,\bm{x})>0\},
\quad \mathfrak{R}_b =\{ \bm{x} \quad|\quad \rho_b(\bm{x}) >0\}.
\end{equation}
Of course
$$\mathfrak{R}(t)=\{ \bm{x}=\bm{\varphi}(t|0,\bar{\bm{x}}||\bm{v}) \quad|\quad \bar{\bm{x}} \in \mathfrak{R}(0) \}, $$
and the equation \eqref{Z8.0} is valid for
$\forall (t,\bar{\bm{x}}) \in [0,T]\times \mathbb{R}^3$.

We suppose
\begin{equation}
\rho^0=\rho_b. \label{Z8.1}
\end{equation}
Then $\mathfrak{R}(0)=\mathfrak{R}_b$ and
$$ \mathfrak{R}(t)=\{ \bm{x}=\bm{\varphi}(t|0,\bar{\bm{x}}||\bm{v}) \quad|\quad \bar{\bm{x}} \in \mathfrak{R}_b \},
$$
that is,
$$\underline{\bm{x}} =\bm{x}+ \bm{\xi}(t,\bm{x}) \in \mathfrak{R}(t) \quad
\Leftrightarrow \quad \bm{x} \in \mathfrak{R}_b. $$
Here we note that $\rho_b(\bm{\varphi}(0|t,\bm{x}||\bm{v}_b))=\rho_b(\bm{x})$ by the axial symmetricity.\\

We claim
\begin{equation}
\Delta \rho+\rho_b\mathrm{div} \bm{\xi}=\mathscr{O}(\varepsilon^2). \label{Z9.1}
\end{equation}

Thanks to \eqref{*7.20}, this implies
\begin{equation}
\delta \rho= -\mathrm{div}(\rho_b\bm{\xi})+\mathscr{O}(\varepsilon^2) \label{7*drho}
\end{equation}\\

Proof of \eqref{Z9.1}. Proof is devided into several steps. 1) Let $(Q, Q_b)$ be a pair of scalar fields such that $\delta Q=\mathscr{O}(\varepsilon)$. Let $D$ be an arbitrary domain in $\mathbb{R}^3$. Put
\begin{align}
& \underline{D}(t)=\{ \underline{\bm{x}}=\bm{x}+ \bm{\xi}(t,\bm{x}) \quad |\quad \bm{x} \in D\}, \nonumber \\
&I=\int_{\underline{D}(t)}Q(t,\underline{\bm{x}})d\underline{\bm{x}},\quad I_b=\int_D Q_b(t,\bm{x})d\bm{x}, \nonumber \\
&\delta I=I- I_b.
\end{align}
Then we have
\begin{equation}
\delta I=\int_D q(t,\bm{x})d\bm{x}, 
\end{equation}
where
\begin{equation}
q(t,\bm{x})=\Delta Q(t,\bm{x})+Q(t,\underline{\bm{x}})(\mathrm{det}(I+D_{\bm{x}} \bm{\xi}(t,\bm{x}))-I).
\end{equation}
Note that
$$q=\Delta Q+Q_b\mathrm{div} \bm{\xi}+ \mathscr{O}(\varepsilon^2).
$$

2) Let $\bm{V}=\bm{v}$ or $\bm{v}_b$. Let $Q$ be a scalar field, and $E$ a domain in $\mathbb{R}^3$. Put
\begin{align}
&E(t)=\{ \bm{z}=\bm{\varphi}(t|0, \bm{y})\quad |\quad \bm{y} \in E \}, \nonumber \\
& I(t)=\int_{E(t)}Q(t,\bm{z})d\bm{z}.
\end{align}
Here we denote $\bm{\varphi}(t|s,\bm{y})=\bm{\varphi}(t|s,\bm{y}||\bm{V})$. 
Then 
\begin{equation}
 I(t)=\int_EQ(t,\bm{z})\mathrm{det}D_{\bm{y}}\bm{\varphi}(t|0,\bm{y}) d\bm{y}. \label{Z9.3}
\end{equation}
Note that
$$\mathrm{det} D_{\bm{y}}\bm{\varphi}(t|0,\bm{y})=
\exp\Big[\int_0^t(\mathrm{div}_{\bm{x}}\bm{V})(t', \bm{\varphi}(t'|0,\bm{y}))dt\Big] >0$$
and
$$
\frac{\partial}{\partial t}\mathrm{det} D_{\bm{y}}\bm{\varphi}(t|0,\bm{y})=(\mathrm{div}_{\bm{x}}\bm{V})(t,\bm{z})\cdot
\mathrm{det} D_{\bm{y}}\bm{\varphi}(t|0,\bm{y}).$$

Differentiating \eqref{Z9.3} by $t$, we get
\begin{equation}
\frac{dI}{dt}=\int_E\Big(\frac{\partial Q}{\partial t}+(\bm{V}|\nabla_{\bm{x}}Q)+Q\mathrm{div}_{\bm{x}}\bm{V}\Big)(t,\bm{z})\mathrm{det}D_{\bm{y}}\bm{\varphi} d\bm{y}
\end{equation}
Therefore, if
\begin{equation}
\frac{\partial Q}{\partial t}+(\bm{V}|\nabla_{\bm{x}} Q)+Q\mathrm{div}\bm{V}=0 \label{Z9.5}
\end{equation}
holds everywhere, then $dI/dt=0$ and
\begin{equation}
I(t)=\int_E Q(0,\bm{y})d\bm{y}.
\end{equation}
 Note that \eqref{Z9.5} holds for $Q=\rho$ or $\rho_b$ respectively with
 $\bm{V}=\bm{v}$ or $\bm{v}_b$, that is, \eqref{EPa}. 
 
 3) Let $D$ be an arbitrary domain in $\mathbb{R}^3$, and $(Q, Q_b)$ be a pair of scalar fields such that 
 \eqref{Z9.5} holds everywhere for $Q$ with $\bm{V}=\bm{v}$ and for $Q_b$ with $\bm{v}_b$.
 
 Then, keeping in mind
 \begin{align*}
 &\underline{\bm{x}}=\bm{x}+ \bm{\xi}(t,\bm{x}) =\bm{\varphi}(t|0,\bar{\bm{x}}|| \bm{v}), \\
 &\bm{x}=\bm{\varphi}(t|0,\bar{\bm{x}}||\bm{v}_b) \quad
 \Leftrightarrow\quad
 \bar{\bm{x}}=\bm{\varphi}(0|t,\bm{x}||\bm{v}_b),
 \end{align*}
 we see
 $$\bm{x} \in D \quad\Leftrightarrow\quad \bar{\bm{x}} \in \bar{D}\quad\Leftrightarrow
 \quad \underline{\bm{x}} \in D(t), $$
 where
 \begin{align*}
 &\bar{D}=\{ \bar{\bm{x}} =\bm{\varphi}(0|t,\bm{x}||\bm{v}_b)\quad|\quad \bm{x} \in D\}, \\
 &\underline{D}(t)=\{ \underline{\bm{x}}=\bm{\varphi}(t|0,\bar{\bm{x}}||\bm{v})\quad |\quad \bar{\bm{x}} \in \bar{D} \}.
 \end{align*}
 
 We consider
 $$I=\int_{\underline{D}}Q(t,\underline{\bm{x}})d\underline{\bm{x}}, \quad I_b=\int_DQ(t,\bm{x})d\bm{x}, $$
 $$\delta I=I -I_b.$$
 
 Thanks to 2) we have
 $$I=\int_{\bar{D}}Q(0,\bar{\bm{x}})d\bar{\bm{x}}, \quad
 I_b=\int_{\bar{D}}Q_b(0, \bar{\bm{x}})d\bar{\bm{x}}.$$
 Therefore
 $$\delta I=\int_{\bar{D}}(Q(0,\bar{\bm{x}})-Q_b(0,\bar{\bm{x}}))d\bar{\bm{x}}.
 $$
 Onthe other hand, thanks to 1), we have
 $$\delta I=\int_Dq(t,\bm{x})d\bm{x}.
 $$
 Hence
 $$
 \int_{\bar{D}}(Q(0,\bar{\bm{x}})-Q_b(0,\bar{\bm{x}}))d\bar{\bm{x}} =\int_Dq(t,\bm{x})d\bm{x}.
 $$
 The change of variable
 $\bar{D}\rightarrow D:\bar{\bm{x}}\mapsto
 \bm{x}=\bm{\varphi}(t|0,\bar{\bm{x}}||\bm{v}_b)$,
 for which we have $\mathrm{det}D_{\bar{\bm{x}}}\bm{x}=1$ since $\mathrm{div} \bm{v}_b=0$, we get 
 $$\int_{\bar{D}}\Big(
 Q(0,\bar{\bm{x}})-Q_b(0,\bar{\bm{x}})-q(t,\bm{x})\Big)d\bar{\bm{x}} =0. $$
 Since $D$, therefore $\bar{D}$ is arbitrary, we can claim
 $$Q(0,\bar{\bm{x}})-Q_b(0,\bar{\bm{x}})=q(t,\bm{x}) \quad \forall \bm{x} \in \mathbb{R}^3. $$
 
 4) Applying the above discussion to $Q=\rho, Q_b=\rho_b$,
 we can claim
 \begin{equation}
 \rho^0(\bar{\bm{x}})-\rho_b(\bar{\bm{x}})=\Delta\rho(t,\bm{x}) +\rho(t,\underline{\bm{x}})\Big(\mathrm{der}(I+D_{\bm{x}} \bm{\xi}(t,\bm{x}))-1\Big).
 \end{equation}
 Supposing $\rho^0=\rho_b$, we can claim
 $$\Delta \rho+\rho_b\mathrm{div} \bm{\xi}=\mathscr{O}(\varepsilon^2).
 $$ Proof is done.      $\square$\\
 
 \textbullet\  As for the equation of motion \eqref{EPb}, we devide it by $\rho$ and $\rho_b$:
 \begin{align}
 &\frac{\partial\bm{v}}{\partial t} +\sum_kv^k\frac{\partial\bm{v}}{\partial x^k}+2\bm{\Omega}\times\bm{v}+
 \bm{\Omega}\times(\bm{\Omega}\times \bm{x})+\frac{1}{\rho}\nabla P+\nabla \Phi=0, \label{Z15.1}
 \\
 &\frac{\partial\bm{v}_b}{\partial t} +\sum_kv_b^k\frac{\partial\bm{v}_b}{\partial x^k}+2\bm{\Omega}\times\bm{v}_b+
 \bm{\Omega}\times(\bm{\Omega}\times \bm{x})+\frac{1}{\rho_b}\nabla P_b+\nabla \Phi_b=0.  \label{Z15.2}
 \end{align}
 Here \eqref{Z15.1} is valid on $ \mathfrak{R}(t)$ and \eqref{Z15.2} is valid on $\mathfrak{R}_b$. So we can apply the operator
 $\Delta$ on $\mathfrak{R}_b$, for which $\underline{\bm{x}}=\bm{x}+ \bm{\xi}\in \mathfrak{R}(t)$ while $\bm{x} \in \mathfrak{R}_b$, to consider
 \begin{equation}
 \Delta\Big[
 \frac{\partial\bm{v}}{\partial t} +\sum_kv^k\frac{\partial\bm{v}}{\partial x^k}+2\bm{\Omega}\times\bm{v}+
 \bm{\Omega}\times(\bm{\Omega}\times \bm{x})+\frac{1}{\rho}\nabla P+\nabla \Phi
 \Big]=0, \label{Z15.3}
 \end{equation}
 
 In order to analyze 
 $\displaystyle \Delta\Big[\frac{\partial\bm{v}}{\partial t}+\sum_k v^k\frac{\partial \bm{v}}{\partial x^k}\Big]$,
 we introduce the operator
 $\displaystyle \frac{D}{Dt}\Big|_{\bm{V}}, \bm{V}=\bm{v}, \bm{v}_b$. Namely, 
 
 \begin{Definition} 
For a scalar field $Q$, we define
 \begin{equation}
 \frac{D}{Dt}\Big|_{\bm{V}}Q(t,\bm{x})=\partial_t Q(t,\bm{x})+(\bm{V}(t,\bm{x})|\nabla_{\bm{x}}Q(t,\bm{x})).
 \end{equation}
 Here $\partial_t=\partial/\partial t$. For a vector field $\bm{Q}$, we define
 \begin{equation}
 \frac{D}{Dt}\Big|_{\bm{V}}\bm{Q}(t,\bm{x})=\partial_t\bm{Q}(t,\bm{x})+
 (D_{\bm{x}}\bm{Q}(t,\bm{x}))\bm{V}(t,\bm{x}). \label{z40}
 \end{equation}
 \end{Definition}
 
 Let us note that 
 $$
 \Big(\frac{D}{Dt}\Big|_{\bm{V}}\bm{Q}\Big)^j=\frac{D}{Dt}\Big|_{\bm{V}}Q^j,
 $$
 since $((D_{\bm{x}}\bm{Q})\bm{V})^j=(\bm{V}|\nabla_{\bm{x}}Q^j)$.\\
 
 We claim the formula:\\
 
 {\bf Formula}: {\it For scalar field $Q$ :
 \begin{equation}
 \Delta \frac{D}{Dt}\Big|_{\bm{v}} Q=\frac{D}{Dt}\Big|_{\bm{v}_b}\Delta Q, \label{Z18.1}
 \end{equation}
 and for vector field $\bm{Q}$:
 \begin{equation}
 \Delta \frac{D}{Dt}\Big|_{\bm{v}} \bm{Q}=\frac{D}{Dt}\Big|_{\bm{v}_b}\Delta \bm{Q}. \label{Z18.2}
 \end{equation}
 }
 
 Here we mean 
 $$\Big( \frac{D}{Dt}\Big|_{\bm{v}}Q \Big)_b=\frac{D}{Dt}\Big|_{\bm{v}_b}Q_b,$$
 so we mean
 $$ \Delta \frac{D}{Dt}\Big|_{\bm{v}}Q(t,\bm{x})=
 \Big(\frac{D}{Dt}\Big|_{\bm{v}}Q\Big)(t,\bm{x}+ \bm{\xi})-\frac{D}{Dt}\Big|_{\bm{v}_b}Q_b(t,\bm{x}).$$ 
 \\
 
 Proof. It is sufficient to verify \eqref{Z18.1}. By definition we see
 \begin{align*}
 \clubsuit&:=\Delta\frac{D}{Dt}\Big|_{\bm{v}}Q(t,\bm{x})-\frac{D}{Dt}\Big|_{\bm{v}_b}\Delta Q_b(t,\bm{x}) \\
 &=(\bm{v}(t,\bm{x}+ \bm{\xi})|\nabla_{\bm{x}}Q(t, \bm{x}+ \bm{\xi})-\Big(\nabla_{\bm{x}}Q(t,\bm{x}+ \bm{\xi})\Big|\frac{\partial \bm{\xi}}{\partial t}\Big) \\
 &-(\bm{v}_b(\bm{x})|(D_{\bm{x}} \bm{\xi})^{\top}\nabla_{\bm{x}}Q(t,\bm{x}+ \bm{\xi}))
 -(\bm{v}_b(\bm{x})|\nabla_{\bm{x}}Q(t,\bm{x}+ \bm{\xi}))
 \end{align*}
 Using \eqref{Z4.2}, we have
 \begin{equation}
 \frac{\partial \bm{\xi}}{\partial t}=\bm{v}(t,\bm{x}+ \bm{\xi})-
 D_{\bm{y}}\bm{\varphi}(t|0,\bar{\bm{x}})D_{\bm{y}}\bm{\varphi}_b(0|t,\bm{x})\bm{v}_b(x). \label{Z19.1}
 \end{equation}
 Here we denote
 \begin{align*}
 &\bm{\varphi}(t|s,\bm{y})=\bm{\varphi}(t|s,\bm{y}||\bm{v}), \quad
 \bm{\varphi}_b(t|s,\bm{y})=\bm{\varphi}(t|s,\bm{y}||\bm{v}_b), \\
 &\bar{\bm{x}}=\bm{\varphi}_b(0|t,\bm{x}).
 \end{align*}
 On the other hand, we have
 \begin{equation}
 D_{\bm{x}} \bm{\xi}(t,\bm{x})=D_{\bm{y}}\bm{\varphi}(t|0,\bar{\bm{x}})D_{\bm{y}}\bm{\varphi}_b(0|t,\bm{x}) - I. \label{Z19.2}
 \end{equation}
 Hence
 \begin{align*}
 \clubsuit&=(\nabla_{\bm{x}}Q(t,\bm{x}+ \bm{\xi})|D_{\bm{y}}\bm{\varphi}(t|0,\bar{\bm{x}})D_{\bm{y}}\bm{\varphi}_b(0|t,x)\bm{v}_b(t,\bm{x})) \\
& -\Big(\bm{v}_b(\bm{x})\Big|\Big(D_{\bm{y}}\bm{\varphi}(t|0,\bar{\bm{x}})D_{\bm{y}}\bm{\varphi}_b(0|t,\bm{x})\Big)^{\top}\nabla_{\bm{x}}Q(t,\bm{x}+ \bm{\xi})\Big) \\
 &=0.
 \end{align*}
 $\square$\\
 
 { \it We have
 \begin{equation}
 \Delta \bm{v}=\frac{D}{Dt}\Big|_{\bm{v}_b} \bm{\xi}. \label{Z20.1}
 \end{equation}
 }\\
 
 Proof. Consider
 $$ \clubsuit:=\frac{D}{Dt}\Big|_{\bm{v}_b} \bm{\xi}(t,\bm{x})=
 \partial_t \bm{\xi}(t,\bm{x})+D_{\bm{x}} \bm{\xi}(t,\bm{x})\bm{v}_b(\bm{x}).
 $$
 Recalling \eqref{Z19.1}, \eqref{Z19.2}, we see
 \begin{align*}
 \clubsuit&=\bm{v}(t,\bm{x}+ \bm{\xi})- D_{\bm{y}}\bm{\varphi}(t|0,\bar{\bm{x}})D_{\bm{y}}\bm{\varphi}_b(0|t,\bm{x})\bm{v}_b(\bm{x}) \\
 &+D_{\bm{x}}\bm{\varphi}(t|0,\bar{\bm{x}})D_{\bm{y}}\bm{\varphi}_b(0|t,\bm{x})\bm{v}_b(\bm{x})-\bm{v}_b(\bm{x}) \\
 &=\bm{v}(t,\bm{x}+ \bm{\xi})-\bm{v}_b(\bm{x}) =\Delta\bm{v} (t,\bm{x}).
 \end{align*}
 $\square$\\
 
 Consequently we have
 \begin{align}
 \Delta\frac{D}{Dt}\Big|_{\bm{v}}\bm{v}&=\frac{D}{Dt}\Big|_{\bm{v}_b}\Delta\bm{v} = \Big(\frac{D}{Dt}\Big|_{\bm{v}_b}\Big)^2 \bm{\xi}\nonumber \\
 &=\frac{\partial^2 \bm{\xi}}{\partial t^2}+
 2\Big(D_{\bm{x}}\frac{\partial \bm{\xi}}{\partial t}\Big)\bm{v}_b
 +D_{\bm{x}}((D_{\bm{x}} \bm{\xi})\bm{v}_b)\bm{v}_b
 \end{align}
 and
 \begin{align}
 \Delta[2\bm{\Omega} \times \bm{v}]&=2\bm{\Omega} \times \frac{D}{Dt}\Big|_{\bm{v}_b} \bm{\xi}\nonumber \\
 &=2\bm{\Omega} \times \frac{\partial \bm{\xi}}{\partial t}
 +2\bm{\Omega} \times ((D_{\bm{x}} \bm{\xi})\bm{v}_b).
 \end{align}\\
 
 Since $\Delta \bm{x}= \bm{\xi}$ clearly, we have
 \begin{equation}
 \Delta [\bm{\Omega}\times(\bm{\Omega} \times \bm{x}]=\bm{\Omega} \times (\bm{\Omega} \times
 \bm{\xi}).
 \end{equation}.
 
 Summing up, we can write
 \begin{align}
 &\Delta\Big[
 \frac{\partial\bm{v}}{\partial t} +\sum_kv^k\frac{\partial\bm{v}}{\partial x^k}+2\bm{\Omega}\times\bm{v}+
 \bm{\Omega}\times(\bm{\Omega}\times \bm{x})\Big] = \nonumber \\
 &=\frac{\partial^2 \bm{\xi}}{\partial t^2}+2\bm{\Omega} \times \frac{\partial \bm{\xi}}{\partial t}
 +\mathcal{B}_{\blacklozenge 1}\frac{\partial \bm{\xi}}{\partial t} +
 \mathcal{L}_{\blacklozenge 2} \bm{\xi}+\mathcal{L}_{\blacklozenge 1} \bm{\xi}
 +L_{\blacklozenge 0} \bm{\xi}, 
 \end{align}
 with
 \begin{subequations}
 \begin{align}
  \mathcal{B}_{\blacklozenge 1}\dot{ \bm{\xi}}&=2(D_{\bm{x}}\dot{ \bm{\xi}})\bm{v}_b, \label{Z21a}\\
 \mathcal{L}_{\blacklozenge 2} \bm{\xi}&=D_{\bm{x}}((D_{\bm{x}} \bm{\xi})\bm{v}_b)\bm{v}_b, \label{Z21b} \\
 \mathcal{L}_{\blacklozenge 1} \bm{\xi}&=2\bm{\Omega} \times ((D_{\bm{x}} \bm{\xi})\bm{v}_b), \label{Z21c} \\
 L_{\blacklozenge 0} \bm{\xi}&=\bm{\Omega} \times(\bm{\Omega} \times \bm{\xi}). \label{Z21d}
 \end{align}
 \end{subequations}\\
 
 \textbullet\  Let us analyze 
 $\displaystyle \Delta\Big(\frac{1}{\rho}\nabla P\Big)$. It is easy to see
 $$
 \Delta\Big(\frac{1}{\rho}\nabla P\Big)=
 -\frac{\Delta\rho}{\rho_b^2}\nabla P_b
 +\frac{1}{\rho_b}\nabla \delta P+
 \frac{1}{\rho_b}(D_{\bm{x}}\nabla P_b) \bm{\xi}+\mathscr{O}(\varepsilon^2) $$
 We can write this as
 $$
 \Delta\Big(\frac{1}{\rho}\nabla P\Big)=-\frac{\delta \rho}{\rho_b^2}\nabla P_b
 -\frac{( \bm{\xi}|\nabla\rho_b)}{\rho_b^2}\nabla P_b
 +\frac{1}{\rho_b}\nabla \delta P + \frac{1}{\rho_b}(D_{\bm{x}}\nabla P_b) \bm{\xi}+ \mathscr{O}(\varepsilon^2),
 $$
 since
 $$\Delta \rho=\delta \rho+( \bm{\xi}|\nabla \rho_b)+\mathscr{O}(\varepsilon^2).$$
 Moreover we can write
 \begin{equation}
 \Delta \Big(\frac{1}{\rho}\nabla P\Big)=-\frac{\delta\rho}{\rho_b^2}\nabla P_b +\frac{1}{\rho_b}
 \nabla \delta P+D_{\bm{x}}\Big(\frac{1}{\rho_b}\nabla P_b\Big) +\mathscr{O}(\varepsilon^2), \label{Z22.1}
 \end{equation}
 since
 $$
 D_{\bm{x}}
\Big(\frac{1}{\rho_b}\nabla P_b\Big) \bm{\xi}=\frac{1}{\rho_b}(D_{\bm{x}}\nabla P_b) \bm{\xi}
 -\frac{( \bm{\xi}|\nabla \rho_b)}{\rho_b^2}\nabla P_b.
 $$\\
 
 Now let us consider $\delta P$. Since $\Delta S=0 $, we have
 \begin{equation}
 \Delta P=\Big(\frac{\partial P}{\partial \rho}\Big)_b\Delta\rho +\mathscr{O}(\varepsilon^2) \label{7*DP}
 \end{equation}
 By \eqref{*7.20}, this means
 \begin{align*}
 \delta P&=\Big(\frac{\partial P}{\partial \rho}\Big)_b\delta\rho+
 \Big(\bm{\xi}\Big|\Big(\frac{\partial P}{\partial\rho}\Big)_b\nabla \rho_b-\nabla P_b\Big) +\mathscr{O}(\varepsilon^2) \\
 &=\Big(\frac{\partial P}{\partial \rho}\Big)_b\delta\rho-
 \Big(\bm{\xi}\Big|\Big(\frac{\partial P}{\partial S}\Big)_b\nabla S_b\Big) +\mathscr{O}(\varepsilon^2) \\
 &=\Big(\frac{\partial P}{\partial \rho}\Big)_b\delta\rho+
 \Big(\frac{\partial P}{\partial S}\Big)_b \delta S +\mathscr{O}(\varepsilon^2).
 \end{align*}
 Here we have used
\eqref{7*dS}. \\

Summing up, we have
 \begin{equation}
 \Delta\Big(\frac{1}{\rho}\nabla P\Big)=\mathcal{L}_0 \bm{\xi}+\Big(D_{\bm{x}}\Big(\frac{\nabla P_b}{\rho_b}\Big)\Big) \bm{\xi}+\mathscr{O}(\varepsilon^2),
 \end{equation}
 with
 \begin{subequations}
 \begin{align}
 \mathcal{L}_0 \bm{\xi}&=-\frac{\nabla  P_b}{\rho_b^2}
 \lceil \delta\rho \rfloor
 +\frac{1}{\rho_b}\nabla \lceil \delta P \rfloor, \\
 \lceil \delta\rho\rfloor &=-\mathrm{div}(\rho_b \bm{\xi}), \\
 \lceil \delta S \rfloor &=-( \bm{\xi}|\nabla S_b), \\ 
 \lceil \delta P \rfloor &=\Big(\frac{\Gamma P}{\rho}\Big)_b\lceil \delta\rho \rfloor +
 \Big(\frac{\partial P}{\partial S}\Big)_b \lceil \delta S \rfloor.
 \end{align}
 \end{subequations} 
 Here we recall $\displaystyle \frac{\partial P}{\partial \rho}=\frac{\Gamma P}{\rho}$ by the definition of $\Gamma$, Definition \ref{Def.G}.\\
 
 Moreover we can avoid to mention $\displaystyle  \Big(\frac{\partial P}{\partial S}\Big)_b $
 in the approximation of $\delta P$. In fact, since
 \begin{align*}
 &\delta P=\Delta P -( \bm{\xi}|\nabla P_b) +\mathscr{O}(\varepsilon^2), \\
 &\Delta \rho=-\rho_b\mathrm{div} \bm{\xi}
  +\mathscr{O}(\varepsilon^2),
 \end{align*}
 \eqref{7*DP} implies
 \begin{equation}
 \delta P=-(\Gamma P)_b\mathrm{div} \bm{\xi}-
 ( \bm{\xi}|\nabla P_b) + \mathscr{O}(\varepsilon^2).
 \end{equation}\\
 
 We note that \cite{FriedmanS} gives another expression of $\displaystyle \Delta \Big(\frac{1}{\rho}\nabla P\Big))$, namely,
 \begin{equation}
 \Delta \Big(\frac{1}{\rho}\nabla P\Big))=\frac{1}{\rho_b}\Big[
 -\nabla(\Gamma_bP_b\mathrm{div} \bm{\xi})-
 ((D_{\bm{x}} \bm{\xi})\nabla P_b)^{\top}
 +(\mathrm{div} \bm{\xi})\nabla P_b\Big]+
 \mathscr{O}(\varepsilon^2). \label{Z24}
 \end{equation}
 
 This expression \eqref{Z24} can be derived by applying the formula
 \begin{align}
 \Delta\nabla Q&=\nabla \Delta Q+(D\nabla Q) \bm{\xi}
 -\nabla( \bm{\xi}|\nabla Q) +\mathscr{O}(\varepsilon^2) \nonumber \\
 &=\nabla\Delta Q -((D \bm{\xi})\nabla Q)^{\top} +\mathscr{O}(\varepsilon^2)
 \end{align}
 to $Q=P$ and noting
 $$
 \Delta P=\Big(\frac{\Gamma P}{\rho}\Big)_b\Delta \rho +\mathscr{O}(\varepsilon^2) =\Gamma_bP_b\mathrm{div} \bm{\xi}+ \mathscr{O}(\varepsilon^2).
$$
and
$$ \frac{\Delta \rho}{\rho_b}=-\mathrm{div} \bm{\xi}+ \mathscr{O}(\varepsilon^2).
$$\\

\textbullet\  Let us look at
$$
\Delta\nabla \Phi(\bm{x})=\nabla \Phi(\bm{x}+ \bm{\xi})-\nabla \Phi_b(\bm{x}), $$
where
\begin{align*}
&\Phi=-4\pi\mathsf{G}\mathcal{K}[\rho], \qquad \Phi_b=-4\pi\mathsf{G}\mathcal{K}[\rho_b], \\
&\mathcal{K}[g](\bm{x})=\frac{1}{4\pi}\int_{\mathbb{R}^3}\frac{g(\bm{x}')}{\|\bm{x}-\bm{x}'\|}d\bm{x}'.
\end{align*}
We introduce the integral operator $\mathcal{K}^{\nabla}$ by
\begin{equation}
\mathcal{K}^{\nabla}[g](\bm{x})=-\frac{1}{4\pi}\int_{\mathbb{R}^3}\frac{\bm{x}-\bm{x}'}{\|\bm{x}-\bm{x}'\|^3}g(\bm{x}')d\bm{x}'.
\end{equation}
When $g \in C^{\alpha}_0(\mathbb{R}^3), 0<\alpha<1$, then 
$\mathcal{K}^{\nabla}[g] \in C^{1+\alpha}(\mathbb{R}^3)$ and
$\nabla \mathcal{K}[g]=\mathcal{K}^{\nabla}[g]$. So
\begin{align*}
\Delta\nabla \Phi(\bm{x})&=-4\pi\mathsf{G}\Big(\mathcal{K}^{\nabla}[\delta\rho](\underline{\bm{x}})
+\nabla\Phi_b(\underline{\bm{x}})-\nabla\Phi_b(\bm{x})\Big) \\
&=-4\pi\mathsf{G}\mathcal{K}^{\nabla}[\delta\rho]+(D_{\bm{x}}\nabla\Phi_b(\bm{x})) \bm{\xi}
+\mathscr{O}(\varepsilon^2).
\end{align*}
Here $\underline{\bm{x}}=\bm{x}+ \bm{\xi}$, and we note 
$\mathcal{K}^{\nabla}[\delta\rho]\in C^{1+\alpha}$ and $\mathscr{O}(\varepsilon)$.

Consequently we can claim
\begin{equation}
\Delta\nabla\Phi=4\pi\mathsf{G}\mathcal{L}_1 \bm{\xi}+
(D_{\bm{x}}\nabla_{\bm{x}}\Phi_b) \bm{\xi}+\mathscr{O}(\varepsilon^2)
\end{equation}
with
\begin{equation}
\mathcal{L}_1 \bm{\xi}=\nabla\mathcal{K}[\mathrm{div}(\rho_b \bm{\xi})].
\end{equation}\\

\textbullet\  Summing up, we have derived the equation used by astrophysicists:

\begin{equation}
\frac{\partial^2 \bm{\xi}}{\partial t^2}+\mathcal{B}_{\blacklozenge}\frac{\partial \bm{\xi}}{\partial t} + \mathcal{L}_{\blacklozenge} \bm{\xi}=0, \label{7.76}
\end{equation}
where
\begin{equation}
\mathcal{B}_{\blacklozenge}\dot{ \bm{\xi}}=2\bm{\Omega} \times \dot{ \bm{\xi}}+
 \mathcal{B}_{\blacklozenge 1}\dot{ \bm{\xi}}, 
 \end{equation}
 with
 \begin{equation}
   \mathcal{B}_{\blacklozenge 1}\dot{ \bm{\xi}}=2(D\dot{ \bm{\xi}})\bm{v}_b,
  \end{equation}
and

  \begin{equation} 
 \mathcal{L}_{\blacklozenge}=\mathcal{L}_{\blacklozenge 2}+\mathcal{L}_{\blacklozenge 1}+
 L_{\blacklozenge 0}+\mathcal{L}_0+L_{0\blacklozenge}+
 4\pi\mathsf{G}\mathcal{L}_1+L_{1\blacklozenge}
 \end{equation}
  with
 \begin{subequations}
 \begin{align} 
 \mathcal{L}_{\blacklozenge 2} \bm{\xi}&=D((D \bm{\xi})\bm{v}_b)\bm{v}_b
 ,   \\
 \mathcal{L}_{\blacklozenge 1} \bm{\xi}&=2\bm{\Omega} \times ((D \bm{\xi})\bm{v}_b),   \\
  L_{\blacklozenge 0} \bm{\xi}&=\bm{\Omega} \times(\bm{\Omega} \times \bm{\xi}), \\
\mathcal{L}_0 \bm{\xi}&=-\frac{\lceil \delta\rho \rfloor}{\rho_b^2} +
\frac{1}{\rho_b}\nabla \lceil \delta P \rfloor, \\
\mbox{with}& \nonumber \\
 \lceil \delta\rho \rfloor&=-\mathrm{div}(\rho_b \bm{\xi}), \quad \lceil \delta S \rfloor =-( \bm{\xi}|\nabla S_b), \nonumber \\
\lceil \delta P \rfloor &=\Big(\frac{\Gamma P}{\rho}\Big)_b
\lceil \delta\rho \rfloor
+\Big(\frac{\partial P}{\partial S}\Big)_b\lceil \delta S \rfloor
=\Big(\frac{\Gamma P}{\rho}\Big)_b \lceil \delta\rho \rfloor+ 
 \Gamma_bP_b( \bm{\xi} |\mathfrak{a}_b)
\\
L_{0\blacklozenge} \bm{\xi}&=\Big(D\Big(\frac{\nabla P_b}{\rho_b}\Big)\Big) \bm{\xi}, \\
\mathcal{L}_1 \bm{\xi}&=-\nabla\mathcal{K}[ \lceil \delta\rho \rfloor   ]=
\nabla \mathcal{K}[\mathrm{div}(\rho_b \bm{\xi})], \\
L_{1\blacklozenge} \bm{\xi}&=(D\nabla \Phi_b) \bm{\xi}.
\end{align}
\end{subequations}

We denote 
$$\mathcal{L}=\mathcal{L}_0+4\pi\mathsf{G}\mathcal{L}_1.$$\\

However we would like to rewrite $\mathcal{L}_{\blacklozenge}$ further. Since $(\rho_b,S_b, \bm{v}_b)$ satisfies the equation of motion
$$(D\bm{v}_b)\bm{v}_b+2\bm{\Omega}\times\bm{v}_b+\bm{\Omega}\times(\bm{\Omega}\times\bm{x})+
\frac{\nabla P_b}{\rho_b}
+\nabla \Phi_b=0 $$
on $\mathfrak{R}_b$, we see
$$L_{0\blacklozenge} \bm{\xi}+L_{1\blacklozenge} \bm{\xi}=-D\Big((D\bm{v}_b)\bm{v}_b+
2\bm{\Omega}\times \bm{v}_b+
\bm{\Omega}\times(\Omega\times\bm{x}\Big) \bm{\xi}.
$$
But we have
$$D\Big(\bm{\Omega}\times(\bm{\Omega}\times \bm{x})\Big) \bm{\xi}=\bm{\Omega}\times(\bm{\Omega}\times \bm{\xi}).$$
Thus

\begin{equation}
L_{0\blacklozenge}+L_{1\blacklozenge}+L_{\blacklozenge 0}=\mathcal{L}_{\blacklozenge 00}+\mathcal{L}_{\blacklozenge 01},
\end{equation}
where
\begin{equation}
\mathcal{L}_{\blacklozenge 00} \bm{\xi}=
-\Big(D((D\bm{v}_b)\bm{v}_b)\Big) \bm{\xi},
\quad \mathcal{L}_{\blacklozenge 01} \bm{\xi}=-\Big(D(2\bm{\Omega} \times \bm{v}_b)\Big) \bm{\xi}.
\end{equation}

We define
\begin{equation}
\mathcal{L}_{\blacklozenge 0}:=\mathcal{L}_{\blacklozenge 00}+\mathcal{L}_{\blacklozenge 01}.
\end{equation}

Consequently we have
 the equation \eqref{ZLASO} with  \eqref{7.73} \eqref{7.74} \eqref{7.75} with \eqref{7.76a} - \eqref{ZLASO7.9f}, as the equation of linear adiabatic stellar oscillations  used by astrophysicists after \cite{Lebovitz} \cite{FriedmanS} \cite{LyndenBO}.\\

We note that the derivation above does not require that $\bm{v}_b$ is small, say, $\mathscr{O}(\varepsilon)$ during the linear approximation procedure, while we supposed 
$\bm{v}_b=\mathscr{O}(\varepsilon)$ to derive {\bf ELASO}  \eqref{Eqxi} in Section 3. 
But we recall that \eqref{Z8.1}: $\rho^0=\rho_b$ and \eqref{Z14.1}: $S^0=S_b$ are supposed, while only \eqref{3.1}: $\{ \rho^0>0\}=\{ \rho_b >0 \}$ was supposed in Section 3. Moreover we note that $ \bm{\xi}(0,\bm{x})=\bm{0}$ by the definition
\eqref{z7}, while the Lagrange displacement
$\bm{u}$ in our equation \eqref{Eqxi} which was defined by
\eqref{1.16} or \eqref{3.8}:
$$ \bm{u}(t, \bar{\bm{x}})=\bm{\varphi}(t|0,\bar{\bm{x}}|| \bm{v})
-\bm{\varphi}(t|0,\bar{\bm{x}} ||\bm{v}_b)+
\bm{u}^0(\bar{\bm{x}}) $$
enjoys $\bm{u}(0,\bar{\bm{x}})=\bm{u}^0(\bar{\bm{x}})$, which $\mathscr{O}(\varepsilon)$ was arbitrary as long as
$$\rho^0-\rho_b=-\mathrm{div}(\rho_b\bm{u}^0),\quad
S^0-S_b=-(\nabla S_b|\bm{u}^0) $$ 
and can be $\not= \bm{0}$.\\

\section{Analysis of {\bf ELASO($\blacklozenge$) } }

\subsection{Trivial displacements}

 By the definition \eqref{z7} it holds that $ \bm{\xi}(0,\bm{x})=0$. But it is possible to consider solutions $\tilde{ \bm{\xi}}= \bm{\xi}+\bm{\eta}$ of ELASO
 ($\blacklozenge$) \eqref{ZLASO}. After \cite{FriedmanS} let us call $\bm{\eta}(t,\bm{x})$ a trivial displacement if it satisfies
\begin{align}
&0=-\mathrm{div}(\rho_b \bm{\eta}), \label{Z1} \\
&0=-(\nabla S_b|\bm{\eta}), \label{Z2} \\
&\bm{0}=\partial_t\bm{\eta} +(D\bm{\eta})\bm{v}_b-(D\bm{v}_b)\bm{\eta}. \label{Z3}
\end{align}
Here $\partial_t=\partial/\partial t$. Then
\begin{align}
&\lceil \delta\rho ||\bm{\eta} \rfloor:=-\mathrm{div}(\rho_b \bm{\eta})=0, \\
&\lceil \delta S || \bm{\eta} \rfloor:=-(\nabla S_b|\bm{\eta})=0, \\
&\lceil \delta P||\bm{\eta} \rfloor:=
\Big(\frac{\Gamma P}{\rho}\Big)_b\lceil \delta\rho || \bm{\eta} \rfloor
-\Big(\frac{\partial P}{\partial S}\Big)_b(\nabla S_b|\bm{\eta})=0, \\
&\lceil \delta\bm{v} ||\bm{\eta} \rfloor:=\partial_t\bm{\eta} +(D\bm{\eta})\bm{v}_b-(D\bm{v}_b)\bm{\eta}=\bm{0}.
\end{align} 
Here we note that
\begin{align*}
\delta\bm{v}&=\Delta \bm{v} -(D\bm{v}_b) \bm{\xi}+\mathscr{O}(\varepsilon^2) \quad (\mbox{see}\quad \eqref{z19}), \\
&=\frac{D}{Dt}\Big|_{\bm{v}_b} \bm{\xi}-(D\bm{v}_b) \bm{\xi}+ \mathscr{O}(\varepsilon^2) \quad( \mbox{see} \quad \eqref{Z20.1}) \\
&=\partial_t \bm{\xi}+ (D \bm{\xi})\bm{v}_b
-(D\bm{v}_b) \bm{\xi}+ \mathscr{O}(\varepsilon^2), \quad(\mbox{see}\quad \eqref{z40}),
\end{align*}
so 
$$
\lceil \delta\bm{v} || \bm{\xi}\rfloor =\partial_t \bm{\xi}+ (D \bm{\xi})\bm{v}_b-(D\bm{v}_b) \bm{\xi}
$$ approximates $\delta \bm{v}$.\\

Let $\bm{\eta}(t,\bm{x})$ be a trivial displacement. Then \eqref{Z1} and \eqref{Z3}  imply $\mathcal{L}_0\bm{\eta}=\bm{0}$ and $\mathcal{L}_1\bm{\eta}=\bm{0}$. Moreover we have
\begin{equation}
\partial_t^2\bm{\eta}+\mathcal{B}_{\blacklozenge 1}\partial_t\bm{\eta}+
(\mathcal{L}_{\blacklozenge 2}+\mathcal{L}_{\blacklozenge 00})\bm{\eta}=\bm{0}, \label{Z4}
\end{equation}
and
\begin{equation}
2\bm{\Omega} \times \partial_t\bm{\eta}+
(\mathcal{L}_{\blacklozenge 1}+\mathcal{L}_{\blacklozenge 01})\bm{\eta}=\bm{0}. \label{Z5}
\end{equation}
Therefore we can claim that $\bm{\eta}$ satisfies the equation \eqref{ZLASO}, so $\tilde{ \bm{\xi}}= \bm{\xi}+\bm{\eta}$ should be a solution of \eqref{ZLASO}. 

Proof of \eqref{Z4}. Eliminating $\eta_t\bm{\eta}, \partial_t^2\bm{\eta}$ by \eqref{Z3}, we see that the left-hand side of \eqref{Z4} reduces to
\begin{align*}
&(D((D\bm{v}_b)\bm{\eta})-(D\bm{v}_b)(D\bm{\eta}))\bm{v}_b
-((D\bm{v}_b)(D\bm{v}_b)-D((D\bm{v}_b)\bm{v}_b))\bm{\eta} \\
&=\langle D^2\bm{v}_b|\bm{\eta},\bm{v}_b\rangle -\langle D^2\bm{v}_b|\bm{\eta},\bm{v}_b\rangle =\bm{0}.
\end{align*}
Here we denote
\begin{equation}
\langle D^2\bm{f} | \bm{X}, \bm{Y}\rangle^i=\sum_{k.\ell}\frac{\partial^2 f^i}{\partial x^k \partial x^{\ell}}X^kY^{\ell}.
\end{equation}
$\square$

Proof of \eqref{Z5}. Eliminating $\partial_t\bm{\eta}$ by \eqref{Z3}, we see that
the left-hand side of \eqref{Z5} reduces to
$$
2\bm{\Omega}\times (D\bm{v}_b)\bm{\eta} -(D(2\bm{\Omega} \times \bm{v}_b))\bm{\eta} =\bm{0}.
$$

$\square$ \\

Of course $\bm{\eta}=\bm{0}$ is a trivial displacement, which could be called the trivial trivial displacement, and we are interested in non-trivial trivial displacements in this sense. 

Let us consider 
\begin{equation}
\bm{\eta}=\frac{1}{\rho_b}\nabla S_b \times \nabla f, \label{Z6}
\end{equation}
where $f$ is a scalar fields wich satisfies
\begin{equation}
\frac{D}{Dt}\Big|_{\bm{v}_b}f=\partial_tf +(\bm{v}_b|f)=0. \label{Z7}
\end{equation}
Then $\eta$ is a trivial displacement.

Proof. Clearly \eqref{Z1}, \eqref{Z2} are satisfied, for any $f$. Usng the properties of the background
\begin{equation}
\sum_k\bm{v}_b\partial_k\rho_b=0, \qquad \sum_k\bm{v}_b\partial_kS_b=0, \label{Z8}
\end{equation}
where $\partial_k=\partial/\partial x^k$, we can verify
\begin{equation}
((D\bm{\eta})\bm{v}_b-(D\bm{v}_b)\bm{\eta})^j=\frac{1}{\rho_b}\Big(
-\sum_{\alpha,\beta}E^{j\alpha\beta}(\partial_{\alpha}S_b)(\partial_{\beta}f)
+\nabla S_b \times \nabla(\bm{v}_b|f) \Big), \label{9}
\end{equation}
where
\begin{equation}
E^{j\alpha\beta}=\sum_k\Big(
\epsilon^{jk\beta}\partial_kv_b^{\alpha}+\epsilon^{j\alpha k}\partial_kv_b^{\beta}+\epsilon^{k\alpha\beta}\partial_kv_b^j \Big),\quad
\epsilon^{\lambda\mu\nu}=
\mathrm{sign}\left(\begin{array}{rrr}
1& 2& 3 \\
\lambda &\mu&\nu
\end{array}\right).
\end{equation}
But we can verify
\begin{equation}
E^{j\alpha \beta}=\epsilon^{j\alpha\beta}\sum_k\partial_k\bm{v}_b^k,
\end{equation}
therefore $E^{j\alpha\beta}=0$ for $\forall j,\alpha, \beta$,
since $\mathrm{div}\bm{v}_b=0$. Thus
$$(D\bm{\eta})\bm{v}_b
-(D\bm{v}_b)\bm{\eta}=\frac{1}{\rho_b}\nabla S_b \times \nabla(\bm{v}_b|f)$$
so that \eqref{Z3} is satisfied, since
$$\partial_t\bm{\eta}=
\frac{1}{\rho_b}\nabla S_b \times \nabla \partial_t f.
$$ 
$\square$

General solution of \eqref{Z7} is
\begin{equation}
f(t,\bm{x})=f^0(\bm{\varphi}(0|t,\bm{x}|| \bm{v}_b)), \label{Z13}
\end{equation}
where $\bar{\bm{x}} \mapsto f^0(\bar{\bm{x}})$ is an arbitrary function on $\mathfrak{R}_b$.Then
\begin{equation}
\bm{\eta}(0, \bm{x})=\frac{1}{\rho_b(\bm{x})}
\nabla S_b \times \nabla_{\bar{\bm{x}}}f^0(\bm{x}).
\end{equation}
When $\nabla S_b \not= \bm{0}$, 
that is, $S_b$ is not a constant, then this construction gives many non-trivial trivial solutions.
However when $\nabla S_b=\bm{0}$, $\bm{\eta}$ given by \eqref{Z6} is trivial, $=\bm{0}$.
If $\bm{v}_b=\mathscr{O}(\varepsilon)$ the displacement given by
\begin{equation}
\bm{\eta}(t,\bm{x})=\frac{1}{\rho_b(\bar{\bm{x}})}\mathrm{rot}_{\bar{\bm{x}}}\bm{f}(\bar{\bm{x}}),\quad \bar{\bm{x}}=\bm{\varphi}(0|t,\bm{x}||\bm{v}_b)
\end{equation}
satisfies \eqref{Z1}, \eqref{Z3} approximately with
errors $\mathscr{O}(\varepsilon)$, provided that the function $\bm{f}$ is small
$\mathscr{O}(\varepsilon)$. 
We compare these observations with the discussion of Theorem \ref{Th.*4} in Section 4.
\\

\subsection{Functional analysis of ELASO($\blacklozenge$) }

  Let us consider the equation \eqref{ZLASO} in the functional space
$\mathfrak{H}=L^2(\mathfrak{R}_b, \rho_bd\bm{x}; \mathbb{C}^3)$. We suppose\\

{\bf Assumption} : {\it The equation of state enjoys {\bf Assumtion} \ref{Ass.0} and the background $(\rho_b, S_b, \bm{v}_b)$ is barotropic.}   \\

The operator $\mathcal{B}:\dot{ \bm{\xi}} \mapsto 2\bm{\Omega} \times \dot{ \bm{\xi}}$ can be considered as a bounded linear operator in $\mathcal{B}(\mathfrak{H})$, and $\mathcal{L}=\mathcal{L}_0+4\pi\mathsf{G}\mathcal{L}_1$ can be extended to a self-adjoint operator bounded from below in $\mathfrak{H}$, by following the discussion of Section 4. We look at the resting terms
$\mathcal{B}_{\blacklozenge 1}, \mathcal{L}_{\blacklozenge2}, \mathcal{L}_{\blacklozenge 1}, 
\mathcal{L}_{\blacklozenge 0} $.\\

First, we see that $\mathcal{B}_{\blacklozenge 1}$ is skew-symmetric considered on $C_0^{\infty}(\mathfrak{R}_b)$. In fact,
for $\bm{v}_1, \bm{v}_2 \in C_0^{\infty}$, we have
\begin{align*}
 &(\mathcal{B}_{\blacklozenge 1}\bm{v}_1|\bm{v}_2)_{\mathfrak{H}}=
 2\int \omega_b\Big(\frac{\partial \bm{v}_1}{\partial\phi}\Big|\bm{v}_2\Big)\rho_bd\bm{x} = \\
 &=-2\int \omega_b\Big(\bm{v}_1\Big|\frac{\partial \bm{v}_2}{\partial\phi}\Big)\rho_bd\bm{x} = -
 (\bm{v}_1|\mathcal{B}_{\blacklozenge 1}\bm{v}_2)_{\mathfrak{H}}.
 \end{align*}

However this operator is not bounded in $\mathfrak{H}$-norm, if $\bm{v}_b\not= \bm{0}$, say, $\omega_b \not=0$. \\

$\mathcal{L}_{\blacklozenge 2}, \mathcal{L}_{\blacklozenge 1}$ are symmetric on $C_0^{\infty}(\mathfrak{R}_b)$. In fact
we have,
for $\bm{u}_1, \bm{u}_2 \in C_0^{\infty}$, 
\begin{align*}
(\mathcal{L}_{\blacklozenge 2}\bm{u}_1|\bm{u}_2)_{\mathfrak{H}}&=
-\int \omega_b^2\Big(\frac{\partial\bm{u}_1}{\partial\phi}\Big|\frac{\partial\bm{u}_2}{\partial\phi}\Big)\rho_bd\bm{x} \\
&=(\bm{u}_1|\mathcal{L}_{\blacklozenge 2}\bm{u}_2)_{\mathfrak{H}}, \\
(\mathcal{L}_{\blacklozenge 1}\bm{u}_1|\bm{u}_2)_{\mathfrak{H}}&=2\Omega\int
\omega_b\Big(\frac{\partial u_1^1}{\partial\phi}(u_2^2)^*+
u_1^2\Big(\frac{\partial u_2^1}{\partial \phi}\Big)^*\Big)\rho_bd\bm{x} \\
&=(\bm{u}_1|\mathcal{L}_{\blacklozenge1}\bm{u}_2)_{\mathfrak{H}}.
\end{align*}
However $\mathcal{L}_{\blacklozenge 2}+\mathcal{L}_{\blacklozenge 1}$ cannot be bounded from below, if $\bm{v}_b\not=0$. In fact,
we see
\begin{align}
(\mathcal{L}_{\blacklozenge 2}\bm{u}|\bm{u})_{\mathfrak{H}}&=-\int \omega_b^2\Big\|\frac{\partial\bm{u}}{\partial\phi}\Big\|^2\rho_bd\bm{x}, \\
(\mathcal{L}_{\blacklozenge 1}\bm{u}|\bm{u})_{\mathfrak{H}}&=
4\Omega\int
\omega_b\mathfrak{Re}\Big[\frac{\partial u^1}{\partial {\phi}}(u^2)^*\Big]\rho_bd\bm{x}.
\end{align}
Therefore
$$((\mathcal{L}_{\blacklozenge 2}+\mathcal{L}_{\blacklozenge 1})\bm{u}|\bm{u})_{\mathfrak{H}}=
-m^2\int \omega_b^2\|\bm{u}\|^2\rho_bd\bm{x},
$$
when
 
$$\bm{u}(\bm{x})=\bm{U}(r, \vartheta)e^{\mathrm{i}m\phi }, \quad
 \bm{U}=
 \begin{bmatrix}
 U \\
  U \\
  0
  \end{bmatrix}
  , \quad \mbox{while}\quad
U \in C_0^{\infty}(\mathfrak{R}_b; \mathbb{R}).$$

Suppose $\omega_b \not=0$, and
 we suppose $\omega_b^2 \geq \kappa^2 >0$, $\kappa$ being a constant, on a domain $\mathcal{V}(\not= \emptyset), \mathrm{supp}[U] \subset 
\mathcal{V}$. Then we have
$$
((\mathcal{L}_{\blacklozenge 2}+\mathcal{L}_{\blacklozenge 1})\bm{u}|\bm{u})_{\mathfrak{H}}\leq -m^2\kappa^2 \|\bm{u}\|_{\mathfrak{H}}^2.
$$

Therefore
$$
\inf_{\bm{u} \in C_0^{\infty}}
\frac{((\mathcal{L}_{\blacklozenge 2}+\mathcal{L}_{\blacklozenge 1})\bm{u}|\bm{u})_{\mathfrak{H}}}{\|\bm{u}\|_{\mathfrak{H}}^2 }
\leq -m^2\kappa^2
$$
for $\forall m$ so that
$$\inf_{\bm{u} \in C_0^{\infty}}
\frac{((\mathcal{L}_{\blacklozenge 2}+\mathcal{L}_{\blacklozenge 1})\bm{u}|\bm{u})_{\mathfrak{H}}}{\|\bm{u}\|_{\mathfrak{H}}^2 }
=-\infty.$$
  Thus
$\mathcal{L}_{\blacklozenge 2}+\mathcal{L}_{\blacklozenge 1} $ is not bounded from below, if $\omega_b\not=0$.\\

Finally, as for the operator $\mathcal{L}_{\blacklozenge 0}=\mathcal{L}_{\blacklozenge 00}+\mathcal{L}_{\blacklozenge 01}$, it is a simple multiplication by the real-valued matrix function $W_b \in 
 C(\mathbb{R}^3; \mathbb{R}^{3\times 3})$,
therefore $\mathcal{L}_{\blacklozenge 0}$
can be considered as a symmetric bounded operator on $\mathfrak{H}$.\\

Summing up, when $\bm{v}_b\not=0$, or at least, when
$\bm{v}_b \not =\mathscr{O}(\varepsilon)$, by which we cannot neglect terms $\mathcal{B}_{\blacklozenge 1}\partial_t \bm{\xi},
\mathcal{L}_{\blacklozenge 2} \bm{\xi}, \mathcal{L}_{\blacklozenge 1} \bm{\xi}, \mathcal{L}_{\blacklozenge 0} \bm{\xi}$ as $\mathscr{O}(\varepsilon^2)$, the application of the Hille-Yosida theory as discussed in Section 4 does not work
for the unboundedness of $\mathcal{B}_{\blacklozenge 1}$, which implies the unboundedness of the total $\mathcal{B}_{\blacklozenge}$, and the unboundedness from below of
$\mathcal{L}_{\blacklozenge 2}
+\mathcal{L}_{\blacklozenge 1}$, which implies the unboundedness from below of the total $\mathcal{L}_{\blacklozenge}$. However,
if we can suppose that $\bm{v}_b =\mathscr{O}(\varepsilon)$, then the terms $ 
\mathcal{B}_{\blacklozenge 1}\partial_t \bm{\xi},
\mathcal{L}_{\blacklozenge 2} \bm{\xi}, \mathcal{L}_{\blacklozenge 1} \bm{\xi}, \mathcal{L}_{\blacklozenge 0} \bm{\xi}$ are of $\mathscr{O}(\varepsilon^2)$ and negligible. In this situation, $\mathcal{B}_{\blacklozenge}, \mathcal{L}_{\blacklozenge}$ reduce
to $\mathcal{B}$ and $\mathcal{L}=\mathcal{L}_0+4\pi\mathsf{G}\mathcal{L}_1$, that is, the astrophysicists' {\bf ELASO ($\blacklozenge$)} \eqref{ZLASO} reduces to our {\bf ELASO} \eqref{Eqxi}.\\

\textbullet\  However we want to say something for the initial value problem to {\bf ELASO ($\blacklozenge$)}:
\begin{align}
& \frac{\partial^2\bm{\xi}}{\partial t^2}+\mathcal{B}_{\blacklozenge}\frac{\partial\bm{\xi}}{\partial t}+\mathcal{L}_{\blacklozenge}\bm{\xi}=\bm{f}(t,\bm{x}), \label{IVPbl.1} \\
&\bm{\xi}|_{t=0}=\bm{\xi}^0,\quad \frac{\partial\bm{\xi}}{\partial t}\Big|_{t=0}=\bm{v}^0, \label{IVPbl.2}
\end{align}
where $\bm{f}, \bm{\xi}^0, \bm{v}^0$ are given data. Suppose that the background is a barotropic admissible stationary solution.

In order to fix the idea to consider the integro-differential operator $\mathcal{L}_{\blacklozenge}$ in the Hilbert space
$\mathfrak{H}=L^2(\mathfrak{R}_b, \rho_bd\bm{x}; \mathbb{C}^3\}$,
we  define $\bm{L}_{\blacklozenge}$ 
 as
\begin{align}
&\mathsf{D}(\bm{L}_{\blacklozenge})=\{ \bm{u} \in (\mathfrak{G}_{\blacklozenge})_0 \quad |\quad \mbox{Both }\quad  \mathcal{L}_{\blacklozenge 2}\bm{u} \quad\mbox{and}\quad
\mathcal{L}\bm{u}
  \nonumber \\
& \mbox{ belong to} \quad  \mathfrak{H} 
\quad\mbox{in distribution sense}\}, \\
&\bm{L}_{\blacklozenge}\bm{u}=\mathcal{L}_{\blacklozenge}\bm{u} \quad\mbox{for}\quad \bm{u} \in 
\mathsf{D}(\bm{L}_{\blacklozenge}),
\end{align}
where $\mathfrak{G}_{\blacklozenge}$ stands for the functional space 
\begin{equation}
\mathfrak{G}_{\blacklozenge}=\Big\{ \bm{u} \in \mathfrak{H} \quad  \Big| \mathrm{div}(\rho_b\bm{u}) \in L^2\Big(\Big(\frac{\Gamma P}{\rho^2}\Big)_bd\bm{x}\Big), \quad \omega_b\frac{\partial\bm{u}}{\partial\phi} \in \mathfrak{H}\Big\}
\end{equation}
endowed with the norm
\begin{equation}
\|\bm{u}\|_{\mathfrak{G}_{\blacklozenge}}=\Big[
\|\bm{u}\|_{\mathfrak{H}}^2+\Big\|\mathrm{div}(\rho_b\bm{u}); L^2\Big(\Big(\frac{\Gamma P}{\rho^2}\Big)_bd\bm{x}\Big)\Big\|^2+
\Big\|\omega_b\frac{\partial \bm{u}}{\partial\phi}\Big\|_{\mathfrak{H}}^2 \Big]^{\frac{1}{2}},
\end{equation}
and 
$(\mathfrak{G}_{\blacklozenge})_0$ stands for the closure of $C_0^{\infty}(\mathfrak{R}_b)$
in $\mathfrak{G}_{\blacklozenge}$.

Note that, if $\bm{u} \in \mathfrak{G}_{\blacklozenge}$, then we have, a priori, $\mathcal{L}_{\blacklozenge 1}\bm{u} \in \mathfrak{H}$ and
$\mathcal{L}_{\blacklozenge 0}\bm{u} \in \mathfrak{H}$. Thus, if
$\bm{u} \in \mathsf{D}(\bm{L}_{\blacklozenge})$, then $\mathcal{L}_{\blacklozenge} \bm{u} \in \mathfrak{H} $ in distribusion sense. But, if $ \bm{v}_b \not=0$, we may not say that $\mathcal{L}_{\blacklozenge}\bm{u} \in \mathfrak{H}$ guarantees
$\mathcal{L}_{\blacklozenge 2}\bm{u}, \mathcal{L}\bm{u} \in \mathfrak{H}$ separately.\\

This formulation comes from the following observation.\\

Put, for $\bm{u} \in \mathfrak{G}_{\blacklozenge}$, 
\begin{align}
Q_{\blacklozenge 2}[\bm{u}]&:=-\int \omega_b^2\Big\|\frac{\partial \bm{u}}{\partial \phi}\Big\|^2\rho_bd\bm{x}
=-\Big\|\omega_b\frac{\partial\bm{u}}{\partial\phi}\Big\|_{\mathfrak{H}}^2
, \\
Q_{\blacklozenge 1}[\bm{u}]&:=4\Omega\int \omega_b\mathfrak{Re}\Big[
\frac{\partial u^1}{\partial\phi}(u^2)^*\Big]\rho_bd\bm{x}, \\
Q_{\blacklozenge 0}[\bm{u}]&:=-\int (W_b\bm{u}|\bm{u})\rho_bd\bm{x}, \\
Q_0[\bm{u}]&:=
 \int \frac{\Gamma P}{\rho^2}|\mathrm{div}(\rho \bm{u})|^2 d\bm{x} +2\mathfrak{Re}\Big[
\int \frac{\Gamma P \mathscr{A} }{\rho}(\bm{u}|\bm{n})\mathrm{div}(\rho \bm{u})^*d\bm{x} \Big] + \nonumber \\
&-\int \frac{\Gamma P \mathscr{A} }{\rho}|(\bm{u}|\bm{n})|^2\|{\nabla}\rho\|d\bm{x}
4\pi\mathsf{G}\int \mathcal{K}[\mathrm{div}(\rho\bm{u})]\mathrm{div}(\rho \bm{u})^*d\bm{x}.
\end{align}

Then we have the estimates

\begin{align}
|Q_{\blacklozenge 1}[\bm{u}]|&\leq 2|\Omega|\Big(\|\bm{u}\|_{\mathfrak{H}}^2+
\Big\|\omega_b\frac{\partial\bm{u}}{\partial\phi}\Big\|^2_{\mathfrak{H}}\Big), \\
|Q_{\blacklozenge 0}[\bm{u}]|&\leq \|W_b\|_{L^{\infty}}\|\bm{u}\|_{\mathfrak{H}}^2.
\end{align}

Put, for $\bm{u} \in \mathfrak{G}_{\blacklozenge}$,
\begin{equation}
Q_0^{\blacklozenge}[\bm{u}]:=\sum_{\nu=0,1,2}Q_{\blacklozenge \nu}[\bm{u}]+Q_0[\bm{u}],
\end{equation}

Then, if $\bm{u} \in C_0^{\infty}$, we have
\begin{align}
Q_{\blacklozenge 2}[\bm{u}]&=(\mathcal{L}_{\blacklozenge 2}\bm{u}|\bm{u})_{\mathfrak{H}}, \\
Q_{\blacklozenge 1}[\bm{u}]&=(\mathcal{L}_{\blacklozenge 1}\bm{u}|\bm{u})_{\mathfrak{H}}, \\
Q_{\blacklozenge 0}[\bm{u}]&=(\mathcal{L}_{\blacklozenge 0}\bm{u}|\bm{u})_{\mathfrak{H}}, \\
Q_0[\bm{u}]&=(\mathcal{L}\bm{u}|\bm{u})_{\mathfrak{H}},
\end{align}
and
\begin{equation}
Q_0^{\blacklozenge}[\bm{u}]=(\mathcal{L}_{\blacklozenge}\bm{u}|\bm{u})_{\mathfrak{H}}.
\end{equation}\\

Let us introduce an equivalent norm $\| \cdot \|_{\breve{\mathfrak{G}}\blacklozenge}$ on $\mathfrak{G}_{\blacklozenge}$ as follows.

 Putting
\begin{align}
K_1&:=2+4|\Omega|, \\
K_0&:=2+4|\Omega|+2\|W_b\|_{L^{\infty}},
\end{align}
we have
\begin{align}
&\Big(\frac{K_0}{2}+1\Big)\|\bm{u}\|_{\mathfrak{H}}^2+
\frac{K_1}{2}\Big\|\omega_b\frac{\partial\bm{u}}{\partial\phi}\Big\|_{\mathfrak{H}}^2 \leq \nonumber \\
& \sum_{\nu=0,1,2}Q_{\blacklozenge \nu}[\bm{u}]+K_0\|\bm{u}\|_{\mathfrak{H}}^2
+K_1\Big\|\omega_b\frac{\partial \bm{u}}{\partial\phi}\Big\|_{\mathfrak{H}}^2 \leq \nonumber \\
&\leq \Big(\frac{3 K_0}{2}-1\Big)\|\bm{u}\|_{\mathfrak{H}}^2+
\frac{3 K_1}{2}\Big\|\omega_b\frac{\partial\bm{u}}{\partial\phi}\Big\|_{\mathfrak{H}}^2.
\end{align}

As for $\mathcal{L}$ and $Q_0$ we know 
\begin{align} 
&\frac{1}{2}\int\Big(\frac{\Gamma P}{\rho^2}\Big)_b|\mathrm{div}(\rho_b \bm{u})|^2 d\bm{x}
\leq Q_0[\bm{u}] +\|\bm{u}\|_{\mathfrak{H}}^2 \leq \nonumber \\
&\leq \frac{3}{2}\int\Big(\frac{\Gamma P}{\rho^2}\Big)_b|\mathrm{div}(\rho_b \bm{u})|^2 d\bm{x}
+2a\|\bm{u}\|_{\mathfrak{H}}^2,
\end{align}
with the positive number $a$  defined by \eqref{7.4.26}.
Therefore, putting
\begin{equation}
Q_{\blacklozenge}[\bm{u}]:=Q_0^{\blacklozenge}[\bm{u}]+(a+K_0)\|\bm{u}\|_{\mathfrak{H}}^2+
K_1\Big\|\omega_b\frac{\partial \bm{u}}{\partial \phi}\Big\|_{\mathfrak{H}}^2,
\end{equation}
when
\begin{align}
Q_{\blacklozenge}[\bm{u}]&=(K_1-1)\Big\|\omega_b\frac{\partial \bm{u}}{\partial \phi}\Big\|_{\mathfrak{H}}^2
+(a+K_0)\|\bm{u}\|_{\mathfrak{H}}^2
+Q_{\blacklozenge 1}[\bm{u}]+
Q_{\blacklozenge 0}[\bm{u}]+Q_0[\bm{u}]  \nonumber \\
&\geq \frac{K_1}{2}\Big\|\omega_b\frac{\partial \bm{u}}{\partial \phi}\Big\|_{\mathfrak{H}}^2
+\frac{K_0-2}{2}\|\bm{u}\|_{\mathfrak{H}}^2 
+\frac{1}{2}\int\Big(\frac{\Gamma P}{\rho^2}\Big)_b|\mathrm{div}(\rho_b \bm{u})|^2 d\bm{x}
\geq 0,
\end{align}
we have
\begin{equation}
\frac{1}{C}\|\bm{u}\|_{\mathfrak{G}_{\blacklozenge}}^2
\leq Q_{\blacklozenge}[\bm{u}]+\|\bm{u}\|_{\mathfrak{H}}^2
\leq C\|\bm{u}\|_{\mathfrak{G}_{\blacklozenge}}^2,
\end{equation}
that is, 
\begin{equation}
\|\cdot\|_{\breve{\mathfrak{G}}\blacklozenge}:=\sqrt{ Q_{\blacklozenge}+\|\cdot\|_{\mathfrak{H}}^2}
\end{equation}
 is equivalent to $\|\cdot\|_{\mathfrak{G}_{\blacklozenge}}$.\\

Although we do not know whether this operator $\bm{L}_{\blacklozenge}$ is a closed operator in
$\mathfrak{H}$ or not, we can claim, at least, the following

\begin{Proposition}
If a sequence $(\bm{u}_n)_{n=0,1,2,\cdots}$ in $\mathsf{D}(\bm{L}_{\blacklozenge})$ converges to a limit $\bm{u}$ in $\mathfrak{H}$ and $\mathcal{L}_{\blacklozenge 2}\bm{u}_n$ and $\mathcal{L}\bm{u}_n$ converge in  $\mathfrak{H}$, then the limit $\bm{u}$ belongs to $\mathsf{D}(\bm{L}_{\blacklozenge})$ and $\bm{L}_{\blacklozenge}\bm{u}_n$ converge to $\bm{L}_{\blacklozenge}\bm{u}$.
\end{Proposition}

Proof. Let $\mathcal{L}_{\blacklozenge 2}\bm{u}_n \rightarrow \bm{f}_{\blacklozenge 2}$,
$\mathcal{L}\bm{u} \rightarrow \bm{f} $ in $\mathfrak{H}$. Of course $\mathcal{L}_{\blacklozenge 0}\bm{u}_n \rightarrow
\mathcal{L}_{\blacklozenge 0}\bm{u}$ in $\mathfrak{H}$. Put $\bm{f}_{\blacklozenge 0}:=\mathcal{L}_{\blacklozenge 0}\bm{u}$. Since 
$\mathcal{L}_{\blacklozenge 2}\bm{u}_n$ are supposed to converge, we have
\begin{align*}
&-\Big\|\omega_b\frac{\partial}{\partial\phi}(\bm{u}_{n+p}-\bm{u}_n)\Big\|_{\mathfrak{H}}^2=
Q_{\blacklozenge 2}[\bm{u}_{n+p}-\bm{u}_n] \\
&=(\mathcal{L}_{\blacklozenge 2}\bm{u}_{n+p}-\mathcal{L}_{\blacklozenge 2}\bm{u}_n|\bm{u}_{n+p}-\bm{u}_n)_{\mathfrak{H}} \rightarrow 0
\end{align*}
as $n+p \geq n \rightarrow \infty$. For any $\bm{\psi} \in C_0^{\infty}$, we have
\begin{align*}
&|(\mathcal{L}_{\blacklozenge 1}\bm{u}_{n+p}-\mathcal{L}_{\blacklozenge 1}\bm{u}_n|\bm{\psi})_{\mathfrak{H}}|
=|Q_{\blacklozenge 1}(\bm{u}_{n+p}-\bm{u}_n, \bm{\psi})| \\
&\leq 2|\Omega|
\Big\|\omega_b\frac{\partial}{\partial\phi}(\bm{u}_{n+p}-\bm{u}_n)\Big\|_{\mathfrak{H}}\|\bm{\psi}\|_{\mathfrak{H}} .
\end{align*}
Therefore
$$\|\mathcal{L}_{\blacklozenge 1}\bm{u}_{n+p}-\mathcal{L}_{\blacklozenge 1}\bm{u}_n\|_{\mathfrak{H}}
\leq 2|\Omega|
\Big\|\omega_b\frac{\partial}{\partial\phi}(\bm{u}_{n+p}-\bm{u}_n)\Big\|_{\mathfrak{H}} \rightarrow 0,
$$ that is, $\mathcal{L}_{\blacklozenge 1}\bm{u}_n$ converge in $\mathfrak{H}$. Put
$\bm{f}_{\blacklozenge 1}:=\mathcal{L}_{\blacklozenge 1}\bm{u}$.
Since $\mathcal{L}\bm{u}_n$ converge in ${\mathfrak{H}}$,
$$Q_0[\bm{u}_{n+p}-\bm{u}_n]=
(\mathcal{L}\bm{u}_{n+p}-\mathcal{L}\bm{u}_n|
\bm{u}_{n+p}-\bm{u}_n)_{\mathfrak{H}} \rightarrow 0.
$$
Summing up, we have that $Q_{\blacklozenge}[\bm{u}_{n+p}-\bm{u}_n] \rightarrow 0$ and $\| \bm{u}_{n+p}-\bm{u}_n \|_{\breve{\mathfrak{G}}\blacklozenge}
\rightarrow 0$. This implies that $\bm{u}_n$ converge in $\mathfrak{G}_{\blacklozenge}$ to the limit
$\bm{u} \in \mathfrak{G}_{\blacklozenge}$. Since
$\bm{u}_n \in (\mathfrak{G}_{\blacklozenge})_0$, $\bm{u} \in (\mathfrak{G}_{\blacklozenge})_0$. We see $\bm{u} \in \mathsf{D}(\bm{L}_{\blacklozenge})$
and $\mathcal{L}_{\blacklozenge}\bm{u}_n \rightarrow
\mathcal{L}\bm{u}=
\bm{f}_{\blacklozenge2}+\bm{f}_{\blacklozenge 1}+
\bm{f}_{\blacklozenge 0}+\bm{f}$. $\square$\\

If the condition could be weakend to that the total $\bm{L}_{\blacklozenge}\bm{u}_n$
converge, then it would be told that $\bm{L}_{\blacklozenge}$ is a closed operator,
which is the case when $\omega_b=0$.\\

\subsection{  Splitting {\bf ELASO($\blacklozenge$)} with respect to tha azimuthal numbers }

Since we do not know whether $\bm{L}_{\blacklozenge}$ is a self-adjoint operator in $\mathfrak{H}$ or not, we try to consider its action on vector fields $\bm{u}$ of the form $\bm{u}(\bm{x})=\check{\bm{u}}(r,\vartheta)e^{\mathrm{i}m\phi}$, where $m$ is a fixed
integer. 

Precisely speaking, we put
\begin{equation}
\check{\mathfrak{R}}_b:=\{ (r,\vartheta)|(r\sin\vartheta,0,r\cos\vartheta) \in \mathfrak{R}_b\} \subset [0,+\infty[\times \mathbb{R}/2\pi\mathbb{Z}.
\end{equation}
Since $\mathfrak{R}_b$ is axially symmetric, it holds that 
$$\check{\mathfrak{R}}_b=\{ (r,\vartheta)|(r\sin\vartheta\cos\phi,
r\sin\vartheta\sin\phi,r\cos\vartheta) \in \mathfrak{R}_b\quad \forall \phi \} .
$$
While $\rho_b(\bm{x})=\rho_b^{\flat}(r,\vartheta)$, $\check{\mathfrak{R}}_b=\{ (r,\vartheta)|\rho_b^{\flat}(r,\vartheta) >0\}$. 

Let $m$ be an integer. 

\begin{Definition}
We denote by $\mathfrak{H}_{|m}=\mathfrak{H}_{|m}(\mathfrak{R}_b; \mathbb{C}^d)$, $d=1,3$, the set of all functions $\bm{u} \in \mathfrak{H}=L^2(\mathfrak{R}_b, \rho_bd\bm{x}; \mathbb{C}^d)$ such that there exists a function $\check{\bm{u}}$ defined on $\check{\mathfrak{R}}_b$ for which it holds that
$\bm{u}(\bm{x})=\check{\bm{u}}(r,\vartheta)e^{\mathrm{i}m\phi}$. 

We denote 
\begin{align}
&\mathfrak{D}_{|m}:=C_0^{\infty}(\mathfrak{R}_b)\cap \mathfrak{H}_{|m}, \\
&\mathfrak{G}_{\blacklozenge|m}:=\mathfrak{G}_{\blacklozenge}\cap \mathfrak{H}_{|m}
\end{align}
\end{Definition}

 Clearly, $u(\bm{x}) =\check{u}(r,\vartheta)e^{\mathrm{i}m\phi}$ belongs to $\mathfrak{D}_{|m}$ if and only if there exists $\check{u}^{\sharp} \in C_0^{\infty}(\check{\mathfrak{R}}_b^{\sharp})$ such that
$\check{u}(r,\vartheta)=
\check{u}^{\sharp}(r\sin\vartheta,0, r\cos\vartheta)$, where $\check{\mathfrak{R}}_b^{\sharp}:=\{(r\sin\vartheta, 0, r\cos\vartheta) | (r,\vartheta \in \check{\mathfrak{R}}_b\}$.

The following calculus can be easily verified:

1) If $q \in \mathfrak{H} $ is axially symmetric, that is, $\in \mathfrak{H}_{|0}$, like $\rho_b, S_b, P_b$, and if
 $\bm{u} \in \mathfrak{H}_{|m}, q\bm{u} \in \mathfrak{H}$, then $q\bm{u} \in \mathfrak{H}_{|m}$ and
${(q\bm{u})}\check \ =q^{\flat}\check{\bm{u}}$.

2) If $\bm{Q} \in \mathfrak{H}_{|m}$ and $\mathrm{div}\bm{Q} \in \mathfrak{H}$, then $\mathrm{div}\bm{Q} \in \mathfrak{H}_{|m}$. In fact
$$(\mathrm{div}\bm{Q})\check \  (r,\vartheta)=
\frac{1}{r^2}\frac{\partial}{\partial r}(r^2\check{Q}^r)+
\frac{1}{r\sin\vartheta}\frac{\partial}{\partial \vartheta}(\sin\vartheta \check{Q}^{\vartheta})+
\frac{\mathrm{i}m}{r\sin\vartheta}\check{Q}^{\phi}$$
for $\check{\bm{Q}}=\check{Q}^r\bm{e}_r+\check{Q}^{\vartheta}
\bm{e}_{\vartheta}+\check{Q}^{\phi}\bm{e}_{\phi}$.

3) If $Q \in \mathfrak{H}_{|m} $ and $\nabla Q \in \mathfrak{H}$, then $\nabla Q \in \mathfrak{H}_{|m}$. In fact
$$(\nabla Q)\check\  =\frac{\partial \check{Q}}{\partial r}\bm{e}_r+
\frac{1}{r}\frac{\partial \check{Q}}{\partial\vartheta}\bm{e}_{\vartheta}+
\frac{\mathrm{i}m}{r\sin\vartheta}\check{Q}\bm{e}_{\phi}.
$$

4) If $\bm{Q} \in \mathfrak{H}_{|m}$, then $(\bm{Q}|\nabla S_b) \in \mathfrak{H}_{|m}$. In fact
$$(\bm{Q}|\nabla S_b)\check\  =\check{Q}^r\frac{\partial S_b}{\partial r}+\check{Q}^{\vartheta}\frac{\partial S_b}{\partial \vartheta}. $$

5) If $ q \in C^{\alpha}(\mathfrak{R}_b\cap \partial\mathfrak{R}_b)\cap \mathfrak{H}_{|m}$, then 
$\mathcal{K}[q] \in \mathfrak{H}_{|m}$ . In fact
$$\mathcal{K}[q]\check\  (r,\vartheta)=
\int\int_{\check{\mathfrak{R}}_b}
\check{q}(r',\vartheta')K_{|m}(r,\vartheta, r',\vartheta')(r')^2\sin\vartheta' dr'd\vartheta',
$$
where
$$K_{|m}(r,\vartheta, r',\vartheta'):=
\frac{1}{2\pi}
\int_0^{\pi}
\frac{\cos(m\phi)d\phi}{\sqrt{r^2+(r')^2-2rr'(\sin\vartheta\sin\vartheta'\cos\phi+\cos\vartheta\cos\vartheta')}}.
$$

and so on.\\

Thanks to these calculus, we can claim that the operators $ \bm{B}_{\blacklozenge | m},
\bm{L}_{\blacklozenge | m} $
can be considered as operators from dense subspace of $\mathfrak{H}_{|m}$ into $\mathfrak{H}_{|m}$, where

\begin{Definition}
We put
\begin{align}
\bm{B}_{\blacklozenge | m} &:= \mathcal{B}_{\blacklozenge}\restriction \mathfrak{H}_{|m}, \\
\bm{L}_{\blacklozenge | m}& :=
\bm{L}_{\blacklozenge} \restriction \mathfrak{H}_{|m}
\end{align}
\end{Definition}

Let us consider $\bm{u} \in \mathfrak{G}_{\blacklozenge|m}$. Then, since
$\displaystyle \frac{\partial\bm{u}}{\partial \phi}=\mathrm{i}m \bm{u}$, we see

\begin{equation}
0 \leq Q_{\blacklozenge}[\bm{u}]\leq Q^{\blacklozenge}_0[\bm{u}] + a(|m|)\|\bm{u}\|_{\mathfrak{H}}^2,
\end{equation}
where
\begin{equation}
a(|m|):=a+ K_0 + K_1\|\omega_b\|_{L^{\infty}}^2|m|^2.
\end{equation}
Thus $Q_{\blacklozenge}\restriction \mathfrak{G}_{\blacklozenge |m}$ is equivalent to the norm of
$$
\mathfrak{G}_{\blacklozenge|m}=\Big\{ \bm{u} \in \mathfrak{H}_{|m} \Big|\  \mathrm{div}(\rho_b\bm{u}) \in L^2\Big(\Big(\frac{\Gamma P}{\rho^2}\Big)_bd\bm{x}\Big) \Big\}.
$$
Of course, for $\bm{u} \in \mathfrak{D}_{|m}$, we have
\begin{equation}
Q^{\blacklozenge}_0[u]=(\mathcal{L}_{\blacklozenge}\bm{u}|\bm{u})_{\mathfrak{H}} \geq -a(|m|)\|\bm{u}\|_{\mathfrak{H}}^2,
\end{equation}
that is, $\mathcal{L}_{\blacklozenge}\restriction \mathfrak{D}_{|m}$ is bounded from below, and symmetric. Therefore
$ \bm{L}_{\blacklozenge |m}:=\bm{L}\restriction \mathfrak{H}_{|m}$ is self-adjoint, while
$$
\mathsf{D}(\bm{L}_{\blacklozenge |m})=\{ \bm{u} \in (\mathfrak{G}_{\blacklozenge |m})_0 \ |\ \mathcal{L}_{\blacklozenge}\bm{u} \in \mathfrak{H} \}.
$$
Here $(\mathfrak{G}_{\blacklozenge|m})_0$ stands for the closure of $\mathfrak{D}_{|m}=C_0^{\infty}(\mathfrak{R}_b)\cap \mathfrak{H}_{|m}$
in $\mathfrak{G}_{\blacklozenge |m}$. Note that it holds
$(\mathfrak{G}_{\blacklozenge|m})_0=(\mathfrak{G}_{\blacklozenge})_0\cap
\mathfrak{H}_{|m}$, and, for $\bm{u} \in \mathfrak{G}_{\blacklozenge |m}$, we have
$$ \mathcal{L}_{\blacklozenge 2}\bm{u}=-m^2\omega_b^2\bm{u},
\quad
\mathcal{L}_{\blacklozenge 1}\bm{u}=-2m\mathrm{i}\Omega\omega_b
\begin{bmatrix}
-u^2 \\
\\
u^1 \\
\\
0
\end{bmatrix},
$$
which belong to $\mathfrak{H}$ automatically.

On the other hand, 
$\bm{B}_{\blacklozenge |m}=\mathcal{B}_{\blacklozenge}\restriction \mathfrak{H}_{|m}$ is an operator in $\mathfrak{H}_{|m}$. Actually, for $\bm{v} \in \mathfrak{H}_{|m}$, it turns out to be that
$$\bm{B}_{\blacklozenge |m}\bm{v}=2\bm{\Omega}\times \bm{v} + 2\mathrm{i}m \omega_b\bm{v}.
$$
therefore $\bm{B}_{\blacklozenge |m}$ is a bounded skew-symmetric linear operator
in $\mathfrak{H}_{|m}$ with
\begin{equation}
\mathfrak{Re}[(\bm{B}_{\blacklozenge |m}  \bm{v}|\bm{v})_{\mathfrak{H}}]=0
\quad\mbox{for}\quad\forall \bm{v} \in \mathfrak{H}_{|m}.
\end{equation}

Consequently we can apply the Hille-Yosida theory to the initial value problem \eqref{IVPbl.1}\eqref{IVPbl.2} to claim:\\

\begin{Proposition} \label{Prop.|m}
 If $\bm{\xi}^0 \in \mathsf{D}(\bm{L}_{\blacklozenge |m}), \bm{v}^0 \in (\mathfrak{G}_{\blacklozenge |m})_0,
\bm{f} \in C([0,+\infty[; \mathfrak{H}_{|m})$, then the problem
\eqref{IVPbl.1}\eqref{IVPbl.2}
admits a unique solution
$\bm{\xi}=\bm{\xi}(t,\bm{x})$
in
$$C^2([0,+\infty[; \mathfrak{H}_{|m})
\cap C^1([0,+\infty[; (\mathfrak{G}_{\blacklozenge|m})_0) \cap 
C([0,+\infty[; \mathsf{D}(\bm{L}_{\blacklozenge |m})).
$$
 Moreover the energy
$$
E(t, \bm{u})=\|\bm{u}\|_{\breve{\mathfrak{G}}\blacklozenge}^2+\Big\|\frac{\partial\bm{u}}{\partial t}\Big\|_{\mathfrak{H}}^2
=
\|\bm{u}\|_{\mathfrak{H}}^2+Q_{\blacklozenge}[\bm{u}]+\Big\|\frac{\partial\bm{u}}{\partial t}\Big\|_{\mathfrak{H}}^2
$$
of the solution enjoys the inequality
\begin{equation}
\sqrt{E(t,\bm{\xi})}\leq
e^{\kappa(|m|)t}\sqrt{ E(0,\bm{\xi})}+
\int_0^te^{\kappa(|m|)(t-s)}
\|\bm{f}(s)\|_{\mathfrak{H}} ds,
\end{equation}
where
\begin{equation}
\kappa(|m|):=1+a(|m|).
\end{equation}
\end{Proposition}

Note that 
$$E(0,\bm{\xi})=(1+a(|m|))\|\bm{\xi}^0\|_{\mathfrak{H}}^2+
(\bm{L}_{\blacklozenge}\bm{\xi}^0|\bm{\xi}^0)_{\mathfrak{H}}+
\|\bm{v}^0\|_{\mathfrak{H}}^2.
$$\\

Now, given a positive integer $m^{\vee}=1,2,\cdots$, we consider the Hilbert space 
$\mathfrak{H}_{||m^{\vee}}$ defined as

\begin{Definition}
The Hilbert space 
$\mathfrak{H}_{||m^{\vee}}$ is the set of all linear combinations $$\bm{u}=\sum_{|m|\leq m^{\vee}} \bm{u}_m
=\sum_{|m|\leq m^{\vee}}\check{\bm{u}}_m(r,\vartheta)e^{\mathrm{i}m\phi},
$$
where $ \bm{u}_m \in \mathfrak{H}_{|m}$. 

Put
\begin{align}
\mathfrak{G}_{\blacklozenge ||m^{\vee}}&:=\Big\{ \bm{u} \in \mathfrak{H}_{||m^{\vee}}\ \Big|\ \mathrm{div}(\rho_b\bm{u}) \in
L^2\Big(\Big(\frac{\Gamma P}{\rho^2}\Big)_bd\bm{x}\Big) \Big\}, \\
(\mathfrak{G}_{\blacklozenge||m^{\vee}})_0&:=\mbox{the closure of }\quad C_0^{\infty}(\mathfrak{R}_b)\cap \mathfrak{H}_{||m^{\vee}}\quad\mbox{in}\quad \mathfrak{G}_{\blacklozenge ||m^{\vee}}.
\end{align}

Denote
\begin{align}
\bm{B}_{\blacklozenge ||m^{\vee}}&:=\mathcal{B}_{\blacklozenge}\restriction \mathfrak{H}_{|| m^{\vee}}, \\
\bm{L}_{\blacklozenge ||m^{\vee}}&:=\bm{L}_{\blacklozenge}\restriction \mathfrak{H}_{||m^{\vee}}
\end{align}

\end{Definition}

Then $\bm{L}_{\blacklozenge ||m^{\vee}}$ turns out to be a self-adjoint operator in $\mathfrak{H}_{||m^{\vee}}$ with
$$\mathsf{D}(\bm{L}_{\blacklozenge||m^{\vee}})=
\{ \bm{u} \in (\mathfrak{G}_{\blacklozenge ||m^{\vee}})_0 \ |\  \mathcal{L}_{\blacklozenge}\bm{u} \in \mathfrak{H} \}
$$
and $\bm{B}_{\blacklozenge ||m^{\vee}}$ turns out to be a bounded linear operator in $\mathfrak{H}_{||m^{\vee}}$.  \\

We can claim the following

\begin{Theorem}
Let $m^{\vee}$ be a positive integer.
 If $\bm{\xi}^0 \in \mathsf{D}(\bm{L}_{\blacklozenge ||m^{\vee}}), \bm{v}^0 \in (\mathfrak{G}_{\blacklozenge ||m^{\vee}})_0,
\bm{f} \in C([0,+\infty[; \mathfrak{H}_{||m^{\vee}})$, then the problem
\eqref{IVPbl.1}\eqref{IVPbl.2}
admits a unique solution
$\bm{\xi}=\bm{\xi}(t,\bm{x})$
in
$$C^2([0,+\infty[; \mathfrak{H}_{||m^{\vee}})
\cap C^1([0,+\infty[; (\mathfrak{G}_{\blacklozenge||m^{\vee}})_0) \cap 
C([0,+\infty[; \mathsf{D}(\bm{L}_{\blacklozenge ||m^{\vee}})).
$$
Moreover the energy
$$
E(t, \bm{u})=
\|\bm{u}\|_{\breve{\mathfrak{G}}\blacklozenge}^2+\Big\|\frac{\partial\bm{u}}{\partial t}\Big\|_{\mathfrak{H}}^2=
\|\bm{u}\|_{\mathfrak{H}}^2+Q_{\blacklozenge}[\bm{u}]+\Big\|\frac{\partial\bm{u}}{\partial t}\Big\|_{\mathfrak{H}}^2
$$
of the solution enjoys the inequality
\begin{equation}
\sqrt{E(t,\bm{\xi})}\leq
e^{\kappa(m^{\vee})t}\sqrt{ E(0,\bm{\xi})}+
\int_0^te^{\kappa(m^{\vee})(t-s)}
\|\bm{f}(s)\|_{\mathfrak{H}} ds.
\end{equation}

\end{Theorem}

Proof. First we note that, if $m\not=m'$, $\bm{u} \in C^1([0,+\infty[;\mathfrak{G}_{\blacklozenge |m})$, $\bm{u}' \in C^1([0,+\infty[; \mathfrak{G}_{\blacklozenge |m'})$, then $(\bm{u}|\bm{u}')_{\mathfrak{H}}=0,  Q_{\blacklozenge}(\bm{u}, \bm{u}')=0$, 
$ (\partial_t\bm{u}|\partial_t\bm{u}')_{\mathfrak{H}}=0$, therefore
$E(t,\bm{u}+\bm{u}')=E(t,\bm{u})+E(t,\bm{u}')$. Therefore the assertion is a direct consequence of Proposition \ref{Prop.|m}, since
 $\kappa(|m|) \leq \kappa(m^{\vee})$ for $|m|\leq m^{\vee}$. 
 $\square$\\

When  $\|\omega_b\|_{L^{\infty}} \not=0$, then we see
$\kappa(m^{\vee})t \geq 2\|\omega_b\|_{L^{\infty}}^2(m^{\vee})^2 t \rightarrow \infty$ as $m^{\vee} \rightarrow \infty$ for $\forall t>0$ Therefore we can not deduce from this theorem the existence of solutions to the problem \eqref{IVPbl.1}\eqref{IVPbl.2}
for non-trivial data $(\bm{\xi}^0,\bm{v}^0)
\in \mathsf{D}(\bm{L}_{\blacklozenge})\times (\mathfrak{G}_{\blacklozenge})_0$ by considering the sequence
of the solutions $\bm{\xi}=\bm{\xi}_{m^{\vee}}(t, \bm{x}) $ to the data
$(\bm{\xi}_{m^{\vee}}^0, \bm{v}_{m^{\vee}}^0) \in \mathsf{D}(\bm{L}_{\blacklozenge ||m^{\vee}})\times
(\mathfrak{G}_{\blacklozenge ||m^{\vee}})_0$
which approximate the data 
$(\bm{\xi}^0,\bm{v}^0)$ and passing to the limit of $m^{\vee}\rightarrow \infty$.\\

\subsection{ Axially symmetric perturbations and Rayleigh's instability}

Let us consider {\bf ELASO($\blacklozenge$) }
in the cylindrical co-ordinates $(\varpi, z, \phi)$:
$$
x^1=\varpi\cos\phi.\quad x^2=\varpi\sin\phi,\quad x^3=z.
$$
by following \cite{Lebovitz1970}.

For a vector field $\displaystyle \bm{Q}=\sum_{k=1}^3Q^k\frac{\partial}{\partial x^k}$ we denote
$$\bm{Q}^{<}=
\begin{bmatrix}
Q^{\varpi} \\
\\
Q^z \\
\\
Q^{\phi}
\end{bmatrix},
$$
where 
\begin{align*}
\bm{Q}=\sum Q^k\frac{\partial}{\partial x^k}&=Q^{\varpi}\frac{\partial}{\partial\varpi}+Q^z\frac{\partial}{\partial z}+\frac{Q^{\phi}}{\varpi}\frac{\partial}{\partial \phi} \\
&=Q^{\varpi}\bm{e}_{\varpi}+Q^z\bm{e}_z+Q^{\phi}\bm{e}_{\phi},
\end{align*}that is,
\begin{align*}
Q^{\varpi}&=Q^1\cos\phi+Q^2\sin\phi \\
Q^z&=Q^3 \\
Q^{\phi}&=-Q^1\cos\phi+Q^2\cos\phi.
\end{align*}
Then {\bf ELASO($\blacklozenge$)} reads
\begin{equation}
\frac{\partial^2\bm{\xi}^{<}}{\partial t^2}+\Big(\mathcal{B}_{\blacklozenge}\frac{\partial \bm{\xi}}{\partial t}\Big)^{<}
+(\mathcal{L}_{\blacklozenge}\bm{\xi})^{<}=0, \label{Cyl.1}
\end{equation}
where
\begin{subequations}
\begin{align}
&(\mathcal{B}_{\blacklozenge}\dot{\bm{\xi}})^{<}=2\omega_b
\begin{bmatrix}
\partial_{\phi}\dot{\xi}^{\phi} \\
\\
\partial_{\phi}\dot{\xi}^z \\
\\
\partial_{\phi}\dot{\xi}^{\varpi}
\end{bmatrix}+
2\Omega_b
\begin{bmatrix}
-\dot{\xi}^{\phi} \\
\\
0 \\
\\
\dot{\xi}^{\varpi}
\end{bmatrix}
, \\
&(\mathcal{L}_{\blacklozenge})^{<}=
((\mathcal{L}_{\blacklozenge 2}+\mathcal{L}_{\blacklozenge 1}+\mathcal{L}_{\blacklozenge 0})\bm{\xi})^{<}+
(\mathcal{L}\bm{\xi})^{<}, \\
&((\mathcal{L}_{\blacklozenge 2}+\mathcal{L}_{\blacklozenge 1}+\mathcal{L}_{\blacklozenge 0})\bm{\xi})^{<}
=
\omega_b^2
\begin{bmatrix}
\partial_{\phi}^2\xi^{\varpi} \\
\\
\partial_{\phi}^2\xi^{z} \\
\\
\partial_{\phi}^2\xi^{\phi}
\end{bmatrix}
-2\Omega \omega_b
\begin{bmatrix}
\partial_{\phi}\xi^{\phi} \\
\\
0 \\
\\
\partial_{\phi}\xi^{\varpi}
\end{bmatrix}
+ \nonumber \\
&+2\Omega_b
\begin{bmatrix}
\varpi\Big(\frac{\partial \Omega_b}{\partial \varpi}\xi^{\varpi}+
\frac{\partial \Omega_b}{\partial z}\xi^z\Big) \\
\\
0 \\
\\
0
\end{bmatrix}, \\
&(\mathcal{L}\bm{\xi})^{<}=
\frac{1}{\rho_b}
\begin{bmatrix}
\frac{\partial}{\partial\varpi} \lceil\delta P \rfloor \\
\\
\frac{\partial}{\partial z}\lceil \delta P \rfloor\\
\\
\frac{\partial}{\partial \phi} \lceil \delta P \rfloor
\end{bmatrix}
-\frac{\lceil \delta \rho \rfloor}{\rho_b^2}
\begin{bmatrix}
\frac{\partial P_b}{\partial \varpi} \\
\\
\frac{\partial P_b}{\partial z} \\
\\
0
\end{bmatrix}
-4\pi\mathsf{G}
\begin{bmatrix}
\frac{\partial}{\partial \varpi}\mathcal{K}[\lceil \delta\rho \rfloor ] \\
\\
\frac{\partial}{\partial z}\mathcal{K}[\lceil \delta\rho \rfloor] \\
\\
\frac{\partial}{\partial \phi}\mathcal{K}[ \lceil\delta\rho \rfloor] 
\end{bmatrix}
, \\
&\lceil \delta\rho \rfloor=-\frac{1}{\varpi}\frac{\partial}{\partial \varpi}(\varpi \rho_b\xi^{\varpi})
-\frac{\partial}{\partial z}(\rho_b\xi^z)-\frac{1}{\varpi}
\frac{\partial}{\partial \phi}(\rho_b\xi^{\phi}), \\
&\lceil \delta P \rfloor=\Big(\frac{\Gamma P}{\rho}\Big)_b\lceil \delta \rho \rfloor+
(\Gamma P)_b(\xi^{\varpi}\mathfrak{a}_b^{\varpi}+\xi^z\mathfrak{a}_b^z), \\
&\mathfrak{a}_b^{\varpi}=\Big(
\frac{1}{\rho}\frac{\partial \rho}{\partial \varpi}-
\frac{1}{\Gamma P}\frac{\partial P}{\partial \varpi}
\Big)_b, \quad
\mathfrak{a}_b^{z}=\Big(
\frac{1}{\rho}\frac{\partial \rho}{\partial z}-
\frac{1}{\Gamma P}\frac{\partial P}{\partial z}
\Big)_b.
\end{align}
\end{subequations}

Here we denote
\begin{equation}
\Omega_b:=\Omega + \omega_b.
\end{equation}\\

Now let us consider a solution $\bm{\xi}(t,\bm{x})$ such that
$\displaystyle \frac{\partial}{\partial\phi}\xi^{\varpi}= \frac{\partial}{\partial\phi}\xi^{z}= \frac{\partial}{\partial\phi}\xi^{\phi}=0$, namely, an axisymmetric perturbation $\bm{\xi}=\bm{\xi}(t, \varpi, z)$. 

Then it holds that 
$$\frac{\partial}{\partial \phi}\delta\rho=0,\quad \frac{\partial}{\partial\phi}\delta P=0, \quad
\frac{\partial}{\partial \phi}\mathcal{K}[\delta \rho]=0,
$$
therefore $(\mathcal{L}\bm{\xi})^{\phi}=0$. Moreover
$(\mathcal{L}\bm{\xi})^{<}$ involve only $\xi^{\varpi}, \xi^z$, and involve neither 
$\xi^{\phi}$ nor its derivatives, while
$$\delta\rho=-\frac{1}{\varpi}\frac{\partial}{\partial \varpi}(\varpi \rho_b\xi^{\varpi})
-\frac{\partial}{\partial z}(\rho_b\xi^z).$$

Therefore the $\phi$-component of \eqref{Cyl.1} reads
\begin{equation}
\frac{\partial^2\xi^{\phi}}{\partial t^2}+2\Omega_b\frac{\partial \xi^{\varpi}}{\partial t}=0.
\end{equation}
Integrating this equation, we have
$$
\frac{\partial \xi^{\phi}}{\partial t}+2\Omega_b\xi^{\varpi}=C(\varpi, z)\quad \forall t \geq 0.
$$
Assuming that the initial data 
$\displaystyle {\bm{\xi}}^0=\bm{\xi}|_{t=0}, {\bm{v}}^0=\frac{\partial \bm{\xi}}{\partial t}\Big|_{t=0}$
satisfy
\begin{equation}
{v}^{\phi 0}+2\Omega_b {\xi}^{\varpi 0}=0,
\end{equation}
we have the identity
\begin{equation}
\frac{\partial \xi^{\phi}}{\partial t}+2\Omega_b\xi^{\varpi}=0 \quad \forall t \geq 0. \label{Cyl.2}
\end{equation}
This identity means that the angular momentum per unit mass with respect to $z$-axis,
$$J=\varpi (\varpi \Omega+v^{\phi})$$
of the perturbed motion coincides with the angular momentum per unit mass
$$J_b=\varpi (\varpi \Omega+v_b^{\phi})=\varpi^2\Omega_b$$
of the background, that is, $\Delta J=0$. In fact, noting that
$$\bm{v}_b^{<}=\omega_b
\begin{bmatrix}
0 \\
0 \\
1
\end{bmatrix},\quad
\Delta\varpi=\xi^{\varpi},\quad
\Delta\bm{v}=\frac{D}{Dt}\Big|_{\bm{v}_b}\bm{\xi},
$$
we see
\begin{align*}
\Delta J&=(\Delta\varpi)(\varpi \Omega +\bm{v}_b^{\phi})+
\varpi\Delta (\varpi+v^{\phi}) \\
&=2\varpi\Omega\xi^{\varpi}+
\varpi\omega_b\xi^{\varpi}+\varpi\frac{D}{Dt}\Big|_{\bm{v}_b}\xi^{\phi} \\
&=2\varpi\Omega \xi^{\varpi}+\varpi \omega_b\xi^{\varpi}+
\varpi\Big(\frac{\partial \xi^{\phi}}{\partial t}+\omega_b\xi^{\varpi}\Big) \\
&=\varpi\Big(\frac{\partial \xi^{\phi}}{\partial t}+2(\Omega+\omega_b)\xi^{\varpi}\Big) \\
&=0.
\end{align*}

The $\varpi$-component of \eqref{Cyl.1} turns out to be
\begin{equation}
\frac{\partial^2\xi^{\varpi}}{\partial t^2}+c_{\varpi}\xi^{\varpi}+c_z\xi^z+
(\mathcal{L}\bm{\xi})^{\varpi}=0, \label{Cyl.3}
\end{equation}
where
\begin{equation}
c_{\varpi}:=\frac{1}{\varpi^3}\frac{\partial}{\partial \varpi}(\varpi^4 \Omega_b^2), \quad
c_{z}:=\frac{1}{\varpi^3}\frac{\partial}{\partial z}(\varpi^4 \Omega_b^2).
\end{equation}
If $\displaystyle \frac{\partial}{\partial z}\omega_b =0$, as in barotropic background, then 
$ c_{\varpi}=\mathscr{C}, c_z=0$, where $\mathscr{C}$ defined by
\begin{equation}
\mathscr{C}:=\frac{1}{\varpi^3}\frac{d}{d \varpi}(\varpi^4 \Omega_b^2) 
\end{equation}
is so called `{\bf Rayleigh stability discriminant }'. 

The $z$-component of \eqref{Cyl.2} reads
\begin{equation}
\frac{\partial^2\xi^z}{\partial t^2}+(\mathcal{L}\bm{\xi})^{z} =0, \label{Cyl.5}
\end{equation}
Since $(\mathcal{L}\bm{\xi})^{<} $ does not involve $\xi^{\phi}$, the equations
\eqref{Cyl.3}\eqref{Cyl.5} form a system for unknowns $\xi^{\varpi}, \xi^z$.
when a solution $(\xi^{\varpi}, \xi^z)$ of this system is given, then, integrating \eqref{Cyl.2}, we get $\xi^{\phi}$ as
$$\xi^{\phi}(t)={\xi}^{\phi 0}-2\Omega_b\int_0^t
\xi^{\varpi}(t')dt'.
$$\\

Let us consider the existence theorem for the system of equations \eqref{Cyl.3}\eqref{Cyl.5} for the unkown $(\xi^{\varpi}(t, \varpi, z), \xi^{z}(t,\varpi, z))$. We assume that the background is barotropic.

We use the following notations:
\begin{align*}
\mathfrak{R}^{<}&=\{ (\varpi, z) \in [0,+\infty[\times \mathbb{R} \  |\  (\varpi,0, z) \in \mathfrak{R}_b \}, \\
\mathfrak{h}^d&=L^2( \mathfrak{R}^{<}, 2\pi \varpi d\varpi dz; \mathbb{C}^d),\quad d=1,2,3, \\
(\bm{q}^{\eqslantless})^{>}(\bm{x})&=
\begin{bmatrix}
q^{\varpi}(\sqrt{(x^1)^2+(x^2)^2 }, x^3) \\
\\
0 \\
\\
q^{z}(\sqrt{(x^1)^2+(x^2)^2 }, x^3) 
\end{bmatrix}
\quad\mbox{for}\quad
\bm{q}^{\eqslantless}=
\begin{bmatrix}
q^{\varpi}(\varpi, z) \\
\\
0 \\
\\
q^z(\varpi, z)
\end{bmatrix}
\in \mathfrak{h}^2, \\
\mathfrak{D}^{\eqslantless}&=\{ \bm{q}^{\eqslantless} \in \mathfrak{h}^2 \  |\  
(\bm{q}^{\eqslantless})^{>} \in C_0^{\infty}(\mathfrak{R}_b; \mathbb{C}^3) \},
\end{align*}\\

The system \eqref{Cyl.3}\eqref{Cyl.5} reads
\begin{equation}
\frac{\partial^2 \bm{\xi}^{\eqslantless}}{\partial t^2}+\mathcal{L}_{\blacklozenge}^{\eqslantless}\bm{\xi}^{\eqslantless}=0,
\end{equation}
where
\begin{equation}
\mathcal{L}_{\blacklozenge}^{\eqslantless}\bm{\xi}^{\eqslantless}=
\Big((\mathcal{B}_{\blacklozenge}\partial_t+ \mathcal{L}_{\blacklozenge}). (\bm{\xi}^{\eqslantless})^{>} \Big)^{\eqslantless}=
\begin{bmatrix}
\Big((\mathcal{B}_{\blacklozenge}\partial_t+ \mathcal{L}_{\blacklozenge}). (\bm{\xi}^{\eqslantless})^{>} \Big)^{\varpi} \\
\\
\Big((\mathcal{B}_{\blacklozenge}\partial_t+ \mathcal{L}_{\blacklozenge}). (\bm{\xi}^{\eqslantless})^{>} \Big)^{z}
\end{bmatrix},
\end{equation}
\begin{subequations}
\begin{align}
&\Big((\mathcal{B}_{\blacklozenge}\partial_t+ \mathcal{L}_{\blacklozenge}). (\bm{\xi}^{\eqslantless})^{>} \Big)^{\varpi} = \\
&=\mathscr{C}\xi^{\varpi}+
\frac{1}{\rho}\frac{\partial}{\partial \varpi}\lceil \delta P \rfloor
-\frac{\lceil \delta \rho \rfloor }{\rho^2}\frac{\partial P}{\partial \varpi}
-4\pi\mathsf{G}\frac{\partial}{\partial\varpi}\mathcal{K}[\lceil \delta \rho \rfloor ], \\
&\Big((\mathcal{B}_{\blacklozenge}\partial_t+ \mathcal{L}_{\blacklozenge}). (\bm{\xi}^{\eqslantless})^{>} \Big)^{z} = \\
&=
\frac{1}{\rho}\frac{\partial}{\partial z}\lceil \delta P \rfloor
-\frac{\lceil \delta \rho \rfloor }{\rho^2}\frac{\partial P}{\partial z}
-4\pi\mathsf{G}\frac{\partial}{\partial z}\mathcal{K}[\lceil \delta \rho \rfloor ].
\end{align}
\end{subequations}
Here 
\begin{subequations}
\begin{align}
\lceil \delta \rho \rfloor &=-\mathrm{div}(\rho (\bm{\xi}^{\eqslantless})^{>}) = \nonumber \\
&=-\Big[
\frac{1}{\varpi}\frac{\partial}{\partial\varpi}(\varpi \rho {\xi}^{\varpi}+\frac{\partial}{\partial z}(\rho \xi^z)\Big], \\
\lceil \delta P \rfloor &=\frac{\Gamma P}{\rho}\lceil \delta \rho \rfloor
+\Gamma P(\xi^{\varpi}\mathfrak{a}^{\varpi}+\xi^z\mathfrak{a}^z), 
\\
&\mathfrak{a}^{\varpi}=
\frac{1}{\rho}\frac{\partial \rho}{\partial \varpi}-
\frac{1}{\Gamma P}\frac{\partial P}{\partial \varpi}
, \quad
\mathfrak{a}^{z}=
\frac{1}{\rho}\frac{\partial \rho}{\partial z}-
\frac{1}{\Gamma P}\frac{\partial P}{\partial z}.
\end{align}
\end{subequations}
and
$\rho, P, \Gamma $ stand for
$\rho_b, P_b, \Gamma_b$.\\

Then, for $\bm{\xi}^{\eqslantless} =(\xi^{\varpi}, \xi^z)^{\top} \in \mathfrak{D}^{\eqslantless}$, it holds that
$$
(\mathcal{L}_{\blacklozenge}^{\eqslantless}\bm{\xi}^{\eqslantless}|\bm{\xi}^{\eqslantless})_{\mathfrak{h}^2}=
Q_{\blacklozenge}^{\eqslantless}[\bm{\xi}^{\eqslantless}],
$$
where
\begin{align*}
Q_{\blacklozenge}^{\eqslantless}[\bm{\xi}^{\eqslantless}]&:=2\pi\int_{\mathfrak{R}^{<}}\mathscr{C}|\xi^{\varpi}|^2\rho \varpi d\varpi dz +
2\pi\int_{\mathfrak{R}^{<}}\frac{\Gamma P}{\rho^2}|\lceil \delta\rho \rfloor|^2\varpi d\varpi dz + \\
&+4\pi\mathfrak{Re}\Big[\int_{\mathfrak{R}^<}\frac{\Gamma P\mathscr{A}}{\rho}(\bm{\xi}|\bm{n})\lceil \delta \rho \rfloor^*\varpi d\varpi dz \Big]
-2\pi\int_{\mathfrak{R}^<}\frac{\Gamma P \mathscr{A}}{\rho}\|\nabla \rho \|^2|(\bm{\xi}|\bm{n})|^2\varpi d\varpi dz \\
&-8\pi\mathsf{G}\int_{\mathfrak{R}^<}\mathcal{K}[\lceil \delta\rho \rfloor] \lceil \delta\rho \rfloor ^
*\varpi d\varpi dz,
\end{align*}
with
\begin{align*}
\lceil \delta \rho \rfloor &=-\Big[
\frac{1}{\varpi}\frac{\partial}{\partial\varpi}(\varpi \rho {\xi}^{\varpi})+
\frac{\partial}{\partial z}(\rho \xi^z)
\Big], 
\\
(\bm{\xi}|\bm{n})&=
\frac{ \xi^{\varpi}\partial_{\varpi}\rho+\xi^z\partial_z\rho }{\sqrt{ (\partial_{\varpi}\rho)^2+(\partial_z\rho)^2 }}.
\end{align*}

Recall that $\mathscr{A}=\mathscr{A}_b$ is such that
$$ \mathfrak{a}_b=-\mathscr{A}\bm{n}=-\mathscr{A}\frac{\nabla \rho_b}{\|\nabla \rho_b\|}. $$

Put
\begin{equation}
\mathfrak{g}:=\Big\{
\bm{\xi}^{\eqslantless}=
\begin{bmatrix}
\xi^{\varpi} \\
 \xi^z
 \end{bmatrix} 
 \in \mathfrak{h}^2 \  \Big|\  \lceil \delta\rho \rfloor \in
 L^2(\mathfrak{R}^<,\frac{\Gamma P}{\rho^2}2\pi\varpi d\varpi dz ; \mathbb{C}) \Big\}
 \end{equation}
 and denote the closure of $\mathfrak{D}^{\eqslantless}$ in $\mathfrak{g}$ by $\mathfrak{g}_0$.
 
 We can consider the quadratic form $Q_{\blacklozenge}^{\eqslantless}$ on $\mathfrak{g}$, and it is bounded from below, since
 $\inf
 \mathscr{C}(\varpi) > -\infty
 $,
 and by Proposition \ref{Prop.8v1}, Proposition \ref{Prop.8v2}. Thus we see that $\mathcal{L}_{\blacklozenge}^{\eqslantless}\restriction \mathfrak{D}^{\eqslantless}$ has the Friedrichs extension 
 $\bm{L}_{\blacklozenge}^{\eqslantless}$ which is a self-adjoint operator in $\mathfrak{h}^2$ with domain
 $$
 \mathsf{D}(\bm{L}_{\blacklozenge}^{\eqslantless})=
 \{ \bm{\xi}^{\eqslantless} \in \mathfrak{g}_0\  |\  \mathcal{L}_{\blacklozenge}^{\eqslantless}\bm{\xi}^{\eqslantless} \in \mathfrak{h}^2 \}.
 $$
 It holds that
 $$(\bm{L}_{\blacklozenge}^{\eqslantless}\bm{\xi}^{\eqslantless}|\bm{\xi}^{\eqslantless})_{\mathfrak{h}^2}=
 Q_{\blacklozenge}^{\eqslantless}[\bm{\xi}^{\eqslantless}]
 $$
 for $\bm{\xi}^{\eqslantless} \in \mathsf{D}(\bm{L}_{\blacklozenge}^{\eqslantless})$.\\
 
 Consequently we can claim
 
 \begin{Proposition} \label{Cyl.Prop1}
 If $\bm{\xi}^{\eqslantless 0} \in \mathsf{D}(\bm{L}_{\blacklozenge}^{\eqslantless}), \bm{v}^{\eqslantless 0} \in \mathfrak{g}_0$, then the initial value problem
\begin{equation}
\frac{\partial^2\bm{\xi}^{\eqslantless}}{\partial t^2}+
\bm{L}_{\blacklozenge}^{\eqslantless}\bm{\xi}^{\eqslantless}=0,\quad
\bm{\xi}^{\eqslantless}|_{t=0}=\bm{\xi}^{\eqslantless 0},
\frac{\partial \bm{\xi}^{\eqslantless}}{\partial t}\Big|_{t=0}=\bm{v}^{\eqslantless 0}
\label{Cyl.16}
\end{equation}
admits a unique solution in
$$C^2([0,+\infty[; \mathfrak{h}^2)
 \cap C^1([0,+\infty[, \mathfrak{g}_0)
\cap C([0,+\infty[; \mathsf{D}(\bm{L}_{\blacklozenge}^{\eqslantless})).$$
\end{Proposition}
 
 Let $\bm{\xi}^{<0}=(\xi^{\varpi 0},\xi^{z 0}, \xi^{\phi 0}) \in \mathsf{D}(\bm{L}_{\blacklozenge}^{\eqslantless})\times \mathfrak{h}^1,
 \bm{v}^{< 0}=(v^{\varpi 0}, v^{z0}, v^{\phi 0}) \in \mathfrak{g}_0\times \mathfrak{h}^1$ satisfy
 $$
 v^{\phi 0}+2\Omega_b (\varpi)\xi^{\varpi 0}=0.
 $$ 
 Let $\bm{\xi}^{\eqslantless} =(\xi^{\varpi}(t,\varpi, z), \xi^z(t,\varpi, z))$ be the solution of Proposition \ref{Cyl.Prop1}. Put
 $$ \xi^{\phi}(t,\varpi, z)=
 \xi^{\phi 0}(\varpi, z)-2\Omega_b(\varpi)
 \int_0^{t}\xi^{\varpi}(t', \varpi, z)dt'.
 $$
 Then $\bm{\xi}^<=(\xi^{\varpi}, \xi^z, \xi^{\phi})$ turns out to be the unique solution
 of the initial value problem
 \begin{equation}
 \frac{\partial^2\bm{\xi}^<}{\partial t^2}+\Big(\mathcal{B}_{\blacklozenge}\frac{\partial\bm{\xi}}{\partial t}+\mathcal{L}_{\blacklozenge}\bm{\xi}\Big)^<=0,
 \quad
 \bm{\xi}^<|_{t=0}=\bm{\xi}^{<0},
 \frac{\partial\bm{\xi}^<}{\partial t}\Big|_{t=0}=\bm{v}^{< 0}
 \end{equation}
 in
 $$ C^2([0,+\infty[; \mathfrak{h}^3)\cap C^1([0,+\infty[; \mathfrak{g}_0\times \mathfrak{h}^1)
 \cap C([0,+\infty[; \mathsf{D}(\bm{L}_{\blacklozenge}^{\eqslantless})\times \mathfrak{h}^1).
 $$
 
 Thus $\bm{\xi}=(\bm{\xi}^<)^>$, namely,
 $$
 \bm{\xi}=
 \begin{bmatrix}
 \xi^1(t, \bm{x}) \\
 \\
 \xi^2(t,\bm{x}) \\
 \\
 \xi^3(t,\bm{x})
 \end{bmatrix}
 =
 \begin{bmatrix}
 \xi^{\varpi}(t,\varpi, z)\cos\phi -\xi^{\phi}(t,\varpi, z)\sin\phi \\
 \\
 \xi^{\varpi}(t,\varpi, z)\sin\phi+\xi^{\phi}(t,\varpi, z)\cos\phi \\
 \\
 \xi^z(t,\varpi, z)
 \end{bmatrix}
 ,
 $$
 is the solution of the equation LASO($\blacklozenge$)
 in
 $$
 C^2([0,+\infty[; \mathfrak{H})\cap C^1([0,+\infty[; \mathfrak{G}_0)  \cap C([0,+\infty[; \mathsf{D}(\bm{L}_{\blacklozenge}).
 $$
 
 Taking the inner product of the equation of \eqref{Cyl.16} with $\partial \bm{\xi}^{\eqslantless}/\partial t$, and integrating it, we can show 
 
 \begin{Proposition} \label{Cyl.Prop2}
 The solution $\bm{\xi}^{\eqslantless}=(\xi^{\varpi}, \xi^z)$ of \eqref{Cyl.16} given by Proposition \ref{Cyl.Prop1} enjoys the energy equality
 \begin{equation}
 \Big\|\frac{\partial \bm{\xi}^{\eqslantless}}{\partial t}\Big\|_{\mathfrak{h}^2}^2+
 Q_{\blacklozenge}^{\eqslantless}[\bm{\xi}^{\eqslantless}(t)]=
 \|\bm{v}^{\eqslantless 0}\|_{\mathfrak{h}^2}^2+Q_{\blacklozenge}^{\eqslantless}[\bm{\xi}^{\eqslantless 0}]. \label{Cyl.EnEq}
 \end{equation}
 \end{Proposition}
 
 On the other hand, if (i) we suppose
 $$-\varepsilon\frac{\|\nabla \rho_b\|^2}{\rho_b} \leq \mathscr{A}_b \leq 0 \quad\mbox{on}\quad \mathfrak{R}_b, \quad 0<\exists \varepsilon <1,
 $$
 and if (ii) we adopt the Cowling approximation, that is, we relpace $\mathcal{L}=\mathcal{L}_0+4\pi\mathsf{G}\mathcal{L}_1$
 by $\mathcal{L}_0$, then we have
 \begin{align*}
 Q_{\blacklozenge}^{\eqslantless}[\bm{\xi}^{\eqslantless}] \geq & 2\pi\int_{\mathfrak{R}^<} \mathscr{C}|\xi^{\varpi}|^2 \varpi d\varpi dz \\
& +2\pi(1-\varepsilon)
 \int_{\mathfrak{R}^<}\frac{\Gamma P}{\rho^2}|\lceil \delta \rho\rfloor|^2
 \varpi d\varpi dz
 \quad \mbox{for}\quad \forall \bm{\xi}^{\eqslantless} \in \mathfrak{g}.
 \end{align*}
 ( Proposition \ref{Prop.8v11}.)
 Moreover if (iii) we suppose `Rayleigh stability condition'
 $$ \mathscr{C} \geq 0 \quad\mbox{on}\quad \mathfrak{R}_b,
  $$
  ( As for the meaning of the Rayleigh criterion, see {\bf Appendix D}.),
  then we have
  $$
 Q_{\blacklozenge}^{\eqslantless}[\bm{\xi}^{\eqslantless}]
 \geq
 2\pi(1-\varepsilon)
 \int_{\mathfrak{R}^<}\frac{\Gamma P}{\rho^2}|\lceil \delta \rho \rfloor|^2
 \varpi d\varpi dz
  \geq 0
 \quad \mbox{for}\quad \forall \bm{\xi}^{\eqslantless} \in \mathfrak{g}.
 $$
  In this case, Proposition \ref{Cyl.Prop2} deduces the stability in the sense that
  \begin{subequations}
  \begin{align}
 &\Big\|\frac{\partial \bm{\xi}^{\eqslantless}}{\partial t}\Big\|_{\mathfrak{h}^2} \leq C_1,
 \quad
 \| \lceil \delta \rho \rfloor \|_{L^2(\mathfrak{R}^<, \frac{\Gamma P}{\rho^2}\cdot 2\pi \varpi d\varpi dz)} 
 \leq \frac{C_1}{ \sqrt{1-\varepsilon} }, \\ 
&  \|\bm{\xi}^{\eqslantless}(t)\|_{\mathfrak{h}^2} \leq C_0+C_1t.
\end{align}
 \end{subequations}
 Here $C_0=\|\bm{\xi}^{\eqslantless 0}\|_{\mathfrak{h}^2} $ and
 $\displaystyle C_1=\Big(\|\bm{v}^{\eqslantless 0}\|_{\mathfrak{h}^2}^2+Q_{\blacklozenge}^{\eqslantless}[\bm{\xi}^{\eqslantless 0}]\Big)^{\frac{1}{2}}$.
 
Consequently we can claim: \\

{\it The associated solution $\bm{\xi}=(\bm{\xi}^<)^>$ of {\bf ELASO($\blacklozenge$) } enjoys
 \begin{subequations}
 \begin{align}
 \Big\|\frac{\partial \bm{\xi}}{\partial t}\Big\|_{\mathfrak{H}}&\leq
 C_1+2\|\Omega_b\|_{L^{\infty}} C_0+2\|\Omega_b\|_{L^{\infty}}C_1t, \\
 \|\bm{\xi}(t)\|_{\mathfrak{H}}&\leq
 \|\bm{\xi}^0\|_{\mathfrak{H}}+
 (C_1+2\|\Omega_b\|_{L^{\infty}}C_0)t+
 \|\Omega_b\|_{L^{\infty}}C_1 t^2.
 \end{align}
 \end{subequations}
 }

%%%%%%%%%%%%%%%%%%%%%%%%%%%%%%%%%%%%%%%%%%%%%%%%%%%%%%%%%%%

\vspace{15mm}

%%%%%%%%%%%%%%%%%%%%%%%%%%%%

{\Large \bf Appendix A}\\

We consider the initial value problem:

\begin{align*}
&\frac{d^2{u}}{dt^2}+{B}\frac{d{u}}{dt}+{L}{u}=0, \label{Apx.1a}\\
&{u}(0)=\overset{\circ}{{u}},\quad \frac{d{u}}{dt}(0)=\overset{\circ}{{v}}. 
\end{align*}

Here the unkown is a function ${u}: t \mapsto {u}(t): [0,+\infty[ \rightarrow \mathsf{X}$,
$\mathsf{X}$ being a Hilbert space endowed with an inner product $(\cdot | \cdot)_{\mathsf{X}}$.

We assume: 

$\mathsf{Y}$ is a Hilbert space endowed with an inner product
$(\cdot | \cdot )_{\mathsf{Y}}$, which is continuously and densely inbedded in $\mathsf{X}$. 
${L}$ is a self-adjoint operator in $\mathsf{X}$ whose domain $\mathsf{D}({L})$
is included in $\mathsf{Y}$ and dense in $\mathsf{Y}$ ( in the $\mathsf{Y}$-norm). 
It holds that
$$ \|u\|_{\mathsf{Y}} \geq \|u\|_{\mathsf{X}} \quad\mbox{for}\quad \forall u \in \mathsf{Y},
$$
$$
( ({L}+a+1){u}_1|{u}_2)_{\mathsf{X}}=({u}_1|{u}_2)_{\mathsf{Y}}
\quad\mbox{for}\quad \forall {u}_1 \in \mathsf{D}({L}), \forall {u}_2 \in \mathsf{Y},$$
where $a$ is a non-negative number.
${B}$ is a bounded linear operator in $\mathsf{X}$ 
with
$$\beta:=\sup\frac{|\mathfrak{Re}[({B}{v}|{v})_{\mathsf{X}}]|}{\|{v}\|_{\mathsf{X}}^2}.
$$ \\

{\bf Remark}:\  {\it If $B$ is skew-symmetric, then we have $\beta=0$. This is the case when $B=\bm{B}=2\bm{\Omega}\times $ on $\mathsf{X}=\mathfrak{H}=L^2(\mathfrak{R},\rho d\bm{x})$. }\\

{\bf Proposition A1}\  {\it For any $ c \in \mathbb{R}$ and $\lambda >a+|c|\beta$ the operator
 ${L} +c{B}+\lambda$ has the bounded linear inverse operator $({L}+c{B}+\lambda)^{-1}$ defined on the whole space $\mathsf{X}$
such that
$$
|\|({L}+c{B}+\lambda)^{-1}\||_{\mathcal{B}(\mathsf{X})}\leq \frac{1}{\lambda-a-|c|\beta}.
$$
}\\

Proof. First we see ${L} +c{B}+\lambda$ is invertible. In fact, if
$$({L} +c{B}+\lambda){u}={f},\quad {u} \in \mathsf{D}({L}),\quad {f} \in \mathsf{X},
$$
then  we see 
\begin{align*}
(\lambda-a-|c|\beta)\|{u}\|_{\mathsf{X}}^2 & \leq 
\|u\|_{\mathsf{Y}}^2+(\lambda - a -1)\|u\|_{\mathsf{X}}^2
+c\mathfrak{Re}[({B}{u}|{u})_{\mathsf{X}}]
 = \\
&=\mathfrak{Re}[ ({u}|{f})_{\mathsf{X}}] \leq \|{u}\|_{\mathsf{X}}\|{f}\|_{\mathsf{X}},
\end{align*}
therefore we have
$$\|{u}\|_{\mathsf{X}}\leq \frac{1}{\lambda-a-|c|\beta}\|{f}\|_{\mathsf{X}}.
$$
Hence the range $\mathsf{R}({L}+c{B}
+\lambda)$ is closed, since ${L}+c{B}+\lambda$ is a closed operator.
We claim that $\mathsf{R}({L}+c{B}
+\lambda) = \mathsf{X}$. In fact, suppose
$$(({L}+c{B}
 +\lambda){u}|{h})_{\mathsf{X}}=0\quad \forall {u} \in \mathsf{D}({L}).
$$
Then 
\begin{align*}
({L}{u}|{h})&=
-((c{B}+\lambda){u}|{h})_{\mathsf{X}} \\
&=({u}|(-c{B}^*-\lambda){h})_{\mathsf{X}}
\end{align*}
for $\forall {u} \in \mathsf{D}({L})$. 
Hence ${h} \in \mathsf{D}({L}^*)$ and
$${L}^*{h}=
-c{B}^*-\lambda{h}.$$
Since ${L}={L}^*$, this means that ${h} \in \mathsf{D}({L})$ and
$$({L}+c{B}^*+\lambda){h}=0.$$
Since ${L}+c{B}^*+\lambda$ is invertible, for
$({B}^*v|v)_{\mathsf{X}}=(Bv|v)_{\mathsf{X}}^*$, it follows that ${h}=0$. 
Therefore
$\mathsf{R}({L}+c{B}+\lambda) =\mathsf{X}$.
Summing up, we have $(L+cB+\lambda)^{-1}\in \mathcal{B}(\mathsf{X})$.
$\square$\\

We are going to apply the Hille-Yosida theory to the initial value problem
{\bf (IVP) }:

\begin{align*}
&\frac{d^2{u}}{d t^2}+
{B}
\frac{d{u}}{d t}+
{L}{u}={f}(t),  \\
&{u}(t)\in \mathsf{D}({L})\quad\mbox{for}\quad \forall t \geq 0, \\
&{u}=\overset{\circ}{{u}},\quad \frac{d{u}}{d t}=\overset{\circ}{{v}}
\quad \mbox{at}\quad t=0, 
 \end{align*}

We put
\begin{align*}
&U=
\begin{bmatrix}
{u} \\
\\
\dot{{u}}
\end{bmatrix},
\quad \dot{{u}}=\frac{d{u}}{d t}, \\
&\mathbf{A}U=
\begin{bmatrix}
O & -I \\
\\
{L} & {B}
\end{bmatrix}
U
=
\begin{bmatrix}
-\dot{{u}} \\
\\
{B}\dot{{u}}+{L}{u}
\end{bmatrix}, \\
&\mathsf{Z}=\mathsf{Y}\times \mathsf{X}  \\
\mbox{with}&  \\
&(U_1|U_2)_{\mathsf{Z}}=
({u}_1|{u}_2)_{\mathsf{Y}}+(\dot{{u}}_1|\dot{{u}}_2)_{\mathsf{X}}, \\
&\mathsf{D}(\mathbf{A})=\mathsf{D}({L})\times \mathsf{Y}.
\end{align*}

Then the initial value problem {\bf (IVP)} can be written as the problem {\bf (IVP\ddag)}:

$$
\frac{dU}{dt}+\mathbf{A}U=F(t),\quad U|_{t=0}=U_0, 
$$
where 
$$
U_0=
\begin{bmatrix}
\overset{\circ}{{u}} \\
\\
\overset{\circ}{{v}}
\end{bmatrix}.
\quad F(t)=
\begin{bmatrix}
0 \\
\\
{f}(t)
\end{bmatrix}.
$$

Applying \cite[Theorem 7.4]{Brezis}, we can claim\\

{\bf Proposition A2}\  {\it If $U_0 \in \mathsf{D}(\mathbf{A})$
and $F \in C([0,+\infty[; \mathsf{Z})$, then there exists a unique solution 
$U \in C^1([0,+\infty[, \mathsf{Z})\cap C([0,+\infty[, \mathsf{D}(\mathbf{A}))
$
to the problem {\bf (IVP\ddag)}. Moreover $E(t)=\|U(t)\|_{\mathsf{Z}}^2$ enjoys
$$
\sqrt{E(t)}\leq e^{\Lambda t}\Big(\sqrt{E(0)} + \int_0^te^{-\Lambda s
}\|F(s)\|_{\mathsf{Z}}ds \Big),
$$
where $\Lambda=1 +a+ \beta$.} \\

Here we consider that $\mathsf{D}(\mathbf{A})$ is equipped with the operator norm
$(\|U\|_{\mathsf{Z}}^2+\|\mathbf{A}U\|_{\mathsf{Z}}^2)^{1/2}$.\\

Proof of  Proposition A2.  Firstly $\mathbf{A}+1+a+ \beta$ is monotone, that is, for $\forall U \in \mathsf{D}(\mathbf{A})$ we have
\begin{align*}
\mathfrak{Re}[(\mathbf{A}U|U)_{\mathsf{Z}}]+(1+a+ \beta)\|U\|_{\mathsf{Z}}^2&=\mathfrak{Re}\Big[-(\dot{{u}}|{u})_{\mathsf{Y}}+({L}{u}|\dot{{u}})_{\mathsf{X}}+
(\mbox{\boldmath$B$}\dot{{u}}|\dot{{u}})_{\mathsf{X}} \Big]+ \\
+&(1+a+ \beta)(\|{u}\|_{\mathsf{Y}}^2+
\|\dot{{u}}\|_{\mathsf{X}}^2) \\
&\geq -(1+a)\mathfrak{Re}[(\dot{{u}}|{u})_{\mathsf{X}}]-\beta \|\dot{{u}}\|_{\mathsf{X}}^2 +\\
+&(1+a+ \beta)\|{u}\|_{\mathsf{X}}^2+
(1+a + \beta)\|\dot{{u}}\|_{\mathsf{X}}^2 \\
&\geq
(1+a)\Big[ \|\dot{{u}}\|_{\mathsf{X}}^2
-\mathfrak{Re}[(\dot{{u}}|{u})_{\mathsf{X}}]+
\|{u}\|_{\mathsf{X}}^2\Big] \\
&\geq 0,
\end{align*}
since $(({L}+a+1){u}|\dot{{u}})_{\mathsf{X}}=({u} | \dot{{u}})_{\mathsf{Y}}=
(\dot{{u}}| {u})_{\mathsf{Y}}^*$ and 
$|\mathfrak{Re}[({B}\dot{{u}}|\dot{{u}})_{\mathsf{X}}]| \leq \beta 
\|\dot{{u}}\|_{\mathsf{X}}^2$.\\

If $\Lambda > \frac{\beta}{2}+\sqrt{\frac{\beta^2}{4}+a}$, then the operator $\mathbf{A}+\Lambda$ has the bounded inverse defined on $\mathsf{Z}$. Actually the equation
$$\mathbf{A}U+\Lambda U=F=
\begin{bmatrix}
{f} \\
\\
\mbox{\boldmath$g$}
\end{bmatrix}
\in \mathsf{Z}
$$
means
$$
\begin{cases}
-\dot{{u}}+\Lambda{u}={f} \in \mathsf{Y} \\
\\
{B}
\dot{{u}}+{L}{u}+\Lambda\dot{{u}}=\mbox{\boldmath$g$} \in \mathsf{X},
\end{cases}
$$
which can be solved as
$$
\begin{cases}
{u}=
({L}+\Lambda
{B}
+\Lambda^2)^{-1}({B}{f}+\Lambda{f}+\mbox{\boldmath$g$})
\in \mathsf{D}({L}),\\
\\
\dot{{u}}=
({L}+\Lambda
{B}
+\Lambda^2)^{-1}({B}{f}+\Lambda{f}+\mbox{\boldmath$g$})
-{f} \in \mathsf{Y},
\end{cases}
$$
thanks to Proposition A1,
since $\Lambda^2 >a+|\Lambda|\beta$ holds for
$\Lambda >\frac{\beta}{2}+\sqrt{\frac{\beta^2}{4}+a}$.  $\square$\\

Therefore,
considering the problem {\bf (IVP)}:

\begin{align*}
&\frac{d^2{u}}{d t^2}+{B}\frac{d{u}}{d t}+
{L}{u}={0},
 \\
&{u}(t)\in \mathsf{D}({L})\quad\mbox{for}\quad \forall t \geq 0,  \\
&{u}=\overset{\circ}{{u}},\quad \frac{d{u}}{d t}=\overset{\circ}{{v}}
\quad \mbox{at}\quad t=0,
\end{align*}
 we can claim \\

{\bf Theorem A}\  {\it
Suppose $\overset{\circ}{{u}} \in \mathsf{D}({L})$, 
$\overset{\circ}{{v}}\in \mathsf{Y}$ and $f \in C([0,+\infty[; \mathsf{X})$.
 Then the initial-boundary value problem {\bf (IVP)}  admits a unique solution
$${u}\in C^2([0,+\infty[, \mathsf{X})\cap
C^1([0,+\infty[,\mathsf{Y})
\cap C([0,+\infty[, \mathsf{D}({L}))$$
and the energy
$$
E(t)=\|{u}\|_{\mathsf{Y}}^2+\|\dot{{u}}\|_{\mathsf{X}}^2 
$$
enjoys the estimate
$$
\sqrt{E(t)}\leq e^{\Lambda t}\Big(\sqrt{E(0)}
+\int_0^te^{-\Lambda s}\|f(s)\|_{\mathsf{X}}ds \Big), $$
where $\Lambda=1+a+ \beta$.} \\

Here $\mathsf{D}({L})$ is equipped with the norm
$(\|{u}\|_{\mathsf{Y}}^2+\|{L}{u}\|_{\mathsf{X}}^2)^{1/2}$.\\

\vspace{10mm}

{\Large \bf Appendix B}\\

Let us imagine the motion $\bm{x}=\bm{x}_p(t)$ of an infenitesimally small percel in the background $(\rho_b, S_b)$. The density $\rho_p(\bm{x}(t))$, pressure $P_p(\bm{x}(t))$,  entropy density $S_p(\bm{x}(t))$ of the percel obey the same EOS: $P=\mathscr{P}(\rho,S)$ as the background. 
Suppose that $P_p(\bm{x}(t))=P_b(\bm{x}(t))$.\\

Here we note the following:
{\it Let $\rho=\mathscr{R}(\rho, S)$ be the inverse function of $P=\mathscr{P}(\rho,S)$, so
$$\rho=\mathscr{R}(\mathscr{P}(\rho,S),S);$$
then it holds
\begin{align*}
\partial_S\mathscr{R}(P, S)&=
-(\partial_{\rho}\mathscr{R}(P,S))(\partial_{S}\mathscr{P}(\mathscr{R}(P,S),S) )\\
&=-\frac{\partial_S\mathscr{P}(\rho,S)}{\partial_{\rho}\mathscr{P}(\rho,S)}\Big|_{\rho=\mathscr{R}(P,S)} ;
\end{align*}
therefore, by Definition \ref{Def.G}, \eqref{Def.aa}, we have
$$
\partial_S\mathscr{R}(P,S)\nabla S=\rho\mathfrak{a}.
$$}\\

Now the equation of motion of the parcel is 
$$
\frac{d^2\bm{x}_p}{dt^2}=-\frac{\nabla P_b}{\rho_p}+\bm{F}(t,\bm{x}_p(t)),
$$while the background satisfies
$$0=-\frac{\nabla P_b}{\rho_b}+\bm{F}(t,\bm{x})
$$ for some vector field of force.
( In our situation, we are concerned with
$$\bm{F}=-\nabla \Phi_b
-(D\bm{v}_b)\bm{v}_b-2\bm{\Omega}\times \bm{v}_b
-\bm{\Omega}\times(\bm{\Omega}\times \bm{x}).
$$
But we need not the specification of $\bm{F}$, wich contains both real and virtual forces.) Therefore the equation of motion is
$$
\frac{d^2\bm{x}_p}{dt^2}=\Big(\frac{1}{\rho_b}-\frac{1}{\rho_p}\Big)\nabla P_b(\bm{x}_p(t)).
$$
On the other hand we have approximately
\begin{align*}
\frac{1}{\rho_b}-\frac{1}{\rho_p}&\approxeq
\frac{1}{\rho_b^2}\partial_S\mathscr{R}(P_b, S_b)
(\nabla S_b|\bm{x}_p(t)-\bm{x}_p(0)) \\
&=-\frac{1}{\rho_b}(\mathfrak{a}_b|\bm{x}_p(t)-\bm{x}_p(0)).
\end{align*}
Since $\nabla P_b$ is parallel to $\bm{n} \approx \bm{n}_0=\bm{n}(\bm{x}_p(0))$,
we suppose that
$$\bm{x}_p(t)=X(t)\bm{n}_0.$$
Then
$$
\frac{1}{\rho_b}-\frac{1}{\rho_p}\approxeq \frac{1}{\rho_b}\mathscr{A}_b\Big|_{\bm{x}=\bm{x}_p(0)}(X(t)-X(0),
$$
and the equation of motion is
\begin{align*}
\frac{d^2X}{dt^2}&\approxeq
\frac{1}{\rho_b}\mathscr{A}_b(\nabla P_b|\bm{n}_0)(X(t)-X(0)) \\
&\approxeq
-\mathscr{N}_b^2\Big|_{\bm{x}=\bm{x}_p(0)}(X(t)-X(0)). 
\end{align*}
\\

The solution of this equation is  
$$
X(t)=X(0)+C_+e^{\mathrm{i}Nt}+C_-e^{-\mathrm{i}Nt} \quad\mbox{where}\quad N:=\mathscr{N}_b(\bm{x}_p(0)).$$ Clearly
$|X(t)-X(0)|=O(1)$ for non-trivial $X(t)$ if and only if $N \in \mathbb{R}$. Therefore the definition \eqref{DefN} is justifiable.

\vspace{15mm}

{\bf\Large Appendix C}\\

We consider the initial value problem {\bf (IV)}:
$$\frac{d^2u}{dt^2}+B\frac{du}{dt}+Lu=f(t), \quad
u|_{t=0}=\overset{\circ}{u}, \quad \frac{du}{dt}\Big|_{t=0}=\overset{\circ}{v}.
$$

The unknown $u(t)$ is valued in a Hilbert space $\mathsf{X}$. $\mathsf{Y}$ is a Hilbert space which is densely imbedded in $\mathsf{X}$. $L$ is a closed linear operator whose domain $\mathsf{D}(L)$ is a dense subset of $\mathsf{Y}$. 
It holds that

\begin{align*}
& \|u\|_{\mathsf{Y}} \geq \|u\|_{\mathsf{X}} \quad\mbox{for}\quad \forall u \in \mathsf{Y} \\
& \mathfrak{Re}[((L+a+1)u|u)_{\mathsf{X}}]=\|u\|_{\mathsf{Y}}^2  \quad
\forall u \in \mathsf{D}(L), 
\end{align*}
where $a$ is a non-negative number. 
Noreover we suppose
$$|(Lu_1|u_2)_{\mathsf{X}}-(u_1|Lu_2)_{\mathsf{X}}| \leq 2\nu\|u_1\|_{\mathsf{X}}\|u_2\|_{\mathsf{X}} \quad \mbox{for}\quad \forall u_1,u_2 \in \mathsf{D}(L),
$$
where $\nu$ is a non-negative number. 
$B$ is a bounded linear operator defined on 
$\mathsf{X}$ into $\mathsf{X}$ with 
$$\beta:=\sup \frac{|\mathfrak{Re}[(Bv|v)_{\mathsf{X}}]|}{\|v\|_{\mathsf{X}}^2}. $$
\\

{\bf Remark}: \  {\it If $L$ is symmetric, we take $\nu=0$. If $B$ is skew-symmetric, then $\beta=0$. }\\

We  claim\\

{\bf Theorem C}\  {\it If $ \overset{\circ}{u} \in \mathsf{D}(L),
\overset{\circ}{v} \in \mathsf{Y}$
and $f \in C([0,+\infty[; \mathsf{X})$, then the problem {\bf (IVP)} admits a unique solution
$u \in C^2([0,+\infty[; \mathsf{X})\cap
C^1([0,+\infty[; \mathsf{Y})
\cap C([0,+\infty[; \mathsf{D}(L))$. The solution enjoys
$$\sqrt{E(t)}
\leq e^{\kappa t}\sqrt{E(0)}+
\int_0^t e^{\kappa (t-s)}\|f(s)\|_{\mathsf{X}}ds.
$$
Here $$\kappa=\kappa(a,\nu,\beta):=
\Big(\frac{5}{4}+\frac{a+\nu}{2}\Big)\vee
\Big(1+\frac{a+\nu}{2}+\beta\Big).
$$
 and }
$$E(t):=\|u(t)\|_{\mathsf{Y}}^2+\Big\|\frac{du}{dt}\Big\|_{\mathsf{X}}^2.
$$
 \\

Let us sketch the proof of {\bf Theorem C}. First we claim\\

{\bf Proposition C1 }\  {\it Let $c\in \mathbb{R}, \lambda >a+|c|\beta$. Then $\lambda+cB+L$ admits the bounded inverse
$(\lambda+cB+L)^{-1}$ with
$$
|\|(\lambda+cB+L)^{-1}\||_{\mathcal{B}(\mathsf{X})}
\leq \frac{1}{\lambda -a -|c|\beta}. \eqno{(C1)}
$$
}\\

Proof. Given $f \in \mathsf{X}$, we have to solve the equation
$$\lambda u+cB u+Lu=f. \eqno{(C2)}
$$
If $u \in \mathsf{D}(L)$ satisfies (C2), then we have
$$\|u\|_{\mathsf{X}} \leq \frac{1}{\lambda-a-|c|\beta}\|f\|_{\mathsf{X}}.
$$
In fact, (C2) implies
$$\mathfrak{Re}[(\lambda u+Lu|u)_{\mathsf{X}}]+c\mathfrak{Re}[(Bu|u)_{\mathsf{X}}]=
\mathfrak{Re}[(f|u)_{\mathsf{X}}],
$$
and
\begin{align*}
\mathfrak{Re}[(\lambda u+Lu|u)_{\mathsf{X}}] &\geq
\mathfrak{Re}[(a+L)u|u)_{\mathsf{X}}]+(\lambda-a)\|u\|_{\mathsf{X}}^2 \\
&\geq (\lambda-a)\|u\|_{\mathsf{X}}^2, \\
c\mathfrak{Re}[(Bu|u)_{\mathsf{X}}]&\geq -|c|\beta\|u\|_{\mathsf{X}}^2 \\
\mathfrak{Re}[(f|u)_{\mathsf{X}}] &\leq \|f\|_{\mathsf{X}}\|u\|_{\mathsf{X}},
\end{align*}
therefore
$$
(\lambda-a-|c|\beta)\|u\|_{\mathsf{X}} \leq \|f\|_{\mathsf{X}}.
$$

Hence the range $\mathsf{R}(\lambda+cB+L)=
\mathsf{D}((\lambda+cB+L)^{-1})$ is a closed subspace of $\mathsf{X}$, since $\lambda+cB+L$ is a closed operator. We claim that 
 $\mathsf{R}(\lambda+cB+L)=\mathsf{X}$. In fact, consider $h \in \mathsf{X}$ which satisfies
$$((\lambda+cB+L)u|h)_{\mathsf{X}}=0 \quad\mbox{for}\quad \forall u \in \mathsf{D}(L).
$$
Then $h \in \mathsf{D}(L^*)$ and $$(u|(\lambda+cB^*+L^*)h)_{\mathsf{X}}=0
\quad\mbox{for}\quad \forall u \in \mathsf{D}(L).
$$
Since $\mathsf{D}(L)$ is dense in $\mathsf{X}$, we have
$$(\lambda+cB^*+L^*)h=0,$$
and
$$\mathfrak{Re}[(\lambda+L^*)h|h)_{\mathsf{X}}]
+c\mathfrak{Re}[(B*h|h)_{\mathsf{X}}]=0.
$$Since
$$((\lambda+L^*)h|h)_{\mathsf{X}}=((\lambda+L)h|h)_{\mathsf{X}}^*,
\quad
(B^*h|h)_{\mathsf{X}}=(Bh|h)_{\mathsf{X}}^*,
$$
we have
$$(\lambda-a)\|h\|_{\mathsf{X}}^2-|c|\beta\|h\|_{\mathsf{X}}^2 \leq 0.$$
Therefore $h=0$, for $\lambda-a-|c|\beta >0$. Therefore
 $\mathsf{R}(\lambda+cB+L)=\mathsf{X}$ and $(\lambda+cB+L)^{-1} \in \mathcal{B}(\mathsf{X})$ and enjoys (C1). $\square$\\

Now we write the problem {\bf (IVP)} as {\bf (IVP\ddag)}:
$$
\frac{dU}{dt}=\mathfrak{A} U+F,\quad
U(0)=U_0 ,
$$
where
\begin{align*}
&U=
\begin{bmatrix}
u \\
\\
v
\end{bmatrix}
=\begin{bmatrix}
u \\
\\
\dot{u}
\end{bmatrix}
=
\begin{bmatrix}
u \\
\\
\frac{du}{dt}
\end{bmatrix},\quad
U_0=\begin{bmatrix}
\overset{\circ}{u} \\
\\
\overset{\circ}{v}
\end{bmatrix}, \\
&\mathfrak{A}=
\begin{bmatrix}
O & I \\
\\
-L& -B
\end{bmatrix}
,
\quad
F=\begin{bmatrix}
0 \\
\\
f
\end{bmatrix}
\end{align*}
and apply Hile-Yosida theory. We consider the problem in the Hilbert space
$\mathsf{Z}=\mathsf{Y}\times \mathsf{X}$ endowed with the norm
$$\|U\|_{\mathsf{Z}}=\Big[\|u\|_{\mathsf{Y}}^2+\|\dot{u}\|_{\mathsf{X}}^2\Big]^{\frac{1}{2}},
$$
while
$\mathsf{D}(\mathfrak{A})=\mathsf{D}(L)\times \mathsf{Y} $.
We claim \\

{\bf Proposition C2 }: \  {\it Put $\kappa=\kappa(a,\nu, \beta)$. Consider $n \in \mathbb{R}$ such that  $|n|$ is large. 

The equation
$$\Big(I-\frac{1}{n}
\mathfrak{A}\Big)U=F,
$$
$F =( g, f)^{\top} \in \mathsf{Z}$ given, admits a unique solution $U \in \mathsf{D}(\mathfrak{A})$, and
$$\|U\|_{\mathsf{Z}} \leq \frac{1}{1-\frac{\kappa}{|n|}}\|F\|_{\mathsf{Z}}.
$$
That is,
$$\Big|\Big\|\Big(I-\frac{1}{n}\mathfrak{A}\Big)^{-1}\Big\|\Big|_{\mathcal{B}(\mathsf{Z})} \leq \frac{1}{1-\frac{\kappa}{|n|}}.
$$
}\\

Proof. We have to solve
\begin{align*}
&u-\frac{1}{n}v =g ( \in \mathsf{Y}), \\
&v+\frac{1}{n}Bv +\frac{1}{n}Lu= f (\in \mathsf{X}).
\end{align*}
Since that $|n| $ is so large that
 $n^2> a+|n|\beta$. Then {\bf Proposition C1} gurantees 
$(n^2+n B+ L)^{-1} \in \mathcal{B}(\mathsf{X})$ and
$$|\|(n^2+nB+L)^{-1}\||_{\mathcal{B}(\mathsf{X})} 
\leq \frac{1}{n^2-a-|n|\beta}.
$$
Keeping in mind this, we can find the solution
\begin{align*}
u&=(n^2+nB+L)^{-1}(nBg +n^2g+nf), \\
v&=n(u-g)=n\Big[(n^2+nB+L)^{-1}(nBg +n^2g+nf)-g\Big].
\end{align*}
Note that $u \in \mathsf{D}(L)$ but $v \in \mathsf{Y}$. Temporally we suppose that
$g \in \mathsf{D}(L)$ so that $v \in \mathsf{D}(L)$. In this situation, we are going to estimate
$\mathfrak{Re}[((\lambda +L)u|u)_{\mathsf{X}}]+\|v\|_{\mathsf{X}}^2$ by
$\mathfrak{Re}[((\lambda +L)g|g)_{\mathsf{X}}]+\|f\|_{\mathsf{X}}^2$, where we put $$\lambda:=1+a.$$ Note that
\begin{align*}
&\mathfrak{Re}[((\lambda +L)u|u)_{\mathsf{X}}]+\|v\|_{\mathsf{X}}^2 =\|u\|_{\mathsf{Y}}^2+\|v\|_{\mathsf{X}}^2, \\
&\mathfrak{Re}[((\lambda +L)g|g)_{\mathsf{X}}]+\|f\|_{\mathsf{X}}^2
=\|g\|_{\mathsf{Y}}^2+\|f\|_{\mathsf{X}}.
\end{align*}

Using 
$$g=u-\frac{1}{n}v, \quad
f=v+\frac{1}{n}Lu+\frac{1}{n}Bv, 
$$
we see
\begin{align*}
((\lambda+L)g|g)_{\mathsf{X}}&=((\lambda+L)u|u)_{\mathsf{X}}-\frac{\lambda}{n}((u|v)_{\mathsf{X}}+(v|u)_{\mathsf{X}}) + \\
&-\frac{1}{n}\Big((Lu|v)_{\mathsf{X}}+(Lv|u)_{\mathsf{X}}\Big)+
\frac{1}{n^2}((\lambda+L)v|v)_{\mathsf{X}}, \\
\|f\|_{\mathsf{X}}^2&=\|v\|_{\mathsf{X}}^2+
\frac{1}{n}\Big((v|Lu)_{\mathsf{X}}+(Lu|v)_{\mathsf{X}}\Big)+
\frac{1}{n^2}\|Lu\|_{\mathsf{X}}^2 + \\
&+\frac{1}{n}((v|Bv)_{\mathsf{X}}+(Bv|v)_{\mathsf{X}})+
\frac{1}{n^2}(Bv|Lu)_{\mathsf{X}}+\frac{1}{n^2}\|Bv\|_{\mathsf{X}}^2.
\end{align*}
Therefore
\begin{align*}
((\lambda+L)g|g)_{\mathsf{X}}+\|f\|_{\mathsf{X}}^2&=((\lambda+L)u|u)_{\mathsf{X}}+
\frac{1}{n^2}\|Lu\|_{\mathsf{X}}^2 + \\
&+\frac{1}{n^2}((\lambda+L)v|v)_{\mathsf{X}}+\|v\|_{\mathsf{X}}^2+\mathfrak{r},
\end{align*}
where
\begin{align*}
\mathfrak{r}&:=-\frac{2\lambda}{n}\mathfrak{Re}[(u|v)_{\mathsf{X}}]+
\frac{1}{n}\Big(-(Lv|u)_{\mathsf{X}}+(v|Lu)_{\mathsf{X}}\Big) + \\
&+\frac{2}{n}\mathfrak{Re}[(Bv|v)_{\mathsf{X}}]+
\frac{1}{n^2}(Bv|Lu)_{\mathsf{X}}+\frac{1}{n^2}\|Bv\|_{\mathsf{X}}^2 \\
&=
-\frac{2\lambda}{n}\mathfrak{Re}[(u|v)_{\mathsf{X}}]+
\frac{1}{n}\Big(-(Lv|u)_{\mathsf{X}}+(v|Lu)_{\mathsf{X}}\Big) + \\
&+\frac{2}{n}\mathfrak{Re}[(Bv|v)_{\mathsf{X}}]+ \\
&+\frac{1}{n}\Big(f-v-\frac{1}{n}Lu\Big|u\Big)_{\mathsf{X}}+
\Big\|f-v-\frac{1}{n}Lu\Big\|_{\mathsf{X}}^2.
\end{align*}

Since
\begin{align*}
&\Big|-\frac{2\lambda}{n}\mathfrak{Re}[(u|v)_{\mathsf{X}}]\Big| \leq \frac{\lambda}{|n|}(\|u\|_{\mathsf{X}}^2+\|v\|_{\mathsf{X}}^2), \\
&\Big|\frac{1}{n}\mathfrak{Re}\Big[-(Lv|u)_{\mathsf{X}}+(v|Lu)_{\mathsf{X}}\Big]
\Big|\leq \frac{\nu}{|n|}(\|u\|_{\mathsf{X}}^2+\|v\|_{\mathsf{X}}^2), \\
&\Big|\frac{2}{n}\mathfrak{Re}[(Bv|v)_{\mathsf{X}}] \Big|\leq \frac{2\beta}{|n|}\|v\|_{\mathsf{X}}^2, \\
&
\Big|\frac{1}{n}\mathfrak{Re}\Big[
\Big(f-v-\frac{1}{n}Lu\Big|u\Big)_{\mathsf{X}} \Big]\Big|
\leq \frac{1}{|n|}
\Big(\frac{1}{2}\|f\|_{\mathsf{X}}^2+
\Big(1+\frac{1}{2|n|}\Big)\|u\|_{\mathsf{X}}^2+
\frac{1}{2}\|v\|_{\mathsf{X}}^2+
\frac{1}{2|n|}\|Lu\|_{\mathsf{X}}^2
\Big),
\end{align*}
we have the estimtae
\begin{align*}
\mathfrak{Re}[\mathfrak{r}]&\geq
-\frac{1}{|n|}\Big(\lambda+\nu+1+\frac{1}{2|n|}\Big)\|u\|_{\mathsf{X}}^2
-\frac{1}{|n|}\Big(\lambda+\nu+2\beta+\frac{1}{2}\Big)\|v\|_{\mathsf{X}}^2 \\
&-\frac{1}{2n^2}\|Lu\|_{\mathsf{X}}^2
-\frac{1}{2|n|}\|f\|_{\mathsf{X}}^2.
\end{align*}
Consequently we have
\begin{align*}
\mathfrak{Re}[((\lambda+L)g|g)_{\mathsf{X}}]+\|f\|_{\mathsf{X}}^2&\geq\mathfrak{Re}[(\lambda+L)u|u)_{\mathsf{X}}] \\
&-
\frac{1}{|n|}\Big(\lambda+\nu+1+\frac{1}{2|n|}\Big)\|u\|_{\mathsf{X}}^2+
\Big(1-\frac{1}{|n|}\Big(\lambda+\nu+2\beta+\frac{1}{2}\Big)\Big)\|v\|_{\mathsf{X}}^2 \\
&
-\frac{1}{2|n|}\|f\|_{\mathsf{X}}^2,
\end{align*}
or
\begin{align*}
\Big(1+\frac{1}{2|n|}\Big)
\Big(\mathfrak{Re}[((\lambda+L)g|g)_{\mathsf{X}}]+\|f\|_{\mathsf{X}}^2\Big)
&\geq
\Big(1-\frac{1}{|n|}\Big(\lambda+\nu+1+\frac{1}{2|n|}\Big)\Big)
\mathfrak{Re}[(\lambda+L)u|u)_{\mathsf{X}}] + \\
&+
\Big(1-\frac{1}{|n|}\Big(\lambda+\nu+2\beta+\frac{1}{2}\Big)\Big)\|v\|_{\mathsf{X}}^2 .
\end{align*}
If we take $\kappa >\kappa(\nu, a, \beta)=\Big(\frac{5}{4}+\frac{a+\nu}{2}\Big)\vee
\Big(1+\frac{a+\nu}{2}+\beta\Big)$, it holds
$$
\mathfrak{Re}[((\lambda+L)g|g)_{\mathsf{X}}]+\|f\|_{\mathsf{X}}^2\geq
\Big(1-\frac{\kappa}{|n|}\Big)^2
\Big(\mathfrak{Re}[((\lambda+L)u|u)_{\mathsf{X}}] + \|v\|_{\mathsf{X}}^2\Big)
$$
for $|n|\geq N_*$, $N_*$ being a sufficiently large positive number.
In fact, we see
\begin{align*}
&
\Big(1+\frac{1}{2|n|}\Big)^{-1}
\Big(1-
\frac{1}{|n|}\Big(\lambda+\nu+1+\frac{1}{2|n|}\Big)\Big)
\Big(1-\frac{\kappa}{|n|}\Big)^{-2} \\
&=\Big(2\kappa -\Big(\frac{5}{2}+a+\nu\Big)\Big)\frac{1}{|n|}+O\Big(\frac{1}{n^2}\Big), \\
&
\Big(1+\frac{1}{2|n|}\Big)^{-1}
\Big(1-\frac{1}{|n|}\Big(\lambda+\nu+2\beta+\frac{1}{2}\Big)\Big)
\Big(1-\frac{\kappa}{|n|}\Big)^{-2} \\
&=\Big(2\kappa -\Big(2+a+\nu+2\beta \Big)\Big)\frac{1}{|n|}+O\Big(\frac{1}{n^2}\Big).
\end{align*}\\

This means
$$
\|F\|_{\mathsf{Z}} \geq \Big(1-\frac{\kappa}{|n|}\Big)\|U\|_{\mathsf{Z}}.$$

Passing to the limit, we can claim this is valid for $f \in \mathsf{Y}, v \in Y_0$, and $\kappa=\kappa(\nu,a,\beta)$
$\square$.\\

Consequently we can apply \cite[Theorem and Corollary 1 of IX.9]{Yosida}
to conclude that \\

{\it $\mathfrak{A}$ is the infenitesimal generator of an equi-continuous group of class ($C_0$)
of operators $\mathscr{S}_t   \in \mathcal{B}(\mathsf{Z}), -\infty <t<+\infty$; If $U_0 \in \mathsf{D}(\mathfrak{A})$ and
$F \in C([0,+\infty[; \mathsf{Z})$, then
$$ U(t)=\mathscr{S}_tU_0+
\int_0^t\mathscr{S}_{t-s}F(s)ds $$
is the  solution of {\bf (IVP\ddag)}
of the class $C^1([0,+\infty[,\mathsf{Z})\cap
C([0,+\infty[; \mathsf{D}(\mathfrak{A}))$. It enjoys
$$\|U(t)\|_{\mathsf{Z}}
\leq e^{\kappa t}\|U_0\|_{\mathsf{Z}}+
\int_0^t e^{\kappa (t-s)}\|F(s)\|_{\mathsf{Z}}ds
$$
for $t \geq 0$, since $|\|\mathscr{S}_t\||_{\mathcal{B}(\mathsf{Z})} \leq e^{\kappa t}$ for $t \geq 0$.}\\

 Interpreting this conclusion to the words of $u$, we completes the proof of {\bf Theorem C}.

%%%%%%%%%%%%%%%%%%%%%%%%%%%%%%%%%%%%

\vspace{15mm}

{\bf\Large Appendix D}\\

Let us consider a motion of infinitesimally small parcel of fluid in the barotropic background $(\rho_b, S_b, \bm{v}_b)$ described in the non-rotating cylindrical coordinate system. So,
$$
\bm{v}_b^<=
\begin{bmatrix}
v^{\varpi} \\
\\
v^z \\
\\
v^{\phi}
\end{bmatrix}
=\begin{bmatrix}
0 \\
\\
0 \\
\\
\varpi^2\Omega_b(\varpi)
\end{bmatrix}
$$
and the motion of the parcell is 
$$\bm{x}=\bm{x}_p(t)=
\begin{bmatrix}
x_p^1(t) \\
\\
x_p^2 (t) \\
\\
x_p^3(t)
\end{bmatrix}
=\begin{bmatrix}
\varpi_p(t) \\
\\
z_p(t) \\
\\
\phi_p(t)
\end{bmatrix}
^>.
$$

We suppose that $\displaystyle v_p^{\phi}(\bm{x}_p(0))=\frac{d}{dt} \phi_p(t)\Big|_{t=0}=0$. 

The equation of motion of the paercell is
$$\frac{d^2 \varpi_p}{dt^2}=K^{\varpi},
$$
where $K^{\varpi}$ is the $\varpi$-component of the force per unit mass, while the background satisfies
$$-\varpi \Omega_b(\varpi)^2=K^{\varpi}.
$$
Here 
$$K^{\varpi}=-\Big[ \frac{1}{\rho}\frac{\partial P}{\partial \varpi}+\frac{\partial \Phi}{\partial \varpi}\Big]_b.
$$
Recall that the $\varpi$-component of the equation of motion reads
$$
\Big[\frac{\partial}{\partial t}+v^{\varpi}\frac{\partial}{\partial \varpi}+ v^z\frac{\partial}{\partial z}+\frac{v^{\phi}}{\varpi}\frac{\partial}{\partial \phi}\Big]v^{\varpi}-\frac{(v^{\phi})^2}{\varpi}
+ \frac{1}{\rho}\frac{\partial P}{\partial \varpi}+\frac{\partial \Phi}{\partial \varpi}=0.
$$

We suppose that the specific angular momentum 
$J_p=\varpi_pv_p^{\phi}$ of the parcel coincides with that $J_b=\varpi^3\Omega_b(\varpi)$ of the background. Then it holds that
$ v_p^{\phi}=\varpi_p^2\Omega_b(\varpi_p)$.

Under these assumptions, the equation of motion of the parcel turns out to be
\begin{align*}
\frac{d^2\varpi_p}{dt^2}&=K^{\varpi}=-\varpi \Omega_b^2=
-\varpi\Big(\frac{v_p^{\phi}}{\varpi^2}\Big)^2
=-\frac{1}{\varpi^3}(v_p^{\phi})^2 \\
&\approxeq -\frac{1}{\varpi^3}\Big[(v_p^{\phi})^2|_{t=0}+
\frac{\partial}{\partial \varpi}(v_p^{\phi})^2 \Big|_{t=0}
(\varpi_p(t)-\varpi_p(0))
\Big] \\
&=-\frac{1}{\varpi^3}\frac{\partial}{\partial\varpi}(\varpi^4\Omega_b^2)\Big|_{t=0}(\varpi_p(t)-\varpi(0)),
\end{align*}
that is, approximately,
$$\frac{d^2\varpi_p}{d t^2}+\mathscr{C}(\varpi_p(0))(\varpi_p(t)-\varpi_p(0))=0.
$$

Now it is clear that the motion $\varpi_p(t)$ is bounded for any initial data if and only if
$\mathscr{C}(\varpi_p(0)) \geq 0$, the Rayleigh's criterion.

\vspace{15mm}

{\bf\Large  Acknowledgment}\\

To do this study has been caused by helpful discussions with Professors Akitaka Matsumura (Osaka University) and Masao Takata (the University of Tokyo), to whom sincere thanks are expressed. This work was supported by Japan Society for the Promotion of Science (JSPS) Grant-in-Aid for Scientific Research (KAKENHI)
grant JP21K03311, and by the Research Institute for Mathematical Sciences, an International Joint Usage/Research Center located in Kyoto University.

\vspace{15mm}

{\bf\Large  The data availability statement}\\

No new data were created or analyzed in this study.

% V\"{a}is\"{a}l\"{a}

\end{document}